# Super Linear Algebra

# SUPER LINEAR ALGEBRA

**W. B. Vasantha Kandasamy**
e-mail: **vasanthakandasamy@gmail.com**
web: **http://mat.iitm.ac.in/~wbv**

**Florentin Smarandache**
e-mail: **smarand@unm.edu**

# CONTENTS

# PREFACE

In this book, the authors introduce the notion of Super linear algebra and super vector spaces using the definition of super matrices defined by Horst (1963). This book expects the readers to be well-versed in linear algebra.

Many theorems on super linear algebra and its properties are proved. Some theorems are left as exercises for the reader. These new class of super linear algebras which can be thought of as a set of linear algebras, following a stipulated condition, will find applications in several fields using computers. The authors feel that such a paradigm shift is essential in this computerized world. Some other structures ought to replace linear algebras which are over a century old.

Super linear algebras that use super matrices can store data not only in a block but in multiple blocks so it is certainty more powerful than the usual matrices.

This book has 3 chapters. Chapter one introduces the notion of super vector spaces and enumerates a number of properties. Chapter two defines the notion of super linear algebra, super inner product spaces and super bilinear forms. Several interesting properties are derived. The main application of these new structures in Markov chains and Leontief economic models are also given in this chapter. The final chapter suggests 161 problems mainly to make the reader understand this new concept and apply them.

The authors deeply acknowledge the unflinching support of Dr.K.Kandasamy, Meena and Kama.

W.B.VASANTHA KANDASAMY
FLORENTIN SMARANDACHE

## Chapter One

# SUPER VECTOR SPACES

This chapter has four sections. In section one a brief introduction about supermatrices is given. Section two defines the notion of super vector spaces and gives their properties. Linear transformation of super vector is described in the third section. Final section deals with linear algebras.

### 1.1 Supermatrices

Though the study of super matrices started in the year 1963 by Paul Horst. His book on matrix algebra speaks about super matrices of different types and their applications to social problems. The general rectangular or square array of numbers such as

$$A = \begin{bmatrix} 2 & 3 & 1 & 4 \\ -5 & 0 & 7 & -8 \end{bmatrix},\ B = \begin{bmatrix} 1 & 2 & 3 \\ -4 & 5 & 6 \\ 7 & -8 & 11 \end{bmatrix},$$

$$C = [3, 1, 0, -1, -2] \text{ and } D = \begin{bmatrix} -7/2 \\ 0 \\ \sqrt{2} \\ 5 \\ -41 \end{bmatrix}$$

are known as matrices.

We shall call them as simple matrices [17]. By a simple matrix we mean a matrix each of whose elements are just an ordinary number or a letter that stands for a number. In other words, the elements of a simple matrix are scalars or scalar quantities.

A supermatrix on the other hand is one whose elements are themselves matrices with elements that can be either scalars or other matrices. In general the kind of supermatrices we shall deal with in this book, the matrix elements which have any scalar for their elements. Suppose we have the four matrices;

$$a_{11} = \begin{bmatrix} 2 & -4 \\ 0 & 1 \end{bmatrix}, \; a_{12} = \begin{bmatrix} 0 & 40 \\ 21 & -12 \end{bmatrix}$$

$$a_{21} = \begin{bmatrix} 3 & -1 \\ 5 & 7 \\ -2 & 9 \end{bmatrix} \text{ and } a_{22} = \begin{bmatrix} 4 & 12 \\ -17 & 6 \\ 3 & 11 \end{bmatrix}.$$

One can observe the change in notation $a_{ij}$ denotes a matrix and not a scalar of a matrix $(1 \leq i, j \leq 2)$.

Let

$$a = \begin{bmatrix} a_{11} & a_{12} \\ a_{21} & a_{22} \end{bmatrix};$$

we can write out the matrix a in terms of the original matrix elements i.e.,

$$a = \left[\begin{array}{cc|cc} 2 & -4 & 0 & 40 \\ 0 & 1 & 21 & -12 \\ \hline 3 & -1 & 4 & 12 \\ 5 & 7 & -17 & 6 \\ -2 & 9 & 3 & 11 \end{array}\right].$$

Here the elements are divided vertically and horizontally by thin lines. If the lines were not used the matrix a would be read as a simple matrix.

Thus far we have referred to the elements in a supermatrix as matrices as elements. It is perhaps more usual to call the elements of a supermatrix as submatrices. We speak of the submatrices within a supermatrix. Now we proceed on to define the order of a supermatrix.

The order of a supermatrix is defined in the same way as that of a simple matrix. The height of a supermatrix is the number of rows of submatrices in it. The width of a supermatrix is the number of columns of submatrices in it.

All submatrices with in a given row must have the same number of rows. Likewise all submatrices with in a given column must have the same number of columns.

A diagrammatic representation is given by the following figure.

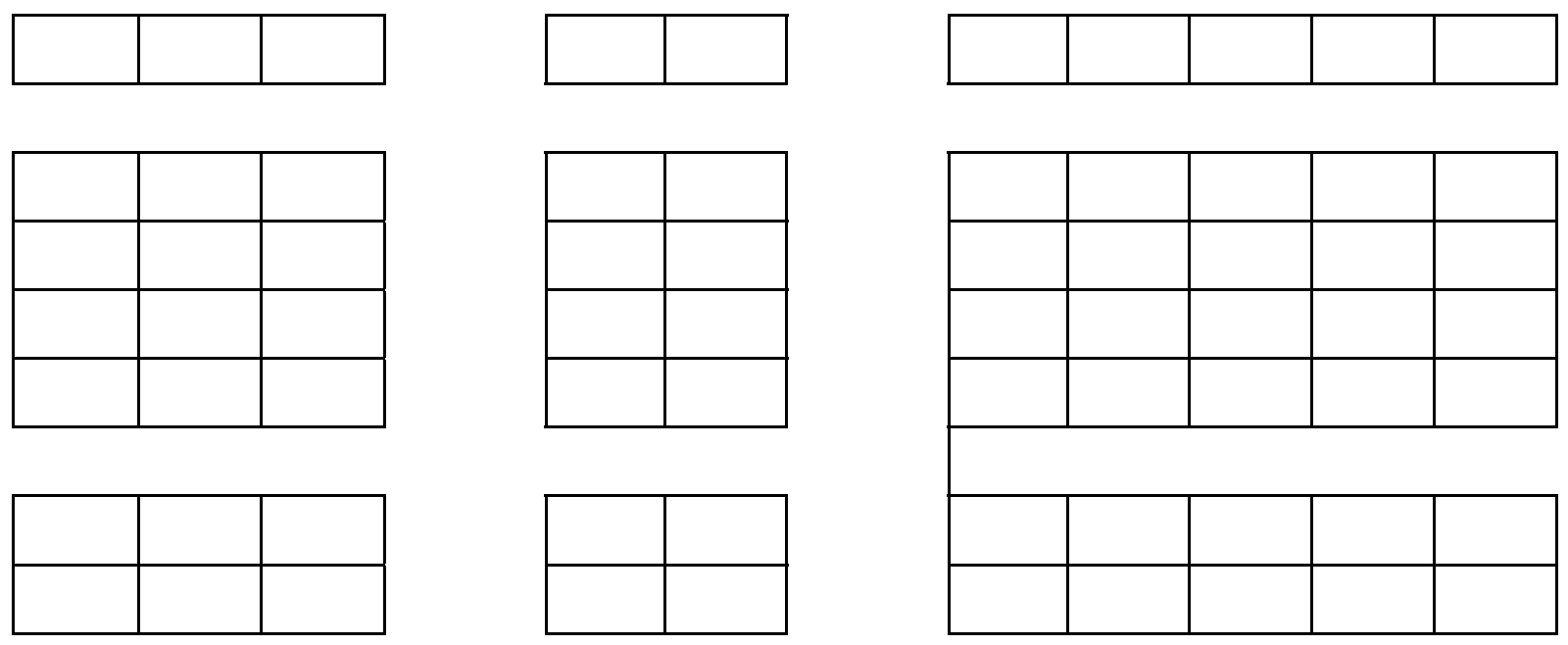

In the first row of rectangles we have one row of a square for each rectangle; in the second row of rectangles we have four rows of squares for each rectangle and in the third row of rectangles we have two rows of squares for each rectangle. Similarly for the first column of rectangles three columns of squares for each rectangle. For the second column of rectangles we have two column of squares for each rectangle, and for the third column of rectangles we have five columns of squares for each rectangle.

Thus we have for this supermatrix 3 rows and 3 columns.

One thing should now be clear from the definition of a supermatrix. The super order of a supermatrix tells us nothing about the simple order of the matrix from which it was obtained by partitioning. Furthermore, the order of supermatrix tells us nothing about the orders of the submatrices within that supermatrix.

Now we illustrate the number of rows and columns of a supermatrix.

***Example 1.1.1:*** Let

$$a = \left[\begin{array}{c|cccc} 3 & 3 & 0 & 1 & 4 \\ -1 & 2 & 1 & -1 & 6 \\ \hline 0 & 3 & 4 & 5 & 6 \\ 1 & 7 & 8 & -9 & 0 \\ 2 & 1 & 2 & 3 & -4 \end{array}\right].$$

a is a supermatrix with two rows and two columns.

Now we proceed on to define the notion of partitioned matrices. It is always possible to construct a supermatrix from any simple matrix that is not a scalar quantity.

The supermatrix can be constructed from a simple matrix this process of constructing supermatrix is called the partitioning.

A simple matrix can be partitioned by dividing or separating the matrix between certain specified rows, or the procedure may be reversed. The division may be made first between rows and then between columns.

We illustrate this by a simple example.

***Example 1.1.2:*** Let

$$A = \begin{bmatrix} 3 & 0 & 1 & 1 & 2 & 0 \\ 1 & 0 & 0 & 3 & 5 & 2 \\ 5 & -1 & 6 & 7 & 8 & 4 \\ 0 & 9 & 1 & 2 & 0 & -1 \\ 2 & 5 & 2 & 3 & 4 & 6 \\ 1 & 6 & 1 & 2 & 3 & 9 \end{bmatrix}$$

is a 6 × 6 simple matrix with real numbers as elements.

$$A_1 = \left[\begin{array}{cc|cccc} 3 & 0 & 1 & 1 & 2 & 0 \\ 1 & 0 & 0 & 3 & 5 & 2 \\ 5 & -1 & 6 & 7 & 8 & 4 \\ 0 & 9 & 1 & 2 & 0 & -1 \\ 2 & 5 & 2 & 3 & 4 & 6 \\ 1 & 6 & 1 & 2 & 3 & 9 \end{array}\right].$$

Now let us draw a thin line between the 2nd and 3rd columns.

This gives us the matrix $A_1$. Actually $A_1$ may be regarded as a supermatrix with two matrix elements forming one row and two columns.

Now consider

$$A_2 = \left[\begin{array}{cccccc} 3 & 0 & 1 & 1 & 2 & 0 \\ 1 & 0 & 0 & 3 & 5 & 2 \\ 5 & -1 & 6 & 7 & 8 & 4 \\ 0 & 9 & 1 & 2 & 0 & -1 \\ \hline 2 & 5 & 2 & 3 & 4 & 6 \\ 1 & 6 & 1 & 2 & 3 & 9 \end{array}\right]$$

Draw a thin line between the rows 4 and 5 which gives us the new matrix $A_2$. $A_2$ is a supermatrix with two rows and one column.

Now consider the matrix

$$A_3 = \left[\begin{array}{cc|cccc} 3 & 0 & 1 & 1 & 2 & 0 \\ 1 & 0 & 0 & 3 & 5 & 2 \\ 5 & -1 & 6 & 7 & 8 & 4 \\ 0 & 9 & 1 & 2 & 0 & -1 \\ \hline 2 & 5 & 2 & 3 & 4 & 6 \\ 1 & 6 & 1 & 2 & 3 & 9 \end{array}\right],$$

$A_3$ is now a second order supermatrix with two rows and two columns. We can simply write $A_3$ as

$$\begin{bmatrix} a_{11} & a_{12} \\ a_{21} & a_{22} \end{bmatrix}$$

where

$$a_{11} = \begin{bmatrix} 3 & 0 \\ 1 & 0 \\ 5 & -1 \\ 0 & 9 \end{bmatrix},$$

$$a_{12} = \begin{bmatrix} 1 & 1 & 2 & 0 \\ 0 & 3 & 5 & 2 \\ 6 & 7 & 8 & 4 \\ 1 & 2 & 0 & -1 \end{bmatrix},$$

$$a_{21} = \begin{bmatrix} 2 & 5 \\ 1 & 6 \end{bmatrix} \text{ and } a_{22} = \begin{bmatrix} 2 & 3 & 4 & 6 \\ 1 & 2 & 3 & 9 \end{bmatrix}.$$

The elements now are the submatrices defined as $a_{11}$, $a_{12}$, $a_{21}$ and $a_{22}$ and therefore $A_3$ is in terms of letters.

According to the methods we have illustrated a simple matrix can be partitioned to obtain a supermatrix in any way that happens to suit our purposes.

The natural order of a supermatrix is usually determined by the natural order of the corresponding simple matrix. Further more we are not usually concerned with natural order of the submatrices within a supermatrix.

Now we proceed on to recall the notion of symmetric partition, for more information about these concepts please refer [17]. By a symmetric partitioning of a matrix we mean that the rows and columns are partitioned in exactly the same way. If the matrix is partitioned between the first and second column and between the third and fourth column, then to be symmetrically partitioning, it must also be partitioned between the first and second rows and third and fourth rows. According to this rule of symmetric partitioning only square simple matrix can be symmetrically partitioned. We give an example of a symmetrically partitioned matrix $a_s$,

***Example 1.1.3:*** Let

$$a_s = \left[\begin{array}{c|cc|c} 2 & 3 & 4 & 1 \\ \hline 5 & 6 & 9 & 2 \\ 0 & 6 & 1 & 9 \\ \hline 5 & 1 & 1 & 5 \end{array}\right].$$

Here we see that the matrix has been partitioned between columns one and two and three and four. It has also been partitioned between rows one and two and rows three and four.

Now we just recall from [17] the method of symmetric partitioning of a symmetric simple matrix.

***Example 1.1.4:*** Let us take a fourth order symmetric matrix and partition it between the second and third rows and also between the second and third columns.

$$a = \left[\begin{array}{cc|cc} 4 & 3 & 2 & 7 \\ 3 & 6 & 1 & 4 \\ \hline 2 & 1 & 5 & 2 \\ 7 & 4 & 2 & 7 \end{array}\right].$$

We can represent this matrix as a supermatrix with letter elements.

$$a_{11} = \begin{bmatrix} 4 & 3 \\ 3 & 6 \end{bmatrix}, \; a_{12} = \begin{bmatrix} 2 & 7 \\ 1 & 4 \end{bmatrix}$$

$$a_{21} = \begin{bmatrix} 2 & 1 \\ 7 & 4 \end{bmatrix} \text{ and } a_{22} = \begin{bmatrix} 5 & 2 \\ 2 & 7 \end{bmatrix},$$

so that

$$a = \begin{bmatrix} a_{11} & a_{12} \\ a_{21} & a_{22} \end{bmatrix}.$$

The diagonal elements of the supermatrix a are $a_{11}$ and $a_{22}$. We also observe the matrices $a_{11}$ and $a_{22}$ are also symmetric matrices.

The non diagonal elements of this supermatrix a are the matrices $a_{12}$ and $a_{21}$. Clearly $a_{21}$ is the transpose of $a_{12}$.

The simple rule about the matrix element of a symmetrically partitioned symmetric simple matrix are (1) The diagonal submatrices of the supermatrix are all symmetric matrices. (2) The matrix elements below the diagonal are the transposes of the corresponding elements above the diagonal.

The forth order supermatrix obtained from a symmetric partitioning of a symmetric simple matrix a is as follows.

$$a = \begin{bmatrix} a_{11} & a_{12} & a_{13} & a_{14} \\ a'_{12} & a_{22} & a_{23} & a_{24} \\ a'_{13} & a'_{23} & a_{33} & a_{34} \\ a'_{14} & a'_{24} & a'_{34} & a_{44} \end{bmatrix}.$$

How to express that a symmetric matrix has been symmetrically partitioned (i) $a_{11}$ and $a^t_{11}$ are equal. (ii) $a^t_{ij}$ $(i \neq j)$; $a^t_{ij} = a_{ji}$ and $a^t_{ji} = a_{ij}$. Thus the general expression for a symmetrically partitioned symmetric matrix;

$$a = \begin{bmatrix} a_{11} & a_{12} & \ldots & a_{1n} \\ a'_{12} & a_{22} & \ldots & a_{2n} \\ \vdots & \vdots & & \vdots \\ a'_{1n} & a'_{2n} & \ldots & a_{nn} \end{bmatrix}.$$

If we want to indicate a symmetrically partitioned simple diagonal matrix we would write

$$D = \begin{bmatrix} D_1 & 0 & \ldots & 0 \\ 0' & D_2 & \ldots & 0 \\ & & & \\ 0' & 0' & \ldots & D_n \end{bmatrix}$$

0' only represents the order is reversed or transformed. We denote $a_{ij}^{t} = a'_{ij}$ just the ' means the transpose.

D will be referred to as the super diagonal matrix. The identity matrix

$$I = \begin{bmatrix} I_s & 0 & 0 \\ 0 & I_t & 0 \\ 0 & 0 & I_r \end{bmatrix}$$

s, t and r denote the number of rows and columns of the first second and third identity matrices respectively (zeros denote matrices with zero as all entries).

***Example 1.1.5:*** We just illustrate a general super diagonal matrix d;

$$d = \left[\begin{array}{ccc|cc} 3 & 1 & 2 & 0 & 0 \\ 5 & 6 & 0 & 0 & 0 \\ \hline 0 & 0 & 0 & 2 & 5 \\ 0 & 0 & 0 & -1 & 3 \\ 0 & 0 & 0 & 9 & 10 \end{array}\right]$$

$$\text{i.e., } d = \begin{bmatrix} m_1 & 0 \\ 0 & m_2 \end{bmatrix}.$$

An example of a super diagonal matrix with vector elements is given, which can be useful in experimental designs.

***Example 1.1.6:*** Let

$$\left[\begin{array}{c|c|c|c} 1 & 0 & 0 & 0 \\ 1 & 0 & 0 & 0 \\ 1 & 0 & 0 & 0 \\ \hline 0 & 1 & 0 & 0 \\ 0 & 1 & 0 & 0 \\ \hline 0 & 0 & 1 & 0 \\ 0 & 0 & 1 & 0 \\ 0 & 0 & 1 & 0 \\ \hline 0 & 0 & 1 & 0 \\ \hline 0 & 0 & 0 & 1 \\ 0 & 0 & 0 & 1 \\ 0 & 0 & 0 & 1 \\ 0 & 0 & 0 & 1 \end{array}\right].$$

Here the diagonal elements are only column unit vectors. In case of supermatrix [17] has defined the notion of partial triangular matrix as a supermatrix.

***Example 1.1.7:*** Let

$$u = \left[\begin{array}{ccc|cc} 2 & 1 & 1 & 3 & 2 \\ 0 & 5 & 2 & 1 & 1 \\ 0 & 0 & 1 & 0 & 2 \end{array}\right]$$

u is a partial upper triangular supermatrix.

***Example 1.1.8:*** Let

$$L = \left[\begin{array}{ccccc} 5 & 0 & 0 & 0 & 0 \\ 7 & 2 & 0 & 0 & 0 \\ 1 & 2 & 3 & 0 & 0 \\ 4 & 5 & 6 & 7 & 0 \\ 1 & 2 & 5 & 2 & 6 \\ \hline 1 & 2 & 3 & 4 & 5 \\ 0 & 1 & 0 & 1 & 0 \end{array}\right];$$

L is partial upper triangular matrix partitioned as a supermatrix.

Thus $T = \left[\begin{array}{c} T \\ \hline a' \end{array}\right]$ where T is the lower triangular submatrix, with

$$T = \begin{bmatrix} 5 & 0 & 0 & 0 & 0 \\ 7 & 2 & 0 & 0 & 0 \\ 1 & 2 & 3 & 0 & 0 \\ 4 & 5 & 6 & 7 & 0 \\ 1 & 2 & 5 & 2 & 6 \end{bmatrix} \text{ and } a' = \begin{bmatrix} 1 & 2 & 3 & 4 & 5 \\ 0 & 1 & 0 & 1 & 0 \end{bmatrix}.$$

We proceed on to define the notion of supervectors i.e., Type I column supervector. A simple vector is a vector each of whose elements is a scalar. It is nice to see the number of different types of supervectors given by [17].

***Example 1.1.9:*** Let

$$v = \left[\begin{array}{c} 1 \\ 3 \\ 4 \\ \hline 5 \\ 7 \end{array}\right].$$

This is a type I i.e., type one column supervector.

$$v = \begin{bmatrix} v_1 \\ v_2 \\ \vdots \\ v_n \end{bmatrix}$$

where each $v_i$ is a column subvectors of the column vector v.

Type I row supervector is given by the following example.

***Example 1.1.10:*** $v^1$ = [2 3 1 | 5 7 8 4] is a type I row supervector. i.e., v' = [$v'_1$, $v'_2$, …, $v'_n$] where each $v'_i$ is a row subvector; $1 \leq i \leq n$.

Next we recall the definition of type II supervectors.

Type II column supervectors.

**DEFINITION 1.1.1:** *Let*

$$a = \begin{bmatrix} a_{11} & a_{12} & \ldots & a_{1m} \\ a_{21} & a_{22} & \ldots & a_{2m} \\ \ldots & \ldots & \ldots & \ldots \\ a_{n1} & a_{n2} & \ldots & a_{nm} \end{bmatrix}$$

$$\begin{aligned} a_1^1 &= [a_{11} \ldots a_{1m}] \\ a_2^1 &= [a_{21} \ldots a_{2m}] \\ &\ldots \\ a_n^1 &= [a_{n1} \ldots a_{nm}] \end{aligned}$$

*i.e.,* $$a = \begin{bmatrix} a_1^1 \\ a_2^1 \\ \vdots \\ a_n^1 \end{bmatrix}_m$$

*is defined to be the type II column supervector.*
*Similarly if*

$$a^1 = \begin{bmatrix} a_{11} \\ a_{21} \\ \vdots \\ a_{n1} \end{bmatrix}, \; a^2 = \begin{bmatrix} a_{12} \\ a_{22} \\ \vdots \\ a_{n2} \end{bmatrix}, \ldots, a^m = \begin{bmatrix} a_{1m} \\ a_{2m} \\ \vdots \\ a_{nm} \end{bmatrix}.$$

*Hence now* $a = [a^1 \; a^2 \ldots a^m]_n$ *is defined to be the type II row supervector.*

*Clearly*

$$a = \begin{bmatrix} a_1^1 \\ a_2^1 \\ \vdots \\ a_n^1 \end{bmatrix}_m = [a^1\ a^2\ \ldots\ a^m]_n$$

*the equality of supermatrices.*

***Example 1.1.11:*** Let

$$A = \begin{bmatrix} 3 & 6 & 0 & 4 & 5 \\ 2 & 1 & 6 & 3 & 0 \\ 1 & 1 & 1 & 2 & 1 \\ 0 & 1 & 0 & 1 & 0 \\ 2 & 0 & 1 & 2 & 1 \end{bmatrix}$$

be a simple matrix. Let a and b the supermatrix made from A.

$$a = \left[\begin{array}{ccc|cc} 3 & 6 & 0 & 4 & 5 \\ 2 & 1 & 6 & 3 & 0 \\ 1 & 1 & 1 & 2 & 1 \\ \hline 0 & 1 & 0 & 1 & 0 \\ 2 & 0 & 1 & 2 & 1 \end{array}\right]$$

where

$$a_{11} = \begin{bmatrix} 3 & 6 & 0 \\ 2 & 1 & 6 \\ 1 & 1 & 1 \end{bmatrix},\ a_{12} = \begin{bmatrix} 4 & 5 \\ 3 & 0 \\ 2 & 1 \end{bmatrix},$$

$$a_{21} = \begin{bmatrix} 0 & 1 & 0 \\ 2 & 0 & 1 \end{bmatrix} \text{ and } a_{22} = \begin{bmatrix} 1 & 0 \\ 2 & 1 \end{bmatrix}.$$

i.e.,
$$a = \begin{bmatrix} a_{11} & a_{12} \\ a_{21} & a_{22} \end{bmatrix}.$$

$$b = \left[\begin{array}{cccc|c} 3 & 6 & 0 & 4 & 5 \\ 2 & 1 & 6 & 3 & 0 \\ 1 & 1 & 1 & 2 & 1 \\ 0 & 1 & 0 & 1 & 0 \\ \hline 2 & 0 & 1 & 2 & 1 \end{array}\right] = \begin{bmatrix} b_{11} & b_{12} \\ b_{21} & b_{22} \end{bmatrix}$$

where

$$b_{11} = \begin{bmatrix} 3 & 6 & 0 & 4 \\ 2 & 1 & 6 & 3 \\ 1 & 1 & 1 & 2 \\ 0 & 1 & 0 & 1 \end{bmatrix}, \; b_{12} = \begin{bmatrix} 5 \\ 0 \\ 1 \\ 0 \end{bmatrix},$$

$b_{21} = [2\ 0\ 1\ 2\ ]$ and $b_{22} = [1]$.

$$a = \left[\begin{array}{ccc|cc} 3 & 6 & 0 & 4 & 5 \\ 2 & 1 & 6 & 3 & 0 \\ 1 & 1 & 1 & 2 & 1 \\ \hline 0 & 1 & 0 & 1 & 0 \\ 2 & 0 & 1 & 2 & 1 \end{array}\right]$$

and

$$b = \left[\begin{array}{cccc|c} 3 & 6 & 0 & 4 & 5 \\ 2 & 1 & 6 & 3 & 0 \\ 1 & 1 & 1 & 2 & 1 \\ 0 & 1 & 0 & 1 & 0 \\ \hline 2 & 0 & 1 & 2 & 1 \end{array}\right].$$

We see that the corresponding scalar elements for matrix a and matrix b are identical. Thus two supermatrices are equal if and only if their corresponding simple forms are equal.

Now we give examples of type III supervector for more refer [17].

***Example 1.1.12:***

$$a = \left[\begin{array}{ccc|cc} 3 & 2 & 1 & 7 & 8 \\ 0 & 2 & 1 & 6 & 9 \\ 0 & 0 & 5 & 1 & 2 \end{array}\right] = [T' \mid a']$$

and

$$b = \left[\begin{array}{ccc} 2 & 0 & 0 \\ 9 & 4 & 0 \\ 8 & 3 & 6 \\ \hline 5 & 2 & 9 \\ 4 & 7 & 3 \end{array}\right] = \left[\begin{array}{c} T \\ \hline b' \end{array}\right]$$

are type III supervectors.

One interesting and common example of a type III supervector is a prediction data matrix having both predictor and criterion attributes.

The next interesting notion about supermatrix is its transpose. First we illustrate this by an example before we give the general case.

***Example 1.1.13:*** Let

$$a = \left[\begin{array}{ccc|cc} 2 & 1 & 3 & 5 & 6 \\ 0 & 2 & 0 & 1 & 1 \\ 1 & 1 & 1 & 0 & 2 \\ \hline 2 & 2 & 0 & 1 & 1 \\ 5 & 6 & 1 & 0 & 1 \\ \hline 2 & 0 & 0 & 0 & 4 \\ 1 & 0 & 1 & 1 & 5 \end{array}\right]$$

$$= \begin{bmatrix} a_{11} & a_{12} \\ a_{21} & a_{22} \\ a_{31} & a_{32} \end{bmatrix}$$

where

$$a_{11} = \begin{bmatrix} 2 & 1 & 3 \\ 0 & 2 & 0 \\ 1 & 1 & 1 \end{bmatrix}, \; a_{12} = \begin{bmatrix} 5 & 6 \\ 1 & 1 \\ 0 & 2 \end{bmatrix},$$

$$a_{21} = \begin{bmatrix} 2 & 2 & 0 \\ 5 & 6 & 1 \end{bmatrix}, \; a_{22} = \begin{bmatrix} 1 & 1 \\ 0 & 1 \end{bmatrix},$$

$$a_{31} = \begin{bmatrix} 2 & 0 & 0 \\ 1 & 0 & 1 \end{bmatrix} \text{ and } a_{32} = \begin{bmatrix} 0 & 4 \\ 1 & 5 \end{bmatrix}.$$

The transpose of a

$$a^t = a' = \left[\begin{array}{ccc|cc|cc} 2 & 0 & 1 & 2 & 5 & 2 & 1 \\ 1 & 2 & 1 & 2 & 6 & 0 & 0 \\ 3 & 0 & 1 & 0 & 1 & 0 & 1 \\ \hline 5 & 1 & 0 & 1 & 0 & 0 & 1 \\ 6 & 1 & 2 & 1 & 1 & 4 & 5 \end{array}\right].$$

Let us consider the transposes of $a_{11}$, $a_{12}$, $a_{21}$, $a_{22}$, $a_{31}$ and $a_{32}$.

$$a'_{11} = a^t_{11} = \begin{bmatrix} 2 & 0 & 1 \\ 1 & 2 & 1 \\ 3 & 0 & 1 \end{bmatrix}$$

$$a'_{12} = a^t_{12} = \begin{bmatrix} 5 & 1 & 0 \\ 6 & 1 & 2 \end{bmatrix}$$

$$a'_{21} = a^t_{21} = \begin{bmatrix} 2 & 5 \\ 2 & 6 \\ 0 & 1 \end{bmatrix}$$

$$a'_{31} = a_{31}^t = \begin{bmatrix} 2 & 1 \\ 0 & 0 \\ 0 & 1 \end{bmatrix}$$

$$a'_{22} = a_{22}^t = \begin{bmatrix} 1 & 0 \\ 1 & 1 \end{bmatrix}$$

$$a'_{32} = a_{32}^t = \begin{bmatrix} 0 & 1 \\ 4 & 5 \end{bmatrix}.$$

$$a' = \begin{bmatrix} a'_{11} & a'_{21} & a'_{31} \\ a'_{12} & a'_{22} & a'_{32} \end{bmatrix}.$$

Now we describe the general case. Let

$$a = \begin{bmatrix} a_{11} & a_{12} & \cdots & a_{1m} \\ a_{21} & a_{22} & \cdots & a_{2m} \\ \vdots & \vdots & & \vdots \\ a_{n1} & a_{n2} & \cdots & a_{nm} \end{bmatrix}$$

be a n × m supermatrix. The transpose of the supermatrix a denoted by

$$a' = \begin{bmatrix} a'_{11} & a'_{21} & \cdots & a'_{n1} \\ a'_{12} & a'_{22} & \cdots & a'_{n2} \\ \vdots & \vdots & & \vdots \\ a'_{1m} & a'_{2m} & \cdots & a'_{nm} \end{bmatrix}.$$

a' is a m by n supermatrix obtained by taking the transpose of each element i.e., the submatrices of a.

Now we will find the transpose of a symmetrically partitioned symmetric simple matrix. Let a be the symmetrically partitioned symmetric simple matrix.

Let a be a m × m symmetric supermatrix i.e.,

$$a = \begin{bmatrix} a_{11} & a_{21} & \cdots & a_{m1} \\ a_{12} & a_{22} & \cdots & a_{m2} \\ \vdots & \vdots & & \vdots \\ a_{1m} & a_{2m} & \cdots & a_{mm} \end{bmatrix}$$

the transpose of the supermatrix is given by a'

$$a' = \begin{bmatrix} a'_{11} & (a'_{12})' & \cdots & (a'_{1m})' \\ a'_{12} & a'_{22} & \cdots & (a'_{2m})' \\ \vdots & \vdots & & \vdots \\ a'_{1m} & a'_{2m} & \cdots & a'_{mm} \end{bmatrix}$$

The diagonal matrix $a_{11}$ are symmetric matrices so are unaltered by transposition. Hence

$$a'_{11} = a_{11},\ a'_{22} = a_{22},\ \ldots,\ a'_{mm} = a_{mm}.$$

Recall also the transpose of a transpose is the original matrix. Therefore

$$(a'_{12})' = a_{12},\ (a'_{13})' = a_{13},\ \ldots,\ (a'_{ij})' = a_{ij}.$$

Thus the transpose of supermatrix constructed by symmetrically partitioned symmetric simple matrix a of a' is given by

$$a' = \begin{bmatrix} a_{11} & a_{12} & \cdots & a_{1m} \\ a'_{21} & a_{22} & \cdots & a_{2m} \\ \vdots & \vdots & & \vdots \\ a'_{1m} & a'_{2m} & \cdots & a_{mm} \end{bmatrix}.$$

Thus a = a'.
Similarly transpose of a symmetrically partitioned diagonal matrix is simply the original diagonal supermatrix itself;

i.e., if

$$D = \begin{bmatrix} d_1 & & & \\ & d_2 & & \\ & & \ddots & \\ & & & d_n \end{bmatrix}$$

$$D' = \begin{bmatrix} d_1' & & & \\ & d_2' & & \\ & & \ddots & \\ & & & d_n' \end{bmatrix}$$

$d'_1 = d_1$, $d'_2 = d_2$ etc. Thus D = D'.

Now we see the transpose of a type I supervector.

***Example 1.1.14:*** Let

$$V = \left[\begin{array}{c} 3 \\ 1 \\ 2 \\ \hline 4 \\ 5 \\ 7 \\ \hline 5 \\ 1 \end{array}\right]$$

The transpose of V denoted by V' or $V^t$ is

$$V' = [3\ 1\ 2 \mid 4\ 5\ 7 \mid 5\ 1].$$

If

$$V = \begin{bmatrix} v_1 \\ v_2 \\ v_3 \end{bmatrix}$$

where

$$v_1 = \begin{bmatrix} 3 \\ 1 \\ 2 \end{bmatrix}, v_2 = \begin{bmatrix} 4 \\ 5 \\ 7 \end{bmatrix} \text{ and } v_3 = \begin{bmatrix} 5 \\ 1 \end{bmatrix}$$

$$V' = [v'_1 \; v'_2 \; v'_3].$$

Thus if

$$V = \begin{bmatrix} v_1 \\ v_2 \\ \vdots \\ v_n \end{bmatrix}$$

then

$$V' = [v'_1 \; v'_2 \; \ldots \; v'_n].$$

***Example 1.1.15:*** Let

$$t = \left[\begin{array}{cccc|cc} 3 & 0 & 1 & 1 & 5 & 2 \\ 4 & 2 & 0 & 1 & 3 & 5 \\ 1 & 0 & 1 & 0 & 1 & 6 \end{array}\right]$$

= [T | a ]. The transpose of t

$$\text{i.e., } t' = \left[\begin{array}{ccc} 3 & 4 & 1 \\ 0 & 2 & 0 \\ 1 & 0 & 1 \\ 1 & 1 & 0 \\ \hline 5 & 3 & 1 \\ 2 & 5 & 6 \end{array}\right] = \left[\begin{array}{c} T' \\ \hline a' \end{array}\right].$$

The addition of supermatrices may not be always be defined.

***Example 1.1.16:*** For instance let

$$a = \begin{bmatrix} a_{11} & a_{12} \\ a_{21} & a_{22} \end{bmatrix}$$

and

$$b = \begin{bmatrix} b_{11} & b_{12} \\ b_{21} & b_{22} \end{bmatrix}$$

where

$$a_{11} = \begin{bmatrix} 3 & 0 \\ 1 & 2 \end{bmatrix}, \quad a_{12} = \begin{bmatrix} 1 \\ 7 \end{bmatrix}$$

$$a_{21} = [4 \;\; 3], \quad a_{22} = [6].$$

$$b_{11} = [2], \quad b_{12} = [1 \;\; 3]$$

$$b_{21} = \begin{bmatrix} 5 \\ 2 \end{bmatrix} \text{ and } b_{22} = \begin{bmatrix} 4 & 1 \\ 0 & 2 \end{bmatrix}.$$

It is clear both a and b are second order square supermatrices but here we cannot add together the corresponding matrix elements of a and b because the submatrices do not have the same order.

## 1.2 Super Vector Spaces and their properties

This section for the first time introduces systematically the notion of super vector spaces and analyze the special properties associated with them. Throughout this book F will denote a field in general. R the field of reals, Q the field of rationals and $Z_p$ the field of integers modulo p, p a prime. These fields all are real; whereas C will denote the field of complex numbers.

We recall $X = (x_1\ x_2 \mid x_3\ x_4\ x_5 \mid x_6)$ is a super row vector where $x_i \in F$; F a field; $1 \leq i \leq 6$. Suppose $Y = (y_1\ y_2 \mid y_3\ y_4\ y_5 \mid y_6)$ with $y_i \in F$; $1 \leq i \leq 6$ we say X and Y are super vectors of the same type. Further if $Z = (z_1\ z_2\ z_3\ z_4 \mid z_5\ z_6)$ $z_i \in F$; $1 \leq i \leq 6$ then we don't say Z to be a super vector of same type as X or Y. Further same type super vectors X and Y over the same field are equal if and only if $x_i = y_i$ for i = 1, 2, …, 6. Super vectors of same type can be added the resultant is once again a super vector of the same type. The first important result about the super vectors of same type is the following theorem.

**THEOREM 1.2.1**: *This collection of all super vectors* $S = \{X = (x_1\ x_2 \ldots x_r \mid x_{r+1} \ldots x_i \mid x_{i+1} \ldots x_{t+1} \mid x_{t+2} \ldots x_n) \mid x_i \in F\}$*; F a field,* $1 \leq i \leq n$. $\{1 < 2 < \ldots < r < r + 1 < \ldots < i < i+1 < \ldots < t+1 < \ldots < n\}$ *of this type is an abelian group under component wise addition.*

*Proof:* Let

$$X = (x_1\ x_2 \ldots x_r \mid x_{r+1} \ldots x_i \mid x_{i+1} \ldots x_{t+1} \mid x_{t+2} \ldots x_n)$$

and

$$Y = (y_1\ y_2 \ldots y_r \mid y_{r+1} \ldots y_i \mid y_{i+1} \ldots y_{t+1} \mid y_{t+2} \ldots y_n) \in S.$$

$$X + Y = \{(x_1 + y_1\ x_2 + y_2 \ldots x_r + y_r \mid x_{r+1} + y_{r+1} \ldots x_i + y_i \mid x_{i+1} + y_{i+1} \ldots x_{t+1} + y_{t+1} \mid x_{t+2} + y_{t+2} \ldots x_n + y_n)\}$$

is again a super vector of the same type and is in S as $x_i + y_i \in F$; $1 \leq i \leq n$.

Clearly $(0\ 0 \ldots 0 \mid 0 \ldots 0 \mid 0 \ldots 0 \mid 0 \ldots 0) \in S$ as $0 \in F$.

Now if

$$X = (x_1\ x_2 \ldots x_r \mid x_{r+1} \ldots x_i \mid x_{i+1} \ldots x_{t+1} \mid x_{t+2} \ldots x_n) \in S$$

then

$$-X = (-x_1\ -x_2 \ldots -x_r \mid -x_{r+1} \ldots -x_i \mid -x_{i+1} \ldots -x_{t+1} \mid -x_{t+2} \ldots -x_n) \in S$$

with

$$X + (-X) = (-X) + X = (0\ 0 \ldots 0 \mid 0 \ldots 0 \mid 0 \ldots 0 \mid 0 \ldots 0)$$

Also X + Y = Y + X.

Hence S is an abelian group under addition.

We first illustrate this situation by some simple examples.

***Example 1.2.1:*** Let Q be the field of rationals. Let S = $\{(x_1\ x_2\ x_3 \mid x_4\ x_5) \mid x_1, \ldots, x_5 \in Q\}$. Clearly S is an abelian group under component wise addition of super vectors of S. Take any two super vectors say X = (3 2 1 | –5 3) and Y = (0 2 4 | 1 –2) in S.

We see X + Y = (3 4 5 | – 4 1) and X + Y $\in$ S. Also (0 0 0 0 | 0 0) acts as the super row zero vector which can also be called as super identity or super row zero vector. Further if X = (5 7 – 3| 0 –1) then –X = (–5 –7 3| 0 1) is the inverse of X and we see X + (–X) = (0 0 0 | 0 0). Thus S is an abelian group under componentwise addition of super vectors.

If X' = (3 1 1 4 | 5 6 2) is any super vector. Clearly X' $\notin$ S, given in example 1.2.1 as X' is not the same type of super vector, as X' is different from X = $(x_1\ x_2\ x_3 \mid x_4\ x_5)$.

***Example 1.2.2***: Consider the set S = $\{(x_1 \mid x_2\ x_3 \mid x_4\ x_5) \mid x_i \in Q; 1 \leq i \leq 5\}$. S is an additive abelian group. We call such groups as matrix partition groups.

Every matrix partition group is a group. But every group in general is not a partition group we also call the matrix partition group or super matrix group or super special group.

***Example 1.2.3***: Let S = $\{(x_1\ x_2\ x_3) \mid x_i \in Q; 1 \leq i \leq 3\}$. S is a group under component wise addition of row vectors but S is not a matrix partition group only a group.

***Example 1.2.4:*** Let

$$P = \left\{ \left( \begin{array}{c|cc} x_1 & x_5 & x_6 \\ x_2 & x_2 & x_8 \\ \hline x_3 & x_9 & x_{10} \\ x_4 & x_{11} & x_{12} \end{array} \right) \mid x_i \in Q; i = 1, 2, \ldots, 12 \right\}.$$

Clearly P is a group under matrix addition, which we choose to call as partition matrix addition. P is a partition abelian group or we call them as super groups. Now we proceed on to define super vector space.

**DEFINITION 1.2.1:** *Let V be an abelian super group i.e. an abelian partitioned group under addition, F be a field. We call V a super vector space over F if the following conditions are satisfied*

- *(i)* *for all $v \in V$ and $c \in F$, vc and cv are in V. Further vc = cv we write first the field element as they are termed as scalars over which the vector space is defined.*
- *(ii)* *for all $v_1, v_2 \in V$ and for all $c \in F$ we have $c(v_1 + v_2) = cv_1 + cv_2$.*
- *(iii)* *also $(v_1 + v_2) c = v_1c + v_2c$.*
- *(iv)* *for a, $b \in F$ and $v_1 \in V$ we have $(a + b) v_1 = av_1 + bv_1$ also $v_1 (a + b) = v_1a + v_1b$.*
- *(v)* *for every $v \in V$ and $1 \in F$, $1.v = v.1 = v$*
- *(vi)* *$(c_1 c_2) v = c_1 (c_2v)$ for all $v \in V$ and $c_1, c_2 \in F$.*

*The elements of V are called "super vectors" and elements of F are called "scalars".*

We shall illustrate this by the following examples.

***Example 1.2.5:*** Let $V = \{(x_1\ x_2\ x_3 \mid x_4) \mid x_i \in R; 1 \leq i \leq 4$, the field of reals$\}$. V is an abelian super group under addition. Q be the field of rationals V is a super vector space over Q. For if $10 \in Q$ and $v = (\sqrt{2}\ \ 5\ 1 \mid 3) \in V$; $10v = (10\sqrt{2}\ \ 50\ 10 \mid 30) \in V$.

***Example 1.2.6:*** Let $V = \{(x_1\ x_2\ x_3 \mid x_4) \mid x_i \in R$, the field of reals $1 \leq i \leq 4\}$. V is a super vector space over R. We see there is difference between the super vector spaces mentioned in the example 1.2.5 and here.

We can also have other examples.

***Example 1.2.7:*** Let

$$V = \left\{ \begin{pmatrix} y_1 \\ y_2 \\ y_3 \end{pmatrix} \middle| \; y_1, y_2, y_3 \in Q \right\}.$$

Clearly V is a super group under addition and is an abelian super group. Take Q the field of rationals. V is a super vector space over Q. Take 5 $\in$ Q,

$$v = \begin{pmatrix} -1 \\ 2 \\ 4 \end{pmatrix} \text{ in V.}$$

$$5v = \begin{pmatrix} -5 \\ 10 \\ 20 \end{pmatrix} \in V.$$

As in case of vector space which depends on the field over which it has to be defined so also are super vector space.

The following example makes this more explicit.

***Example 1.2.8:*** Let

$$V = \left\{ \begin{pmatrix} y_1 \\ y_2 \\ y_3 \end{pmatrix} \middle| \; y_1, y_2, y_3 \in Q \text{; the field of rational} \right\};$$

V is an abelian super group under addition. V is a super vector space over Q; but V is not a super vector space over the field of reals R. For $\sqrt{2} \in$ R;

$$v = \begin{pmatrix} 5 \\ 1 \\ 3 \end{pmatrix} \in V.$$

$$\sqrt{2}\, v = \sqrt{2} \begin{pmatrix} 5 \\ 1 \\ 3 \end{pmatrix} = \begin{pmatrix} 5\sqrt{2} \\ \sqrt{2} \\ 3\sqrt{2} \end{pmatrix} \notin V$$

as $5\sqrt{2}$, $\sqrt{2}$ and $3\sqrt{2} \notin Q$. So V is not a super vector space over R.

We can also have V as a super n-tuple space.

***Example 1.2.9:*** Let $V = \{F^{n_1} \mid \ldots \mid F^{n_t}\}$ where F is a field. V is a super abelian group under addition so V is a super vector space over F.

***Example 1.2.10:*** Let $V = \{(Q^3 \mid Q^3 \mid Q^2) = \{(x_1\ x_2\ x_3 \mid y_1\ y_2\ y_3 \mid z_1\ z_2) \mid x_i, y_k, z_j \in Q;\ 1 \le i \le 3;\ 1 \le k \le 3;\ 1 \le j \le 2\}$. V is a super vector space over Q. Clearly V is not a super vector space over the field of reals R.

Now as we have matrices to be vector spaces likewise we have super matrices are super vector spaces.

***Example 1.2.11:*** Let

$$A = \left\{ \left( \begin{array}{cc|ccc} x_1 & x_2 & x_9 & x_{10} & x_{11} \\ x_3 & x_4 & x_{12} & x_{13} & x_{14} \\ x_5 & x_6 & x_{15} & x_{16} & x_{17} \\ x_7 & x_8 & x_{18} & x_{19} & x_{20} \end{array} \right) \middle|\ x_i \in Q; 1 \le i \le 20 \right\}$$

be the collection of super matrices with entries from Q. A is a super vector space over Q.

***Example 1.2.12:*** Let

$$V = \left\{ \left( \begin{array}{cc|ccc} x_1 & x_2 & x_5 & x_6 & x_7 \\ x_3 & x_4 & x_8 & x_9 & x_{10} \end{array} \right) \middle| \; x_i \in R; 1 \le i \le 10 \right\}.$$

V is a super vector space over Q.

***Example 1.2.13:*** Let

$$V = \left\{ \left( \begin{array}{ccc|cc} x_1 & x_2 & x_5 & x_6 & x_7 \\ x_3 & x_4 & x_8 & x_9 & x_{10} \end{array} \right) \middle| \; x_i \in R; 1 \le i \le 10 \right\}$$

V is a super vector space over R. V is also a super vector space over Q. However soon we shall be proving that these two super vector spaces are different.

***Example 1.2.14:*** Let

$$A = \left\{ \left( \begin{array}{cc|cc} a_1 & a_2 & a_5 & a_6 \\ a_3 & a_4 & a_7 & a_8 \\ \hline a_9 & a_{10} & a_{13} & a_{14} \\ a_{11} & a_{12} & a_{15} & a_{16} \end{array} \right) \middle| \; a_i \in Q; \; 1 \le i \le 16 \right\}.$$

V is a super vector space over Q. However V is not a super vector space over R.

We call the elements of the super vector space V to be super vectors and elements of F to be just scalars.

**DEFINITION 1.2.2:** *Let V be a super vector space over the field F. A super vector $\beta$ in V is said be a linear combination of super vectors $\alpha_1, \ldots, \alpha_n$ in V provided there exists scalars $c_1, \ldots, c_n$ in F such that* $\beta = c_1\alpha_1 + \ldots + c_n\alpha_n = \sum_{i=1}^{n} c_i\alpha_i$.

We illustrate this by the following example.

***Example 1.2.15:*** Let V = {($a_1$ $a_2$ | $a_3$ $a_4$ $a_5$ | $a_6$)| $a_i \in$ Q; $1 \leq i \leq 6$}. V is a super vector space over Q. Consider β = (7 5 | 0 2 8 | 9) a super vector in V. Let $\alpha_1$ = (1 1 | 2 0 1 | –1), $\alpha_2$ = (5 –3 | 1 2 5 | 5) and $\alpha_3$ = (0 7 | 3 1 2 | 8) be 3 super vectors in V. We can find a, b, c in Q such that $a\alpha_1 + b\alpha_2 + c\alpha_3 = \beta$.

***Example 1.2.16:*** Let

$$A = \left\{ \left. \left( \begin{array}{c|c} a & c \\ b & d \end{array} \right) \right| a, b, c, d \in Q \right\}.$$

A is a super vector space over Q.
Let

$$\beta = \left( \begin{array}{c|c} 12 & 5 \\ 8 & -1 \end{array} \right) \in A.$$

We have for

$$\left( \begin{array}{c|c} 2 & 1 \\ 1 & -1 \end{array} \right), \left( \begin{array}{c|c} 4 & 1 \\ 4 & 3 \end{array} \right) \in A$$

such that for scalars 4, 1 ∈ Q we have

$$4 \left( \begin{array}{c|c} 2 & 1 \\ 1 & -1 \end{array} \right) + 1 \left( \begin{array}{c|c} 4 & 1 \\ 4 & 3 \end{array} \right)$$

$$= \left( \begin{array}{c|c} 8 & 4 \\ 4 & -4 \end{array} \right) + \left( \begin{array}{c|c} 4 & 1 \\ 4 & 3 \end{array} \right)$$

$$= \left( \begin{array}{c|c} 12 & 5 \\ 8 & -1 \end{array} \right) = \beta.$$

Now we proceed onto define the notion of super subspace of a super vector space V over the field F.

**DEFINITION 1.2.3:** *Let V be a super vector space over the field F. A proper subset W of V is said to be super subspace of V if W itself is a super vector space over F with the operations of super vector addition and scalar multiplication on V.*

**THEOREM 1.2.2:** *A non-empty subset W of V, V a super vector space over the field F is a super subspace of V if and only if for each pair of super vectors α, β in W and each scalar c in F the super vector cα + β is again in W.*

*Proof:* Suppose that W is a non empty subset of V; where V is a super vector space over the field F. Suppose that cα + β belongs to W for all super vectors α, β in W and for all scalars c in F. Since W is non-empty there is a super vector p in W and hence (–1)p + p = 0 is in W. Thus if α is any super vector in W and c any scalar, the super vector cα = cα + 0 is in W. In particular, (–1) α = –α is in W. Finally if α and β are in W then α + β = 1. α + β is in W. Thus W is a super subspace of V.

Conversely if W is a super subspace of V, α and β are in W and c is a scalar certainly cα + β is in W.

*Note:* If V is any super vector space; the subset consisting of the zero super vector alone is a super subspace of V called the zero super subspace of V.

**THEOREM 1.2.3:** *Let V be a super vector space over the field F. The intersection of any collection of super subspaces of V is a super subspace of V.*

*Proof:* Let $\{W_\alpha\}$ be the collection of super subspaces of V and let $W = \bigcap_\alpha W_\alpha$ be the intersection. Recall that W is defined as the set of all elements belonging to every $W_\alpha$ (For if $x \in W = \bigcap W_\alpha$ then x belongs to every $W_\alpha$). Since each $W_\alpha$ is a super subspace each contains the zero super vector. Thus the zero super vector is in the intersection W and W is non empty. Let α and β be super vectors in W and c be any scalar. By definition of W both α and β belong to each $W_\alpha$ and because each $W_\alpha$ is a super subspace, the super vector $c\alpha + \beta$ is in every $W_\alpha$. Thus $c\alpha + \beta$ is again in W. By the theorem just proved; W is a super subspace of V.

**DEFINITION 1.2.4**: *Let S be a set of super vectors in a super vector space V. The super subspace spanned by S is defined to be the intersection W of all super subspaces of V which contain S. When S is a finite set of super vectors, that is S = $\{\alpha_1, \ldots, \alpha_n\}$ we shall simply call W, the super subspace spanned by the super vectors $\{\alpha_1, \ldots, \alpha_n\}$.*

**THEOREM 1.2.4:** *The super subspace spanned by a non empty subset S of a super vector space V is the set of all linear combinations of super vectors in S.*

*Proof:* Given V is a super vector space over the field F. W be a super subspace of V spanned by S. Then each linear combination $\alpha = x_1\alpha_1 + \ldots + x_n\alpha_n$ of super vectors $\alpha_1, \ldots, \alpha_n$ in S is clearly in W. Thus W contains the set L of all linear combinations of super vectors in S. The set L, on the other hand, contains S and is non-empty. If $\alpha$, $\beta$ belong to L then $\alpha$ is a linear combination.

$$\alpha = x_1\alpha_1 + \ldots + x_m\alpha_m$$

of super vectors $\alpha_1, \ldots, \alpha_m$ in S and $\beta$ is a linear combination.

$$\beta = y_1\beta_1 + \ldots + y_m\beta_m$$

of super vectors $\beta_j$ in S; $1 \leq j \leq m$. For each scalar,

$$c\alpha + \beta = \sum_{i=1}^{m} (cx_i)\alpha_i + \sum_{j=1}^{m} y_j\beta_j$$

$x_i, y_i \in F$; $1 \leq i, j \leq m$.
Hence $c\alpha + \beta$ belongs to L. Thus L is a super subspace of V.

Now we have proved that L is a super subspace of V which contains S, and also that any subspace which contains S contains L. It follows that L is the intersection of all super subspaces containing S, i.e. that L is the super subspace spanned by the set S.

Now we proceed onto define the sum of subsets.

**DEFINITION 1.2.5**: *If $S_1, \ldots, S_K$ are subsets of a super vector space V, the set of all sums $\alpha_1 + \ldots + \alpha_K$ of super vectors $\alpha_i$ in $S_i$ is called the sum of the subsets $S_1, S_2, \ldots, S_K$ and is denoted by $S_1 + \ldots + S_K$ or by $\sum_{i=1}^{K} S_i$.*

*If $W_1, \ldots, W_K$ are super subspaces of the super vector space V, then the sum $W = W_1 + W_2 + \ldots + W_K$ is easily seen to be a super subspace of V which contains each of super subspace $W_i$. i.e. W is the super subspace spanned by the union of $W_1, W_2, \ldots, W_K$, $1 \leq i \leq K$.*

***Example 1.2.17:*** Let

$$A = \left\{ \left( \begin{array}{c|cccc} x_1 & x_2 & x_9 & x_{10} & x_{11} \\ x_3 & x_4 & x_{12} & x_{13} & x_{14} \\ \hline x_5 & x_6 & x_{15} & x_{16} & x_{17} \\ x_7 & x_8 & x_{18} & x_{19} & x_{20} \end{array} \right) \middle| \, x_i \in Q; 1 \leq i \leq 16 \right\}$$

be a super vector subspace of V over Q.
Let

$$W_1 = \left\{ \left( \begin{array}{c|ccc} x_1 & 0 & 0 & 0 \\ x_3 & 0 & 0 & 0 \\ \hline 0 & x_6 & x_{13} & x_{14} \\ 0 & x_8 & x_{15} & x_{16} \end{array} \right) \middle| \, x_1, x_3, x_6, x_8, x_{13}, x_{14}, x_8, x_{15}, x_{16} \in Q \right\}$$

$W_1$ is clearly a super subspace of V.

Let

$$W_2 = \left\{ \left( \begin{array}{c|ccc} 0 & 0 & 0 & 0 \\ 0 & 0 & 0 & 0 \\ \hline x_5 & 0 & 0 & 0 \\ x_7 & 0 & 0 & 0 \end{array} \right) \middle| \; x_5, x_6 \in Q \right\},$$

$W_2$ is a super subspace of V.
Take

$$W_3 = \left\{ \left( \begin{array}{c|ccc} 0 & x_2 & x_9 & x_{10} \\ 0 & x_4 & x_{11} & x_{12} \\ \hline 0 & 0 & 0 & 0 \\ 0 & 0 & 0 & 0 \end{array} \right) \middle| \; x_2, x_9, x_4, x_{10}, x_{11}, x_{12} \in Q \right\}$$

a proper super vector subspace of V.
Clearly $V = W_1 + W_2 + W_3$ i.e.,

$$\left( \begin{array}{c|ccc} x_1 & x_2 & x_9 & x_{10} \\ x_3 & x_4 & x_{11} & x_{12} \\ \hline x_5 & x_6 & x_{13} & x_{14} \\ x_7 & x_8 & x_{15} & x_{16} \end{array} \right) = \left( \begin{array}{c|ccc} x_1 & 0 & 0 & 0 \\ x_3 & 0 & 0 & 0 \\ \hline 0 & x_6 & x_{11} & x_{14} \\ 0 & x_8 & x_{15} & x_{16} \end{array} \right) +$$

$$\left( \begin{array}{c|ccc} 0 & 0 & 0 & 0 \\ 0 & 0 & 0 & 0 \\ \hline x_5 & 0 & 0 & 0 \\ x_7 & 0 & 0 & 0 \end{array} \right) + \left( \begin{array}{c|ccc} 0 & x_2 & x_9 & x_{10} \\ 0 & x_4 & x_{11} & x_{12} \\ \hline 0 & 0 & 0 & 0 \\ 0 & 0 & 0 & 0 \end{array} \right).$$

The super subspace

$$W_i \bigcap W_j = \left( \begin{array}{c|ccc} 0 & 0 & 0 & 0 \\ 0 & 0 & 0 & 0 \\ \hline 0 & 0 & 0 & 0 \\ 0 & 0 & 0 & 0 \end{array} \right); \; i \neq j; \; 1 \leq i, j \leq 3 .$$

***Example 1.2.18:*** Let V = {(a b c | d e | f g h) | a, b, c d, e, f, g, h $\in$ Q} be a super vector space over Q. Let $W_1$ = {(a b c | 0 e | 0 0 0 0) | a, b, c, e $\in$ Q}, $W_1$ is a super space of V. Take $W_2$ = {(0 0 c | 0 0 | f g h) | f, g, h, c $\in$ Q}; $W_2$ is a super subspace of V.

Clearly V = $W_1 + W_2$ and $W_1 \cap W_2$ = {(0 0 c | 0 0 | 0 0 0) | c $\in$ Q} is a super subspace of V. In fact $W_1 \cap W_2$ is also a super subspace of both $W_1$ and $W_2$.

***Example 1.2.19:*** Let

$$V = \left\{ \left( \begin{array}{c} a \\ b \\ \hline c \\ d \\ e \\ f \\ \hline g \end{array} \right) \middle| \; a, b, c, d, e, f, g \in R \right\}.$$

V is a super vector space over Q. Take

$$W_1 = \left\{ \left( \begin{array}{c} 0 \\ 0 \\ \hline c \\ d \\ e \\ f \\ \hline 0 \end{array} \right) \middle| \; c, d, e, f \in R \right\},$$

$W_1$ is a super subspace of V.

Let

$$W_2 = \left\{ \left( \begin{array}{c} a \\ b \\ \hline 0 \\ 0 \\ 0 \\ 0 \\ \hline g \end{array} \right) \middle| \; a, b, g \in R \right\},$$

$W_2$ is a super subspace of V. In fact $V = W_1 + W_2$ and

$$W = W_1 \cap W_2 = \left( \begin{array}{c} 0 \\ 0 \\ \hline 0 \\ 0 \\ 0 \\ 0 \\ \hline 0 \end{array} \right)$$

is the super zero subspace of V.

***Example 1.2.20:*** Let

$$V = \left\{ \left( \begin{array}{cc|ccc|c} x_1 & x_2 & x_3 & x_4 & x_5 & x_6 \\ x_7 & x_4 & x_9 & x_{10} & x_{11} & x_{12} \end{array} \right) \right.$$

such that $x_i \in Q$; $1 \le i \le 12\}$, be the super vector space over Q. Let

$$W_1 = \left\{ \left( \begin{array}{cc|ccc|c} x_1 & x_2 & 0 & 0 & 0 & 0 \\ x_7 & x_8 & 0 & 0 & 0 & 0 \end{array} \right) \middle| \; x_{1,} x_2, x_7, x_8 \in Q \right\}$$

be the super subspace of the super vector space V.

$$W_2 = \left\{ \left( \begin{array}{cc|ccc|c} 0 & 0 & 0 & 0 & 0 & x_6 \\ 0 & 0 & x_9 & x_{10} & x_{11} & x_{12} \end{array} \right) \middle| \, x_6, x_9, x_{10}, x_{11}, x_{12} \in Q \right\}$$

be a super subspace of the super vector space V. Clearly V ≠ $W_1$ + $W_2$. But

$$W_2 \cap W_1 = \left( \begin{array}{cc|ccc|c} 0 & 0 & 0 & 0 & 0 & 0 \\ 0 & 0 & 0 & 0 & 0 & 0 \end{array} \right)$$

the zero super matrix of V.

Now we proceed onto define the notion of basis and dimension of a super vector space V.

**DEFINITION 1.2.6:** *Let V be a super vector space over the field F. A subset S of V is said to be linearly dependent (or simply dependent) if there exists distinct super vectors $\alpha_1, \alpha_2, \ldots, \alpha_n$ in S and scalars $c_1, c_2, \ldots, c_n$ in F, not all of which are zero such that $c_1\alpha_1 + c_2\alpha_2 + \ldots + c_n\alpha_n = 0$. A set which is not linearly dependent is called linearly independent. If the set S contains only a finitely many vectors $\alpha_1, \alpha_2, \ldots, \alpha_n$ we some times say that $\alpha_1, \alpha_2, \ldots, \alpha_n$ are dependent (or independent) instead of saying S is dependent (or independent).*

***Example 1.2.21:*** Let V = {($x_1$ $x_2$ | $x_3$ $x_4$ $x_5$ $x_6$ | $x_7$) | $x_i \in Q$; $1 \le i \le 7$} be a super vector space over Q. Consider the super vectors $\alpha_1, \alpha_2, \ldots, \alpha_8$ of V given by

$$\alpha_1 = (1\ 2 \mid 3\ 5\ 6 \mid 7)$$
$$\alpha_2 = (5\ 6 \mid -1\ 2\ 0\ 1 \mid 8)$$
$$\alpha_3 = (2\ 1 \mid 8\ 0\ 1\ 2 \mid 0)$$
$$\alpha_4 = (1\ 1 \mid 1\ 1\ 0\ 3 \mid 2)$$
$$\alpha_5 = (3\ -1 \mid 8\ 1\ 0\ -1 \mid -4)$$
$$\alpha_6 = (8\ 1 \mid 0\ 1\ 1\ 1 \mid -2)$$
$$\alpha_7 = (1\ 2 \mid 2\ 0\ 0\ 1 \mid 0)$$

and

$$\alpha_8 = (3\ 1 \mid 2\ 3\ 4\ 5 \mid 6).$$

Clearly $\alpha_1, \alpha_2, \ldots, \alpha_8$ forms a linearly dependent set of super vectors of V.

***Example 1.2.22:*** Let V = $\{(x_1\ x_2 \mid x_3\ x_4) \mid x_i \in Q\}$ be a super vector space over the field Q.
Consider the super vector

$$\alpha_1 = (1\ 0 \mid 0\ 0),$$
$$\alpha_2 = (0\ 1 \mid 0\ 0),$$
$$\alpha_3 = (0\ 0 \mid 1\ 0)$$

and

$$\alpha_4 = (0\ 0 \mid 0\ 1).$$

Clearly the super vectors $\alpha_1, \alpha_2, \alpha_3, \alpha_4$ form a linearly independent set of V. If we take the super vectors (1 0 | 0 0), (2 1 | 0 0) and (1 4 | 0 0) they clearly form a linearly dependent set of super vectors in V.

**DEFINITION 1.2.7:** *Let V be a super vector space over the field F. A super basis or simply a basis for V is clearly a dependent set of super vectors V which spans the space V. The super space V is finite dimensional if it has a finite basis.*

*Let* $V = \{(x_1 \ldots x_r \mid x_{r+1} \ldots x_k \mid \mid x_{k+1} \ldots x_n)\}$ *be a super vector space over a field F; i.e.* $x_i \in F$; $1 \leq i \leq n$.. *Suppose*

$$W_1 = \{(x_1 \ldots x_r \mid 0 \ldots 0 \mid 0 \ldots 0 \mid 0 \ldots 0)\} \subseteq V$$

*then we call* $W_1$ *a special super subspace of V.*

$$W_2 = \{(0 \ldots 0 \mid x_{r+1} \ldots x_t \mid 0 \ldots 0 \mid 0 \ldots 0) \mid x_{r+1}, \ldots, x_t \in F\}$$

*is again a special super subspace of V.*

$$W_3 = \{(0 \ldots 0 \mid 0 \ldots 0 \mid x_{t+1} \ldots x_k \mid 0 \ldots 0) \mid x_{t+1}, \ldots, x_k \in F\}$$

*is again a special super subspace of V.*

We now illustrate thus situation by the following examples.

***Example 1.2.23:*** Let

$$V = \{(x_1 \mid x_2\ x_3\ x_4 \mid x_5\ x_6) \mid x_i \in Q;\ 1 \le i \le 6\}$$

be a super vector space over Q. The special super subspaces of V are

$$W_1 = \{(x_1 \mid 0\ 0\ 0 \mid 0\ 0) \mid x_1 \in Q\}$$

is a special super subspace of V.

$$W_2 = \{(0 \mid x_2\ x_3\ x_4 \mid 0\ 0) \mid x_2, x_3, x_4 \in Q\}$$

is a special super subspace of V.

$$W_3 = \{(0 \mid 0\ 0\ 0 \mid x_5\ x_6) \mid x_5, x_6 \in Q\}$$

is also a special super subspace of V.

$$W_4 = \{(x_1 | x_2\ x_3\ x_4 | 0\ 0)\}$$

is a special super subspace of V.

$$W_5 = \{(x_1 | 0\ 0\ 0 |\ x_5\ x_6) | x_1\ x_5\ x_6 \in Q\}$$

is a special super subspace of V and

$$W_6 = \{(0 \mid x_2\ x_3\ x_4 | x_5\ x_6) | \ x_2, x_3, x_4, x_5, x_6 \in Q\}$$

is a special super subspace of V. Thus V has only 6 special super subspaces. However if

$$P = \{(0 \mid x_2\ 0\ x_4 \mid 0\ 0) | x_2,\ x_4 \in Q\}$$

is only a super subspace of V and not a special super subspace of V. Likewise

$$T = \{(x_1 | \ 0\ x_3\ 0 | x_5\ 0 \ | x_1, x_3, x_5 \in Q\}$$

is only a super subspace of V and not a special super subspace of V.

***Example 1.2.24:*** Let

$$V = \left\{ \left( \begin{array}{c|cc|cc} x_1 & x_6 & x_{11} & x_{17} & x_{18} \\ x_2 & x_7 & x_{12} & x_{14} & x_{20} \\ x_3 & x_8 & x_{13} & x_{21} & x_{22} \\ \hline x_4 & x_9 & x_{14} & x_{23} & x_{24} \\ x_5 & x_{10} & x_{15} & x_{25} & x_{26} \end{array} \right) \middle| \ x_i \in Q; 1 \le i \le 26 \right\}$$

be a super vector space over Q. The special super subspaces of V are as follows.

$$W_1 = \left\{ \left( \begin{array}{c|cc|cc} x_1 & 0 & 0 & 0 & 0 \\ x_2 & 0 & 0 & 0 & 0 \\ x_3 & 0 & 0 & 0 & 0 \\ \hline 0 & 0 & 0 & 0 & 0 \\ 0 & 0 & 0 & 0 & 0 \end{array} \right) \middle| \, x_1, x_2, x_3 \in Q \right\}$$

is a special super subspace of V

$$W_2 = \left\{ \left( \begin{array}{c|cc|cc} 0 & 0 & 0 & 0 & 0 \\ 0 & 0 & 0 & 0 & 0 \\ 0 & 0 & 0 & 0 & 0 \\ \hline x_4 & 0 & 0 & 0 & 0 \\ x_5 & 0 & 0 & 0 & 0 \end{array} \right) \middle| \, x_4, x_5 \in Q \right\}$$

is a special super subspace of V.

$$W_3 = \left\{ \left[ \begin{array}{c|cc|cc} 0 & x_6 & x_{11} & 0 & 0 \\ 0 & x_7 & x_{12} & 0 & 0 \\ 0 & x_8 & x_{13} & 0 & 0 \\ \hline 0 & 0 & 0 & 0 & 0 \\ 0 & 0 & 0 & 0 & 0 \end{array} \right] \, x_6, x_{11}, x_7, x_8, x_{12}, x_{13} \in Q \right\}$$

is a special super subspace of V.

$$W_4 = \left\{ \left( \begin{array}{c|cc|cc} 0 & 0 & 0 & 0 & 0 \\ 0 & 0 & 0 & 0 & 0 \\ 0 & 0 & 0 & 0 & 0 \\ \hline 0 & x_9 & x_{14} & 0 & 0 \\ 0 & x_{10} & x_{15} & 0 & 0 \end{array} \right) \middle| \, x_9, x_{10}, x_{14} \text{ and } x_{15} \in Q \right\}$$

is a special super subspace of V.

$$W_5 = \left\{ \left( \begin{array}{c|cc|cc} 0 & 0 & 0 & x_{17} & x_{18} \\ 0 & 0 & 0 & x_{19} & x_{20} \\ 0 & 0 & 0 & x_{21} & x_{22} \\ \hline 0 & 0 & 0 & 0 & 0 \\ 0 & 0 & 0 & 0 & 0 \end{array} \right) \middle| \; x_{17}, x_{18}, x_{19}, x_{20}, x_{21} \text{ and } x_{22} \in Q \right\}$$

is a special super subspace of V.

$$W_6 = \left\{ \left( \begin{array}{c|cc|cc} 0 & 0 & 0 & 0 & 0 \\ 0 & 0 & 0 & 0 & 0 \\ 0 & 0 & 0 & 0 & 0 \\ \hline 0 & 0 & 0 & x_{23} & x_{24} \\ 0 & 0 & 0 & x_{25} & x_{26} \end{array} \right) \middle| \; x_{23}, x_{24}, x_{25} \text{ and } x_{26} \in Q \right\}$$

is a special super subspace of V.

$$W_7 = \left\{ \left( \begin{array}{c|cc|cc} x_1 & 0 & 0 & 0 & 0 \\ x_2 & 0 & 0 & 0 & 0 \\ x_3 & 0 & 0 & 0 & 0 \\ \hline x_4 & 0 & 0 & 0 & 0 \\ x_5 & 0 & 0 & 0 & 0 \end{array} \right) \middle| \; x_1 \text{ to } x_5 \in Q \right\}$$

is also a special super subspace of V.

$$W_8 = \left\{ \left( \begin{array}{c|cc|cc} x_1 & x_6 & x_{11} & 0 & 0 \\ x_2 & x_7 & x_{12} & 0 & 0 \\ x_3 & x_8 & x_{13} & 0 & 0 \\ \hline 0 & 0 & 0 & 0 & 0 \\ 0 & 0 & 0 & 0 & 0 \end{array} \right) \middle| \; x_1, x_2, x_3, x_6, x_{11}, x_7, x_8, x_{12} \text{ and } x_{13} \in Q \right\}$$

is also a special super subspace of V and so on,

$$W_t = \left\{ \left( \begin{array}{c|cc|cc} 0 & x_6 & x_{11} & x_{17} & x_{18} \\ 0 & x_7 & x_{12} & x_{19} & x_{20} \\ 0 & x_8 & x_{13} & x_{21} & x_{22} \\ \hline x_4 & x_9 & x_{14} & x_{23} & x_{24} \\ x_5 & x_{10} & x_{15} & x_{25} & x_{26} \end{array} \right) \middle| \; x_i \in Q;\ 4 \leq i \leq 26 \right\}$$

is also a special super subspace of V.

Now we have seen the definition and examples of special super subspace of a super vector space V. We now proceed onto define the standard basis or super standard basis of V.

Let F be a field $V = (F^{n_1} \mid F^{n_2} \mid \ldots \mid F^{n_n})$ be a super vector space over F. The super vectors $\in_1, \ldots, \in_{n_1}, \in_{n_1+1}, \ldots, \in_{n_n}$ given by

$$\begin{aligned} \in_1 &= (1\ 0 \ldots 0 | 0 \ldots 0 | 0 \ldots | 00 \ldots 0) \\ \in_2 &= (0\ 1 \ldots 0 | 0 \ldots | 0 \ldots | 0 \ldots 0) \\ &\vdots \\ \in_{n_1} &= (0\ \ldots\ 1 |\ 0 \ldots 0 | 0 \ldots | 0 \ldots 0) \\ \in_{n_1+1} &= (0\ \ldots\ 0 | 1\ 0 \ldots 0 | \ldots | 0 \ldots 0) \\ &\vdots \\ \in_{n_2} &= (0\ \ldots\ 0 |\ 0 \ldots 1 | 0 \ldots | 0 \ldots 0) \\ &\vdots \\ \in_{n_n} &= (0\ \ldots\ 0 |\ 0 \ldots 0 | \ldots | 0 \ldots 01) \end{aligned}$$

forms a linearly independent set and it spans V; so these super vectors form a basis of V known as the super standard basis of V.

We will illustrate this by the following example.

***Example 1.2.25:*** Let $V = \{(x_1\ x_2\ x_3 \mid x_4\ x_5) \mid x_i \in Q;\ 1 \leq i \leq 5\}$ be a super vector space over Q. The standard basis of V is given by

$$\in_1 = (1\ 0\ 0 \mid 0\ 0),$$
$$\in_2 = (0\ 1\ 0 \mid 0\ 0),$$
$$\in_3 = (0\ 0\ 1 \mid 0\ 0),$$
$$\in_4 = (0\ 0\ 0 \mid 1\ 0),$$

and

$$\in_5 = (0\ 0\ 0 \mid 0\ 1),$$

***Example 1.2.26:*** Let V = $\{(x_1\ x_2\ x_3\ x_4\ x_5 \mid x_6\ x_7\ x_8) \mid x_i \in Q; 1 \leq i \leq 8\}$ be a super vector space over Q. The standard basis for V is given by

$$\in_1 = (1 \mid 0\ 0\ 0\ 0 \mid 0\ 0\ 0),$$
$$\in_2 = (0 \mid 1\ 0\ 0\ 0 \mid 0\ 0\ 0),$$
$$\in_3 = (0 \mid 0\ 1\ 0\ 0 \mid 0\ 0\ 0),$$
$$\in_4 = (0 \mid 0\ 0\ 1\ 0 \mid 0\ 0\ 0),$$
$$\in_5 = (0 \mid 0\ 0\ 0\ 1 \mid 0\ 0\ 0),$$
$$\in_6 = (0 \mid 0\ 0\ 0\ 0 \mid 1\ 0\ 0),$$
$$\in_7 = (0 \mid 0\ 0\ 0\ 0 \mid 0\ 1\ 0),$$

and

$$\in_8 = (0 \mid 0\ 0\ 0\ 0 \mid 0\ 0\ 1),$$

Clearly it can be checked by the reader $\in_1, \in_2, \ldots, \in_8$ forms a super standard basis of V.

***Example 1.2.27:*** Let

$$V = \left\{ \left( \begin{array}{c|cc} x_1 & x_5 & x_6 \\ x_2 & x_7 & x_{18} \\ x_3 & x_9 & x_{10} \\ \hline x_4 & x_{11} & x_{12} \end{array} \right) \middle| \; x_i \in Q; 1 \leq i \leq 12 \right\}$$

be a super vector space over Q. The standard basis for V is ;

$$\in_1 = \left(\begin{array}{c|cc} 1 & 0 & 0 \\ 0 & 0 & 0 \\ 0 & 0 & 0 \\ \hline 0 & 0 & 0 \end{array}\right), \in_2 = \left(\begin{array}{c|cc} 0 & 0 & 0 \\ 1 & 0 & 0 \\ 0 & 0 & 0 \\ \hline 0 & 0 & 0 \end{array}\right), \in_3 = \left(\begin{array}{c|cc} 0 & 0 & 0 \\ 0 & 0 & 0 \\ 0 & 0 & 0 \\ \hline 1 & 0 & 0 \end{array}\right),$$

$$\in_4 = \left(\begin{array}{c|cc} 0 & 1 & 0 \\ 0 & 0 & 0 \\ 0 & 0 & 0 \\ \hline 0 & 0 & 0 \end{array}\right), \in_5 = \left(\begin{array}{c|cc} 0 & 0 & 1 \\ 0 & 0 & 0 \\ 0 & 0 & 0 \\ \hline 0 & 0 & 0 \end{array}\right), \in_6 = \left(\begin{array}{c|cc} 0 & 0 & 0 \\ 0 & 1 & 0 \\ 0 & 0 & 0 \\ \hline 0 & 0 & 0 \end{array}\right),$$

$$\in_7 = \left(\begin{array}{c|cc} 0 & 0 & 0 \\ 0 & 0 & 1 \\ 0 & 0 & 0 \\ \hline 0 & 0 & 0 \end{array}\right), \in_8 = \left(\begin{array}{c|cc} 0 & 0 & 0 \\ 0 & 0 & 0 \\ 0 & 1 & 0 \\ \hline 0 & 0 & 0 \end{array}\right), \in_9 = \left(\begin{array}{c|cc} 0 & 0 & 0 \\ 0 & 0 & 0 \\ 0 & 0 & 1 \\ \hline 0 & 0 & 0 \end{array}\right),$$

$$\in_{10} = \left(\begin{array}{c|cc} 0 & 0 & 0 \\ 0 & 0 & 0 \\ 0 & 0 & 0 \\ \hline 0 & 1 & 0 \end{array}\right), \in_{11} = \left(\begin{array}{c|cc} 0 & 0 & 0 \\ 0 & 0 & 0 \\ 0 & 0 & 0 \\ \hline 0 & 0 & 1 \end{array}\right) \text{ and } \in_{12} = \left(\begin{array}{c|cc} 0 & 0 & 0 \\ 0 & 0 & 0 \\ 1 & 0 & 0 \\ \hline 0 & 0 & 0 \end{array}\right).$$

The reader is expected to verify that $\in_1, \in_2, \ldots, \in_{12}$ forms a super standard basis of V.

Now we are going to give a special notation for the super row vectors which forms a super vector space and the super matrices which also form a super vector space. Let X = $(x_1 \ldots x_t \mid x_{t+1} \ldots x_k \mid \ldots \mid x_{r+1} \ldots x_n)$ be a super row vector with entries from Q.

Define X = $(A_1 \mid A_2 \mid \ldots \mid A_m)$ where each $A_i$ is a row vector $A_1$ corresponds to the row vectors $(x_1 \ldots x_t)$, the set of row vectors $(x_{t+1} \ldots x_k)$ to $A_2$ and so on. Clearly $m \leq n$.

Likewise a super matrix is also given a special representation.

Suppose

$$A = \left(\begin{array}{ccc|cc} x_1 & x_2 & x_3 & x_{16} & x_{17} \\ x_4 & x_5 & x_6 & x_{18} & x_{19} \\ x_7 & x_5 & x_9 & x_{20} & x_{21} \\ \hline x_{10} & x_{11} & x_{12} & x_{22} & x_{23} \\ x_{13} & x_{14} & x_{15} & x_{24} & x_{25} \end{array}\right) = \left(\begin{array}{c|c} A_1 & A_2 \\ \hline A_3 & A_4 \end{array}\right)$$

where $A_1$ is a 3 × 3 matrix given by

$$A_1 = \begin{pmatrix} x_1 & x_2 & x_3 \\ x_4 & x_5 & x_6 \\ x_7 & x_5 & x_9 \end{pmatrix}, A_2 = \begin{pmatrix} x_{16} & x_{17} \\ x_{15} & x_{19} \\ x_{20} & x_{21} \end{pmatrix}$$

is a 3 × 2 rectangular matrix

$$A_3 = \begin{pmatrix} x_{10} & x_{11} & x_{12} \\ x_{15} & x_{14} & x_{15} \end{pmatrix}$$

is again a rectangular 2 × 3 matrix with entries from Q and

$$A_4 = \begin{pmatrix} x_{22} & x_{23} \\ x_{24} & x_{25} \end{pmatrix}$$

is again a 2 × 2 square matrix.

We see the components of a super row vector are row vectors where as the components of a super matrix are just matrices.

Now we proceed onto prove the following theorem.

**THEOREM 1.2.5:** *Let V be a super vector space which is spanned by a finite set of super vectors $\beta_1$, …, $\beta_m$. Then any independent set of super vectors in V is finite and contains no more than m elements.*

*Proof:* Given V is a super vector space. To prove the theorem it suffices to show that every subset S of V which contains more than m super vectors is linearly dependent. Let S be such a set. In S there are distinct super vectors $\alpha_1, \ldots, \alpha_n$ where $n > m$. Since $\beta_1, \beta_2, \ldots, \beta_m$ span V their exists scalars $A_{ij}$ in F such that

$$\alpha_j = \sum_{i=1}^{m} A_{ij}\beta_i .$$

For any n-scalars $x_1, \ldots, x_n$ we have

$$\begin{aligned} x\alpha_1 + \ldots x_n\alpha_n &= \sum_{j=1}^{n} x_j\alpha_j \\ &= \sum_{j=1}^{n} x_j \sum_{i=1}^{m} A_{ij}\beta_i \\ &= \sum_{j=1}^{n}\sum_{i=1}^{m} (A_{ij}x_j)\beta_i \\ &= \sum_{i=1}^{m}\left(\sum_{j=1}^{n} A_{ij}x_j\right)\beta_i . \end{aligned}$$

Since $n > m$ we see there exists scalars $x_1, \ldots, x_n$ not all zero such that

$$\sum_{j=1}^{n} A_{ij}x_j = 0;\ 1 \le i \le m$$

Hence $x_1\alpha_1 + \ldots + x_n\alpha_n = 0$ which proves S is a linearly dependent set.

The immediate consequence of this theorem is that any two basis of a finite dimensional super vector space have same number of elements.

As in case of usual vector space when we say a supervector space is finite dimensional it has finite number of elements in its basis.

We illustrate this situation by a simple example.

***Example 1.2.28:*** Let $V = \{(x_1\ x_2\ x_3 \mid x_4) \mid x_i \in Q;\ 1 \le i \le 4\}$ be a super vector space over Q. It is very clear that V is finite dimensional and has only four elements in its basis. Consider a set

S = {(1 0 1 | 0), (1 2 3 | 4), (4 0 0 | 3), (0 1 2 | 1) and (1 2 0 | 2)} $\{x_1, x_2, x_3, x_4, x_5\} \subseteq$ V, to S is a linearly dependent subset of V;

i.e. to show this we can find scalars $c_1, c_2, c_3, c_4$ and $c_5$ in Q not all zero such that

$\sum c_i x_i = 0$. $c_1$ (1 0 1 | 0) + $c_2$ (1 2 3 | 4) + $c_3$ (4 0 0 | 3) + $c_4$ ( 0 1 2 | 1) + $c_5$ (1 2 0 | 2) = (0 0 0 | 0)
gives

$$c_1 + c_2 + 4c_3 + c_5 = 0$$
$$2c_2 + c_4 + 2c_5 = 0$$
$$c_1 + 3c_2 + 2c_4 = 0$$
$$4c_2 + 3c_3 + c_4 + 2c_5 = 0.$$

It is easily verified we have non zero values for $c_1, \ldots, c_5$ hence the set of 5 super vectors forms a linearly dependent set.

It is left as an exercise for the reader to prove the following simple lemma.

**LEMMA 1.2.1:** *Let S be a linearly independent subset of a super vector space V. Suppose β is a vector in V and not in the super subspace spanned by S, then the set obtained by adjoining β to S is linearly independent.*

We state the following interesting theorem.

**THEOREM 1.2.6:** *If W is a super subspace of a finite dimensional super vector space V, every linearly independent subset of W is finite and is part of a (finite basis for W).*

Since super vectors are also vectors and they would be contributing more elements while doing further operations. The above theorem can be given a proof analogous to usual vector spaces.

Suppose $S_0$ is a linearly independent subset of W. If S is a linearly independent subset of W containing $S_0$ then S is also a linearly independent subset of V; since V is finite dimensional, S contains no more than dim V elements.

We extend $S_0$ to a basis for W as follows: $S_0$ spans W, then $S_0$ is a basis for W and we are done. If $S_0$ does not span W we use the preceding lemma to find a super vector $\beta_1$ in W such that the set $S_1 = S_0 \cup \{\beta_1\}$ is independent. If $S_1$ spans W, fine. If not, we apply the lemma to obtain a super vector $\beta_2$ in W such that $S_2 = S_1 \cup \{\beta_2\}$ is independent.

If we continue in this way then (in not more than dim V steps) we reach at a set $S_m = S_0 \cup \{\beta_1, \ldots, \beta_m\}$ which is a basis for W.

The following two corollaries are direct and is left as an exercise for the reader.

**COROLLARY 1.2.1:** *If W is a proper super subspace of a finite dimensional super vector space V, then W is finite dimensional and dim W < dim V.*

**COROLLARY 1.2.2:** *In a finite dimensional super vector space V every non empty linearly independent set of super vectors is part of a basis.*

However the following theorem is simple and is left for the reader to prove.

**THEOREM 1.2.7:** *If $W_1$ and $W_2$ are finite dimensional super subspaces of a super vector space V then $W_1 + W_2$ is finite dimensional and dim $W_1$ + dim $W_2$ = dim ($W_1 \cap W_2$) + dim ($W_1 + W_2$).*

We have seen in case of super vector spaces we can define the elements of them as n × m super matrices or as super row vectors or as super column vectors.

So how to define linear transformations of super vector spaces. Can we have linear transformations from a super vector space to a super vector space when both are defined over the same field F?

## 1.3 Linear Transformation of Super Vector Spaces

For us to have a meaningful linear transformation, if V is a super vector space, super row vectors having n components ($A_1$, …, $A_n$) where each $A_i$ a is row vector of same length then we should have W also to be a super vector space with super row vectors having only n components of some length, need not be of identical length. When we say two super vector have same components we mean that both the row vector must have same number of partitions. For instance $X = (x_1, x_2, \ldots, x_n)$ and $Y = (y_1, y_2, \ldots, y_m)$, $m \neq n$, the number of partitions in both of them must be the same if $X = (A_1 \mid \ldots \mid A_t)$ then $Y = (B_1 \mid \ldots \mid B_t)$ where $A_i$'s and $B_j$'s are row vectors $1 \leq i, j \leq t$.

Let $$X = (2\ 1\ 0\ 5\ 6\ -1)$$
and
$$Y = (1\ 0\ 2\ 3\ 4\ 5\ 7\ 8).$$
If X is partitioned as
$$X = (2 \mid 1\ 0\ 5 \mid 6\ -1)$$
and
$$Y = (1\ 0\ 2 \mid 3\ 4\ 5 \mid 7\ 8\ 1).$$

$X = (A_1 \mid A_2 \mid A_3)$ and $Y = (B_1 \mid B_2 \mid B_3)$ where $A_1 = 2$, $B_1 = (1\ 0\ 2)$; $A_2$ $(1\ 0\ 5)$, $B_2 = (3\ 4\ 5)$, $A_3 = (6\ -1)$ and $B_3 = (7\ 8\ 1)$.

We say the row vectors X and Y have same number of partitions or to be more precise we say the super vectors have same number of partitions. We can define linear transformation between two super vector spaces. Super vectors with same number of elements or with same number of partition of the row vectors; otherwise we cannot define linear transformation.

Let V be a super vectors space over the field F with super vector $X \in V$ then $X = (A_1 \mid \ldots \mid A_n)$ where each $A_i$ is a row vector. Suppose W is a super vector space over the same field F if for a super row vector, $Y \in W$ and if $Y = (B_1 \mid \ldots \mid B_n)$ then we say V and W are super vector spaces with same type of super row vectors or the number of partitions of the row vectors in both V and W are equal or the same.

We call such super vector spaces as same type of super vector spaces.

**DEFINITION 1.3.1:** *Let V and W be super vector spaces of the same type over the same field F. A linear transformation from V into W is a function T from V into W such that $T(c\alpha + \beta) = cT\alpha + \beta$ for all scalars c in F and the super vectors $\alpha, \beta \in V$;*

*i.e. if $\alpha = (A_1 \mid \ldots \mid A_n) \in V$ then $T\alpha = (B_1 \mid \ldots \mid B_n) \in W$,*

*i.e. T acts on $A_1$ in such a way that it is mapped to $B_1$ i.e. first row vector of α i.e. $A_1$ is mapped into the first row vector $B_1$ of Tα. This is true for $A_2$ and so on.*

Unless this is maintained the map T will not be a linear transformation preserving the number of partitions. We first illustrate it by an example. As our main aim of introducing any notion is not for giving nice definition but our aim is to make the reader understand it by simple examples as the very concept of super vectors happen to be little abstract but very useful in practical problems.

***Example 1.3.1:*** Let V and W be two super vector spaces of same type defined over the field Q. Let $V = \{(x_1\ x_2\ x_3 \mid x_4\ x_5 \mid x_6) \mid x_i \in Q;\ 1 \leq i \leq 6\}$ and $W = \{(x_1\ x_2 \mid x_3\ x_4\ x_5 \mid x_6\ x_7\ x_8) \mid x_i \in Q;\ 1 \leq i \leq 8\}$.

We see both of them have same number of partitions and we do not demand the length of the vectors in V and W to be the same but we demand only the length of the super vectors to be the same, for here we see in both the super vector spaces V and W super vectors are of length 3 only but as vectors V has natural length 6 and W has natural length 8.

Let $T : V \to W$

$T(x_1\ x_2\ x_3 \mid x_4\ x_5 \mid x_6)$
$= (x_1 + x_3\quad x_2 + x_3 \mid x_4\ \ x_4 + x_5\quad x_5 \mid x_6\ 0\ {-x_6})$.

It is easily verified that T is a linear transformation from V into W.

***Example 1.3.2:*** Suppose $V = \{(x_1\ x_2 \mid x_3\ x_4\ x_5) \mid x_1, \ldots, x_5 \in Q\}$ and $W = \{(x_1 \mid x_2\ x_3 \mid x_4\ x_5) \mid x_i \in Q;\ 1 \leq i \leq 5\}$ both super vector spaces over F. Suppose we define a map $T : V \to W$ by $T[(x_1\ x_2 \mid x_3\ x_4\ x_5)] = (x_1 + x_2 \mid x_3 + x_4,\ x_5 \mid 0\ 0)$.

T is a linear transformation but does not preserve partitions. So such linear transformation also exists on super vector spaces.

***Example 1.3.3:*** Let $V = \{(x_1\ x_2\ x_3 \mid x_4\ x_5 \mid x_6) \mid x_i \in Q;\ 1 \leq i \leq 6\}$ and $W = \{(x_1\ x_2 \mid x_3\ x_4) \mid x_i \in Q,\ 1 \leq i \leq 4\}$. Then we cannot define a linear transformation of the super vector spaces V and W. So we demand if we want to define a linear transformation which is not partition preserving then we demand the number partition in the range space (i.e. the super vector space which is the range of T) must be greater than the number of partitions in the domain space.

Thus with this demand in mind we define the following linear transformation of two super vector spaces.

**DEFINITION 1.3.2:** *Let $V = \{(A_1 \mid A_2 \mid \ldots \mid A_n) \mid A_i$ row vectors with entries from a field F\} be a super vector space over F. Suppose $W = \{(B_1 \mid B_2 \mid \ldots \mid B_m)$, $B_i$ row vectors from the same field F; $i = 1, 2, \ldots, m\}$ be a super vector space over F. Clearly $n \leq m$. Then we call T the linear transformation i.e. $T: V \to W$ where $T(A_i) = B_j$, $1 \leq i \leq n$ and $1 \leq j \leq m$ and entries $B_K$ in W which do not have an associated $A_i$ in V are just put as zero row vectors and if T is a linear transformation from $A_i$ to $B_j$; T is called as the linear transformation which does not preserve partition but T acts more like an embedding. Only when $m = n$ we can define the notion of partition preserving linear transformation of super vector spaces from V into W. But when*

*$n > m$ we will not be in a position to define linear transformation from super vector space V into W.*

With these conditions we will give yet some more examples of linear transformation from a super vector space V into a super vector space W both defined over the same field F.

***Example 1.3.4**:* Let V = $\{(x_1\ x_2\ x_3 \mid x_4\ x_5\ x_6 \mid x_7) \mid x_i \in Q; 1 \leq i \leq 7\}$ be a super vector space over Q. W = $\{(x_1\ x_2 \mid x_3\ x_4 \mid x_5\ x_6 \mid x_7\ x_8 \mid x_9) \mid x_i \in Q; 1 \leq i \leq 9\}$ be a super vector space over Q. Define $T : V \to W$ by $T(x_1\ x_2\ x_3 \mid x_4 \mid x_5\ x_6 \mid x_7) \to (x_1 + x_2\ \ x_2 + x_3 \mid x_3 + x_4 \mid x_5 + x_6\ x_5 \mid 0\ 0 \mid x_9)$. It is easily verified T is a linear transformation from V to W, we can have more number of linear transformations from V to W. Clearly T does not preserve the partitions. We also note that number of partitions in V is less than the number of partitions in W.

We give yet another example.

***Example 1.3.5**:* Let $T : V \to W$ be a linear transformation from V into W; where V = $\{(x_1\ x_2 \mid x_3 \mid x_4\ x_5\ x_6) \mid x_i \in Q; 1 \leq i \leq 6\}$ is a super vector space over Q. Let W = $\{(x_1\ x_2 \mid x_3 \mid x_4\ x_5\ x_6\ x_7) \mid x_i \in Q; 1 \leq i \leq 7\}$ be a super vector space over Q. Define $T((x_1\ x_2\ x_3 \mid x_4\ x_5\ x_6) = (x_1 + x_2 \mid x_2 + x_3\ x_2 \mid x_4 + x_5\ x_5 + x_6\ x_6 + x_4\ x_4 + x_5 + x_6)$

It is easily verified that T is a linear transformation from the super vector space V into the super vector space W.

Now we proceed into define the kernel of T or null space of T.

**DEFINITION 1.3.3:** *Let V and W be two super vector spaces defined over the same field F. Let $T : V \to W$ be a linear transformations from V into W. The null space of T which is a super subspace of V is the set of all super vectors α in V such that $T\alpha = 0$. It is easily verified that null space of T; $N = \{\alpha \in V \mid T(\alpha) = 0\}$ is a super subspace of V. For we know $T(0) = 0$ so N is non empty.*

If

$$T\alpha_1 = T\alpha_2 = 0$$

then

$$\begin{aligned} T(c\alpha_1 + \alpha_2) &= cT\alpha_1 + T\alpha_2 \\ &= c.0 + 0 \\ &= 0. \end{aligned}$$

So that for every $\alpha_1$ $\alpha_2 \in$ N, $c\alpha_1 + \alpha_2 \in$ N. Hence the claim. We see when V is a finite dimensional super vector space then we see some interesting properties relating the dimension can be made as in case of vector spaces.

Now we proceed on to define the notion of super null subspace and the super rank space of a linear transformation from a super vector space V into a super vector space W.

**DEFINITION 1.3.4:** *Let V and W be two super vector spaces over the field F and let T be a linear transformation from V into W. The super null space or null super space of T is the set of all super vectors $\alpha$ in V such that $T\alpha = 0$. If V is finite dimensional, the super rank of T is the dimension of the range of T and nullity of T is the dimension of the null space of T.*

This is true for both linear transformations preserving the partition as well as the linear transformations which does not preserve the partition.

***Example 1.3.6:*** Let V = $\{(x_1\ x_2 \mid x_3\ x_4\ x_5) \mid x_i \in Q; 1 \leq i \leq 5\}$ be a super vector space over Q and W = $\{(x_1\ x_2 \mid x_3\ x_4) \mid x_i \in Q; 1 \leq i \leq 4\}$ be a super vector space over Q. Let T: V $\to$ W defined by $T(x_1\ x_2 \mid x_3\ x_4\ x_5) = (x_1 + x_2, x_2 \mid x_3 + x_4, x_4 + x_5)$. T is easily verified to be a linear transformation.

The null super subspace of T is N = $\{(0\ 0 \mid k, k, -k) \mid k \in Q\}$ which is a super subspace of V. Now dim V = 5 and dim W = 4. Find dim N and prove rank T + nullity T = 5.

Suppose V is a finite dimensional super vector space over a field F. We call B = $\{x_1, \ldots, x_n\}$ to be a basis of V if each of the $x_i$'s are super vectors from V and they form a linearly independent set and span V. Suppose V = $\{(x_1 \mid \ldots \mid \ldots \mid \ldots \mid \ldots \mid$ $x_n) \mid x_i \in Q; 1 \leq i \leq n\}$ then dimension of V is n and V has B to be its basis then B has only n-linearly independent elements in it which are super vectors.

So in case of super vector spaces the basis B forms a set which contains only supervectors.

**THEOREM 1.3.1:** *Let V be a finite dimensional super vector space over the field F i.e., $V = \{(x_1 \mid x_2 \mid \ldots \mid \ldots \mid \ldots \mid x_n) \mid x_i \in F; 1 \leq i \leq n\}$ Let $(\alpha_1, \alpha_2, \ldots, \alpha_n)$ be a basis of V i.e., each $\alpha_i$ is a super vector; i = 1, 2, …, n. Let $W = \{(x_1 \mid \ldots \mid \ldots \mid x_m) \mid x_i \in F; 1 \leq i \leq m\}$ be a super vector space over the same field F and let $\beta_1, \ldots, \beta_n$ be super vectors in W. Then there is precisely one linear transformation T from V into W such that $T\alpha_j = \beta_j$; j = 1, 2, …, n.*

*Proof:* To prove that there exists some linear transformation T from V into W with $T\alpha_j = \beta_j$ we proceed as follows:

Given $\alpha$ in V, a super vector there is a unique n-tuple of scalars in F such that $\alpha = x_1\alpha_1 + \ldots + x_n \alpha_n$ where each $\alpha_i$ is a super vector and $\{\alpha_1, \ldots, \alpha_n\}$ is a basis of V; $(1 \leq i \leq n)$. For this $\alpha$ we define $T\alpha = x_1\beta_1 + \ldots + x_n\beta_n$.

Then T is well defined rule for associating with each super vector $\alpha$ in V a super vector $T\alpha$ in W. From the definition it is clear that $T\alpha_j = \beta_j$ for each j. To show T is linear let $\beta = y_1\alpha_1 + \ldots + y_n\alpha_n$ be in V for any scalar $c \in F$. We have $c\alpha + \beta = (cx_1 + y_1)\, \alpha_1 + \ldots + (cx_n + y_n)\, \alpha_n$ and so by definition $T(c\alpha + \beta) = (cx_1 + y_1)\, \beta_1 + \ldots + (cx_n + y_n)\, \beta_n$.

On the other hand

$$c(T\alpha) + T\beta = c \sum_{i=1}^{n} x_i\beta_i + \sum_{i=1}^{n} y_i\beta_i = \sum_{i=1}^{n} (cx_i + y_i)\, \beta_i$$

and thus

$$T(c\alpha + \beta) = c(T\alpha) + T\beta.$$

If U is a linear transformation from V into W with $U\alpha_j = \beta_j$; j = 1, 2, …, n then for the super vector $\alpha = \sum_{i=1}^{n} x_i\alpha_i$ we have $U\alpha =$ $U(\sum_{i=1}^{n} x_i\alpha_i) = \sum_{i=1}^{n} x_i U\alpha_i = \sum_{i=1}^{n} x_i\beta_i$, so that U is exactly the rule T which we have just defined above. This proves that the linear transformation with $T\alpha_j = \beta_j$ is unique.

Now we prove a theorem relating rank and nullity.

**THEOREM 1.3.2:** *Let V and W be super vector spaces over the field F of same type and let T be a linear transformation from V into W. Suppose that V is finite-dimensional. Then super rank T + super nullity T = dim V.*

*Proof:* Let V and W be super vector spaces of the same type over the field F and let T be a linear transformation from V into W. Suppose the super vector space V is finite dimensional with $\{\alpha_1, \ldots, \alpha_k\}$ a basis for the super subspace which is the null super space N of V under the linear transformation T. There are super vectors $\{\alpha_{k+1}, \ldots, \alpha_n\}$ in V such that $\{\alpha_1, \ldots, \alpha_n\}$ is a basis for V.

We shall prove $\{T\alpha_{k+1}, \ldots, T\alpha_n\}$ is a basis for the range of T. The super vectors $\{T\alpha_{k+1}, \ldots, T\alpha_n\}$ certainly span the range of T and since $T\alpha_j = 0$ for $j \leq k$, we see $T\alpha_{k+1}, \ldots, T\alpha_n$ span the range. To prove that these super vectors are linearly independent; suppose we have scalars $c_i$ such that

$$\sum_{i=k+1}^{n} c_i(T\alpha_i) = 0.$$

This says that $T(\sum_{i=k+1}^{n} c_i\alpha_i) = 0$ and accordingly the super vector $\alpha = \sum_{i=k+1}^{n} c_i\alpha_i$ is in the null super space of T. Since $\alpha_1, \ldots, \alpha_k$ form a basis of the null super space N there must be scalars $b_1, \ldots, b_k$ such that $\alpha = \sum_{i=1}^{k} b_i\alpha_i$. Thus

$$\sum_{i=1}^{k} b_i\alpha_i - \sum_{j=k+1}^{n} c_j\alpha_i = 0$$

Since $\alpha_1, \ldots, \alpha_n$ are linearly independent we must have $b_1 = b_2 = \ldots = b_k = c_{k+1} = \ldots = c_n = 0$.
If r is the rank of T, the fact that $T\alpha_{k+1}, \ldots, T\alpha_n$ form a basis for the range of T tells us that r = n – k. Since k is the nullity of T and n is the dimension of V, we have the required result.

Now we want to distinguish the linear transformation T of usual vector spaces from the linear transformation T of the super vector spaces.

To this end we shall from here onwards denote by $T_s$ the linear transformation of a super vector space V into a super vector space W.

Further if $V = \{(A_1 \mid \ldots \mid A_n) \mid A_i$ are row vectors with entries from F, a field$\}$ and V a super vector space over the field F and $W = \{(B_1 \mid \ldots \mid B_n) \mid B_i$ are row vectors with entries from the same field F$\}$ and W is also a super vector space over F. We say $T_s$ is a linear transformation of a super vector space V into W if $T = (T_1 \mid \ldots \mid T_n)$ where $T_i$ is a linear transformation from $A_i$ to $B_i$; i = 1, 2, …, n. Since $A_i$ is a row vector and $B_i$ is a row vector $T_i(A_i) = B_i$ is a linear transformation of the vector space with collection of row vectors $A_i = (x_1 \ldots x_i)$ with entries from F into the vector space of row vectors $B_i$ with entries from F. This is true for each and every i; i = 1, 2, …, n.

Thus a linear transformation $T_s$ from a super vector space V into W can itself be realized as a super linear transformation as $T_s = (T_1 \mid \ldots \mid T_n)$.

From here on words we shall denote the linear transformation of finite dimensional super vector spaces by $T_s = (T_1 \mid T_2 \mid \ldots \mid T_n)$ when the linear transformation is partition preserving in case of linear transformation which do not preserve partition will also be denoted only by $T_s = (T_1 \mid T_2 \mid \ldots \mid T_n)$. Now if $(A_1 \mid \ldots \mid A_n) \in V$ then $T(A_1 \mid \ldots \mid A_n) = (T_1A_1 \mid T_2A_2 \mid \ldots \mid T_nA_n) = (B_1 \mid B_2 \mid \ldots \mid B_n) \in W$ in case $T_s$ is a partition preserving linear transformation.

If $T_s$ is not a partition preserving transformation and if $(B_1 \mid \ldots \mid B_m) \in W$ we know m > n so $T(A_1 \mid \ldots \mid A_n) = (T_1A_1 \mid \ldots \mid T_nA_n \mid 0\ 0 \mid \ldots \mid 0\ 0\ 0) = (T_1A_1 \mid 0\ 0 \ldots \mid T_2A_2 \mid 0\ 0 \mid 0\ 0 \mid \ldots \mid T_nA_n)$ in whichever manner the linear transformation has been defined.

**THEOREM 1.3.3:** *Let $V = \{(A_1 \mid \ldots \mid A_n) \mid A_i$'s are row vectors with entries from F; $1 \le i \le n\}$ a super vector space over F. $W = \{(B_1 \mid \ldots \mid B_n) \mid B_i$'s are row vectors with entries from F; $1 \le i \le n\}$ a super vector space over F. Let $T_s = (T_1 \mid \ldots \mid T_n)$ and $U_s = (U_1 \mid \ldots \mid U_n)$ be linear transformations from V into W. The function $T_s + U_s = (T_1 + U_1 \mid \ldots \mid T_n + U_n)$ defined by $(T_s + U_s)(\alpha) = (T_s + U_s)(A_1 \mid \ldots \mid A_n)$ (where $\alpha \in V$ is such that $\alpha = (A_1 \mid \ldots \mid A_n) = (T_1A_1 + U_1A_1 \mid T_2A_2 + U_2A_2 \mid \ldots \mid T_nA_n + U_nA_n)$ is a linear, transformation from V into W. If d is any element of F, the function $dT = (dT_1 \mid \ldots \mid dT_n)$ defined by $(dT)(\alpha) = d(T\alpha) = d(T_1A_1 \mid \ldots \mid T_nA_n)$ is a linear transformation from V into W, the set of all linear transformations from V into W together with addition and scalar multiplication defined above is a super vector space over the field F.*

*Proof:* Suppose $T_s = (T_1 \mid \ldots \mid T_n)$ and $U_s = (U_1 \mid \ldots \mid U_n)$ are linear transformations of the super vector space V into the super vector space W and that we define $(T_s + U_s)$ as above then

$$(T_s + U_s)(d\alpha + \beta) = T_s(d\alpha + \beta) + U_s(d\alpha + \beta)$$

where $\alpha = (A_1 \mid A_2 \mid \ldots \mid A_n)$ and $\beta = (C_1 \mid \ldots \mid C_n) \in V$ and $d \in F$.

$(T_s = U_s)(d\alpha + \beta)$

$= (T_1 + U_1 \mid \ldots \mid T_n + U_n)(dA_1 + C_1 \mid dA_2 + C_2 \mid \ldots \mid dA_n + C_n)$

$= (T_1(dA_1 + C_1) \mid T_2(dA_2 + C_2) \mid \ldots \mid T_n(dA_n + C_n)) + (U_1(dA_1 + C_1) \mid \ldots \mid U_n(dA_n + C_n))$

$= T_1(dA_1 + C_1) + U_1(dA_1 + C_1) \mid \ldots \mid T_n(dA_n + C_n) + U_n(dA_n+C_n))$

$= (dT_1A_1 + T_1C_1 \mid \ldots \mid dT_nA_n + T_nC_n) + (dU_1A_1 + U_1C_1) \mid \ldots \mid dU_nA_n + U_nC_n)$

$= (dT_1A_1 \mid \ldots \mid dT_nA_n) + (T_1C_1 \mid \ldots \mid T_nC_n) + (dU_1A_1 \mid \ldots \mid dU_nA_n) + (U_1C_1 \mid \ldots \mid U_nC_n)$

$= (dT_1A_1 \mid \ldots \mid dT_nA_n) + (dU_1A_1 \mid \ldots \mid dU_nA_n) + (T_1C_1 \mid \ldots \mid T_nC_n) + (U_1C_1 \mid \ldots \mid U_nC_n)$

$= (d(T_1+U_1)A_1 \mid \ldots \mid d(T_n + U_n) A_n) + ((T_1+U_1) C_1 \mid \ldots \mid (T_n + U_n) C_n)$

which shows $(T_s + U_s)$ is a linear transformation. Similarly

$$\begin{aligned} &(eT_s)(d\alpha + \beta) \\ =\ & e(T_s(d\alpha + \beta)) \\ =\ & e(T_s d\alpha + T_s\beta) \\ =\ & e[d(T_s\alpha)] + eT_s\beta \\ =\ & ed[T_1A_1 \mid \ldots \mid T_nA_n] + e[T_1C_1 \mid \ldots \mid T_nC_n] \\ =\ & edT_s\alpha + eT_s\beta \\ =\ & d(eT_s)\alpha + eT_s\beta \end{aligned}$$

which shows $eT_s$ is a linear transformation.

We see the elements $T_s$, $U_s$ which are linear transformations from super vector spaces are also super vectors as $T_s = (T_1 \mid \ldots \mid T_n)$ and $U_s = (U_1 \mid \ldots \mid U_n)$. Thus the collection of linear transformations $T_s$ from a super vector space V into a super vector space W is a vector space over F. Since each of the linear transformation are super vectors we can say the collection of linear transformation from super vector spaces is again a super vector space over the same field.

Clearly the zero linear transformation of V into W denoted by $0_s = (0 \mid \ldots \mid 0)$ will serve as the zero super vector of linear transformations. We shall denote the collection of linear transformations from the super vector space V into the super vector space W by SL(V, W) which is a super vector space over F, called the linear transformations of the super vector space V into the super vector space W.

Now we have already said the natural dimension of a super vector space is its usual dimension i.e., if $X = (x_1 \mid x_2 \mid \ldots \mid \ldots \mid \ldots \mid x_n)$ then dimension of X is n. So if $X = (A_1 \mid \ldots \mid A_k)$ then $k \leq n$ and if $k < n$ we do not call the natural dimension of X to be k but only as n.

However we cannot say if the super vector space V is of natural dimension n and the super vector space W is of natural dimension m then SL (V, W) is of natural dimension mn.

For we shall first describe how $T_s$ looks like and the way the dimension of SL (V, W) is determined by a simple example.

***Example 1.3.7:*** Let V = $\{(x_1\ x_2\ x_3 \mid x_4\ x_5 \mid x_6\ x_7) \mid x_i \in Q; \mid 1 \leq i \leq 7\}$ be a super vector space over Q. Suppose W = $\{(x_1\ x_2 \mid x_3\ x_4 \mid x_5\ x_6\ x_7\ x_8\ x_9) \mid x_i \in Q; \mid 1 \leq i \leq 9\}$ be a super vector space over Q. Clearly the natural dimension of V is 7 and that of W is 9. Let SL (V, W) denote the super space of all linear transformations from V into W.

Let $T_s : V \rightarrow W$;

$$T_s\ (x_1\ x_2\ x_3 \mid x_4\ x_5 \mid x_6\ x_7) = (x_1 + x_2\ x_2 + x_3 \mid x_4\ x_5 \mid x_6\ 0x_7\ 0x_6)$$

i.e., $T_s = (T_1 \mid T_2 \mid T_3)$ such that $T_1(x_1\ x_2\ x_3) = (x_1 + x_2, x_2 + x_3)$, $T_2\ (x_4, x_5) = (x_4, x_5)$ and $T_3\ (x_6, x_7) = (x_6, 0, x_7, 0, x_6)$.

Clearly natural dimension of $T_1$ is 6, the natural dimension of $T_2$ is 4 and that of $T_3$ is 10. Thus the natural dimension of SL (V, W) is 20. But we see the natural dimension of V is n = 7 and that of W is 9 and the natural dimension of L (V, W) is 63, when V and W are just vector spaces. But when V and W are super vector spaces of natural dimension 7 and 9, the dimension of SL(V, W) is 20. Thus we see the linear transformation of super vector spaces lessens the dimension of SL (V, W).

We also see that the super dimension of SL (V, W) is not unique even if the natural dimension of V and W are fixed, They vary according to the length of the row vectors in the super vector $\alpha = (A_1 \mid \ldots \mid A_k)$; k < n, i.e., they are dependent on the partition of the row vectors.

This is also explained by the following example.

***Example 1.3.8:*** Let V = $\{(x_1\ x_2\ x_3\ x_4 \mid x_5 \mid x_6\ x_7) \mid x_i \in Q; 1 \leq i \leq 7\}$ be a super vector space over Q and W = $\{(x_1\ x_2\ x_3 \mid x_4\ x_5 \mid x_6\ x_7\ x_8\ x_9) \mid x_i \in Q; \mid 1 \leq i \leq 9\}$ be a super vector space over Q. Clearly the natural dimension of V is 7 and that of W is

9. Now let SL(V, W) be the set of all linear transformation of V into W.

Now if $T_s \in$ SL (V, W) then $T_s = (T_1 \mid T_2 \mid T_3)$, where dimension of $T_1$ is 12, dimension of $T_2$ is 2 and dimension of $T_3$ is 8. The super dimension of SL(V, W) is 12 + 2 + 8 = 22. Thus it is not 63 and this dimension is different from that given in example 1.3.7 which is just 20.

***Example 1.3.9:*** Let V = $\{(x_1\ x_2\ x_3\ x_4\ x_5 \mid x_6 \mid x_7) \mid x_i \in Q;\ 1 \le i \le 7\}$ be a super vector space over Q.

Let W = $\{(x_1\ x_2\ x_3 \mid x_4\ x_5\ x_6 \mid x_7\ x_8\ x_9) \mid x_i \in Q; \mid 1 \le i \le 9\}$ be a super vector space over Q. Clearly the natural dimension of V is 7 and that of W is 9.

Now let SL(V, W) be the super vector space of linear transformations of V into W. Let $T_s = (T_1 \mid T_2 \mid T_3) \in$ SL(V, W) dimension of $T_1$ is 15 dimension of $T_2$ is 3 and that of $T_3$ is 3. Thus the super dimension of SL(V, W) is 21.

Now we can by using number theoretic techniques find the minimal dimension of SL(V, W) and the maximal dimension of SL(V, W). Also one can find how many distinct super vector spaces of varied dimension is possible given the natural dimension of V and W.

These are proposed as open problems is the last chapter of this book.

***Example 1.3.10*:** Given V = $\{(x_1\ x_2 \mid x_3\ x_4 \mid x_5\ x_6) \mid x_i \in Q;\ 1 \le i \le 6\}$ is a super vector space over Q. W = $\{(x_1\ x_2\ x_3 \mid x_4\ x_5 \mid x_6) \mid x_i \in Q;\ 1 \le i \le 6\}$ is also a super vector space over Q. Both have the natural dimension to be 6. SL(V, W) be the super vector space of all linear transformation from V into W. The super dimension of SL(V, W) is 12.

Suppose in the same example V = $\{(x_1\ x_2 \mid x_3\ x_4\ x_5 \mid x_6) \mid x_i \in Q;\ 1 \le i \le 6\}$ a super vector space over Q and W = $\{(x_1\ x_2 \mid x_3 \mid x_4\ x_5\ x_6) \mid x_i \in Q;\ 1 \le i \le 6\}$ a super vector space over Q. Let SL(V, W) be the super vector space of linear transformations from V into W. The super dimension of SL(V, W) is 10.

Suppose V = {($x_1$ | $x_2$ $x_3$ $x_4$ | $x_5$ $x_6$) | $x_i \in Q$; $1 \leq i \leq 6$} a super vector space over Q and W = {($x_1$ $x_2$ $x_3$ | $x_4$ $x_5$ | $x_6$) | $x_i \in Q$; $1 \leq i \leq 6$} a super vector space over Q. Let SL(V, W) be the super vector space of all linear transformation from V into W.

The natural dimension of SL(V, W) is 11. Thus we have seen that SL(V, W) is highly dependent on the way the row vectors are partitioned and we have different natural dimensions for different partitions.

So we make some more additions in the definitions of super vector spaces.

Let V = {($x_1$ $x_2$ | … | … | $x_n$) | $x_i \in F$; F a field; $1 \leq i \leq n$} be a super vector space over F. If V = {($A_1$ | $A_2$ | … | $A_k$) | $A_i$ are row vectors with entries from the field F; i = 1, 2, …, k, $k \leq n$} Suppose the number of elements in $A_i$ is $n_i$; i = 1, 2, …, n then we see natural dimension of V is $n = n_1 + \ldots + n_k$ and is denoted by ($n_1$, …, $n_k$).

Let W = {($x_1$ | … | … | $x_m$) | $x_i \in F$, F a field $i \leq i \leq m$) be a super vector space over the field F of natural dimension m. Let W = {($B_1$ | … | $B_k$), $k \leq m$; $B_i$'s row vectors with entries from F; i = 1, 2, …, k}. Then natural dimension of W is $m = m_1 + \ldots + m_k$ where $m_i$ is the number of elements in the row vector $B_i$, $1 \leq i \leq k$.

Now the collection of all linear transformations from V into W be denoted by SL(V, W) which is again a super vector space over F. Now the natural dimension of SL(V, W) = $m_1n_1 + \ldots + m_kn_k$ clearly $m_1n_1 + \ldots + m_kn_k \leq mn$.

Now we state this in the following theorem.

**THEOREM 1.3.4:** *Let V = {($x_1$ | … | $x_n$) / $x_i \in F$; $i \leq i \leq n$} be a super vector space over F of natural dimension n, where V = {($A_1$ | … | $A_k$) | $A_i$'s are row vectors of length $n_i$ and entries of $A_i$ are from F, $i \leq i \leq k$, $k \leq n$ | $n_1 + \ldots + n_k = n$}. W = {($x_1$ | … | … | … | $x_m$) / $x_i \in F$, $1 \leq i \leq m$} is a super vector space of natural dimension m over the field F, where W = {($B_1$ | … | $B_k$) | $B_i$'s are row vectors of length $m_i$ with entries from F, $i \leq i \leq k$, $k < m$ such that $m_1 + \ldots + m_k = m$}. Then the super vector space SL(V,*

*W) of all linear transformations from V into W is finite dimensional and has dimension $m_1n_1 + m_2n_2 + \ldots + m_kn_k \leq mn$.*

*Proof:* Let B = $\{\alpha_1, \ldots, \alpha_n\}$ and $B^1 = \{\beta_1, \ldots, \beta_m\}$ be a basis for V and W respectively where each $\alpha_i$ and $\beta_j$ are super row vectors in V and W respectively; $1 \leq i \leq n$ and $1 \leq j \leq m$. For each pair of integers $(p_i, q_i)$ with $1 \leq p_i \leq m_i$ and $1 \leq q_i \leq n_i$; i =1, 2, …, k.

We define a linear transformation

$$E^{p_i,q_i}(\alpha_t) = \begin{cases} 0 & \text{if } t \neq q_i \\ \beta_{p_i} & \text{if } t = q_i \end{cases}$$

$$= \delta_{t_{q_i}} \beta_{p_i} \text{ for } i = 1, 2, 3, \ldots, k.$$

Thus $E^{p,q} = [E^{p_1,q_1} \mid \ldots \mid E^{p_k,q_k}] \in SL(V, W)$. $E^{p,q}$ is a linear transformation from V into W. From earlier results each $E^{p_i,q_i}$ is unique so $E^{p,q}$ is unique and by properties for vector spaces that the $m_in_i$ transformations $E^{p_i,q_i}$ form a basis for $L(A_i, B_i)$. So dimension of SL(V, W) is $m_1n_1 + \ldots + m_kn_k$.

Now we proceed on to define the new notion of linear operator on a super vector space V i.e., a linear transformation from V into V.

**DEFINITION 1.3.5:** *Let $V = \{(A_1 \mid \ldots \mid A_k) \mid A_i$ is a row vector with entries from a field F with number of elements in $A_i$ to be $n_i$; $i = 1, 2, \ldots, k\} = \{(x_1 \mid \ldots \mid \ldots \mid \ldots \mid x_n) \mid x_i \in F$; $i = 1, 2, \ldots, n\}$; $k \leq n$ and $n_1 + \ldots + n_k = n$; be a super vector space over the field F. A linear transformation $T = (T_1 \mid T_2 \mid \ldots \mid T_k)$ from V into V is called the linear operator on V.*

*Let SL(V, V) denote the set of all linear operators from V to V, the dimension of SL(V, V) = $n_1^2 + \ldots + n_k^2 \leq n^2$.*

**LEMMA 1.3.1:** *Let V be a super vector space over the field F, let $U_1^s$ $T_1^s$ and $T_2^s$ be linear operators on V; let c be an element of F.*

*(a)* $I^sU^s = U^sI^s = U^s$

*(b)* $U^s(T_1^s + T_2^s) = U^sT_1^s + U^sT_2^s$;

$(T_1^s + T_2^s)U^s = T_1^sU^s + T_2^sU^s$

*(c)* $c(U^sT_1^s) = (cU^s)\,T_1^s = U(cT_1^s)$.

*Proof:* Given V = {$(A_1 \mid \ldots \mid A_k) \mid A_i$ are row vectors with entries from the field F}. Let $U^s = (U_1 \mid \ldots \mid U_k)$, $T_1^s = (T_1^1 \mid \ldots \mid T_1^k)$ and $T_2^s = (T_2^1 \mid \ldots \mid T_2^k)$ and $I_s = (I_1 \mid \ldots \mid I_k)$ ($I = I^s = I_s$, the identity operator) be linear operators from V into V. Now

$$\begin{aligned} IU^s &= (I_1 \mid \ldots \mid I_k)\,(U_1 \mid \ldots \mid U_k) \\ &= (I_1U_1 \mid \ldots \mid I_kU_k) \\ &= (U_1I_1 \mid \ldots \mid U_kI_k) \\ &= (U_1 \mid \ldots \mid U_k)\,(I_1 \mid I_2 \mid \ldots \mid I_k). \end{aligned}$$

$$\begin{aligned} &U^s\,(T_2^s + T_2^s) \\ &= (U_1 \mid \ldots \mid U_k)\,[(T_2^1 \mid \ldots \mid T_1^k) \quad + (T_2^1 \mid \ldots \mid T_2^k)] \\ &= [U_1\,(T_1^1 + T_2^1) \mid \ldots \mid U_k\,(T_1^k + T_2^k)] \\ &= [(U_1T_1^1 + U_1\,T_2^1) \mid \ldots \mid U_k\,T_1^k + U_kT_2^k] \\ &= (U_1T_1^1 \mid \ldots \mid U_kT_1^k) + (U_1T_2^1 \mid \ldots \mid U_kT_2^k) \\ &= U^sT_1^s + U^sT_2^s. \end{aligned}$$

On similar lines one can prove $(T_1^s + T_2^s)\,U^s = T_1^sU^s + T_2^sU^s$.

(c) To prove $c(U^sT_1^s)$

$$= (cU^s)(T_1^s) = U^s(cT_1^s).$$

i.e., $c[(U_1 \mid \ldots \mid U_k)\,(T_1^1 \mid \ldots \mid T_1^k)]$

$= c(U_1T_1^1 | \ldots | U_kT_1^k) \quad = c(U_1T_1^1 | \ldots | cU_kT_1^k)$.

Now

$$(cU^s)T_1^s = (cU_1 | \ldots | (cU_k)\,(T_1^1 | \ldots | T_1^k) = (c_1U_1T_1^1 | \ldots | cU_kT_1^k]$$

So

$$c(U^sT_1^s) = (cU^s)\,T_1^s \ldots \text{ I}$$

Consider

$$\begin{aligned} U^s(cT_1^s) &= (U_1 | \ldots | U_k)\,(cT_1^1 | \ldots | cT_1^k) \\ &= (U_1c_1T_1^1 | \ldots | Uc_kT_1^k) \\ &= (cU_1T_1^1 | \ldots | cU_kT_1^k); \end{aligned}$$

from I we see

$$c(U^sT_1^s) = c(U^s)T_1^s = U^s(cT_1^s).$$

We call SL(V, V) a super linear algebra. However we will define this concept elaborately.

***Example 1.3.11:*** Let V = {($x_1$ $x_2$ $x_3$ | $x_4$ | $x_5$ $x_6$ | $x_7$ $x_8$) | $x_i \in$ Q; $1 \le i \le 8$} be a super vector space over Q. Let SL(V, V) denote the collection of all linear operators from V into V.

We see the super dimension of SL(V, V) is 18, not 64 as in case of L($V_1$, $V_1$) where $V_1$ = {($x_1$, …, $x_8$) | $x_i \in$ Q; $1 \le i \le 8$} is a vector space over Q of dimension 8. When V is a super vector space of natural dimension 8 but SL(V, V) is of dimension 18.

These concepts now leads us to define the notion of general super vector spaces. For all the while we were only defining super vector spaces specifically only when the elements were super row vectors or super matrices and we have only studied their properties now we proceed on to define the notion of general super vector spaces.

The super vector spaces using super row vectors and super matrices were first introduced mainly to make the reader how they function. The functioning of them was also illustrated by examples and further many of the properties were derived when the super vector spaces were formed using the super row vectors. However when the problem of linear transformation of super vector spaces was to be carried out one faced with some simple problems however one can also define linear transformation of super vector spaces by not disturbing the partitions or by preserving the partition but the elements with in the partition which are distinct had to be changed or defined depending on the elements in the partitions of the range space.

However the way of defining them in case of super row vectors remain the same only changes come when we want to speak of SL(V, W) and SL(V, V).

They are super vector spaces in that case how the elements should look like only at this point we have to make necessary changes, with which the super vector space status is maintained however it affects the natural dimension which have to be explained.

**DEFINITION 1.3.6:** *Let $V_1, \ldots, V_n$ be n vector spaces of finite dimensions defined over a field F. $V = (V_1 \mid V_2 \mid \ldots \mid \ldots \mid V_n)$ is called the super vector space over F. Since we know if $V_i$ is any vector space over F of dimension say $n_i$ then $V \cong F^{n_i}$. = $\{(x_1, \ldots, x_{n_i}); x_i \in F; 1 \leq i \leq n_i\}$ Thus any vector space of any finite dimension can always be realized as a row vector with the number of elements in that row vector being the dimension of the vector space under consideration. Thus if $n_1, \ldots, n_n$ are the dimensions of vector spaces $V_1, \ldots, V_n$ over the field F then $V \cong (F^{n_1} \mid F^{n_2} \mid \ldots \mid F^{n_n})$ which is a collection of super row vectors, hence V is nothing but a super vector space over F.*

Thus this definition is in keeping with the definition of super vector spaces.

Thus without loss of generality we will for the convenience of notations identify a super vector space elements only by a super row vector.

Now we can give examples of a super vector spaces.

**DEFINITION 1.3.7:** *Let $V = (V_1 \mid \ldots \mid V_n)$ be a super vector space over the field F. Let $n_i$ be the dimension of the vector space $V_i$ over F, i = 1, 2, …, n; then the dimension of V is $n_1 + \ldots + n_n$, we call this as the natural dimension of the super vector space V. Thus if $V = (V_1 \mid \ldots \mid V_n)$ is a super vector space of dimension $n_1 + \ldots + n_n$ over the field F, then we can say $V = (F^{n_1} \mid \ldots \mid F^{n_n})$.*

***Example 1.3.12:*** Let $V = (V_1 \mid V_2 \mid V_3)$ be a super vector space over Q, where $V_1$ = {set of all 2 × 2 matrices with entries from Q}. $V_1$ is a vector space of dimension 4 over Q.

$V_2$ = {All polynomials of degree less than or equal to 5 with coefficients from Q}; $V_2$ is a vector space of dimension 6 over Q and $V_3$ = {set of all 3 × 4 matrices with entries from Q}; $V_3$ is a vector space of dimension 12 over Q. Clearly

$$V = (V_1 \mid V_2 \mid V_3) \cong (Q^4 \mid Q^5 \mid Q^6)$$
$$= \{(x_1\ x_2\ x_3\ x_4 \mid x_5 \ldots x_{10} \mid x_{11} \ldots x_{22}) \mid x_i \in Q \mid 1 \leq i \leq 22\}$$

is nothing but a collection of super row vectors, with natural dimension 22.

Now we proceed on to give a representation of transformations from finite dimensional super vector space V into W by super matrices.

Let V be a super vector space of natural dimension n given by $V = \{(A_1 \mid \ldots \mid A_k) \mid A_i$ is a row vector with entries from the field F of length $n_i$; i = 1, 2, …, k and $n_1 + n_2 + \ldots + n_k = n\}$and let $W = \{(B_1 \mid \ldots \mid B_k) \mid B_i$ is a row vector with entries from the field F of length $m_i$, i = 1, 2, …, k and $m_i + \ldots + m_k = m\}$, where W is a super vector space of natural dimension m over F.

Let $B = \{\alpha_1, \ldots, \alpha_n\}$ be a basis for V where $\alpha_j$ is a super row vector and $B' = \{\beta_1, \ldots, \beta_m\}$ be a basis for W where $\beta_i$ is a super row vector. Let $T_s$ be any linear transformation from V into W, then $T_s$ is determined by its action on the super vectors $\alpha_j$ each of the n super vectors $T_s\, \alpha_j$ is uniquely expressible as a linear combination $T_s\, \alpha_j = \sum_{i=1}^{m} A_{ij}\beta_i$

Here $T_s = (T_1 \mid \ldots \mid T_k)$, k < n.

$$(T_1\alpha_j \mid \ldots \mid T_k\alpha_j) = \left(\sum_{i=1}^{m_1} A_{ij}^{n_1}\beta_i \mid \ldots \mid \alpha_i \sum_{i=1}^{m_k} A_{ij}^{n_k}\beta_i\right)$$ of the super row vector $\beta_i$; $1 \leq i \leq m$; the scalars $A_{ij}^{n_i} \ldots A_{m_i j}^{n_i}$ being the coordinates of $T_i\alpha_j$ in the ordered basis B'; true for i = 1, 2, …, k. Thus the transformation $T_i$ is determined by the $m_i n_i$ scalars $A_{ij}^{n_i}$. The $m_i \times n_i$ matrix $A^i$ defined by A (i, j) = $A_{ij}^{n_i}$ is called the matrix of $T_i$ relative to the pair in ordered basis B and B'. This is true for every i.

Thus the transformation super matrix is a m × n super matrix A given by

$$A = \left(\begin{array}{c|c|c|c} A^1_{m_1 x n_1} & 0 & 0 & 0 \\ \hline 0 & A^2_{m_2 x n_2} & 0 & 0 \\ \hline 0 & 0 & \ddots & 0 \\ \hline 0 & 0 & 0 & A^k_{m_k x n_k} \end{array}\right).$$

With this related super matrix with entries from the field F; one can understand how the transformation takes place.

This will be explicitly described by examples. Clearly the natural order of this m × n matrix is $m_n \times n_1 + \ldots + m_k \times n_k$.

***Example 1..3.13:*** Let

$V = \{(x_1\ x_2\ x_3 \mid x_4\ x_5 \mid x_6\ x_7\ x_8\ x_9) \mid x_i \in Q \mid 1 \leq i \leq 9\}$ be a super vector space over Q.

$W = \{(x_1\ x_2 \mid x_3\ x_4\ x_5 \mid x_6\ x_7) \mid x_i \in Q \mid 1 \leq i \leq 7\}$ be a super vector space over Q.

Let SL(V, W) denote the set of all linear transformations from V into W.

Consider the 7 × 9 super matrix

$$A = \left(\begin{array}{ccc|cc|cccc} 1 & 0 & 0 & 0 & 0 & 0 & 0 & 0 & 0 \\ 1 & 0 & 2 & 0 & 0 & 0 & 0 & 0 & 0 \\ \hline 0 & 0 & 0 & 2 & 1 & 0 & 0 & 0 & 0 \\ 0 & 0 & 0 & 0 & -1 & 0 & 0 & 0 & 0 \\ 0 & 0 & 0 & 0 & 1 & 0 & 0 & 0 & 0 \\ \hline 0 & 0 & 0 & 0 & 0 & 1 & 0 & 1 & 0 \\ 0 & 0 & 0 & 0 & 0 & 0 & 2 & 0 & 1 \end{array}\right)$$

gives the associated linear transformation $T_s(x_1\ x_2\ x_3 |\ x_4\ x_5\ |\ x_6\ x_7\ x_8\ x_9)$

$$= (T_1 \mid T_2 \mid T_3)\ \ [x_1\ x_2\ x_3 |\ x_4\ x_5\ |\ x_6\ x_7\ x_8\ x_9]$$

$$= [T_1\ (x_1\ x_2\ x_3) |\ T_2\ (x_4\ x_5)\ |\ T_3\ (x_6\ x_7\ x_8\ x_9)]$$

$$= [x_1, x_2 + 2x_3 |\ 2x_4 + x_5, -x_5, x_5\ |\ x_6 + x_8, 2x_7 + x_9] \in W.$$

Thus we see as incase of usual vector spaces to every linear transformation from V into W, we have an associated super matrix whose non diagonal terms are zero and diagonal matrices give the components of the transformation $T_s$. Here also ',' is put in the super vector for the readers to understand the transformation, by a default of notation.

We give yet another example so that the reader does not find it very difficult to understand when this notion is described abstractly.

***Example 1.3.14:*** Let

$V = \{(x_1\ x_2\ |\ x_3\ x_4\ |\ x_5\ x_6\ x_7\ |\ x_8\ x_9) \mid x_i \in Q; 1 \leq i \leq 9\}$

be a super vector space over Q and

$W = \{(x_1\ x_2\ x_3\ |\ x_4\ |\ x_5\ x_6\ x_7\ |\ x_8\ x_9\ x_{10}) \mid x_i \in Q; 1 \leq i \leq 10\}$

be another super vector over Q. Let SL (V, W) be the super vector space over Q.

Consider the 10 × 9 super matrix

$$A = \left(\begin{array}{cc|cc|ccc|cc}
1 & 0 & 0 & 0 & 0 & 0 & 0 & 0 & 0 \\
1 & 2 & 0 & 0 & 0 & 0 & 0 & 0 & 0 \\
0 & 1 & 0 & 0 & 0 & 0 & 0 & 0 & 0 \\
\hline
0 & 0 & 1 & -1 & 0 & 0 & 0 & 0 & 0 \\
\hline
0 & 0 & 0 & 0 & 1 & 0 & -1 & 0 & 0 \\
0 & 0 & 0 & 0 & 0 & 1 & 0 & 0 & 0 \\
0 & 0 & 0 & 0 & 0 & 1 & 1 & 0 & 0 \\
\hline
0 & 0 & 0 & 0 & 0 & 0 & 0 & 1 & 1 \\
0 & 0 & 0 & 0 & 0 & 0 & 0 & 0 & 1 \\
0 & 0 & 0 & 0 & 0 & 0 & 0 & 1 & 0
\end{array}\right).$$

The transformation $T_s$: V → W associated with A is given by
$T_s\ (x_1\ x_2 \mid x_3\ x_4 \mid x_5\ x_6\ x_7 \mid x_8\ x_9)$
$= (x_1,\ 1 + 2x_2,\ x_2 \mid x_3 - x_4 \mid x_5 - x_7,\ x_6,\ x_6 + x_7 \mid x_8 + x_9,\ x_9,\ x_8)$.
Thus to every linear transformation $T_s$ of V into W we have a super matrix associated with it and conversely with every appropriate super matrix A we have a linear transformation $T_s$ associated with it.
Thus SL (V, W) can be described as

$$\left\{\left(\begin{array}{cc|cc|ccc|cc}
a_1 & a_2 & 0 & 0 & 0 & 0 & 0 & 0 & 0 \\
a_3 & a_4 & 0 & 0 & 0 & 0 & 0 & 0 & 0 \\
a_5 & a_6 & 0 & 0 & 0 & 0 & 0 & 0 & 0 \\
\hline
0 & 0 & a_7 & a_8 & 0 & 0 & 0 & 0 & 0 \\
\hline
0 & 0 & 0 & 0 & a_9 & a_{10} & a_{11} & 0 & 0 \\
0 & 0 & 0 & 0 & a_{12} & a_{13} & a_{14} & 0 & 0 \\
0 & 0 & 0 & 0 & a_{15} & a_{16} & a_{17} & 0 & 0 \\
\hline
0 & 0 & 0 & 0 & 0 & 0 & 0 & a_{18} & a_{19} \\
0 & 0 & 0 & 0 & 0 & 0 & 0 & a_{20} & a_{21} \\
0 & 0 & 0 & 0 & 0 & 0 & 0 & a_{22} & a_{23}
\end{array}\right)\right.$$

such that $a_i \in Q$; $1 \leq i \leq 23$}.
Thus the dimension of SL(V, W) is $3 \times 2 + 1 \times 2 + 3 \times 3 + 3 \times 2 = 6 + 2 + 9 + 6 = 23$.

Now we give the general working for the fact SL (V,W) is isomorphic to diagonal m × n super matrices. Before we go for deep analysis we just give a few examples of what we mean by a super diagonal matrix.

***Example 1.3.15:*** Let

$$\left(\begin{array}{cc|cc|ccc|cc} 8 & 1 & 0 & 0 & 0 & 0 & 0 & 0 & 0 \\ 6 & 7 & 0 & 0 & 0 & 0 & 0 & 0 & 0 \\ \hline 0 & 0 & 5 & 6 & 0 & 0 & 0 & 0 & 0 \\ \hline 0 & 0 & 0 & 0 & 7 & 1 & 0 & 0 & 0 \\ 0 & 0 & 0 & 0 & 6 & 8 & -1 & 0 & 0 \\ \hline 0 & 0 & 0 & 0 & 0 & 0 & 0 & 3 & 1 \\ 0 & 0 & 0 & 0 & 0 & 0 & 0 & 6 & 7 \\ 0 & 0 & 0 & 0 & 0 & 0 & 0 & 6 & 0 \end{array}\right)$$

be a 8 × 9 super matrix.
We call A the super diagonal matrix as only the diagonal matrices are non zero and rest of the matrices and zero. It is important to mention here that in a super diagonal matrix we do not need the super matrix to be a square matrix; it can be any matrix expect a super row matrix or super column matrix.

Thus we can say if $A = \left(\begin{array}{c|c|c|c} A_{11} & A_{12} & \cdots & A_{1n} \\ \hline A_{21} & A_{22} & \cdots & A_{2n} \\ \hline & & \cdots & \\ \hline A_{n1} & A_{n2} & \cdots & A_{nn} \end{array}\right)$

where $A_{ij}$ are simple matrices we say A is a super diagonal matrix if $A_{11}, A_{22}, \ldots, A_{nn}$ are non zero matrices and $A_{ij}$ is a zero matrix if $i \neq j$.

The only demand we place is that the number of row partitions of A is equal to the number of column partitions of A.

***Example 1.3.16:*** Let A be a super diagonal matrix given by

$$A = \left(\begin{array}{cccc|c|ccc|ccccc}
9 & 0 & 1 & 2 & 0 & 0 & 0 & 0 & 0 & 0 & 0 & 0 & 0 \\
2 & 1 & 0 & 0 & 0 & 0 & 0 & 0 & 0 & 0 & 0 & 0 & 0 \\
5 & 6 & 1 & 0 & 0 & 0 & 0 & 0 & 0 & 0 & 0 & 0 & 0 \\
\hline
0 & 0 & 0 & 0 & 5 & 0 & 0 & 0 & 0 & 0 & 0 & 0 & 0 \\
\hline
0 & 0 & 0 & 0 & 0 & 1 & 1 & 1 & 0 & 0 & 0 & 0 & 0 \\
0 & 0 & 0 & 0 & 0 & 2 & 0 & 0 & 0 & 0 & 0 & 0 & 0 \\
0 & 0 & 0 & 0 & 0 & 5 & 0 & 7 & 0 & 0 & 0 & 0 & 0 \\
\hline
0 & 0 & 0 & 0 & 0 & 0 & 0 & 0 & 1 & 0 & 1 & 1 & 2 \\
0 & 0 & 0 & 0 & 0 & 0 & 0 & 0 & 0 & 1 & 0 & 1 & 0 \\
0 & 0 & 0 & 0 & 0 & 0 & 0 & 0 & 0 & 2 & 1 & 0 & 7
\end{array}\right).$$

We see A is a $10 \times 13$ matrix which is a super diagonal matrix. Only the number of row partitions equals to the number of column partitions equal to 4.

**THEOREM 1.3.5:** *Let $V = \{(x_1\ x_2 \mid \ldots \mid \ldots \mid x_t \ldots x_n) = (A_1 \mid \ldots \mid A_k) \mid x_i \in F$ and $A_i$ is a row vector with entries from the field F; $1 \leq i \leq n$ and $1 \leq t \leq k$; $k \leq n\}$ be a super vector space over F.*

Let $W = \{(x_1\ x_2 \mid \ldots \mid \ldots \mid \ldots \mid x_t \ldots\ x_n) = (B_1 \mid \ldots \mid B_k) \mid x_i \in F$ and $B_i$ is a row vector with entries from the field F with $1 \leq i \leq m$ and $k \leq m$ $1 \leq t \leq k\}$ be a super vector space of same type as V. Let SL (V, W) be the collection of all linear transformations from V into W, SL (V, W) is a super vector space over F and for a set of basis $B = \{\alpha_1, \ldots, \alpha_n\}$and $B^1 = \{\beta_1, \ldots, \beta_m\}$ of V and W respectively. For each linear transformation $T_s$ from V into W there is a $m \times n$ super diagonal matrix A with entries from F such that $T_s \to A$ is a one to one correspondence between the set of all linear transformations from V into W and the set of all m × n super diagonal matrices over the field F.

*Proof:* The super diagonal matrix A associated with $T_s$ is called the super diagonal matrix of $T_s$ relative to the basis B and $B^1$. We know $T_s : (A_1 \mid \ldots \mid A_k) \to (B_1 \mid \ldots \mid B_k)$ where $T_s = (T_1 \mid \ldots \mid T_k)$ and each $T_i$ is a linear transformation from $A_i \to B_i$ where $A_i$ is of dimension $n_i$ and $B_i$ is of dimension $m_i$; i = 1, 2, …, k.
So we have matrix $M^i_j = [T_i\alpha_j]_{[C_i]}$; j = 1, 2, …, $n_i$. $C_i$ a component basis from $B^1$; this is true for i = 1, 2, …, k.
So for any $T_s = (T_1 \mid \ldots \mid T_k)$ and $U_s = (U_1 \mid \ldots \mid U_k)$ in SL (V, W), $cT_s + U_s$ is SL(V, W) for any scalar c in F .
Now $T_i : V \to W$ is such that $T_i(A_j) = (0)$ if $i \neq j$ and $T_i(A_i) = B_i$ and this is true for i = 1, 2, …, k.
Thus the related matrix of $T_s$ is a super diagonal matrix where

$$A = \left(\begin{array}{c|c|c|c} (M_1)_{m_1\times n_1} & 0 & \ldots & 0 \\ \hline 0 & (M_2)_{m_2\times n_2} & 0 & 0 \\ \hline 0 & 0 & \ldots & 0 \\ \hline 0 & 0 & 0 & (M_k)_{m_k\times n_k} \end{array}\right)$$

Thus $M_i$ is a $m_i \times m_i$ matrix associated with the linear transformation $T_i$: $A_i \to B_i$ true for i = 1, 2, …, k. Hence the claim. Likewise we can say that in case of a super vector space $V = \{(x_1 \ldots \mid \ldots \mid \ldots \mid x_t \ldots x_n) \mid x_i \in F$, F a field; $1 \leq i \leq n\} = \{(A_1 \mid \ldots \mid A_k) \mid A_i$ row vectors with entries from the field F; $1 \leq i \leq k\}$over F. We have SL (V, V) is such that there is a one to one correspondence between the n × n super diagonal square matrix with entries from F i.e., SL(V, V) is also a super vector space over F. Further the marked difference between SL(V, W) and SL (V, V) is that SL (V, W) is isomorphic to class of all m × n rectangular super diagonal matrices with entries from F and the diagonal matrices of these super diagonal matrices need not be square matrices but in case of the super vector space SL (V,

V), we have this space to be isomorphic to the collection of all n × n super square diagonal matrices where each of the diagonal matrices are also square matrices.

We will illustrate this situation by a simple example.

***Example 1.3.17:*** Let

$$V = \{(x_1\ x_2\ x_3\ x_4 \mid x_5\ x_6 \mid x_7\ x_8\ x_9 \mid x_{10}) \mid x_i \in Q;\ 1 \le i \le 10\}$$

be a super vector space over Q.

Let SL (V, V) denote the set of all linear operators from V into V. Let $T_s$ be a linear operator on V. Then let A be the super diagonal square matrix associated with $T_s$,

$$A = \left(\begin{array}{cccc|cc|ccc|c}
1 & 0 & 1 & 1 & 0 & 0 & 0 & 0 & 0 & 0 \\
0 & 1 & 1 & 0 & 0 & 0 & 0 & 0 & 0 & 0 \\
1 & 0 & 0 & 1 & 0 & 0 & 0 & 0 & 0 & 0 \\
1 & 0 & 1 & 0 & 0 & 0 & 0 & 0 & 0 & 0 \\
\hline
0 & 0 & 0 & 0 & 1 & 5 & 0 & 0 & 0 & 0 \\
0 & 0 & 0 & 0 & 0 & 2 & 0 & 0 & 0 & 0 \\
\hline
0 & 0 & 0 & 0 & 0 & 0 & 1 & 0 & 2 & 0 \\
0 & 0 & 0 & 0 & 0 & 0 & 0 & 1 & 2 & 0 \\
0 & 0 & 0 & 0 & 0 & 0 & 2 & 0 & 1 & 0 \\
\hline
0 & 0 & 0 & 0 & 0 & 0 & 0 & 0 & 0 & 3
\end{array}\right).$$

Clearly A is a 10 × 10 super square matrix. The diagonal matrices are also square matrices.

Now $T_s\ (x_1\ x_2\ x_3\ x_4 \mid x_5\ x_6 \mid x_7\ x_8\ x_9 \mid x_{10})$

$$= (x_1 + x_3 + x_4,\ x_2 + x_3,\ x_1 + x_4,\ x_1 + x_3 \mid x_5 + 5x_6,\ 2x_6 \mid x_7 + 2x_9,\ x_8 + 2x_9,\ 2x_7 + x_9 \mid 3x_{10}) \in V.$$

Thus we see in case of linear operators $T_s$ of super vector spaces the associated super matrices of $T_s$ is a square super diagonal matrix whose diagonal matrices are also square matrices.

We give yet another example before we proceed on to work with more properties.

***Example 1.3.18:*** Let
$V = (x_1 x_2 \mid x_3 x_4 x_5 \mid x_6 x_7 \mid x_8\ x_9 \mid x_{10} x_{11} x_{12}) \mid x_i \in Q; 1 \leq i \leq 12\}$
be a super vector space over Q. V = {($A_1$ | $A_2$ | $A_3$ | $A_4$ | $A_5$) | $A_i$ are row vectors with entries from Q; 1 ≤ i ≤ 5}. Let us consider a super diagonal 12 × 12 square matrix A with (2 × 2, 3 × 3, 2 × 2, 2 × 2, 3 × 3) ordered diagonal matrices with entries from Q.

$$A = \left(\begin{array}{cc|ccc|cc|cc|ccc}
1 & 1 & 0 & 0 & 0 & 0 & 0 & 0 & 0 & 0 & 0 & 0 \\
2 & 1 & 0 & 0 & 0 & 0 & 0 & 0 & 0 & 0 & 0 & 0 \\ \hline
0 & 0 & 1 & 2 & 0 & 0 & 0 & 0 & 0 & 0 & 0 & 0 \\
0 & 0 & 3 & 0 & 1 & 0 & 0 & 0 & 0 & 0 & 0 & 0 \\
0 & 0 & 0 & 1 & -4 & 0 & 0 & 0 & 0 & 0 & 0 & 0 \\ \hline
0 & 0 & 0 & 0 & 0 & 2 & 0 & 0 & 0 & 0 & 0 & 0 \\
0 & 0 & 0 & 0 & 0 & 1 & 5 & 0 & 0 & 0 & 0 & 0 \\ \hline
0 & 0 & 0 & 0 & 0 & 0 & 0 & 1 & 2 & 0 & 0 & 0 \\
0 & 0 & 0 & 0 & 0 & 0 & 0 & 2 & 1 & 0 & 0 & 0 \\ \hline
0 & 0 & 0 & 0 & 0 & 0 & 0 & 0 & 0 & 1 & 2 & 3 \\
0 & 0 & 0 & 0 & 0 & 0 & 0 & 0 & 0 & 3 & 1 & 2 \\
0 & 0 & 0 & 0 & 0 & 0 & 0 & 0 & 0 & 2 & 3 & 1
\end{array}\right).$$

The linear transformation associated with A is given by

$T_s\ (x_1 x_2 \mid x_3 x_4 x_5 \mid x_6 x_7 \mid x_8\ x_9 \mid x_{10} x_{11} x_{12}) =$
$(x_1 + x_2, 2x_1 + x_2 \mid x_3 + 2x_4, 3x_3 + x_5, x_4 - 4x_5 \mid 2x_6, x_6 + 5x_7 \mid$
$x_8 + 2x_9, 2x_8 + x_9 \mid x_{10} + 2x_{11} + 3x_{12}, 3x_{10} + x_{11} + 2x_{12},\quad 2x_{10} + 3x_{11} + x_{12})$.
Thus we can say given an appropriate super diagonal square matrix with entries from Q we have a linear transformation $T_s$ from V into V and conversely given any $T_s \in$ SL(V, V) we have a square super diagonal matrix associated with it. Hence we can say SL(V, V) =

$$
\left\{A = \left(\begin{array}{cc|ccc|cc|cc|ccc}
a_1 & a_2 & 0 & 0 & 0 & 0 & 0 & 0 & 0 & 0 & 0 & 0 \\
a_3 & a_4 & 0 & 0 & 0 & 0 & 0 & 0 & 0 & 0 & 0 & 0 \\
\hline
0 & 0 & a_5 & a_6 & a_7 & 0 & 0 & 0 & 0 & 0 & 0 & 0 \\
0 & 0 & a_8 & a_9 & a_{10} & 0 & 0 & 0 & 0 & 0 & 0 & 0 \\
0 & 0 & a_{11} & a_{12} & a_{13} & 0 & 0 & 0 & 0 & 0 & 0 & 0 \\
\hline
0 & 0 & 0 & 0 & 0 & a_{14} & a_{15} & 0 & 0 & 0 & 0 & 0 \\
0 & 0 & 0 & 0 & 0 & a_{16} & a_{17} & 0 & 0 & 0 & 0 & 0 \\
\hline
0 & 0 & 0 & 0 & 0 & 0 & 0 & a_{18} & a_{19} & 0 & 0 & 0 \\
0 & 0 & 0 & 0 & 0 & 0 & 0 & a_{20} & a_{21} & 0 & 0 & 0 \\
\hline
0 & 0 & 0 & 0 & 0 & 0 & 0 & 0 & 0 & a_{22} & a_{23} & a_{24} \\
0 & 0 & 0 & 0 & 0 & 0 & 0 & 0 & 0 & a_{25} & a_{26} & a_{27} \\
0 & 0 & 0 & 0 & 0 & 0 & 0 & 0 & 0 & a_{28} & a_{29} & a_{30}
\end{array}\right)\right.
$$

such that $a_i \in Q$; $1 \le i \le 30\}$.
We see the dimension of SL (V, V) = $2^2 + 3^2 + 2^2 + 2^2 + 3^2 = 30$.
Thus we can say if V = $\{(A_1 \mid \ldots \mid A_k) \mid A_i$ is a row vector with entries from a field F and each $A_i$ is of length $n_i$, $1 \le i \le k\}$; V is a super vector space over F; then if $T_s \in SL(V, V)$ then we have an associated A, where A is a $(n_1 + n_2 + \ldots + n_k) \times (n_1 + n_2 + \ldots + n_k)$ square diagonal matrix and dimension of SL (V, V) is $n_1^2 + n_2^2 + \ldots + n_k^2$.
i.e.,

$$
A = \left(\begin{array}{c|c|c|c}
(A_1)_{n_1 \times n_1} & 0 & 0 & 0 \\
\hline
0 & (A_2)_{n_2 \times n_2} & & 0 \\
\hline
0 & 0 & & 0 \\
\hline
0 & 0 & 0 & (A_k)_{n_k \times n_k}
\end{array}\right).
$$

Now we define when is a linear operator from V into V invertible before we proceed onto define the notion of super linear algebras.

**DEFINITION 1.3.8:** *Let $V = \{(A_1 \mid \ldots \mid A_k) \mid A_i$ are row vectors with entries from the field F with length of each $A_i$ to be $n_i$; $i = 1, 2, \ldots, k\}$ be a super vector space over F of dimension $n_1 + \ldots + n_k = n$. Let $T_s$: $V \to V$ be a linear operator on V. We say $T_s = (T_1 \mid \ldots \mid T_k)$ is invertible if their exists a linear operator $U_s$ from V into V such that UT is the identity function of V and TU is also the identity function on V.*

*If $T_s = (T_1 \mid \ldots \mid T_k)$ is invertible implies each $T_i : A_i \to A_i$ is also invertible and $U_s = (U_1 \mid \ldots \mid U_k)$ is denoted by $T_s^{-1} = (T_1^{-1} \mid \ldots \mid T_k^{-1})$.*

*Thus we can say $T_s$ is invertible if and only if $T_s$ is one to one i.e., $T_s\alpha = T_s\beta$ implies $\alpha = \beta$.*

*$T_s$ is onto that is range of $T_s$ is all of V.*

Now if $T_s$ is an invertible linear operator on V and if A is the associated square super diagonal matrix with entries from F then each of the diagonal matrix $M_1, \ldots, M_k$ are invertible matrices i.e., if

$$A = \left(\begin{array}{c|c|c|c} M_1 & 0 & \ldots & 0 \\ \hline 0 & M_2 & 0 & 0 \\ \hline 0 & 0 & \ldots & 0 \\ \hline 0 & 0 & 0 & M_k \end{array}\right)$$

then each $M_i$ is an invertible matrix. So we say A is also an invertible super square diagonal matrix.

It is pertinent to make a mention here that every $T_s$ in SL(V, V) need not be an invertible linear transformation from V into V.

Now we proceed on to define the notion of super linear algebra.

## 1.4 Super Linear Algebra

In this section for the first time we define the notion of super linear algebra and give some of its properties.

**DEFINITION 1.4.1:** *Let $V = (V_1 \mid \ldots \mid V_n)$ be a super vector space over a field F. We say V is a super linear algebra over F if and only if for every pair of super row vectors $\alpha$, $\beta$ in V the product of $\alpha$ and $\beta$ denoted by $\alpha\beta$ is defined in V in such a way that*

(a) *multiplication of super vector in V is associative i.e., if $\alpha$, $\beta$ and $\gamma \in V$ then $\alpha(\beta\gamma) = (\alpha\beta)\gamma$.*
(b) *multiplication is distributive $(\alpha + \beta)\gamma = \alpha\gamma + \beta\gamma$ and $\alpha(\beta + \gamma) = \alpha\beta + \alpha\gamma$ for every $\alpha, \beta, \gamma \in V$.*
(c) *for each scalar c in F $c(\alpha\beta) = (c\alpha)\beta = \alpha(c\beta)$.*

*If there is an element $1_e$ in V such that $1_e\alpha = \alpha 1_e$ for every $\alpha \in V$ we call the super linear algebra V to be a super linear algebra with identity over F. The super linear algebra V is called commutative if $\alpha\beta = \beta\alpha$ for all $\alpha$ and $\beta$ in V.*

We give examples of super linear algebras.

***Example 1.4.1:*** Let

$$V = \{(x_1 x_2 \mid x_3 x_4 x_5 \mid x_6)\; x_i \in Q;\; 1 \le i \le 6\}$$

be a super vector space over Q. Define for $\alpha, \beta \in V$;

$$\alpha = (x_1 x_2 \mid x_3 x_4 x_5 \mid x_6)$$

and

$$\beta = (y_1 y_2 \mid y_3 y_4 y_5 \mid y_6),$$

$$\alpha\beta = (x_1 y_1\; x_2 y_2 \mid x_3 y_3\; x_4 y_4\; x_5 y_5 \mid x_6 y_6).$$

where $x_i, y_j \in Q$; $1 \le i, j \le 6$.
Clearly $\alpha\beta \in V$; so V is a super linear algebra, it can be easily checked that the product is associative. Also it is easily verified the operation is distributive $\alpha(\beta + \gamma) = \alpha\beta + \alpha\gamma$ and $(\alpha + \beta)\gamma = \alpha\gamma + \beta\gamma$ for all $\alpha, \beta, \gamma \in V$.
This V is a super linear algebra. Now the very natural question is that, "is every super vector space a super linear algebra?" The truth is as in case of usual linear algebra, every super linear algebra is a super vector space but in general every super vector space need not be a super linear algebra.

We prove this only by examples.

***Example 1.4.2:*** Let

$$V = \left\{ \left( \begin{array}{ccc|cc} a_1 & a_2 & a_3 & a_{10} & a_{11} \\ a_4 & a_5 & a_6 & a_{12} & a_{13} \\ \hline a_7 & a_8 & a_9 & a_{14} & a_{15} \end{array} \right) \middle| a_i \in Q_i \;\; 1 \le i \le 15 \right\}.$$

Clearly V is a super vector space over Q but is not a super linear algebra.

***Example 1.4.3:*** Let

$$V = \left\{ \left( \begin{array}{cc|c} a_1 & a_2 & a_7 \\ a_3 & a_4 & a_8 \\ \hline a_5 & a_6 & a_9 \end{array} \right) | a_i \in Q; 1 \le i \le 9 \right\};$$

V is a super vector space over Q. V is a super linear algebra for multiplication is defined in V. Let

$$A = \left( \begin{array}{cc|c} 1 & 0 & 1 \\ 2 & 1 & 0 \\ \hline 0 & 1 & 2 \end{array} \right)$$

and

$$B = \left( \begin{array}{cc|c} 0 & 1 & 0 \\ 1 & 2 & 2 \\ \hline 1 & 0 & 1 \end{array} \right)$$

$$A\,B = \left( \begin{array}{cc|c} 1 & 1 & 1 \\ 1 & 4 & 2 \\ \hline 3 & 2 & 2 \end{array} \right) \in V.$$

Now we have seen that in general all super vector spaces need not be super linear algebras.

We proceed on to define the notion of super characteristic values or we may call it as characteristic super values. We have also just now seen that the collection of all linear operators of a super vector space to itself is a super linear algebra.

**DEFINITION 1.4.2:** *Let* $V = \{(A_1 \mid \ldots \mid A_k) \mid A_i$ *are row vectors with entries from a field F}* *and let* $T_s$ *be a linear operator on V.*

*i.e.,* $T_s: V \rightarrow V$ *i.e.,* $T_s: (A_l \mid \ldots \mid A_k) \rightarrow (A_1 \mid \ldots \mid A_k)$

*i.e.,* $T = T_s = (T_1 \mid \ldots \mid T_k)$ *with* $T_i: A_i \rightarrow A_i$; $i = 1, 2, \ldots, k$. *A characteristic super value is* $c = (c_1 \, c_2 \ldots c_k)$ *in F (i.e., each* $c_i \in F$*) such that there is a non zero super vector* $\alpha$ *in V with* $T\alpha = c\alpha$ *i.e.,* $T_i\alpha_i = c_i\alpha_i$, $\alpha_i \in A_i$ *true for each i.*

*i.e.,* $T\alpha = c\alpha$.
*i.e.,* $(T_1\alpha_1 \mid \ldots \mid T_k\alpha_k) = (c_1\alpha_1 \mid \ldots \mid c_k\alpha_k)$.

*The k-tuple* $(c_1 \ldots c_k)$ *is a characteristic super value of* $T = (T_1 \mid \ldots \mid T_k)$,

(a) *We have for any* $\alpha$ *such that* $T\alpha = c\alpha$, *then* $\alpha$ *is called the characteristic super vector of T associated with the characteristic super value* $c = (c_1, \ldots, c_k)$

(b) *The collection of all super vectors* $\alpha$ *such that* $T\alpha = c\alpha$ *is called the characteristic super vector space associated with c.*

*Characteristic super values are often called characteristic super vectors, latent super roots, eigen super values, proper super values or spectral super values.*

We shall use in this book mainly the terminology characteristic super values.

It is left as an exercise for the reader to prove later.

**THEOREM 1.4.1:** *Let $T_s$ be a linear operator on a finite dimensional super vector space V and let c be a scalar of n tuple. Then the following are equivalent.*

i. *c is a characteristic super value of $T_s$*
ii. *The operator $(T_s – cI)$ is singular.*
iii. *det $(T_s – cI) = (0)$.*

We now proceed on to define characteristic super values and characteristic super vectors for any square super diagonal matrix A. We cannot as in case of other matrices define the notion of characteristic super values of any square super matrix as at the first instance we do not have the definition of determinant in case of super matrices. As the concept of characteristic values are defined in terms of the determinant of matrices so also the characteristic super values can only be defined in terms of the determinant of super matrices. So we just define the determinant value in case of only square super diagonal matrix whose diagonal elements are also squares.

We first give one or two examples of square super diagonal matrix.

***Example 1.4.4:*** Let A be a square super diagonal matrix where

$$A = \left(\begin{array}{cccc|ccc|ccc} 0 & 1 & 2 & 3 & 0 & 0 & 0 & 0 & 0 & 0 \\ 5 & 0 & 1 & 1 & 0 & 0 & 0 & 0 & 0 & 0 \\ 2 & 0 & 0 & 1 & 0 & 0 & 0 & 0 & 0 & 0 \\ \hline 0 & 0 & 0 & 0 & 9 & 2 & 1 & 0 & 0 & 0 \\ 0 & 0 & 0 & 0 & 0 & 1 & 2 & 0 & 0 & 0 \\ \hline 0 & 0 & 0 & 0 & 0 & 0 & 0 & 6 & 1 & 0 \\ 0 & 0 & 0 & 0 & 0 & 0 & 0 & 1 & 1 & 0 \\ 0 & 0 & 0 & 0 & 0 & 0 & 0 & 0 & 1 & 1 \\ 0 & 0 & 0 & 0 & 0 & 0 & 0 & 1 & 0 & 1 \\ 0 & 0 & 0 & 0 & 0 & 0 & 0 & 1 & 0 & 2 \end{array}\right).$$

This is a square super diagonal matrix whose diagonal terms are not square matrices. So for such type of matrices we cannot define the notion of determinant of A.

***Example 1.4.5:*** Consider the super square diagonal matrix A given by

$$A = \left(\begin{array}{cccc|ccc|ccccc}
3 & 1 & 0 & 2 & 0 & 0 & 0 & 0 & 0 & 0 & 0 & 0 \\
1 & 0 & 1 & 1 & 0 & 0 & 0 & 0 & 0 & 0 & 0 & 0 \\
5 & 0 & 0 & 1 & 0 & 0 & 0 & 0 & 0 & 0 & 0 & 0 \\
0 & 2 & 1 & 0 & 0 & 0 & 0 & 0 & 0 & 0 & 0 & 0 \\
\hline
0 & 0 & 0 & 0 & 5 & 0 & 1 & 0 & 0 & 0 & 0 & 0 \\
0 & 0 & 0 & 0 & 0 & 1 & 1 & 0 & 0 & 0 & 0 & 0 \\
0 & 0 & 0 & 0 & 0 & 2 & 2 & 0 & 0 & 0 & 0 & 0 \\
\hline
0 & 0 & 0 & 0 & 0 & 0 & 0 & 1 & 0 & -1 & 2 & 0 \\
0 & 0 & 0 & 0 & 0 & 0 & 0 & 0 & 1 & 3 & 1 & 6 \\
0 & 0 & 0 & 0 & 0 & 0 & 0 & 1 & 0 & 0 & 2 & -1 \\
0 & 0 & 0 & 0 & 0 & 0 & 0 & 0 & 2 & 5 & 1 & 0 \\
0 & 0 & 0 & 0 & 0 & 0 & 0 & 1 & 3 & 1 & 0 & 2
\end{array}\right)$$

Clearly A is the square super diagonal matrix which diagonal elements are also square matrices, these super matrices we venture to define as square super square diagonal matrix or strong square super diagonal matrix.

***Example 1.4.6:*** Let A be a 10 × 12 super diagonal matrix.

$$A = \left(\begin{array}{cccc|ccc|ccccc}
6 & 3 & 1 & 2 & 0 & 0 & 0 & 0 & 0 & 0 & 0 & 0 \\
1 & 1 & 0 & 1 & 0 & 0 & 0 & 0 & 0 & 0 & 0 & 0 \\
3 & 0 & 0 & 7 & 0 & 0 & 0 & 0 & 0 & 0 & 0 & 0 \\
\hline
0 & 0 & 0 & 0 & 7 & 1 & 0 & 0 & 0 & 0 & 0 & 0 \\
0 & 0 & 0 & 0 & 0 & 1 & 2 & 0 & 0 & 0 & 0 & 0 \\
0 & 0 & 0 & 0 & 0 & 1 & 1 & 0 & 0 & 0 & 0 & 0 \\
\hline
0 & 0 & 0 & 0 & 0 & 0 & 0 & 4 & 0 & 1 & 1 & -1 \\
0 & 0 & 0 & 0 & 0 & 0 & 0 & 0 & 2 & 1 & 0 & 0 \\
0 & 0 & 0 & 0 & 0 & 0 & 0 & 3 & 5 & 0 & 0 & 1 \\
0 & 0 & 0 & 0 & 0 & 0 & 0 & 1 & 2 & 3 & 4 & 0
\end{array}\right).$$

We see the diagonal elements are not square matrices hence A is not a square matrix but yet A is a super diagonal matrix. Thus unlike in usual matrices where we cannot define the notion of diagonal if the matrix is not a square matrix in case of super matrices which are not square super matrices we can define the concept of super diagonal even if the super matrix is not a square matrix.

So we can call a rectangular super matrix to be a super diagonal matrix if in that super matrix all submatrices are zero except the diagonal matrices.

***Example 1.4.7:*** Let A be a square super diagonal matrix given below

$$A = \left(\begin{array}{ccccc|ccc|ccccc}
9 & 0 & 1 & 1 & 1 & 0 & 0 & 0 & 0 & 0 & 0 & 0 & 0 \\
0 & 1 & 0 & 1 & 0 & 0 & 0 & 0 & 0 & 0 & 0 & 0 & 0 \\
2 & 0 & 1 & 0 & 0 & 0 & 0 & 0 & 0 & 0 & 0 & 0 & 0 \\
0 & 3 & 0 & 0 & 4 & 0 & 0 & 0 & 0 & 0 & 0 & 0 & 0 \\
\hline
0 & 0 & 0 & 0 & 0 & 2 & 0 & 1 & 0 & 0 & 0 & 0 & 0 \\
0 & 0 & 0 & 0 & 0 & 0 & 3 & 2 & 0 & 0 & 0 & 0 & 0 \\
0 & 0 & 0 & 0 & 0 & 1 & 1 & 0 & 0 & 0 & 0 & 0 & 0 \\
\hline
0 & 0 & 0 & 0 & 0 & 0 & 0 & 0 & 2 & 1 & 0 & 1 & 1 \\
0 & 0 & 0 & 0 & 0 & 0 & 0 & 0 & 0 & 1 & 0 & 2 & 0 \\
0 & 0 & 0 & 0 & 0 & 0 & 0 & 0 & 1 & 1 & 0 & 1 & 2 \\
0 & 0 & 0 & 0 & 0 & 0 & 0 & 0 & 2 & 0 & 1 & 1 & 1 \\
0 & 0 & 0 & 0 & 0 & 0 & 0 & 0 & 0 & 1 & 1 & 2 & 1 \\
0 & 0 & 0 & 0 & 0 & 0 & 0 & 0 & 2 & 0 & 1 & 1 & 1
\end{array}\right)$$

This super square matrix A is a diagonal super square matrix as the main diagonal are matrices. Hence A is only a super square diagonal matrix but is not a super square diagonal square matrix as the diagonal matrices are not square matrices.

**DEFINITION 1.4.2:** *Let A be a square super diagonal matrix whose diagonal matrices are also square matrices then the super determinant of A is defined as*

$$|A| = \left[\begin{array}{c|c|c|c|c|c}
|A_1| & 0 & 0 & 0 & 0 & 0 \\
\hline
0 & |A_2| & 0 & 0 & 0 & 0 \\
\hline
0 & 0 & |A_3| & 0 & 0 & 0 \\
\hline
0 & 0 & 0 & 0 & 0 & 0 \\
\hline
0 & 0 & 0 & 0 & |A_{n-1}| & 0 \\
\hline
0 & 0 & 0 & 0 & 0 & |A_n|
\end{array}\right]$$

$$= (|A_1|, |A_2|, \ldots, |A_n|).$$

*where each submatrix $A_i$ of A is a square matrix and $|A_i|$ denotes the determinant of $A_i$, i = 1, 2, …, n.*

***Example 1.4.8:*** Let A be a super square diagonal matrix;

$$\left(\begin{array}{cc|ccc|cc|cccc} 2 & 1 & & 0 & & & 0 & & & 0 & \\ 0 & 1 & & & & & & & & & \\ \hline & & 3 & 1 & 2 & & & & & & \\ & 0 & 0 & 0 & 1 & & 0 & & & 0 & \\ & & 1 & 2 & 0 & & & & & & \\ \hline & 0 & & 0 & & 3 & 1 & & & 0 & \\ & & & & & 0 & 1 & & & & \\ \hline & & & & & & & 0 & 1 & 2 & 0 \\ & 0 & & 0 & & & 0 & 0 & 0 & 3 & 4 \\ & & & & & & & 0 & 1 & 0 & 0 \\ & & & & & & & 1 & 0 & 0 & 0 \end{array}\right).$$

Now the super determinant of

$$A = |A| = \left[\begin{array}{c|c|c|c} \begin{vmatrix} 2 & 1 \\ 0 & 1 \end{vmatrix} & \begin{vmatrix} 3 & 1 & 2 \\ 0 & 0 & 1 \\ 1 & 2 & 0 \end{vmatrix} & \begin{vmatrix} 3 & 1 \\ 0 & 1 \end{vmatrix} & \begin{vmatrix} 0 & 1 & 2 & 0 \\ 0 & 0 & 3 & 4 \\ 0 & 1 & 0 & 0 \\ 1 & 0 & 0 & 0 \end{vmatrix} \end{array}\right]$$

$$= [\ 2 \mid -5 \mid 3 \mid -8 \mid].$$

We see the resultant is a super vector. Thus the super determinant of a square super diagonal square matrix which we define as a super determinant is always a super vector. Further if the square super diagonal matrix has n components then we have the super determinant to have a super row vector with n partition and the natural length of the super vector is also n.

Thus the super determinant of a super matrix is defined if and only if the super matrix is a square super diagonal square matrix.

Now having defined the determinant of a square super diagonal matrix, we proceed on to define super characteristic value associated with a square super diagonal square matrix. At this point it has become pertinent to mention here that all linear operators $T_s$ can be associated with a super matrix A, where A is a super square diagonal square matrix.

Now we first illustrate it by an example. We have already defined the notion of super polynomial $p(x) = [p_1(x) \mid p_2(x) \mid \ldots \mid p_n(x)]$.

Now we will be making use of this definition also.

***Example 1.4.9:*** Let V = {(Q[x] | Q[x] | Q[x] | Q[x]) | Q[x] are polynomials with coefficients from the rational field Q}. V is a super vector space of infinite dimension called the super vector space of polynomials of infinite dimension over Q. Any element $p(x) = (p_1(x) \mid p_2(x) \mid p_3(x) \mid p_4(x))$ such $p_i(x) \in Q[x]$; $1 \leq i \leq 4$ or more non abstractly $p(x) = [x^3+1 \mid 2x^2 – 3x+1 \mid 5x^7 + 3x^2 + 3x + 1 \mid x^5 – 2x + 1] \in V$ is a super polynomial of V.

This polynomial p(x) can also be given the super row vector representation by p(x) = (1 0 0 1 | 1 –3 2 | 1 3 3 0 0 0 0 5 | 1 –2 0 0 0 1].

Here it is pertinent to mention that the super row vectors will not be of the same type. Still in interesting to note that V = {(Q[x] | … | Q[x]) | Q[x] are polynomial rings} over the field Q is a super linear algebra, for if $p(x) = (p_1(x) \mid \ldots \mid p_n(x))$ and $q(x) = (q_1 (x) \mid \ldots \mid q_n(x)) \in V$ then $p(x)\, q(x) = (p_1(x)\, q_1(x) \mid \ldots \mid p_n(x)q_n(x)) \in V$.

Thus the super vector space of polynomials of infinite dimension is a super linear algebra over the field over which they are defined.

***Example 1.4.10:*** Let

$V = \{(Q^5[x] \mid Q^3[x] \mid Q^6[x] \mid Q^2[x] \mid Q^3[x]) \mid Q^i[x]$

is a polynomial of degree less than or equal i; i = 3, 5, 6 and 2}, V is a super vector space over Q and V is a finite dimensional super vector space over Q.

For the dimension of V is 6 + 4 + 7 + 3 + 4 = 24. Thus

$$V \cong \{(x_1 x_2 x_3 x_4 x_5 x_6 \mid x_7 x_8 x_9 x_{10} \mid x_{11} x_{12} x_{13} x_{14} x_{15} x_{16} x_{17} \mid x_{18} x_{19} x_{20} \mid x_{21} x_{22} x_{23} x_{24}) \mid x_i \in Q; 1 \leq i \leq 24\}$$

is a super vector space of dimension 24 over Q.

Clearly V is a super vector space of super polynomials of finite degree. Further V is not a super linear algebra.

So any element p(x) = $\{(x^3+1 \mid x^2+4 \mid x^5+3x^4 + x^2+1 \mid x+1 \mid x^2+3x-1)\}$ = (1 0 0 1 | 4 0 1 | 1 0 1 0 3 1 | 1 1 | -1 3 1) is the super row vector representation of p(x).

How having illustrated by example the super determinant and super polynomials now we proceed on to define the notion of super characteristic values and super characteristic polynomial associated with a square super diagonal square matrix with entries from a field F.

**DEFINITION 1.4.4:** *Let*

$$A = \left(\begin{array}{c|c|c|c} A_1 & 0 & & 0 \\ \hline 0 & A_2 & & 0 \\ \hline 0 & 0 & & 0 \\ \hline 0 & 0 & & A_n \end{array}\right)$$

*be a square super diagonal square matrix with entries from a field F, where each $A_i$ is also a square matrix i = 1, 2, …, n. A super characteristic value of A or characteristic super value of A (both mean the same) in F is a scalar n-tuple c = ($c_1$ | … | $c_n$) in F such that the super matrix | A – cI | is singular ie non invertible ie [A – cI] is again a square super diagonal super square matrix given as follows.*

$$A - cI = \left(\begin{array}{c|c|c|c} A_1 - c_1 I & 0 & & 0 \\ \hline 0 & A_2 - c_2 I & & 0 \\ \hline 0 & 0 & & 0 \\ \hline 0 & 0 & 0 & A_n - c_n I \end{array}\right)$$

*where c = ($c_1$ | … | $c_n$) as mentioned earlier $c_i \in F$; $1 \le i \le n$. c is the super characteristic value of A if and only if super det (A–cI) = (det ($A_1$-$c_1$I) | … | det ($A_n$-$c_n$I))*
*= (0 | 0 | … | 0) or equivalently if and only of super det [cI – A] = (det ($A_1$ – $c_1$I) | … | det ($A_n$ – $c_n$I)) = (0 | … | 0), we form the super matrix (xI-A) = ((xI – $A_1$) | … | (xI – $A_n$)) with super polynomial entries and consider the super polynomial f = det (xI – A) = (det (xI – $A_1$) | … | det (xI – $A_n$)) = [$f_1$ | … | $f_n$].*

*Clearly the characteristic super value of A in F are just the super scalars c in F such that f(c) = ($f_1$($c_1$) | $f_2$($c_2$) | … | $f_n$($c_n$)) = (0 | … | 0). For this reason f is called the super characteristic polynomial (characteristic super polynomial) of A. It is important to note that f is a super monic polynomial which has super deg exactly ($n_1$ | … | $n_n$) where $n_i$ is the order of the square matrix $A_i$ of A for i = 1, 2, …, n.*

We say a super polynomial p(x) = [$p_1$(x) | … | $p_n$(x)] to be a super monic polynomial if every polynomial $p_i$(x) of p (x) is monic for i = 1, 2, …, n.

Based on this we can define the new notion of similarly square super diagonal square matrices.

**DEFINITION 1.4.5:** *Let A be a square super diagonal square matrix with entries from a field F.*

$$A = \left(\begin{array}{c|c|c|c} A_1 & 0 & & 0 \\ \hline 0 & A_2 & & 0 \\ \hline & & & \\ \hline 0 & 0 & & A_n \end{array}\right)$$

*where each $A_i$ is a square matrix of order $n_i \times n_i$, i = 1, 2, …, n. Let B be another square super diagonal square matrix of same order ie let*

$$B = \left(\begin{array}{c|c|c|c} B_1 & 0 & & 0 \\ \hline 0 & B_2 & & 0 \\ \hline & & & \\ \hline 0 & 0 & & B_n \end{array}\right)$$

*where each $B_i$ is a $n_i \times n_i$ matrix for i = 1, 2, …, n.*

We say A and B are similar super matrices if there exists an invertible square super diagonal square matrix P;

$$P = \left(\begin{array}{c|c|c|c} P_1 & 0 & & 0 \\ \hline 0 & P_2 & & 0 \\ \hline & & & \\ \hline 0 & 0 & & P_n \end{array}\right)$$

where each $P_i$ is a $n_i \times n_i$ matrix for i = 1, 2, …, n such that each $P_i$ is invertible i.e., $P_i^{-1}$ exists for each i = 1, 2, …, n; and is such that

$$B = P^{-1} A P = \left(\begin{array}{c|c|c|c} P_1^{-1} & 0 & & 0 \\ \hline 0 & P_2^{-1} & & 0 \\ \hline & & & \\ \hline 0 & & & P_n^{-1} \end{array}\right) \times$$

$$\left(\begin{array}{c|c|c|c} A_1 & 0 & & 0 \\ \hline 0 & A_2 & & 0 \\ \hline & & & \\ \hline 0 & & & A_n \end{array}\right) \times \left(\begin{array}{c|c|c|c} P_1 & 0 & & 0 \\ \hline 0 & P_2 & & 0 \\ \hline & & & \\ \hline 0 & & & P_n \end{array}\right)$$

$$= \left(\begin{array}{c|c|c|c} B_1 & 0 & & 0 \\ \hline 0 & B_2 & & 0 \\ \hline & & & \\ \hline 0 & & & B_n \end{array}\right)$$

$$= \left(\begin{array}{c|c|c|c} P_1^{-1}A_1P_1 & 0 & & 0 \\ \hline 0 & P_2^{-1}A_2P_2 & & 0 \\ \hline & & & \\ \hline 0 & & & P_n^{-1}A_nP_n \end{array}\right).$$

If $B = P^{-1} A P$ then super determinant of $(xI – B)$

= super determinant of $(xI – P^{-1}A P)$ i.e., $(\det (xI-B_1) \mid \ldots \mid \det (xI – B_n))$

$= (\det (xI - P_1^{-1}A_1P_1) \mid \ldots \mid \det (xI – P_n^{-1}A_nP_n))$

$= (\det(P_1^{-1}(xI – A_1)P_1) \mid \ldots \mid \det P_n^{-1}(xI – A_n)P_n)$

$= (\det P_1^{-1} \det (xI – A_1) \det P_1 \mid \ldots \mid \det P_n^{-1} \det (xI – A_n) \det P_n)$

$= (\det (xI – A_1) \mid \ldots \mid \det (xI – A_n))$.

Thus this result enables one to define the characteristic super polynomial of the operator $T_s$ as the characteristic super polynomial of any $(n_1 \times n_1 \mid \ldots \mid n_n \times n_n)$ square super diagonal square matrix which represents $T_s$ in some super basis for V.

Just as for square super diagonal matrices the characteristic super values of $T_s$ will be the roots of the characteristic super polynomial for $T_s$. In particular this shows us that $T_s$ cannot have more than $n_1 + \ldots + n_n$ characteristic super values.

It is pertinent to point out that $T_s$ may not have any super characteristic values. This is shown by the following example.

***Example 1.4.11:*** Let $T_s$ be a linear operator on $V = \{(x_1x_2 \mid x_3x_4 \mid x_5x_6) \mid x_i \in Q; 1 \leq i \leq 6\}$ the super vector space over Q, which is represented by a square super diagonal square matrix

$$A = \left(\begin{array}{cc|cc|cc} 0 & -1 & 0 & 0 & 0 & 0 \\ 1 & 0 & 0 & 0 & 0 & 0 \\ \hline 0 & 0 & 0 & -1 & 0 & 0 \\ 0 & 0 & 1 & 0 & 0 & 0 \\ \hline 0 & 0 & 0 & 0 & 0 & 1 \\ 0 & 0 & 0 & 0 & -1 & 0 \end{array}\right) = \left(\begin{array}{c|c|c} A_1 & 0 & 0 \\ \hline 0 & A_2 & 0 \\ \hline 0 & 0 & A_3 \end{array}\right).$$

The characteristic super polynomial for $T_s$ or for A is super determinant of

$$(xI - A) = \left\|\left(\begin{array}{cc|cc|cc} x & 1 & 0 & 0 & 0 & 0 \\ -1 & x & 0 & 0 & 0 & 0 \\ \hline 0 & 0 & x & 1 & 0 & 0 \\ 0 & 0 & -1 & x & 0 & 0 \\ \hline 0 & 0 & 0 & 0 & x & -1 \\ 0 & 0 & 0 & 0 & 1 & x \end{array}\right)\right\|$$

i.e., super det $(xI - A) = [x^2+1 \mid x^2+1 \mid x^2+1]$.
This super polynomial has no real roots, $T_s$ has no characteristic super values.

Now we proceed on to discuss about when a super linear operator $T_s$ on a finite dimensional super vector space V is super diagonalizable.

**DEFINITION 1.4.6:** *Let $V = (V_1 \mid \ldots \mid V_n)$ be a super vector space over the field F of super dimension $(n_1 \mid \ldots \mid n_n)$, ie each vector space $V_i$ over the field F is of dimension $n_i$ over F, i = 1, 2, …, n. We say a linear operator $T_s$ on V is super diagonalizable if there is a super basis for V, each super vector of which is a characteristic super vector of $T_s$.*

Recall as in case of usual matrices or usual linear operators of a vector space we in case of super vector spaces using a linear operator $T_s = (T_1 \mid \ldots \mid T_n)$ on V, to have the characteristic super vector $\alpha = (\alpha_1 \mid \ldots \mid \alpha_n)$ as the characteristic super vector, if $T_s\alpha = c\alpha$ i.e., $(T_1\alpha_1 \mid \ldots \mid T_n\alpha_n) = (c_1\alpha_1 \mid \ldots \mid c_n\alpha_n)$ where $c = (c_1 \mid \ldots \mid c_n))$ is the characteristic super value associated with $T_s$.

If the super characteristic value are denoted by

$$(c_1^1 \ldots c_1^{n_1} \mid c_2^1 \ldots c_2^{n_2} \mid \ldots \mid c_n^1 \ldots c_n^{n_n})$$

and for a super basis

$$B = (\alpha_1^1 \ldots \alpha_{n_1}^1 \mid \ldots \mid \alpha_n^1 \;\ldots\; \alpha_{n_n}^n)$$

for V.

$$T_s\alpha = c\alpha \text{ i.e., } T_i\ \alpha_t^i = c_i^t\alpha_t^i$$

for $t = 1, 2, \ldots, n_i$ and $i = 1, 2, \ldots, n$.

$$[T_s]_B = \left(\begin{array}{ccc|ccc|c|ccc}
c_1^1 & 0 & 0 & & & & & & & \\
0 & c_1^2 & 0 & & 0 & & 0 & & 0 & \\
0 & 0 & c_1^{n_1} & & & & & & & \\
\hline
& & & c_2^1 & 0 & 0 & & & & \\
& 0 & & 0 & c_2^2 & 0 & 0 & & 0 & \\
& & & 0 & 0 & c_2^{n_2} & & & & \\
\hline
& 0 & & & & & & & 0 & \\
\hline
& & & & & & & c_n^1 & 0 & 0 \\
& 0 & & & 0 & & 0 & 0 & c_n^2 & 0 \\
& & & & & & & 0 & 0 & c_n^{n_n}
\end{array}\right)$$

We certainly require the super scalars

$$[c_1^1 c_1^2 \ldots c_i^{n_1} \mid c_1^1 \ldots c_2^{n_2} \mid \ldots \mid c_n^1 \ldots c_n^{n_n}]$$

must be distinct for $i = 1, 2, \ldots, n$.

The scalars can be identical when each $T_i$ is a scalar multiple of the identity operator. But in general we may not have them to be distinct suppose $T_s$ is a super diagonalizable operator.

Let $(c_1^1, c_1^2, \ldots, c_i^{k_1})$, $(c_2^1, c_2^2, \ldots, c_2^{k_2})$, $\ldots (c_n^1, c_n^2, \ldots, c_n^{k_n})$ be the distinct characteristic values of $T_i$ of $T_s$ for i = 1, 2, …, n where $T_s = [T_1 \mid \ldots \mid T_n]$. Then we have a basis B for which $T_s$ is represented by a super diagonal matrix for which its diagonal entries are $c_1^i, c_2^i, \ldots, c_n^i$ each repeated a certain number of times. If $c_t^i$ is represented $d_t^i$ times then the super matrix has super block form,

ie $[T_s]_B = [T_s]B =$

$$\left(\begin{array}{c|c|c|c}
\begin{matrix} c_1^1 I_1^1 & 0 & 0 \\ 0 & c_1^2 I_2^1 & 0 \\ 0 & 0 & c_1^{k_1} I_{k_1}^1 \end{matrix} & 0 & 0 & 0 \\ \hline
0 & \begin{matrix} c_2^1 I_1^2 & 0 & 0 \\ 0 & c_2^2 I_1^2 & 0 \\ 0 & 0 & c_2^{k_2} I_{k_2}^2 \end{matrix} & 0 & 0 \\ \hline
0 & & & 0 \\ \hline
0 & 0 & 0 & \begin{matrix} c_n^1 I_1^n & 0 & 0 \\ 0 & c_n^2 I_2^n & 0 \\ 0 & 0 & c_n^{k_n} I_{k_n}^n \end{matrix}
\end{array}\right)$$

where $I_j^t$ is a $d_j^t \times d_j^t$ identity matrix, j = 1, 2, …, $k_t$ and t = 1, 2, …, n. Thus we see the characteristic super polynomial for $T_s$ is the product of linear factors

$$f = ((x - c_1^1)^{d_1^1} \ldots (x - c_1^{k_1})^{d_{k_1}^1} \mid (x - c_2^1)^{d_1^2} (x - c_2^2)^{d_2^2} \ldots$$
$$(x - c_2^{k_2})^{d_{k_2}^2} \mid \ldots \mid (x - c_n^1)^{d_1^n} \ldots (x - c_n^{k_n})^{d_{k_n}^n}) = (f_1 \mid \ldots \mid f_n).$$

We leave the following lemma as an exercise for the reader.

**LEMMA 1.4.2:** *Suppose $T_s\alpha = c\alpha$. If $f = (f_1 \mid \ldots \mid f_n)$ is any super polynomial then $f(T_s) = (f_1(T_1) \mid \ldots \mid f_n(T_n))((\alpha_1 \mid \ldots \mid \alpha_n)) = (f_1(T_1)\,\alpha_1 \mid \ldots \mid f_n(T_n)\,\alpha_n) = (f_1(c_1)\,\alpha_1 \mid \ldots \mid f_n(c_n)\,\alpha_n)$.*

**LEMMA 1.4.3:** *Let $T_s = (T_1 \mid \ldots \mid T_n)$ be a linear operator on a finite dimensional super vector space $V = (V_1 \mid \ldots \mid V_n)$. $\{(c_1^1 \ldots c_{k_1}^1), \ldots, (c_1^n \ldots c_{k_n}^n)\}$ be the distinct set of super characteristic values of $T_s$ and let $(W_{i_1}^1 \mid \ldots \mid W_{i_n}^n)$ be the super subspace of the characteristic super vectors associated with characteristic super values $(c_{i_1}^1, \ldots, c_{i_n}^n)$; $1 \le i_t \le k_t$; $t = 1, 2, \ldots, n$.*
*If*

$$W = (W_1^1 + \ldots + W_{k_1}^1 \mid W_1^2 + \ldots + W_{k_2}^2 \mid \ldots \mid W_1^n + \ldots + W_{k_n}^n)$$
$$= (W^1 \mid \ldots \mid W^n);$$

*then super dimension of*

$$W = (\dim W_1^1 + \ldots + \dim W_{k_1}^1 \mid \ldots \mid \dim W_1^n + \ldots + \dim W_{k_n}^n).$$

*In fact if*

$$B = (B_1^1 \ldots B_{k_1}^1 \mid \ldots \mid B_1^n \ldots B_{k_n}^n)$$

*where $(B_{i_1}^1 \ldots B_{i_n}^n)$ is the ordered basis for $(W_{i_1}^i \mid \ldots \mid W_{i_n}^i)$ then B is ordered basis of W.*

*Proof:* We prove the result for one subspace $W^t = W_1^t + \ldots + W_{k_t}^t$, this being true for every t, t = 1, 2, …, n thus we see it is true for the super subspace $W = (W^1 \mid \ldots \mid W^n)$.

The space $W^t = W_1^t + \ldots + W_{k_t}^t$ is the subspace spanned by all the characteristic vectors of $T_t$ where $T_s = (T_1 \mid \ldots \mid T_n)$ and $1 \le t \le n$. Usually when one forms the sum $W^t$ of subspaces $W_i^t$; $1 \le i \le k_t$ one expects $\dim W^t < \dim W_1^t + \ldots + \dim W_{k_t}^t$ because of linear relations which may exist between vectors in the various spaces. From the above lemma the characteristic spaces associated with different characteristic values are independent of one another.

Suppose that for each $i_t$ we have a vector $\beta_i^t$ in $W_i^t$ and assume that $\beta_i^t + \ldots + \beta_{k_t}^t = 0$ we shall show that $\beta_i^t = 0$ for each i, i = 1, 2, …, $k_t$. Let $f_t$ be any polynomial of the super polynomial $f = (f_1 \mid \ldots \mid f_n)$; $1 \le t \le n$.

Since $T_t\beta_i^t = c_i^t\, \beta_i^t$ the proceeding lemma tells us that 0 = $f_t(T_t); 0 = f_t(T_t)\beta_1^t + \ldots + f_t(T_t)\, \beta_{k_t}^t$

$$= f_t(c_1^t)\beta_1^t + \ldots + f_t(c_{k_t}^t)\beta_{k_t}^t .$$

Choose polynomials $f_t^1, f_t^2, \ldots, f_t^{k_t}$ such that

$$f_t^{i_t}(c_{j_t}^t) = \delta_{i_t j_t} = \begin{cases} 0 & i_t \neq j_t \\ 1 & i_t = j_t \end{cases},$$

for t = 1, 2, …, n. Then 0 = $f_t^i(T_t); 0 = \sum \delta_{i_t j_t}\, \beta_{j_t}^t = \beta_{i_t}^t$.

Now let $B_{i_t}^t$ be a basis of $W_{i_t}^t$ and let $B^t = (B_1^t \ldots B_{k_t}^t)$

Then $B^t$ spans $W^t = W_1^t + \ldots + W_{k_t}^t$, this is true for every t = 1, 2, …, n. Also $B^t$ is a linearly independent sequence of vectors. Any linear relation between the vectors in $B^t$ will have the form $\beta_1^t + \ldots + \beta_{k_t}^t = 0$; where $\beta_{i_t}^t$ is some linear combination of vectors in $\beta_{i_t}^t$ We have just shown $\beta_{i_t}^t$ = 0 for each $i_t$ = 1, 2, .., $k_t$ and for each t = 1, 2, …, n. Since each $B_{i_t}^t$ is linearly independent we see that we have only the trivial relation between the vectors in $B^t$; since this is true for each t we have only trivial relation between the super vectors in

$$B = (B^1 \mid \ldots \mid B^n) = (B_1^1 \ldots B_{k_1}^1 \mid B_1^2 \ldots B_{k_2}^2 \mid \ldots \mid B_1^n \ldots B_{k_n}^n).$$

Hence B is the ordered super basis for

$$W = (W_1^1 + \ldots + W_{k_1}^1 \mid \ldots \mid W_1^n + \ldots + W_{k_n}^n).$$

**THEOREM 1.4.2:** *Let $T_s = (T_1 \mid \ldots \mid T_n)$ be a linear operator on a finite dimensional super space $V = (V_1 \mid \ldots \mid V_n)$ of dimension $(n_1, \ldots, n_n)$ over the field F. Let $((c_1^1 \ldots c_{k_1}^1), \ldots, (c_1^n \ldots c_{k_n}^n))$ be the distinct characteristic super values of $T_s$ and $W_i = (W_{i_1}^1 \mid \ldots \mid W_{i_n}^n)$ be the null super space of*

$$(T - c_i I) = ((T_1 - c_{i_1}^1 I^1) \mid \ldots \mid (T_n - c_{i_n}^n I^n)).$$

*Then the following are equivalent*

*(i) $T_s$ is super diagonalizable*

*(ii) The characteristic super polynomial for $T_s$ is $f = (f_1 \mid \ldots \mid f_n)$*
$= ((x - c_1^1)^{d_1^1} \ldots (x - c_{k_1}^1)^{d_{k_1}^1} \mid \ldots \mid (x - c_1^n)^{d_1^n} \ldots (x - c_{k_n}^n)^{d_{k_n}^n})$
*and dim $W_{i_t}^t = d_i^t$; $1 \leq t \leq k$; $t = 1,2, \ldots, n$.*

*(iii) dim V = $(\dim W_1^1 + \ldots + \dim W_{k_1}^1 \mid \ldots \mid \dim W_1^n + \ldots + \dim_{k_n}^n)$*
$= (\dim V_1 \mid \ldots \mid \dim V_n) = (n_1, \ldots, n_n)$.

*Proof:* We see that (i) always implies (ii).
If the characteristic super polynomial $f = (f_1 \mid \ldots \mid f_n)$ is the product of linear factors as in (ii) then $(d_1^1 + \ldots + d_{k_1}^1 \mid \ldots \mid d_1^n + \ldots + d_{k_n}^n) = (\dim V_1 \mid \ldots \mid \dim V_n)$.
Therefore (ii) implies (iii) holds. By the lemma just proved we must have $V = (V_1 \mid \ldots \mid V_n)$
$= (W_1^1 + \ldots + W_{k_1}^1 \mid \ldots \mid W_1^n + \ldots + W_{k_n}^n)$,
i.e., the characteristic super vectors of $T_s$ span V.

Next we proceed on to define some more properties for super polynomials we have just proved how the super diagonalization of a linear operator $T_s$ works and the associated super polynomial.

Suppose (F[x] | … | F[x]) is a super vector space of polynomials over the field F. V = ($V_1$ | … | $V_n$) be a super vector space over the field F. Let $T_s$ be a linear operator on V. Now we are interested in studying the class of super polynomial, which annihilate $T_s$. Specifically suppose $T_s$ is a linear operator on V, a super vector space V over the field F. If p = ($p_1$ | … | $p_n$) is a super polynomial over F and q = ($q_1$ | … | $q_n$) another super polynomial over F; then

$$(p + q)\, T_s = ((p_1 + q_1)\, T_1 \mid \ldots \mid (p_n + q_n)\, T_n)$$
$$= (p_1 (T_1) \mid \ldots \mid p_n (T_n)) + (q_1(T_1) \mid \ldots \mid q_n (T_n)).$$
$$(pq)\, (T_s) = (p_1(T_1)q_1(T_1) \mid \ldots \mid p_n(T_n)q_n(T_n)].$$

Therefore the collection of super polynomials p which super annihilate $T_s$ in the sence that $p(T_s) = (p_1(T_1) \mid \ldots \mid p_n(T_n)) = (0 \mid \ldots \mid 0)$, is a super ideal of the super polynomial algebra (F[x] | … | F[x]). Now if A = (F[x] | … | F[x]) is a super polynomial algebra we can define a super ideal I of A as I = ($I_1$ | … | $I_n$) where each $I_t$ is an ideal of F[x]. Now we know if F[x] is the polynomial algebra any polynomial $p_t(x)$ in F[x] will generate an ideal $I_t$ of F[x]. In the same way for any super polynomial p(x) = ($p_1(x)$ | … | $p_n(x)$) of A = [F[x] | … | F[x]], we can associate a super ideal I = ($I_1$ | … | $I_n$) of A.

Suppose $T_s$ = ($T_1$ | … | $T_n$) is a linear operator on V, a ($n_1$ | … | $n_n$) dimensional super vector space. We see the first I,

$$T_1, \ldots, T_1^{n_1^2}, I, T_2, \ldots, T_2^{n_2^2}, \ldots, I, T_n, \ldots, T_n^{n_n^2}$$

has $(n_1^2 + 1, n_2^2 + 1, \ldots, n_n^2 + 1)$ powers of $T_s$ i.e., $n_t^2 + 1$ powers of $T_t$ for t = 1, 2, …, n.

The sequence of $(n_1^2 + 1, \ldots, n_n^2 + 1)$ of super operators in SL (V, V), the super space of linear operators on V. The space SL(V,V) is of dimension $(n_1^2, \ldots, n_n^2)$. Therefore the sequence of $(n_1^2 + 1, \ldots, n_n^2 + 1)$ operators in $T_1$, …, $T_n$ must be linearly dependent as each sequence I, $T_t$, …, $T_t^{n_t}$ is linearly dependent; t = 1, 2, …, n i.e., we have $c_0^t\, I_t + c_1^t\, T_t + \ldots + c_{n_1^2}^t\, T_t^{n_t^2} = 0$ true for each t, t = 1, 2, …, n; i.e., $(c_0^1 I_1 + c_1^1 T_t + \ldots + c_{n_1^2}^1 T_1^{n_1^2} \mid \ldots \mid c_0^n I_n + c_1^n T_n + \ldots + c_{n_n^2}^n T_n^{n_n^2})$ = (0 | … | 0) for some scalar $c_{i_n}^1 \ldots c_{i_n}^n$ not all zero. So the super ideal of super polynomial which annihilate $T_s$ contains non-zero super polynomials of degree $(n_1^2, \ldots, n_n^2)$ or less.

In view of this we now proceed on to define the notion of minimal super polynomial for the linear operator $T_s = (T_1 \mid \ldots \mid T_n)$.

**DEFINITION 1.4.7:** *Let $T_s = (T_1 \mid \ldots \mid T_n)$ be a linear operator on a finite dimensional super n vector space $V = (V_1 \mid V_2 \mid \ldots \mid V_n)$; of dimension $(n_1, \ldots, n_n)$ over the field F. The minimal super polynomial for T is the unique monic super generator of the super ideal of super polynomials over F which super annihilate T.*

The super minimal polynomial is the generator of the super polynomial super ideal is characterized by being the monic super polynomial of minimum degree in the super ideal. That means that the minimal super polynomial $p = (p_1(x) \mid \ldots \mid p_n(x))$ for the linear operator $T_s$ is uniquely determined by the following properties.

(1) $p = (p_1 \mid \ldots \mid p_n)$ is a monic super polynomial over the scalar field F.
(2) $p(T_s) = (p_1(T_1) \mid \ldots \mid p_n(T_n)) = (0 \mid \ldots \mid 0)$.
(3) No super polynomial over F which annihilate $T_s$ has smaller degree than p.

In case of square super diagonal square matrices A we define the minimal super polynomial as follows:
If A is a $(n_1 \times n_1, \ldots, n_n \times n_n)$ square super diagonal super matrix over F, we define the minimal super polynomial for A in an analogous way as the unique monic super generator of the super ideal of all super polynomials over F which super annihilate A, i.e.; which annihilate each of the diagonal matrices $A_t$ ; t = 1, 2, …, n.

If the operator $T_s$ represented in some ordered super basis by the super square diagonal square matrix then $T_s$ and A have the same minimal super polynomial. That is because $f(T_s) = (f_1(T_1) \mid \ldots \mid (f_n(T_n))$ is represented in the super basis by the super diagonal square matrix $f(A) = (f_1(A) \mid \ldots \mid f_n(A_n))$ so that $f(T) = (0 \mid \ldots \mid 0)$ if and only if $f(A) = (0 \mid . . \mid 0)$, i.e.; $f_i(T_i) = 0$ if and only if $f_i(A_i) = 0$.

In fact from the earlier properties mentioned in this book similar square super diagonal square matrices have the same minimal super polynomial. That fact is also clear from the definition because $f(P^{-1}AP) = P^{-1}f(A)P$ i.e., $(f_1(P_1^{-1}A_1P_1) \mid \ldots \mid f_n(P_n^{-1}A_nP_n)) = (P_1^{-1}f_1(A_1)P_1 \mid \ldots \mid P_n^{-1}f_n(A_n)P_n)$ for every super polynomial $f = (f_1 \mid \ldots \mid f_n)$.

Yet another basic remark which we should make about minimal super polynomials of square super diagonal square matrices is that if A is a n × n square super diagonal square matrix of orders $n_1 \times n_1, \ldots, n_n \times n_n$ with entries in the field F. Suppose $F_1$ is a field which contains F as a subfield we may regard A as a square super diagonal square matrix either over F or over $F_1$ it may so appear that we obtain two different super minimal polynomials for A.

Fortunately that is not the case and the reason for it is if we find out what is the definition of super minimal polynomial for A, regarded as a square super diagonal square matrix over the field F. We consider all monic super polynomials with coefficients in F which super annihilate A, and we choose one with the least super degree.

If $f = (f_1 \mid \ldots \mid f_n)$ is a monic super polynomial over F.

$$\begin{aligned} f &= (f_1 \mid \ldots \mid f_n) \\ &= (x^{k_1} + \sum_{j_1=0}^{k_1-1} a_{j_1}^1 x^{j_1} \mid \ldots \mid x^{k_n} + \sum_{j_n=0}^{k_n-1} a_{j_n}^n x^{j_n}) \end{aligned}$$

then $f(A) = (f_1(A_1) \mid \ldots \mid f_n(A_n)) = (0 \mid 0 \mid \ldots \mid 0)$ merely say that we have a linear super relation between the power of A i.e.,

$$(A_1^{k_1} + a_{k_1-1}^1 A_1^{k_1-1} + \ldots + a_1^1 A_1 + a_0^1 I_1 \mid$$
$$\ldots \mid A_n^{k_n} + a_{k_n-1}^n A_n^{k_n-1} + \ldots + a_1^n A_n + a_0^n I_n) = (0 \mid 0 \mid \ldots \mid 0) \ldots \text{I}$$

The super degree of the minimal super polynomial is the least super positive degree ($k_1$ | … | $k_n$) such that there is a linear super relation of the above form I in $(I_1 A_1, \ldots, A_1^{k_1-1}; \ldots; I_n, A_n, \ldots, A_n^{k_n-1})$. Furthermore by the uniquiness of the minimal super polynomial there is for that ($k_1$, …, $k_n$) one and only one relation mentioned in I, i.e., once the minimal ($k_1$, …, $k_n$) is determined there are unique set of scalars $(a_0^1 \ldots a_{k_1-1}^1, \ldots, a_0^n \ldots a_{k_n-1}^n)$ in F such that I holds good. They are the coefficients of the minimal super polynomial. Now for each n-tuple ($k_1$, …, $k_n$) we have in I a system of $(n_1^2, \ldots, n_n^2)$ linear equations for the unknowns $(a_0^1, \ldots, a_{k_1-1}^1, \ldots, a_0^n, \ldots, a_{k_n-1}^n)$. Since the entries of A lie in F the coefficients of the system of equations in I lie in F. Therefore if the system has a super solution with $a_0^1, \ldots, a_{k_1-1}^1, \ldots, a_0^n, \ldots, a_{k_n-1}^n$ in $F_1$ it has a solution with $a_0^1, \ldots, a_{k_1-1}^1, \ldots, a_0^n, \ldots, a_{k_n-1}^n$ in F. Thus it must be now clear that the two super minimal polynomials are the same.

Now we prove an interesting theorem about the linear operator $T_s$.

**THEOREM 1.4.3:** *Let $T_s$ be a linear operator on an ($n_1$, …, $n_n$) dimensional super vector space V = ($V_1$ | … | $V_n$) or [let A be an ($n_1 \times n_1$, …, $n_n \times n_n$) square super diagonal square matrix]. The characteristic super polynomial and minimal super polynomials for $T_s$ (for A) have the same super roots expect for multiplicities.*

*Proof:* Let p = ($p_1$ | … | $p_n$) be the minimal super polynomial for $T_s$ = ($T_1$ | … | $T_n$) i.e., $p_i$ is the minimal polynomial of $T_i$; i = 1, 2, …, n. Let c = ($c_1$ | … | $c_n$) be a scalar. We want to prove p(c) = ($p_1(c_1)$ | … | $p_n$ ($c_n$)) = (0 | 0 | … | 0) if and only if c is the characteristic super value for $T_s$.

First we suppose $p(c) = (p_1(c_1) \mid \ldots \mid p_n(c_n)) = (0 \mid 0 \mid \ldots \mid 0)$. Then $p = (p_1 \mid \ldots \mid p_n) = (x - c)\ q = ((x - c_1)q_1 \mid \ldots \mid (x - c_n)q_n)$ where $q = (q_1 \mid \ldots \mid q_n)$ is a super polynomial. Since super degree q < super degree p ie $(\deg q_1 \mid \ldots \mid \deg q_n) < (\deg p_1 \mid \ldots \mid \deg p_n)$, the definition of the minimal super polynomial p tells us that $q(T_s) = (q_1(T_1) \mid \ldots \mid q_n(T_n)) \neq (0 \mid \ldots \mid 0)$. Choose a super vector $\beta = (\beta_1 \mid \ldots \mid \beta_n)$ such that

$$q(T_s)\ \beta = (q_1(T_1)\ \beta_1 \mid \ldots \mid q_n(T_n)\beta_n) \neq (0 \mid 0 \mid \ldots \mid 0).$$

Let $\alpha = (\alpha_1 \mid \ldots \mid \alpha_n)$

$= (q_1(T_1)\ \beta_1 \mid \ldots \mid q_n(T_n)\ \beta_n)$.

Then

$$\begin{aligned}(0 \mid 0 \mid \ldots \mid 0) &= p(T_s)\ \beta \\ &= (p_1(T_1)\ \beta_1 \mid \ldots \mid p_n(T_n)\ \beta_n) \\ &= (T_s - cI)\ q(T_s)\ \beta \\ &= ((T_1 - c_1I_1)\ q_1\ (T_1)\ \beta_1 \mid \ldots \mid (T_n - c_nI_n)\ q_n(T_n)\ \beta_n) \\ &= ((T_1 - c_1I_1)\ \alpha_1 \mid \ldots \mid (T_n - c_nI_n)\ \alpha_n) \\ &= (T - cI)\ \alpha\end{aligned}$$

and this c is a characteristic super value of $T_s$. Now suppose that c is a characteristic super value of T say $T\alpha = c\alpha$ ie $(T_1\alpha_1 \mid \ldots \mid T_n\alpha_n) = (c_1\alpha_1 \mid \ldots \mid c_n\alpha_n)$ with $\alpha = (\alpha_1 \mid \ldots \mid \alpha_n) \neq (0 \mid \ldots \mid 0)$.

As noted in earlier lemma $p(T_s)\ \alpha = p(c)\ \alpha$

i.e., $(p_1(T_1)\ \alpha_1 \mid \ldots \mid p_n\ (T_n)\ \alpha_n) = (p_1(c_1)\ \alpha_1 \mid \ldots \mid p_n(c_n)\ \alpha_n)$.
Since

$$p(T_s) = (p_1(T_1) \mid \ldots \mid p_n(T_n)) = (0 \mid \ldots \mid 0)$$

and $\alpha = (\alpha_1 \mid \ldots \mid \alpha_n) \neq (0 \mid \ldots \mid 0)$ we have

$$p(c) = (p_1(c_1) \mid \ldots \mid p_n(c_n)) = (0 \mid \ldots \mid 0).$$

Let $T_s = (T_1 \mid \ldots \mid T_n)$ be a diagonalizable linear operator and let $(c_1^1 \ldots c_{k_1}^1)$, $(c_1^2 \ldots c_{k_2}^2)$, …, $(c_1^n \ldots c_{k_n}^n)$ be the distinct characteristic super values of $T_s = (T_1 \mid \ldots \mid T_n)$. Then it is easy to see that the minimal super polynomial for $T_s$ is the minimal polynomial.

$P = (p_1 \mid \ldots \mid p_n)$
$= ((x - c_1^1) \ldots (x - c_{k_1}^1) \mid \ldots \mid (x - c_1^n) \ldots (x - c_{k_n}^n))$.

If $\alpha = (\alpha_1 \mid \ldots \mid \alpha_n)$ is the characteristic super vector, then one of the operators $T_s - c_1 I, \ldots, T_s - c_k I$ sends $\alpha$ into $(0 \mid \ldots \mid 0)$ i.e.,
$(T_1 - c_1^1 I_1 \ldots T_1 - c_{k_1}^1 I_1) \mid \ldots \mid (T_n - c_1^n I_n \ldots T_n - c_{k_n}^n I_n)$
sends $\alpha = (\alpha_1 \mid \ldots \mid \alpha_n)$ into $(0 \mid \ldots \mid 0)$.
Therefore,

$$(T_1 - c_1^1 I_1) \ldots (T_1 - c_{k_1}^1 I_1)\, \alpha_1 = 0$$
$$(T_2 - c_1^2 I_2) \ldots (T_2 - c_{k_2}^2 I_2)\, \alpha_2 = 0$$

so on

$$(T_n - c_1^n I_n) \ldots (T_n - c_{k_n}^n I_n)\, \alpha_n = 0;$$

for every characteristic super vector $\alpha = (\alpha_1 \mid \ldots \mid \alpha_n)$.
There is a super basis for the underlying super space which consists of characteristic super vectors of $T_s$, hence

$$\begin{aligned} p(T_s) &= (p_1(T_1) \mid \ldots \mid p_n(T_n)) \\ &= ((T_1 - c_1^1 I_1) \ldots (T_1 - c_{k_1}^1 I_1) \mid \ldots \mid (T_n - c_1^n I_n) \ldots (T_n - c_{k_n}^n I_n)) \\ &= (0 \mid \ldots \mid 0). \end{aligned}$$

Thus we have concluded if $T_s$ is a diagonalizable operator then the minimal super polynomial for $T_s$ is a product of distinct linear factors.

Now we will indicate the proof of the Cayley Hamilton theorem for linear operators $T_s$ on a super vector space V.

**THEOREM 1.4.4: (CAYLEY HAMILTON):** *Let $T_s$ be a linear operator on a finite dimensional vector space $V = (V_1 \mid \ldots \mid V_n)$. If $f = (f_1 \mid \ldots \mid f_n)$ is the characteristic super polynomial for $T_s = (T_1 \mid \ldots \mid T_n)$ ($f_i$ the characteristic polynomial for $T_i$, i = 1, 2, …,* *n.) then f(T) = (f₁(T₁) | … | fₙ(Tₙ)) = (0 | 0 | … | 0); in other words, the minimal super polynomial divides the characteristic super polynomial for T.*

*Proof:* The proof is only indicated. Let $K[T_s] = (K[T_1] \mid \ldots \mid K[T_n])$ be the super commutative ring with identity consisting of all polynomials in $T_1, \ldots, T_n$ of $T_s$. i.e., $K[T_s]$ can be visualized as a commutative super algebra with identity over the scalar field F. Choose a super basis

$$\{(\alpha_1^1 \ldots \alpha_{n_1}^1) \mid \ldots \mid (\alpha_1^n \ldots \alpha_{n_n}^n)\}$$

for the super vector space $V = (V_1 \mid \ldots \mid V_n)$ and let A be the super diagonal square matrix which represents $T_s$ in the given basis. Then $T_s \alpha_i$

$$= (T_1(\alpha_{i_1}^1) \mid \ldots \mid T_n(\alpha_{i_n}^n))$$

$$= \left( \sum_{j_1=1}^{n_1} A_{j_1 i_1}^1 \alpha_{j_1}^1 \mid \ldots \mid \sum_{j_n=1}^{n_n} A_{j_n i_n}^n \alpha_{j_n}^n \right);$$

$1 \leq j_t \leq n_t$; $t = 1, 2, \ldots, n$.

These equations may be written as in the equivalent form

$$\left( \sum_{j_1=1}^{n_1} (\delta_{i_1 j_1} T_1 - A_{j_1 i_1}^1 I_1)\, \alpha_{j_1}^1 \mid \ldots \mid \ldots \mid \sum_{j_n=1}^{n_n} (\delta_{i_n j_n} T_n - A_{j_n i_n}^n I_n)\, \alpha_{j_n}^n \right)$$

$$= (0 \mid 0 \mid \ldots \mid 0),\ 1 \leq i_t \leq n_t;\ t = 1, 2, \ldots, n.$$

Let $B = (B^1 \mid \ldots \mid B^n)$, we may call as notational blunder and yet denote the element of $(K^{n_1 \times n_1} \mid \ldots \mid K^{n_n \times n_n})$ with entries $B_{ij} = (B_{i_1 j_1}^1 \mid \ldots \mid B_{i_n j_n}^n)$

$$= ((\delta_{i_1 j_1} T_1 - A_{j_1 i_1}^1 I_1) \mid \ldots \mid (\delta_{i_n j_n} T_n - A_{j_n i_n}^n I_n)) \text{ when } n = 2.$$

$$B = \left(\begin{pmatrix} T_1 - A^1_{11}I_1 & -A^1_{21}I_1 \\ -A^1_{12}I_1 & T_1 - A^1_{22}I_1 \end{pmatrix}\right) | \ldots | \left(\begin{pmatrix} T_n - A^n_{11}I_n & -A^n_{21}I_n \\ -A^n_{12}I_n & T_n - A^n_{22}I_n \end{pmatrix}\right)$$

(notational blunder)
and super det B =

$$\begin{aligned} B &= ((T_1 - A^1_{11}I_1)(T_1 - A^1_{22}I_1) - A^1_{12}A^1_{21}I_1)) \\ &\ldots | \; ([(T_n - A^n_{11}I_n)(T_n - A^n_{22}I_n) - A^n_{12}A^n_{21}I_n]) \\ &= [(T_1^2 - (A^1_{11} + A^1_{22})\, T_1 + (A^1_{11}A^1_{22} - A^1_{12}A^1_{21})I_1) | \ldots | \\ &\quad (T_n^2 - (A^n_{11} + A^n_{22})\, T_n + (A^n_{11}A^n_{22} - A^n_{12}A^n_{21})I_n)] \\ &= [f_1(T_1) \mid \ldots \mid f_n(T_n)] \\ &= f(T), \end{aligned}$$

where $f = (f_1 \mid \ldots \mid f_n)$ is the characteristic super polynomial, $f = (f_1 \mid \ldots \mid f_n) = ((x^2 - (\text{trace } A_1)\, x + \det A_1) \mid \ldots \mid (x^2 - (\text{trace } A_n)\, x + \det A_n))$.
For $n > 2$ it is also clear that $f(T) = (f_1(T_1) \mid \ldots \mid f_n(T_n))$ = super det B = $(\det B_1 \mid \ldots \mid \det B_n)$, since $f = (f_1 \mid \ldots \mid f_n)$ is the super determinant of the super diagonal square matrix $xI - A = ((xI_1 - A_1) \mid \ldots \mid (xI_n - A_n))$ whose entries are the super polynomials.

$$\begin{aligned} &(xI - A)_{ij} = ((xI_1 - A_1)_{i_1j_1} | \ldots | (xI_n - A_n)_{i_nj_n}) \\ &= ((\delta_{i_1j_1} x - A^1_{j_1i_1}) | \ldots | (\delta_{i_nj_n} x - A^n_{j_ni_n})); \end{aligned}$$

we wish to show that $f(T) = (f_1(T_1) \mid \ldots \mid f_n(T_n)) = (0 \mid \ldots \mid 0)$.
In order that $f(T) = (f_1(T_1) \mid \ldots \mid f_n(T_n))$ is the zero super operator it is necessary and sufficient that (super det B) $\alpha_k$

$$= ((\det B_1)_{\alpha^1_{k_1}} | \ldots | (\det B_n)_{\alpha^n_{k_n}})$$

$= (0 \mid \ldots \mid 0)$ for $k_t = 1, \ldots, n_t$; $t = 1, 2, \ldots, n$.

By the definition of B the super vectors

$$(\alpha^1_1, \ldots, \alpha^1_{n_1}), \ldots, (\alpha^n_1, \ldots, \alpha^n_{n_n})$$

satisfy the equations;

$$\left(\sum_{j_1=1}^{n_1} B^1_{i_1 j_1}\alpha^1_{j_1} \mid \ldots \mid \sum_{j_n=1}^{n_n} B^n_{i_n j_n}\alpha^n_{j_n}\right)$$

$= ( 0 \mid \ldots \mid 0), 1 \le i_t \le n_t; t = 1, 2, \ldots, n.$

When n = 2 it is suggestive

$$\left(\begin{pmatrix} T_1 - A^1_{11}I_1 & -A^1_{21}I_1 \\ -A^1_{12}I_1 & T_1 - A^1_{22}I_1 \end{pmatrix}\begin{pmatrix}\alpha^1_1 \\ \alpha^1_2\end{pmatrix}\right.$$

$$\left.\mid\ldots\mid\begin{pmatrix} T_n - A^n_{11}I_n & -A^n_{21}I_n \\ -A^n_{12}I_n & T_n - A^n_{22}I_n \end{pmatrix}\begin{pmatrix}\alpha^n_1 \\ \alpha^n_2\end{pmatrix}\right) = \left(\begin{pmatrix}0\\0\end{pmatrix}\mid\ldots\mid\begin{pmatrix}0\\0\end{pmatrix}\right).$$

In this case the classical super adjoint B is the super diagonal matrix $\tilde{B} = [\tilde{B}_1 \mid \ldots \mid \tilde{B}_n]$

$$= \left(\begin{array}{c|c|c|c} \tilde{B}_1 & 0 & & 0 \\ \hline 0 & \tilde{B}_2 & & 0 \\ \hline & & & \\ \hline 0 & 0 & & \tilde{B}_n \end{array}\right)$$

(once again with notational blunder!)

$$\tilde{B} =$$

$$\left(\begin{array}{c|c|c|c} \begin{matrix} T_1 - A^1_{22}I_{11} & A^1_{21}I \\ A^1_{12} & T_1 - A^1_{11}I_1 \end{matrix} & 0 & & 0 \\ \hline 0 & \begin{matrix} T_2 - A^2_{22}I_2 & A^2_{21}I_2 \\ A^2_{12} & T_2 - A^2_{11}I_2 \end{matrix} & & 0 \\ \hline & & & \\ \hline 0 & 0 & & \begin{matrix} T_n - A^n_{22}I_n & A^n_{21}I_n \\ A^n_{12} & T_n - A^n_{11}I_n \end{matrix} \end{array}\right)$$

and

$$\tilde{B}\, B = \begin{pmatrix} \text{super det } B & 0 \\ 0 & \text{superdet } B \end{pmatrix}.$$

$$= \left(\begin{array}{cc|c|c|c} \det B_1 & 0 & 0 & 0 & 0 \\ 0 & \det B_1 & & & \\ \hline 0 & & \begin{matrix} \det B_2 & 0 \\ 0 & \det B_2 \end{matrix} & & 0 \\ \hline 0 & & & & \\ \hline 0 & & 0 & 0 & \begin{matrix} \det B_n & 0 \\ 0 & \det B_n \end{matrix} \end{array}\right).$$

Hence we have

$$\text{super } B \begin{pmatrix} \alpha_1 \\ \alpha_2 \end{pmatrix} = \left(\begin{array}{c|c|c|c} \det B_1 \begin{pmatrix} \alpha_1^1 \\ \alpha_2^1 \end{pmatrix} & 0 & & 0 \\ \hline 0 & \det B_2 \begin{pmatrix} \alpha_1^2 \\ \alpha_2^2 \end{pmatrix} & & 0 \\ \hline & & & \\ \hline 0 & 0 & & \det B_n \begin{pmatrix} \alpha_1^n \\ \alpha_2^n \end{pmatrix} \end{array}\right)$$

$$= \tilde{B}\ B \begin{pmatrix} \alpha_1 \\ \alpha_2 \end{pmatrix}.$$

$$= \left(\begin{array}{c|c|c|c} \tilde{B}_1 B_1 \begin{pmatrix} \alpha_1^1 \\ \alpha_2^1 \end{pmatrix} & 0 & & 0 \\ \hline 0 & \tilde{B}_2 B_2 \begin{pmatrix} \alpha_1^2 \\ \alpha_2^2 \end{pmatrix} & & 0 \\ \hline & & & \\ \hline 0 & 0 & & \tilde{B}_n B_n \begin{pmatrix} \alpha_1^n \\ \alpha_2^n \end{pmatrix} \end{array}\right)$$

$$= \tilde{B}\ B \begin{pmatrix} \alpha_1 \\ \alpha_2 \end{pmatrix}$$

$$= \left(\begin{array}{c|c|c|c} \tilde{B}_1 B_1 \begin{pmatrix} \alpha_1^1 \\ \alpha_2^1 \end{pmatrix} & 0 & & 0 \\ \hline 0 & \tilde{B}_2 B_2 \begin{pmatrix} \alpha_1^2 \\ \alpha_2^2 \end{pmatrix} & & 0 \\ \hline 0 & & & 0 \\ \hline 0 & 0 & & \tilde{B}_n B_n \begin{pmatrix} \alpha_1^2 \\ \alpha_2^2 \end{pmatrix} \end{array}\right) = \left( \begin{pmatrix} 0 \\ 0 \end{pmatrix} | \ldots | \begin{pmatrix} 0 \\ 0 \end{pmatrix} \right).$$

$\tilde{B}$ = super adj B

$$\text{i.e., } \left(\begin{array}{c|c|c|c} B_1 & 0 & 0 & 0 \\ \hline 0 & B_2 & 0 & 0 \\ \hline 0 & 0 & & 0 \\ \hline 0 & 0 & 0 & B_n \end{array}\right) =$$

$$\left(\begin{array}{c|c|c|c} \text{adj } B_1 & 0 & & 0 \\ \hline 0 & \text{adj } B_2 & & 0 \\ \hline 0 & & & 0 \\ \hline 0 & 0 & & \text{adj } B_n \end{array}\right)$$

$$\sum_{j=1}^{n} \tilde{B}_{ki} B_{ij} \alpha^i = \left( \sum_{j_1=1}^{n_1} \tilde{B}_{1_{k_1 i_1}} B_{1 i_j j_1} \alpha_{i_1}^1 \ | \ldots | \ \sum_{j_n=1}^{n_n} \tilde{B}_{n_{k_n i_n}} B_{n\, i_n j_n} \alpha_{i_n}^n \right)$$

$$= (0 \mid \ldots \mid 0)$$

for each pair $k_t$, $i_t$ ; $1 \le t \le n_t$; $t = 1, 2, \ldots, n$.

Thus we can prove as in case of usual vector spaces super $(\det B)_{\alpha_k} = ((\det B_1)_{\alpha_1} \mid \ldots \mid (\det B_n)_{\alpha_n}) = (0 \mid \ldots \mid 0)$.

As in case of ordinary matrices the main use of Cayley Hamilton theorem for super diagonal square matrices is to search for minimal super polynomial of various operators.

We know a super diagonal square matrix A which represents $T_s = (T_1 \mid \ldots \mid T_n)$ in some ordered super basis, then we can calculate the characteristic super polynomial, $f = (f_1 \mid \ldots \mid f_n)$. We know that the minimal super polynomial $p = (p_1 \mid \ldots \mid p_n)$ super divides $f = (f_1 \mid \ldots \mid f_n)$ (i.e., each $p_i$ divides $f_i$; i = 1, 2, …, n then we say $p = (p_1 \mid \ldots \mid p_n)$ super divides $f = (f_1 \mid \ldots \mid f_n)$).

We know when a polynomial $p_i$ divides $f_i$ the two polynomials have same roots for i = 1, 2, …, n. There is no method of finding precisely the roots of a polynomial more so it is still a open problem to find precisely the roots of a super polynomial (unless its super degrees are small) however $f = (f_1 \mid f_2 \mid \ldots \mid f_n)$ factors as

$$((x-c_1^1)^{d_1^1} \ldots (x-c_{k_1}^1)^{d_{k_1}^1} \mid \ldots \mid (x-c_1^n)^{d_1^n} \ldots (x-c_{k_n}^n)^{d_{k_n}^n}).$$

$$(c_1^1 \ldots c_{k_1}^1 \mid \ldots \mid c_1^n, \ldots c_{k_n}^n)$$

are super distinct i.e., we demand only $(c_1^t \ldots c_{k_t}^t)$ to be distinct and $d_{i_t}^t \geq 1$ for every t, t = 1, 2, …, n, then

$$p = (p_1 \mid \ldots \mid p_n)$$
$$= ((x-c_1^1)^{r_1^1} \ldots (x-c_{k_1}^1)^{r_{k_1}^1} \mid \ldots \mid (x-c_1^n)^{r_1^n} \ldots (x-c_{k_n}^n)^{r_{k_n}^n})$$

$1 \leq r_{j_t}^t \leq d_{j_t}^t$, for every t = 1, 2, …, n. That is all we can say in general.

If $f = (f_1 \mid \ldots \mid f_n)$ is a super polynomial given above has super degree $(n_1 \mid \ldots \mid n_n)$ then every super polynomial p, given; we can find an $(n_1 \times n_1, \ldots, n_n \times n_n)$ super diagonal square matrix A =

$$\left(\begin{array}{c|c|c|c} A_1 & 0 & \cdots & 0 \\ \hline 0 & A_2 & 0 & 0 \\ \hline \vdots & & & 0 \\ \hline 0 & 0 & 0 & A_n \end{array}\right)$$

with $A_i$, a $n_i \times n_i$ matrix; i = 1, 2, …, n, which has $f_i$ as its characteristic polynomial and $p_i$ as its minimal polynomial. Now we proceed onto define the notion of super invariant subspaces or an invariant super subspaces.

**DEFINITION 1.4.8**: *Let $V = (V_1| \ldots |V_n)$ be a super vector space and $T_s = (T_1 | \ldots | T_n)$ be a linear operator in V. If $W = (W_1 | \ldots | W_n)$ be a super subspace of V; we say that $W = (W_1 | \ldots | W_n)$ is super invariant under T if for each super vector $\alpha = (\alpha_1 | \ldots |\alpha_n)$ in $W = (W_1 | \ldots | W_n)$ the super vector $T_s(\alpha)$ is in $W = (W_1 | \ldots | W_n)$ i.e. if $T_s(W)$ is contained in W. When the super subspace $W = (W_1 | \ldots | W_n)$ is super invariant under the operator $T_s = (T_1 | \ldots | T_n)$ then $T_s$ induces a linear operator $(T_s)_W$ on the super subspace $W = (W_1 | \ldots | W_n)$. The linear operator $(T_s)_W$ is defined by $(T_s)_W (\alpha) = T_s(\alpha)$ for α in $W = (W_1 | \ldots | W_n)$ but $(T_s)_W$ is a different object from $T_s = (T_1 | \ldots | T_n)$ since its domain is W not V.*

When $V = (V_1 | \ldots | V_n)$ is finite $(n_1, \ldots, n_n)$ dimensional, the invariance of $W = (W_1 | \ldots | W_n)$ under $T_s = (T_1 | \ldots | T_n)$ has a simple super matrix interpretation and perhaps we should mention it at this point. Suppose we choose an ordered basis B $= (B_1 | \ldots | B_n) = (\alpha_1^1 \ldots \alpha_{n_1}^1 | \ldots | \alpha_1^n \ldots \alpha_{n_n}^n)$ for $V = (V_1 | \ldots | V_n)$ such that $B' = (\alpha_1^1 \ldots \alpha_{r_1}^1 | \ldots | \alpha_1^n \ldots \alpha_{r_n}^n)$ is an ordered basis for $W = (W_1 | \ldots | W_n)$; super dim $W = (r_1, \ldots, r_n)$. Let $A = [T_s]_B$ so that

$$T_s \alpha_j = \left[ \sum_{i_1=1}^{n_1} A_{i_1 j_1}^1 \alpha_{i_1}^1 \;|\; \ldots \;|\; \sum_{i_n=1}^{n_n} A_{i_n j_n}^n \alpha_{i_n}^n \right].$$

Since $W = (W_1 \mid \ldots \mid W_n)$ is super invariant under $T_s = (T_1 \mid \ldots \mid T_n)$ and the vector $T_s\alpha_j = (T_1\alpha^1_{j_1} \mid \ldots \mid T_n\alpha^n_{j_n})$ belongs to $W = (W_1 \mid \ldots \mid W_n)$ for $j_t \le r_t$. This means that

$$T_s\alpha_j = \left[\sum_{i_1=1}^{r_1} A^1_{i_1 j_1}\alpha^1_{i_1} \mid \ldots \mid \sum_{i_n=1}^{r_n} A^n_{i_n j_n}\alpha^n_{i_n}\right]$$

$j_t \le r_t$; $t = 1, 2, \ldots, n$. In other words $A^t_{i_t j_t} = (A^1_{i_1 j_1} \mid \ldots \mid A^n_{i_n j_n}) = (0 \mid \ldots \mid 0)$ if $j_t \le r_t$ and $i_t > r_t$.

$$A = \left(\begin{array}{c|c|c|c} A_1 & 0 & \ldots & 0 \\ \hline 0 & A_2 & 0 & 0 \\ \hline \vdots & 0 & & 0 \\ \hline 0 & 0 & 0 & A_n \end{array}\right)$$

$$= \left(\begin{array}{c|c|c|c} \begin{matrix} B_1 & C_1 \\ 0 & D_1 \end{matrix} & 0 & 0 & 0 \\ \hline 0 & \begin{matrix} B_2 & C_2 \\ 0 & D_2 \end{matrix} & & 0 \\ \hline 0 & & & 0 \\ \hline 0 & 0 & 0 & \begin{matrix} B_n & C_n \\ 0 & D_n \end{matrix} \end{array}\right)$$

where $B_t$ is an $r_t \times r_t$ matrix, $C_t$ is a $r_t \times (n_t - r_t)$ matrix and $D_t$ is an $(n_t - r_t) \times (n_t - r_t)$ matrix $t = 1, 2, \ldots, n$.

In view of this we prove the following interesting lemma.

**LEMMA 1.4.4:** *Let $W = (W_1 \mid \ldots \mid W_n)$ be an invariant super subspace for $T_s = (T_1 \mid \ldots \mid T_n)$. The characteristic super polynomial for the restriction operator $(T_s)_W = ((T_1)_{W_1} \mid \ldots \mid (T_n)_{W_n})$ divides the characteristic super polynomial for $T_s$. The minimal super polynomial for*

$(T_s)_w = ((T_1)_{W_1} \mid \ldots \mid (T_n)_{w_n})$ *divides the minimal super polynomial for $T_s$.*

*Proof:* We have $[T_s]_B = A$ where $B = \{B_1 \ldots B_n\}$ is a super basis for $V = (V_1 \mid \ldots \mid V_n)$; with $B_i = \{\alpha_1^i \ldots \alpha_{n_i}^i\}$ a basis for $V_i$, this is true for each i, i = 1, 2, …, n. A is a super diagonal square matrix of the form

$$A = \left(\begin{array}{c|c|c|c} A_1 & 0 & \ldots & 0 \\ \hline 0 & A_2 & 0 & 0 \\ \hline \ldots & & & 0 \\ \hline 0 & 0 & 0 & A_n \end{array}\right)$$

where each

$$A_i = \begin{pmatrix} B_i & C_i \\ 0 & D_i \end{pmatrix}$$

for i = 1, 2, …, n; i.e.

$$A = \left(\begin{array}{c|c|c|c} \begin{matrix} B_1 & C_1 \\ 0 & D_1 \end{matrix} & 0 & 0 & 0 \\ \hline 0 & \begin{matrix} B_2 & C_2 \\ 0 & D_2 \end{matrix} & 0 & 0 \\ \hline 0 & & & 0 \\ \hline 0 & 0 & 0 & \begin{matrix} B_n & C_n \\ 0 & D_n \end{matrix} \end{array}\right)$$

and $B = [(T_s)_w]_{B'}$ where B′ is a basis for the super vector subspace $W = (W_1 \mid \ldots \mid W_n)$ and B is a super diagonal square matrix; i.e.

$$B = \left(\begin{array}{c|c|c|c} B_1 & 0 & \ldots & 0 \\ \hline 0 & B_2 & & 0 \\ \hline \vdots & & & 0 \\ \hline 0 & 0 & 0 & B_n \end{array}\right).$$

Now using the block form of the super diagonal square matrix we have super det (xI – A) = super det (xI – B) × super det (xI – D)

i.e. $(\det (xI_1 - A_1) \mid \ldots \mid \det (xI_n - A_n))$

$$= (\det (xI_1' - B_1) \det (xI_1'' - D_1) \mid \ldots \mid \det (xI_n' - B_n) \det (xI_n'' - D_n)).$$

This proves the restriction operator $(T_s)_W$ super divides the characteristic super polynomial for $T_s$. The minimal super polynomial for $(T_s)_W$ super divides the minimal super polynomial for $T_s$.

It is pertinent to observe that $I_1', I_1'', I_1, \ldots, I_n$ represents different identities i.e. of different order.

The $K^{th}$ row of A has the form

$$A^K = \left(\begin{array}{cc|cc|c|cc}
B_1^{K_1} & C_1^{K_1} & \multicolumn{2}{c|}{0} & & \multicolumn{2}{c}{0} \\
0 & D_1^{K_1} & & & & & \\
\hline
\multicolumn{2}{c|}{0} & B_2^{K_2} & C_2^{K_2} & & \multicolumn{2}{c}{0} \\
& & 0 & D_2^{K_2} & & & \\
\hline
\multicolumn{2}{c|}{\vdots} & & & & & \\
\hline
\multicolumn{2}{c|}{0} & \multicolumn{2}{c|}{0} & \cdots & B_n^{K_n} & C_n^{K_n} \\
& & & & & 0 & D_n^{K_n}
\end{array}\right)$$

where $C_t^{K_t}$ is some $r_t \times (n_t - r_t)$ matrix; true for t = 1, 2, …, n. Thus any super polynomial which super annihilates A also super annihilates D. Thus our claim made earlier that, the minimal super polynomial for B super divides the minimal super polynomial for A is established.

Thus we say a super subspace $W = (W_1 \mid \ldots \mid W_n)$ of the super vector space $V = (V_1 \mid \ldots \mid V_n)$ is super invariant under $T_s = (T_1 \mid \ldots \mid T_n)$ if $T_s(W) \subseteq W$ i.e. each $T_i(W_i) \subseteq W_i$; for i = 1, 2, …, n i.e. if $\alpha = (\alpha_1 \mid \ldots \mid \alpha_n) \in W$

then $T_s\alpha = (T_1\alpha_1 \mid \ldots \mid T_n\alpha_n)$ where

$$\alpha_1 = x_1^1\alpha_1^1 + \ldots + x_{r_1}^1\alpha_{r_1}^1 \;;$$
$$\alpha_2 = x_1^2\alpha_1^2 + \ldots + x_{r_2}^2\alpha_{r_2}^1$$

and so on

$$\alpha_n = x_1^n\alpha_1^n + \ldots + x_{r_n}^n\alpha_{r_n}^n \,.$$
$$T_s\alpha = (t_1^1 x_1^1\alpha_1^1 + \ldots + t_{r_1}^1 x_{r_1}^1\alpha_{r_1}^1 \mid \ldots \mid t_1^n x_1^n\alpha_1^n + \ldots + t_{r_n}^n x_{r_n}^n\alpha_{r_n}^n)\,.$$

Now B described in the above theorem is a super diagonal matrix given by

$$B = \left(\begin{array}{c|c|c|c}
\begin{matrix} t_1^1 & 0 & \ldots & 0 \\ 0 & t_2^1 & \ldots & 0 \\ \vdots & \vdots & & \vdots \\ 0 & 0 & \ldots & t_{r_1}^1 \end{matrix} & 0 & \ldots & 0 \\
\hline
0 & \begin{matrix} t_1^2 & 0 & \ldots & 0 \\ 0 & t_2^2 & \ldots & 0 \\ \vdots & \vdots & & \vdots \\ 0 & 0 & \ldots & t_{r_2}^2 \end{matrix} & & 0 \\
\hline
 & & & \\
\hline
0 & 0 & & \begin{matrix} t_1^n & 0 & \ldots & 0 \\ 0 & t_2^n & \ldots & 0 \\ \vdots & \vdots & & \vdots \\ 0 & 0 & \ldots & t_{r_n}^n \end{matrix}
\end{array}\right)$$

Thus the characteristic super polynomial of B i.e. $(T_s)_W)$ is

$$g = (g_1 \mid \ldots \mid g_n) =$$
$$((x - c_1^1)^{e_1^1} \ldots (x - c_{K_1}^1)^{e_{K_1}^1} \mid \ldots \mid (x - c_1^n)^{e_1^n} \ldots (x - c_{K_n}^n)^{e_{K_n}^n})$$

where $e_i^t = \dim W_i^t$ for i = 1, 2, …, $K_t$ and t = 1, 2, …, n.

Now we proceed onto define $T_s$ super conductor of any $\alpha$ into $W = (W_1 \mid \ldots \mid W_n)$.

**DEFINITION 1.4.6:** *Let $V = (V_1 \mid \dots \mid V_n)$ be a super vector space over the field F. $W = (W_1 \mid \dots \mid W_n)$ be an invariant super subspace of V for the linear operator $T_s = (T_1 \mid \dots \mid T_n)$ of V. Let $\alpha = (\alpha_1 \mid \dots \mid \alpha_n)$ be a super vector in V. The T-super conductor of α into W is the set $S_{T_s}(\alpha;W) = (S_{T_1}(\alpha_1;W_1) \mid \dots \mid S_{T_n}(\alpha_n;W_n))$ which consist of all super polynomials $g = (g_1 \mid \dots \mid g_n)$ (over the scalar field F) such that $g(T_s)\alpha$ is in W, i.e. $(g_1(T_1)\alpha_1 \mid \dots \mid g_n(T_n)\alpha_n) \in W = (W_1 \mid \dots \mid W_n)$. i.e. $g_i(T_i)\alpha_i) \in W_i$ for every i. Or we can equivalently define the $T_s$ – super conductor of α in W is a $T_i$ conductor of $\alpha_i$ in $W_i$ for every i = 1, 2, ..., n. Without loss in meaning we can for convenience drop $T_s$ and write the super conductor of α into W as $S(\alpha;W) = (S(\alpha_1;W_1) \mid \dots \mid S(\alpha_n;W_n))$.*

The collection of polynomials will be defined as super stuffer this implies that the super conductor, the simple super operator $g(T_s) = (g_1(T_1) \mid \dots \mid g_n(T_n))$ leads the super vector α into W. In the special case W = (0 | … | 0), the super conductor is called the $T_s$ super annihilator of α. The following important and interesting theorem is proved.

**THEOREM 1.4.5:** *Let $V = (V_1 \mid \dots \mid V_n)$ be a finite dimensional super vector space over the field F and let $T_s$ be a linear operator on V. Then $T_s$ is super diagonalizable if and only if the minimal super polynomial for $T_s$ has the form*

$$p = (p_1 \mid \dots \mid p_n) = ((x - c_1^1)\dots(x - c_{K_1}^1) \mid \dots \mid (x - c_1^n)\dots(x - c_{K_n}^n)]$$

*where $(c_1^1 \dots c_{K_1}^1 \mid \dots \mid c_1^n \dots c_{K_n}^n)$ are such that each set $c_1^t, \dots, c_{K_t}^t$ are distinct elements of F for t = 1, 2, ..., n.*

*Proof:* We have noted that if $T_s$ is super diagonalizable, its minimal super polynomial is a product of distinct linear factors. To prove the converse let $W = (W_1 \mid \dots \mid W_n)$ be the super subspace spanned by all of the characteristic super vectors of $T_s$ and suppose $W = (W_1 \mid \dots \mid W_n) \neq (V_1 \mid \dots \mid V_n)$ i.e. each $W_i \neq V_i$. By the earlier results proved there is a super vector α not in $W = (W_1 \mid \ldots \mid W_n)$ and a characteristic super value $c_j = (c^1_{j_1}, \ldots c^n_{j_n})$ of $T_s$ such that the super vector $\beta = (T - c_j I)\alpha$ i.e. $(\beta_1 \mid \ldots \mid \beta_n) = ((T_1 - c^1_{j_1} I_1)\alpha_1 \mid \ldots \mid (T_n - c^n_{j_n} I_n)\alpha_n)$ lies in $W = (W_1 \mid \ldots \mid W_n)$. Since $(\beta_1 \mid \ldots \mid \beta_n)$ is in W,

$$\beta = (\beta^1_1 + \ldots + \beta^1_{K_1} \mid \beta^2_1 + \ldots + \beta^2_{K_2} \mid \ldots \mid \beta^n_1 + \ldots + \beta^n_{K_n})$$

where $\beta_t = \beta^t_1 + \ldots + \beta^t_{K_t}$ for t = 1, 2, …, n with $T_s\beta_i = c_i\beta_i$; $1 \leq i \leq K$ i.e. $(T_1\beta^1_{i_1} \mid \ldots \mid T_n\beta^n_{i_n}) = (c^1_{i_1} \beta^1_{i_1} \mid \ldots \mid c^n_{i_n}\beta^n_{i_n})$; $(1 \leq i_t \leq K_t)$ and therefore the super vector

$$\begin{aligned} h(T_s)\beta &= (h_1(c^1_1)\beta^1_1 + \ldots + h_1(c^1_{K_1})\beta^1_{K_1} \mid \ldots \mid \\ &\qquad h_n(c^n_1)\beta^n_1 + \ldots + h_n(c^n_{K_n})\beta^n_{K_n}) \\ &= (h_1(T_1)\beta_1 \mid \ldots \mid h_n(T_n)\beta_n) \end{aligned}$$

is in $W = (W_1 \mid \ldots \mid W_n)$ for every super polynomial $h = (h_1 \mid \ldots \mid h_n)$.

Now $(x - c_j)$ q for some super polynomial q, where $p = (p_1 \mid \ldots \mid p_n)$ and $q = (q_1 \mid \ldots \mid q_n)$.

Thus $p = (x - c_j)$ q implies

$$\begin{aligned} p &= (p_1 \mid \ldots \mid p_n) \\ &= ((x - c^1_{j_1})q_1 \mid \ldots \mid (x - c^n_{j_n})q_n) \end{aligned}$$

i.e. $(q_1 - q_1(c^1_{j_1}) \mid \ldots \mid q_n - q_n(c^n_{j_n})) = ((x - c^1_{j_1})h_1 \mid \ldots \mid (x - c^n_{j_n})h_n)$.

We have

$$\begin{aligned} q(T_s)\alpha - q(c_j)\alpha &= (q_1(T_1)\alpha_1 - q_1(c^1_{j_1})\alpha_1 \mid \ldots \mid \\ &\qquad q_n(T_n)\alpha_n - q_n(c^n_{j_n})\alpha_n) \\ &= h(T_s)(T_s - c_j I)\alpha = h(T_s)\beta \\ &= (h_1(T_1)(T_1 - c^1_{j_1} I_1)\alpha_1 \mid \ldots \mid h_n(T_n)(T_n - c^n_{j_n} I_n)\alpha_n) \\ &= (h_1(T_1)\beta_1 \mid \ldots \mid h_n(T_n)\beta_n). \end{aligned}$$

But $h(T_s)\beta$ is in $W = (W_1 \mid \ldots \mid W_n)$ and since

$$0 = p(T_s)\alpha = (p_1(T_1)\alpha_1 \mid \ldots \mid p_n(T_n)\alpha_n)$$

$= (T_s - c_j I) q(T_s)\alpha$
$= ((T_1 - c^1_{j_1} I_1)\, q_1(T_1)\alpha_1 | \ldots | (T_n - c^n_{j_n} I_n)\, q_n(T_n)\alpha_n)$;

the vector $q(T_s)\alpha$ is in W. Therefore $q(c_j)\alpha$ is in W. Since $\alpha$ is not in W we have $q(c_j) = (q_1(c^1_{j_1}) | \ldots | q_n(c^n_{j_n})) = (0| \ldots |0)$. Thus contradicts the fact that $p = (p_1 | \ldots | p_n)$ has distinct roots.

If $T_s$ is represented by a super diagonal square matrix A in some super basis and we wish to know if $T_s$ is super diagonalizable. We compute the characteristic super polynomial $f = (f_1 | \ldots | f_n)$. If we can factor

$$f = (f_1 | \ldots | f_n)$$
$$= ((x - c^1_1)^{d^1_1} \ldots (x - c^1_{K_1})^{d^1_{K_1}} | \ldots | (x - c^n_1)^{d^n_1} \ldots (x - c^n_{K_n})^{d^n_{K_n}})$$

we have two different methods for determining whether or not T is super diagonalizable. One method is to see whether for each i = 1, 2, …, n we can find $d^t_{i_t}$ independent characteristic super vectors associated with the characteristic super values $c_{i_t}$. The other method is to check whether or not

$$(T_s - c_1 I) \ldots (T_s - c_k I) \text{ i.e. } ((T_1 - c^1_1 I_1) \ldots (T_1 - c^1_{K_1} I_1)$$
$$| \ldots | (T_n - c^n_1) I_n \ldots (T_n - c^n_{K_n} I_n))$$

is the super zero operator.

Several other interesting results in this direction can be derived. Now we proceed onto define the new notion of super independent subsuper spaces of a super vector space V.

**DEFINITION 1.4.10**: *Let $V = (V_1 | \ldots | V_n)$ be a super vector space over F. Let $W_1 = (W_1^1 | \ldots | W_1^n)$, $W_2 = (W_2^n | \ldots | W_2^n) \ldots$ $W_K = (W_K^1 | \ldots | W_K^n)$ be K super subspaces of V. We say $W_1, \ldots, W_K$ are super independent if $\alpha_1 + \ldots + \alpha_K = 0$; $\alpha_i \in W_i$ implies each $\alpha_i = 0$.*

$$\alpha_i = (\alpha_1^i \mid \ldots \mid \alpha_n^i) \in W_i = (W_i^1 | \ldots | W_i^n);$$

*true for i = 1, 2, ..., K. If $W_1$ and $W_2$ are any two super vector subspaces of V = ($V_1$ | ... | $V_n$), we say $W_1 = (W_1^1 \mid \ldots | W_1^n)$ and $W_2 = (W_2^1 | \ldots | W_2^n)$ are super independent if and only if $W_1 \cap W_2 = (W_1^1 \cap W_2^1 | \ldots | W_1^n \cap W_2^n)$ = (0 | 0 | ... |0). If $W_1$, $W_2$, ..., $W_K$ are K super subspaces of V we say $W_1$, $W_2$, ..., $W_K$ are independent if $W_1 \cap W_2 \cap \ldots \cap W_K$ = $(W_1^1 \cap W_2^1 \cap \ldots W_K^1 | \ldots | W_1^n \cap W_2^n \cap \ldots \cap W_K^n)$ = (0 | ... | 0) . The importance of super independence in super subspaces is mentioned below. Let*

$$\begin{aligned} W' &= W'_1 + \ldots + W'_k \\ &= (W_1^1 + \ldots + W_K^1 | \ldots | W_1^n + \ldots + W_K^n) \\ &= (W_1' \ | \ldots | W_n') \end{aligned}$$

*$W_i'$ is a subspace $V_i$ and $W_i' = W_1' + \ldots + W_K'$ true for i = 1, 2, ..., n. Each super vector α in W can be expressed as a sum $\alpha = (\alpha_1' \mid \ldots \mid \alpha_n') = ((\alpha_1^1 + \ldots + \alpha_K^1) | \ldots | \alpha_1^n + \ldots + \alpha_K^n))$ i.e. each $\alpha^t = \alpha_1^t + \ldots + \alpha_K^t$ ; $\alpha^t \in W_t$. If $W_1$, $W_2$, ..., $W_K$ are super independent, then that expression for α is unique; for if*

$$\alpha = (\beta_1 + \ldots + \beta_K) = (\beta_1^1 + \ldots + \beta_K^1 \mid \ldots | \beta_1^n + \ldots + \beta_K^n)$$

*$\beta_i \in W_i$; i = 1, 2, ..., K. $\beta_i = \beta_1^i + \ldots + \beta_n^i$ then*

$$\begin{aligned} \alpha - \alpha = (0|\ldots|0) = ((\alpha_1^1 - \beta_1^1) + \ldots + (\alpha_K^1 - \beta_K^1) | \ldots | \\ (\alpha_1^n - \beta_1^n) + \ldots + (\beta_K^n - \alpha_K^n)) \end{aligned}$$

*hence each $\alpha_i^t - \beta_i^t = 0$ ; 1 ≤ i ≤ K; t = 1, 2, ..., n. Thus $W_1$, $W_2$, ..., $W_K$ are super independent so we can operate with super vectors in W as K-tuples $((\alpha_1^1, \ldots, \alpha_K^1); \ldots, (\alpha_1^n, \ldots, \alpha_K^n)$; $\alpha_i^t \in W_t$ ; 1 ≤ i ≤ K; t = 1, 2, ..., n. in the same way we operate with $R^K$ as K-tuples of real numbers.*

**LEMMA 1.4.5**: *Let $V = (V_1 \mid \ldots \mid V_n)$ be a finite $(n_1, \ldots, n_n)$ dimensional super vector space. Let $W_1, \ldots, W_K$ be super subspaces of V and let $W = (W_1^1 + \ldots W_K^1 \mid \ldots \mid W_1^n + \ldots + W_K^n)$. The following are equivalent*

(a) *$W_1, \ldots, W_K$ are super independent.*
(b) *For each j; $2 \leq j \leq K$, we have $W_j \cap (W_1 + \ldots + W_{j-1}) = \{(0 \mid \ldots \mid 0)\}$*
(c) *If $B_i$ is a super basis of $W_i$, $1 \leq i \leq K$, then the sequence $B = (B_1 \ldots B_K)$ is a super basis for W.*

The proof is left as an exercise for the reader. In any or all of the conditions of the above stated lemma is true then the supersum $W = W_1 + \ldots + W_K = (W_1^1 + \ldots + W_K^1 \mid \ldots \mid W_1^n + \ldots + W_K^n)$ where $W_t = (W_t^1 \mid \ldots \mid W_t^n)$ is super direct or that W is a super direct sum of $W_1, \ldots, W_K$ i.e. $W = W_1 \oplus \ldots \oplus W_K$ i.e. $(W_1^1 \oplus \ldots \oplus W_K^1 \mid \ldots \mid W_1^n \oplus \ldots \oplus W_K^n)$. If each of the $W_i$ is $(1, \ldots, 1)$ dimensional then $W = W_1 \oplus \ldots \oplus W_n = (W_1^1 \oplus \ldots \oplus W_n^1 \mid \ldots \mid W_1^n \oplus \ldots \oplus W_n^n)$.

**DEFINITION 1.4.11:** *Let $V = (V_1 \mid \ldots \mid V_n)$ be a super vector space over the field F; a super projection of V is a linear operator $E_s$ on V such that $E_s^2 = E_s$ i.e. $E_s = (E_1 \mid \ldots \mid E_n)$ then $E_s^2 = (E_1^2 \mid \ldots \mid E_n^2) = (E_1 \mid \ldots \mid E_n)$ i.e. each $E_i$ is a projection on $V_i$; i = 1, 2, ..., n.. Suppose $E_s$ is a projection on V and $R = (R_1 \mid \ldots \mid R_n)$ is the super range of $E_s$ and $N = (N_1 \mid \ldots \mid N_n)$ the super null space or null super space of $E_s$. The super vector $\beta = (\beta_1 \mid \ldots \mid \beta_n)$ is in the super range $R = (R_1 \mid \ldots \mid R_n)$ if and only if $E_s\beta = \beta$ i.e. if and only if $(E_1\beta_1 \mid \ldots \mid E_n\beta_n) = (\beta_1 \mid \ldots \mid \beta_n)$ i.e. each $E_i\beta_i = \beta_i$ for i = 1, 2, ..., n.. If $\beta = E_s\alpha$ i.e. $\beta = (\beta_1 \mid \ldots \mid \beta_n) = (E_1\alpha_1 \mid \ldots \mid E_n\alpha_n)$ where the super vector $\alpha = (\alpha_1 \mid \ldots \mid \alpha_n)$ then $E_s\beta = E_s^2\alpha = E_s\alpha = \beta$. Conversely if $\beta = (\beta_1 \mid \ldots \mid \beta_n) = E_s\beta = (E_1\beta_1 \mid \ldots \mid E_n\beta_n)$ then $\beta = (\beta_1 \mid \ldots \mid \beta_n)$ is in the super range of $E_s$. Thus $V = R \oplus N$ i.e. $V = (V_1 \mid \ldots \mid V_n) = (R_1 \oplus N_1 \mid \ldots \mid R_n \oplus N_n)$.*

*Further the unique expression for = $(\alpha_1 \mid \ldots \mid \alpha_n)$ as a sum of super vectors in R and N is $\alpha = E_s\alpha + (\alpha - E_s\alpha)$ i.e. $\alpha_i = E_i\alpha_1 + (\alpha_i - E_i\alpha_i)$ for i = 1, 2, ..., n. From what we have stated it easily follows that if R and N are super subspace of V such that V = R $\oplus$N i.e. $V = (V_1 \mid \ldots \mid V_n) = (R_1 \oplus N_1 \mid \ldots \mid R_n \oplus N_n)$ then there is one and only one super projection operator $E_s$ which has super range $R = (R_1 \mid \ldots \mid R_n)$ and null super space $N = (N_1 \mid \ldots \mid N_n)$. That operator is called the super projection on R along N.*

*Any super projection $E_s$ is super diagonalizable. If $\{(\alpha_1^1 \ldots \alpha_{r_1}^1 \mid \ldots \mid \alpha_1^n \ldots \alpha_{r_n}^n)$ is a super basis for $R = (R_1 \mid \ldots \mid R_n)$ and $(\alpha_{r_1+1}^1 \ldots \alpha_{n_1}^1 \mid \ldots \mid \alpha_{r_n+1}^n \ldots \alpha_{n_n}^n)$ is a super basis for $N = (N_1 \mid \ldots \mid N_n)$ then the basis $B = (\alpha_1^1 \ldots \alpha_{n_1}^1) \mid \ldots \mid \alpha_1^n \ldots \alpha_{n_n}^n) = (B_1 \mid \ldots \mid B_n)$ super diagonalizes $E_s$.*

$$(E_s)_B = \left(\begin{array}{c|c|c|c} \begin{matrix} I_1 & 0 \\ 0 & 0 \end{matrix} & 0 & \cdots & 0 \\ \hline 0 & \begin{matrix} I_2 & 0 \\ 0 & 0 \end{matrix} & \cdots & 0 \\ \hline 0 & 0 & & 0 \\ \hline 0 & 0 & & \begin{matrix} I_n & 0 \\ 0 & 0 \end{matrix} \end{array}\right)$$

$$= ([E_1] \mid \ldots \mid [E_n]_{B_n})$$

*where $I_t$ is a $r_t \times r_t$ identity matrix; t = 1, 2, ..., n. Thus super projections can be used to describe super direct sum decompositions of the super vector space $V = (V_1 \mid \ldots \mid V_n)$.*

## Chapter Two

# SUPER INNER PRODUCT SUPER SPACES

This chapter has three sections. In section one we for the first time define the new notion of super inner product super spaces. Several properties about super inner products are derived. Further the notion of superbilinear form is introduced in section two. Section three gives brief applications of these new concepts.

### 2.1 Super Inner Product Spaces and their Properties

In this section we introduce the notion of super inner products on super vector spaces which we call as super inner product spaces.

**DEFINITION 2.1.1**: *Let $V = (V_1 \mid \ldots \mid V_n)$ be a super vector space over the field of real numbers or the field of complex numbers. A super inner product on V is a super function which assigns to each ordered pair of super vectors $\alpha = (\alpha_1 \mid \ldots \mid \alpha_n)$ and $\beta = (\beta_1 \mid \ldots \mid \beta_n)$ in V a n-tuple scalar $(\alpha \mid \beta) = ((\alpha_1 \mid \beta_1), \ldots, (\alpha_n \mid \beta_n))$ in F in such a way that for all $\alpha = (\alpha_1 \mid \ldots \mid \alpha_n)$, $\beta = (\beta_1 \mid \ldots \mid \beta_n)$ and $\gamma = (\gamma_1 \mid \ldots \mid \gamma_n)$ in V and for all n-tuple of scalars $c = (c_1, \ldots, c_n)$ in F*

*(a)* $(\alpha + \beta \mid \gamma) = (\alpha \mid \gamma) + (\beta \mid \gamma)$ *i.e.,*
$((\alpha_1 + \beta_1 \mid \gamma_1) \mid \ldots \mid (\alpha_n + \beta_n) \mid \gamma_n)$
$= ((\alpha_1 \mid \gamma_1) + ((\beta_1 \mid \gamma_1) \mid \ldots \mid (\alpha_n \mid \gamma_n) + (\beta_n \mid \gamma_n))$

*(b)* $(c\alpha|\beta) = c(\alpha|\beta)$ *i.e.,* $((c_1\alpha_1 \mid \beta_1) \mid \ldots \mid (c_n\alpha_n \mid \beta_n))$
$= (c_1 (\alpha_1 \mid \beta_1) \mid \ldots \mid c_n (\alpha_n \mid \beta_n))$

*(c)* $(\beta \mid \alpha) = \overline{(\alpha \mid \beta)}$ *i.e.,*
$$((\beta_1 \mid \alpha_1)|\ldots|(\beta_n \mid \alpha_n)) = (\overline{(\alpha_1 \mid \beta_1)}|\ldots|\overline{(\alpha_n \mid \beta_n})$$

*(d)* $(\alpha \mid \alpha) > (0 \mid \ldots \mid 0)$ *if* $\alpha \neq 0$ *i.e.,*
$((\alpha_1 \mid \alpha_1) \mid \ldots \mid (\alpha_n \mid \alpha_n)) > (0 \mid \ldots \mid 0)$.

*All the above conditions can be consolidated to imply a single equation*

$$(\alpha | c\beta + \gamma) = \overline{c}\,(\alpha \mid \beta) + (\alpha \mid \gamma) \text{ i.e.}$$
$$((\alpha_1 \mid c_1\beta_1 + \gamma_1)|\ldots|(\alpha_n \mid c_n\beta_n \mid \gamma_n)) =$$
$$(c_1(\alpha_1 \mid \beta_1) + (\alpha_1 \mid \gamma_1)|\ldots|c_n(\alpha_n \mid \beta_n) + (\alpha_n \mid \gamma_n)).$$

***Example 2.1.1****:* Suppose $V = (F^{n_1} |\ldots| F^{n_n})$ be a super inner product space over the field F. Then for $\alpha \in V$ with

$$\alpha = (\alpha^1_1 \ldots \alpha^1_{n_1} |\ldots| \alpha^n_1 \ldots \alpha^n_{n_n})$$

and $\beta \in V$ where

$$\beta = (\beta^1_1 \ldots \beta^1_{n_1} |\ldots| \beta^n_1 \ldots \beta^n_{n_n})$$
$$(\alpha \mid \beta) = \left( \sum_{j_1} \alpha_{j_1} \overline{\beta}_{j_1} |\ldots| \sum_{j_n} \alpha_{j_n} \overline{\beta}_{j_n} \right).$$

This super inner product is called as the standard super inner product on V or super dot product denoted by $\alpha \cdot \beta = (\alpha \mid \beta)$.

We define super norm of $\alpha = (\alpha_1 \mid \ldots \mid \alpha_n) \in V = (V_1 \mid \ldots \mid V_n)$. Super square root of

$$\sqrt{(\alpha \mid \alpha)} = (\sqrt{(\alpha_1 \mid \alpha_1)}|\ldots|\sqrt{(\alpha_n \mid \alpha_n)}),$$

so super square root of a n-tuple $(x_1 \mid \ldots \mid x_n)$ is $\left(\sqrt{x_1}|\ldots|\sqrt{x_n}\right)$.

We call this super square root of $(\alpha \mid \alpha)$, the super norm viz.

$$\sqrt{(\alpha \mid \alpha)} = \left(\sqrt{(\alpha_1 \mid \alpha_1)} \mid \ldots \mid \sqrt{(\alpha_n \mid \alpha_n)}\right)$$
$$= \left(\|\alpha_1\| \mid \ldots \mid \|\alpha_n\|\right)$$

$$= \quad \| \alpha \|.$$

The super quadratic form determined by the inner product is the function that assigns to each super vector α the scalar n-tuple $\|\alpha\|^2 = \left(\|\alpha_1\|^2 \mid \ldots \mid \|\alpha_n\|^2\right)$. Hence just like an inner product space the super inner product space is a real or complex super vector space together with a super inner product on that space.

We have the following interesting theorem for super inner product space.

**THEOREM 2.1.1:** *Let $V = (V_1 \mid \ldots \mid V_n)$ be a super inner product space over a field F, then for super vectors $\alpha = (\alpha_1 \mid \ldots \mid \alpha_n)$ and $\beta = (\beta_1 \mid \ldots \mid \beta_n)$ in V and any scalar c*

*(i)* $||c\alpha|| = |c|\ ||\alpha||$

*(ii)* $||\alpha|| > (0 \mid \ldots \mid 0)$ *for* $\alpha \neq (0 \mid \ldots \mid 0)$ $(||\alpha_1|| \mid \ldots \mid ||\alpha_n||) > (0 \mid \ldots \mid 0)$; $\alpha = (\alpha_1 \mid \ldots \mid \alpha_n) \neq (0 \mid \ldots \mid 0)$ *i.e.* $\alpha_i \neq 0$; $i = 1, 2, \ldots, n$.

*(iii)* $|(\alpha|\beta)| \leq ||\alpha||\ ||\beta||$ *i.e.,*
$$|((\alpha_1|\beta_1)|\ldots|(\alpha_n|\beta_n))| \underset{s}{\leq} (||\alpha_1||\ ||\beta_1||\ |\ldots|\ ||\alpha_n||\ ||\beta_n||)$$
*i.e. each* $|(\alpha_i|\beta_i)| < ||\alpha_i||\ ||\beta_i||$ *for i = 1, 2, …, n.*

*(iv)* $||\alpha + \beta|| \underset{s}{\leq} ||\alpha|| + ||\beta||$ *i.e.*
$$\left(||\alpha_1 + \beta_1||\ |\ldots|\ ||\alpha_n + \beta_n||\right)$$
$$\underset{s}{\leq} \left(||\alpha_1|| + ||\beta_1||\ |\ldots|\ ||\alpha_n|| + ||\beta_n||\right),$$ $\underset{s}{\leq}$ *denotes each* $||\alpha_i + \beta_i|| \leq ||\alpha_i|| + ||\beta_i||$ *for i = 1, 2, …, n.*

*Proof:* Statements (1) and (2) follows immediately from the various definitions involved. The inequality (iii) is true when $\alpha \neq 0$. If $\alpha \neq (0 \mid \ldots \mid 0)$ i.e. $(\alpha_1 \mid \ldots \mid \alpha_n) \neq (0 \mid \ldots \mid 0)$ i.e. $\alpha_i \neq 0$ for i = 1, 2, …, n, put
$$\gamma = \beta - \frac{(\beta|\alpha)}{\|\alpha\|^2}\alpha$$

where $\gamma = (\gamma_1 \mid \ldots \mid \gamma_n)$, $\beta = (\beta_1 \mid \ldots \mid \beta_n)$ and $\alpha = (\alpha_1 \mid \ldots \mid \alpha_n)$,

$$\begin{aligned}(\gamma_1 \mid \ldots \mid \gamma_n) &= \beta - \frac{(\beta \mid \alpha)}{\|\alpha\|^2}\alpha \\ &= \left(\beta_1 - \frac{(\beta_1 \mid \alpha_1)}{\|\alpha_1\|^2}\alpha_1 \mid \ldots \mid \beta_n - \frac{(\beta_n \mid \alpha_n)}{\|\alpha_n\|^2}\alpha_n\right).\end{aligned}$$

Then $(\gamma|\alpha) = (0 \mid \ldots \mid 0)$ and

$$\begin{aligned}(0 \mid \ldots \mid 0) &\leq \left(\|\gamma_1\|^2 \mid \ldots \mid \|\gamma_n\|^2\right) \\ &= \|\gamma\|^2 = \left(\beta - \frac{(\beta \mid \alpha)}{\|\alpha\|^2}\alpha \,\middle|\, \beta - \frac{(\beta \mid \alpha)}{\|\alpha\|^2}\alpha\right) \\ &= (\beta \mid \beta) - \frac{(\beta \mid \alpha)(\alpha \mid \beta)}{\|\alpha\|^2} \\ &= ((\beta_1 \mid \beta_1) \mid \ldots \mid (\beta_n \mid \beta_n)) \\ &\quad - \left(\frac{(\beta_1 \mid \alpha_1)(\alpha_1 \mid \beta_1)}{\|\alpha_1\|^2} \mid \ldots \mid \frac{(\beta_n \mid \alpha_n)(\alpha_n \mid \beta_n)}{\|\alpha_n\|^2}\right)\end{aligned}$$

$$= \left((\beta_1 \mid \beta_1) - \frac{(\beta_1 \mid \alpha_1)(\alpha_1 \mid \beta_1)}{\|\alpha_1\|^2} \mid \ldots \mid (\beta_n \mid \beta_n) - \frac{(\beta_n \mid \alpha_n)(\alpha_n \mid \beta_n)}{\|\alpha\|^2}\right)$$

$$= \left(\|\beta_1\|^2 - \frac{|(\alpha_1 \mid \beta_1)|^2}{\|\alpha_1\|^2} \mid \ldots \mid \|\beta_n\|^2 - \frac{|(\alpha_n \mid \beta_n)|^2}{\|\alpha_n\|^2}\right).$$

Hence $|(\alpha \mid \beta)|^2 \underset{s}{\leq} \|\alpha\|^2 \, \|\beta\|^2$ i.e. $|(\alpha_i \mid \beta_i)|^2 \leq \|\alpha_i\|^2 \, \|\beta_i\|^2$; i = 1, 2, …, n.

$$\|\alpha + \beta\|^2 =_s \|\alpha\|^2 + (\alpha \mid \beta) + (\beta \mid \alpha) + \|\beta\|^2$$

$$\text{i.e. } \left(\|\alpha_1 + \beta_1\|^2 \mid \ldots \mid \|\alpha_n + \beta_n\|^2\right) =$$

$$\left(\left(\|\alpha_1\|^2 + (\alpha_1 \mid \beta_1) + (\beta_1 \mid \alpha_1) + \|\beta_1\|^2 \mid \ldots \mid\right.\right.$$

$$\|\alpha_n\|^2 (\alpha_n \mid \beta_n) + (\beta_n \mid \alpha_n) + \|\beta_n\|^2 ) )$$

$$= \left( \|\alpha_1\|^2 + 2\mathrm{Re}(\alpha_1 \mid \beta_1) + \|\beta_1\|^2 \mid \ldots \mid \right.$$

$$\left. \|\alpha_n\|^2 + 2\mathrm{Re}(\alpha_n \mid \beta_n) + \|\beta_n\|^2 \right)$$

$$\leq \left( \|\alpha_1\|^2 + 2\|\alpha_1\| \, \|\beta_1\| + \|\beta_1\|^2 \mid \ldots \mid \right.$$

$$\left. \|\alpha_n\|^2 + 2\|\alpha_n\| \, \|\beta_n\| + \|\beta_n\|^2 \right)$$

$$= \left( (\|\alpha_1\| + \|\beta_1\|)^2 \mid \ldots \mid (\|\alpha_n\| + \|\beta_n\|)^2 \right);$$

since each $\|\alpha_i + \beta_i\| \leq \|\alpha_i\| + \|\beta_i\|$ for i = 1, 2, …, n we have $\|\alpha + \beta\| \underset{s}{\leq} \|\alpha\| + \|\beta\|$ '$\underset{s}{\leq}$' indicates that the inequality is super inequality i.e. inequality is true componentwise.

Now we proceed onto define the notion of super orthogonal set, super orthogonal supervectors and super orthonormal set.

**DEFINITION 2.1.2:** *Let α = ($\alpha_1$ | … | $\alpha_n$) and β = ($\beta_1$ | … | $\beta_n$) be super vectors in a super inner product space V = ($V_1$ | … | $V_n$). Then α is super orthogonal to β if (α | β) = (($\alpha_1$ | $\beta_1$) | … | ($\alpha_n$ | $\beta_n$)) = (0 | … | 0) since this implies β is super orthogonal to α, we often simply say α and β are super orthogonal. If S = ($S_1$ | … | $S_n$) is a supersubset of super vectors in V = ($V_1$ | … | $V_n$), S is called a super orthogonal super set provided all pairs of distinct super vectors in S are super orthogonal i.e. by the super orthogonal subset we mean every set $S_i$ in S is an orthogonal set for every i = 1, 2, …, n. i.e. ($\alpha_i$ | $\beta_i$) = 0 for all $\alpha_i$, $\beta_i$ ∈ S; i = 1, 2, …, n. A super orthonormal super set is a super orthogonal set with additional property ||α|| = (||$\alpha_1$|| |…| ||$\alpha_n$||) = (1 | … | 1), for every α in S and every $\alpha_i$ in $S_i$ is such that || $\alpha_i$ || = 1.*

The reader is expected to prove the following simple results.

**THEOREM 2.1.2**: *A super orthogonal super set of nonzero super vector is linearly super independent.*

The following corollary is direct.

**COROLLARY 2.1.1**: If a super vector β = ($\beta_1$ | … | $\beta_n$) is a linear super combination of orthogonal sequence of non-zero super vectors, $\alpha_1$, …, $\alpha_m$ then β in particular is a super linear combination,

$$\beta = \left( \sum_{K_1=1}^{m_1} \frac{(\beta_1 \mid \alpha_{K_1}^1)}{\|\alpha_{K_1}^1\|^2} \alpha_{K_1}^1 \;\middle|\; \cdots \;\middle|\; \sum_{K_n=1}^{m_n} \frac{(\beta_n \mid \alpha_{K_n}^n)}{\|\alpha_{K_n}^n\|^2} \alpha_{K_n}^n \right).$$

We can on similar lines as in case of usual vector spaces derive Gram Schmidt super orthogonalization process for super inner product space V.

**THEOREM 2.1.3:** *Let V = ($V_1$ | ... | $V_n$) be a super inner product space and let* $(\beta_1^1 \ldots \beta_{n_1}^1),\ldots,(\beta_1^n \ldots \beta_{n_n}^n)$ *be any independent super vector in V. Then one may construct orthogonal super vector* $(\alpha_1^1 \ldots \alpha_{n_1}^1),\ldots,(\alpha_1^n \ldots \alpha_{n_n}^n)$ *in V such that for each K = ($K_1$ | ... | $K_n$), the set* $\{(\alpha_1^1,\ldots \alpha_{K_1}^1),\ldots,(\alpha_1^n,\ldots,\alpha_{K_n}^n)\}$ *is a super basis for the super subspace spanned by* $(\beta_1^1 \ldots \beta_{K_1}^1),\ldots,(\beta_1^n \ldots \beta_{K_n}^n)$.

Just we indicate how we can prove, for the proof is similar to usual vector spaces with the only change in case of super vector spaces they occur in n-tuples.

It is left for the reader to prove "Every finite ($n_1$, …, $n_n$) dimensional super inner product superspace has an orthonormal super basis". We can as in case of vector space define the notion of best super approximation for super vector spaces. Let V = ($V_1$ | … | $V_n$) be a super vector space over a field F. W = ($W_1$ | … | $W_n$) be a super subspace of V = ($V_1$ | … | $V_n$). A best super approximation to β = ($\beta_1$ | … | $\beta_n$) by super vectors in W = ($W_1$ | … | $W_n$) is a super vector α = ($\alpha_1$ | … | $\alpha_n$) in W such that

$$\begin{aligned} ||\beta - \alpha|| &= (||\beta_1 - \alpha_1|| \,|\, \ldots \,|\, ||\beta_n - \alpha_n||) \\ &\leq (||\beta_1 - \gamma_1|| \,|\, \ldots \,|\, ||\beta_n - \gamma_n|) = ||\beta - \gamma|| \end{aligned}$$

for every super vector $\gamma = (\gamma_1 \mid \ldots \mid \gamma_n)$ in W.

The reader is expected to prove the following theorem.

**THEOREM 2.1.4:** *Let $W = (W_1 \mid \ldots \mid W_n)$ be a super subspace of a super inner product space $V = (V_1 \mid \ldots \mid V_n)$ and let $\beta = (\beta_1 \mid \ldots \mid \beta_n)$ be a super vector in V.*

*(i)* *The super vector $\alpha = (\alpha_1 \mid \ldots \mid \alpha_n)$ in W is a best super approximation to $\beta = (\beta_1 \mid \ldots \mid \beta_n)$ by super vectors in $W = (W_1 \mid \ldots \mid W_n)$ if and only if $\beta - \alpha = (\beta_1 - \alpha_1 \mid \ldots \mid \beta_n - \alpha_n)$ is super orthogonal to every super vector in W.*

*(ii)* *If a best super approximation to $\beta = (\beta_1 \mid \ldots \mid \beta_n)$ by super vector in W exists, it is unique.*

*(iii)* *If $W = (W_1 \mid \ldots \mid W_n)$ is finite dimension super subspace of V and $\{(\alpha_1^1 \ldots \alpha_{K_1}^1), \ldots (\alpha_1^n \ldots \alpha_{K_n}^n)\}$ is any orthonormal super basis for W then the super vector $\alpha = (\alpha_1 \mid \ldots \mid \alpha_n) =$ $\left( \sum_{K_1} \frac{(\beta_1 \mid \alpha_{K_1}^1)\alpha_{K_1}^1}{||\alpha_{K_1}^1||^2} \middle| \ldots \middle| \sum_{K_n} \frac{(\beta_n \mid \alpha_{K_n}^n)\alpha_{K_n}^n}{||\alpha_{K_n}^n||^2} \right)$ is the best super approximation to $\beta$ by super vectors in W.*

Now we proceed onto define the notion of orthogonal complement of a super subset S of a super vector space V.

Let $V = (V_1 \mid \ldots \mid V_n)$ be an inner product super space and S any set of super vectors in V. The super orthogonal complement of S is the superset $S^{\perp}$ of all super vectors in V which are super orthogonal to every super vector in S.

Let $V = (V_1 \mid \ldots \mid V_n)$ be a super vector space over the field F. Let $W = (W_1 \mid \ldots \mid W_n)$ be a super subspace of a super inner product super space V and let $\beta = (\beta_1 \mid \ldots \mid \beta_n)$ be a super vector in V. $\alpha = (\alpha_1 \mid \ldots \mid \alpha_n)$ in W is called the orthogonal super projection to $\beta = (\beta_1 \mid \ldots \mid \beta_n)$ on $W = (W_1 \mid \ldots \mid W_n)$. If every super vector in V has an orthogonal super projection of $\beta = (\beta_1 \mid \ldots \mid \beta_n)$ on W, the mapping that assigns to each super vector in V its orthogonal super projection on $W = (W_1 \mid \ldots \mid W_n)$ is called the orthogonal super projection of V on W. Suppose $E_s = (E_1 \mid \ldots \mid E_n)$ is the orthogonal super projection of V on W. Then the super mapping $\beta \to \beta - E_s\beta$ i.e., $\beta = (\beta_1 \mid \ldots \mid \beta_n) \to (\beta_1 - E_1\beta_1 \mid \ldots \mid \beta_n - E_n\beta_n)$ is the orthogonal super projection of V on $W^{\perp} = (W_1^{\perp} | \ldots | W_n^{\perp})$.

The following theorem can be easily proved and hence left for the reader.

**THEOREM 2.1.5:** *Let $W = (W_1 \mid \ldots \mid W_n)$ be a finite dimensional super subspace of a super inner product space $V = (V_1 \mid \ldots \mid V_n)$ and let $E_s = (E_1 \mid \ldots \mid E_n)$ be the orthogonal super projection of V on W. Then $E_s$ is an idempotent linear transformation of V onto W; $W^{\perp}$ is the null super subspace of $E_s$ and $V = W \oplus W^{\perp}$ i.e. $V = (V_1 \mid \ldots \mid V_n) = (W_1 \oplus W_1^{\perp} | \ldots | W_n \oplus W_n^{\perp})$.*

Consequent of this one can prove $I - E_s = ((I_1 - E_1 \mid \ldots \mid I_n - E_n)$ is the orthogonal super projection of V on $W^{\perp}$. It is a super idempotent linear transformation of V onto $W^{\perp}$ with null super space W.

We can also prove the following theorem.

**THEOREM 2.1.6:** *Let $\{\alpha_1^1 \ldots \alpha_{n_1}^1 | \ldots | \alpha_1^n \ldots \alpha_{n_n}^n\}$ be a super orthogonal superset of nonzero super vectors in an inner product super space $V = (V_1 \mid \ldots \mid V_n)$. If $\beta = (\beta_1 \mid \ldots \mid \beta_n)$ be a super vector in V, then*

$$\left( \sum_{K_1} \frac{|(\beta_1 | \alpha_{K_1}^1)|^2}{\|\alpha_{K_1}^1\|^2} \middle| \ldots \middle| \sum_{K_n} \frac{|(\beta_n | \alpha_{K_n}^n)|^2}{\|\alpha_{K_n}^n\|^2} \right) \leq (\|\beta_1\|^2 | \ldots | \|\beta_n\|^2)$$

*and equality holds if and only if*

$$\beta = (\beta_1 \mid \ldots \mid \beta_n) = \left( \sum_{K_1} \frac{(\beta_1 \mid \alpha^1_{K_1})}{\|\alpha^1_{K_1}\|^2} \alpha^1_{K_1} \middle| \ldots \middle| \sum_{K_n} \frac{(\beta_n \mid \alpha^n_{K_n})}{\|\alpha^n_{K_n}\|^2} \alpha^n_{K_n} \right).$$

Next we proceed onto define the notion of super linear functional on a super vector space V.

**DEFINITION 2.1.3:** *Let $V = (V_1 \mid \ldots \mid V_n)$ be a super vector space over the field F, a super linear functional $f = (f_1 \mid \ldots \mid f_n)$ from V into the scalar field F is also called a super linear functional on V. i.e. $f: V \to (F \mid \ldots \mid F)$ where V is a super vector space defined over the field F, i.e. $f = (f_1 \mid \ldots \mid f_n): V = (V_1 \mid \ldots \mid V_n) \to (F \mid \ldots \mid F)$ by*

$$f(c\alpha + \beta) = cf(\alpha) + f(\beta)$$

$$\text{i.e. } (f_1(c_1\alpha_1 + \beta_1) \mid \ldots \mid f_n(c_n\alpha_n + \beta_n))$$

$$= (c_1f_1(\alpha_1) + f_1(\beta_1) \mid \ldots \mid c_nf_n(\alpha_n) + f_n(\beta_n))$$

*for all super vector α and β in V where $\alpha = (\alpha_1 \mid \ldots \mid \alpha_n)$ and $\beta = (\beta_1 \mid \ldots \mid \beta_n)$.*

$$SL(V, F) = (L(V_1, F) \mid \ldots \mid L(V_n, F)).$$

*$(V_1^* \mid \ldots \mid V_n^*) = V^*$ denotes the collection of all super linear functionals from the super vector space $V = (V_1 \mid \ldots \mid V_n)$ into $(F \mid \ldots \mid F)$ where $V_i$'s are vector spaces defined over the same field F. $V^* = (V_1^* \mid \ldots \mid V_n^*)$ is called the dual super space or super dual space of V.*

*If $V = (V_1 \mid \ldots \mid V_n)$ is finite $(n_1, \ldots, n_n)$ dimensional over F then $V^* = (V_1^* \mid \ldots \mid V_n^*)$ is also finite $(n_1, \ldots, n_n)$ dimensional over F.*

$$\text{i.e. super dim } V^* = \text{super dim } V.$$

*Suppose $B = (B_1 \mid \ldots \mid B_n) = \{\alpha^1_1 \ldots \alpha^1_{n_1} \mid \ldots \mid \alpha^n_1 \ldots \alpha^n_{n_n}\}$ be a super basis for $V = (V_1 \mid \ldots \mid V_n)$ of dimension $(n_1, \ldots, n_n)$ then for each $i_t$ a unique linear function $f^t_{i_t}(\alpha^t_{j_t}) = \delta_{i_t j_t}$; t = 1, 2, …, n.*

*In this way we obtain from B a set of distinct linear functional* $B^* = (f_1^1 \ldots f_{n_1}^1 | \cdots | f_1^n \ldots f_{n_n}^n)$ *on V. These linear functionals are also linearly super independent as each of the set* $\{f_1^t, \ldots, f_{n_t}^t\}$ *is a linearly independent set for t = 1, 2, ..., n.* $B^*$ *forms a super basis for* $V^* = (V_1^* | \ldots | V_n^*)$ .

The following theorem is left as an exercise for the reader to prove.

**THEOREM 2.1.7**: *Let V = ($V_1$ | ... | $V_n$) be a super vector space over the field F and let* $B = \{\alpha_1^1 \ldots \alpha_{n_1}^1 | \ldots | \alpha_1^n \ldots \alpha_{n_n}^n\}$ *be a super basis for V = ($V_1$ | ... | $V_n$). Then there exists a unique dual super basis* $B^* = \{f_1^1 \ldots f_{n_1}^1 | \cdots | f_1^n \ldots f_{n_n}^n\}$ *for* $V^* = (V_1^* | \ldots | V_n^*)$ *such that* $f_{i_t}^t(\alpha_{j_t}^t) = \delta_{i_t j_t}; 1 \le i_t, j_t \le n_t$ *and for every t = 1, 2,..., n.*

For each linear super functional f = ($f_1$ | … | $f_n$) on V we have

$$f = \left( \sum_{i_1=1}^{n_1} f_1(\alpha_{i_1}^1) f_{i_1}^1 \;| \ldots |\; \sum_{i_n=1}^{n_n} f_n(\alpha_{i_n}^n) f_{i_n}^n \right) = f = (f_1 | \ldots | f_n)$$

and for each super vector α = ($\alpha_1$ | … | $\alpha_n$) in V we have

$$\alpha = \left( \sum_{i_1=1}^{n_1} f_{i_1}^1(\alpha_1)\, \alpha_{i_1}^1 \;| \ldots |\; \sum_{i_n=1}^{n_n} f_{i_n}^n(\alpha_n)\, \alpha_{i_n}^n \right).$$

*Note*: We call f = ($f_1$ | … | $f_n$) to be a super linear functional if
f:V = ($V_1$ | … | $V_n$)→ (F| … |F)
i.e., $f_1 : V_1 \to F$ , …, $f_n : V_n \to F$.
This concept of super linear functional leads us to define the notion of hyper super spaces.

*Let V = ($V_1$ | ... | $V_n$) be super vector space over the field F. f = ($f_1$ | ... | $f_n$) be a super linear functional from V = ($V_1$ | ... | $V_n$) into (F| ... |F). Suppose V = ($V_1$ | ... | $V_n$) is finite ($n_1$, ..., $n_n$) dimensional over F.*

*Let* $N = (N^1_{f_1} | \ldots | N^n_{f_n})$ *be the super null space of f = (f*$_1$ *| ... | f*$_n$*). Then super dimension of*

$$N_f = (dim\, N^1_{f_1}, ..., dim\, N^n_{f_n}) = (dim\, V_1 - 1, ..., dim\, V_n - 1)$$
$$= (n_1 - 1, \ldots, n_n - 1)\,.$$

*In a super vector space* $(n_1, \ldots, n_n)$ *a super subspace of super dimension* $(n_1 - 1, \ldots, n_n - 1)$ *is called a super hyper space or hyper super space.*

Is every hyper super space the null super subspace of a super linear functional. The answer is yes.

**DEFINITION 2.1.4**: *If* $V = (V_1 | \ldots | V_n)$ *be a super vector space over the field F and* $S = (S_1 | \ldots | S_n)$ *be a super subset of V, the super annihilator of* $S = (S_1 | \ldots | S_n)$ *is the super set* $S^0 = (S^0_1 | \ldots | S^0_n)$ *of super linear functionals* $f = (f_1 | \ldots | f_n)$ *on V such that* $f(\alpha) = (f_1(\alpha_1) | \ldots | f_n(\alpha_n)) = (0 | \ldots | 0)$ *for every* $\alpha = (\alpha_1 | \ldots | \alpha_n)$ *in* $S = (S_1 | \ldots | S_n)$*. It is easily verified that* $S^0 = (S^0_1 | \ldots | S^0_n)$ *is a subspace of* $V^* = (V^*_1 | \ldots | V^*_n)$*, whether S =* $(S_1 | \ldots | S_n)$ *is super subspace of* $V = (V_1 | \ldots | V_n)$ *or not. If S =* $(S_1 | \ldots | S_n)$ *is the super set consisting of the zero super vector alone then* $S^0 = V^*$ *i.e.* $(S^0_1 | \ldots | S^0_n) = (V^*_1 | \ldots | V^*_n)$*. If S =* $(S_1 | \ldots | S_n) = V = (V_1 | \ldots | V_n)$ *then* $S^0 = (0 | \ldots | 0)$ *of* $V^* = (V^*_1 | \ldots | V^*_n)$ *.*

The following theorem and the two corollaries is an easy consequence of the definition.

**THEOREM 2.1.8:** *Let* $V = (V_1 | \ldots | V_n)$ *be a finite* $(n_1, \ldots, n_n)$ *dimensional super vector space over the field F, and let* $W = (W_1 | \ldots | W_n)$ *be a super subspace of V.*

*Then*

$$super\ dim\ W + super\ dim\ W^0 = super\ dim\ V.$$
$$i.e.\ (dim\ W_1 + dim\ W^0_1, \ldots, dim\, W_n + dim\, W^0_n)$$
$$= (n_1, \ldots, n_n).$$

**COROLLARY 2.1.2**: *$W = (W_1 \mid \ldots \mid W_n)$ is a $(k_1, \ldots, k_n)$ dimensional super subspace of a $(n_1, \ldots, n_n)$ dimensional super vector space $V = (V_1 \mid \ldots \mid V_n)$ then $W = (W_1 \mid \ldots \mid W_n)$ is the super intersection of $(n_1 - k_1, \ldots, n_n - k_n)$ hyper super spaces in V.*

**COROLLARY 2.1.3:** *If $W_1 = (W_1^1 \mid \ldots \mid W_1^n)$ and $W_2 = (W_2^1 \mid \ldots \mid W_2^n)$ are super subspaces of a finite $(n_1, \ldots, n_n)$ dimensional super vector space then $W_1 = W_2$ if and only if $W_1^0 = W_2^0$ i.e. $W_1^t = W_2^t$ if and only if $(W_1^t)^0 = (W_2^t)^0$ for every t = 1, 2, …, n.*

Now we proceed onto prove the notion of super double dual or double super dual (both mean the same). To consider $V^{**}$ the super dual of $V^*$. If $\alpha = (\alpha_1 \mid \ldots \mid \alpha_n)$ is a super vector in $V = (V_1 \mid \ldots \mid V_n)$ then α induces a super linear functional $L_\alpha = (L^1_{\alpha_1} \mid \ldots \mid L^n_{\alpha_n})$ on $V^* = (V_1^* \mid \ldots \mid V_n^*)$ defined by

$$L_\alpha(f) = (L^1_{\alpha_1}(f_1) \mid \ldots \mid L^n_{\alpha_n}(f_n)) = (f_1(\alpha_1) \mid \ldots \mid f_n(\alpha_n))$$

where $f = (f_1 \mid \ldots \mid f_n)$ in $V^*$. The fact that $L_\alpha$ is linear is just a reformulation of the definition of linear operators on $V^*$.

$$\begin{aligned} L_\alpha(cf + g) &= (L^1_{\alpha_1}(c_1 f_1 + g_1) \mid \ldots \mid L^n_{\alpha_n}(c_n f_n + g_n)) \\ &= (c_1 f_1 + g_1)(\alpha_1) \mid \ldots \mid (c_n f_n + g_n)(\alpha_n)) \\ &= ((c_1 f_1)\,\alpha_1 + g_1(\alpha_1) \mid \ldots \mid (c_n f_n)\alpha_n + g_n(\alpha_n)) \\ &= (c_1 L^1_{\alpha_1}(f_1) + L^1_{\alpha_1}(g_1) \mid \ldots \mid c_n L^n_{\alpha_n}(f_n) + L^n_{\alpha_n}(g_n)). \end{aligned}$$

If $V = (V_1 \mid \ldots \mid V_n)$ is finite $(n_1, \ldots, n_n)$ dimensional and $\alpha = (\alpha_1 \mid \ldots \mid \alpha_n) \neq (0 \mid \ldots \mid 0)$ then $L_\alpha = (L^1_{\alpha_1} \mid \ldots \mid L^n_{\alpha_n}) \neq (0 \mid \ldots \mid 0)$, in other words there exists a linear super functional $f = (f_1 \mid \ldots \mid f_n)$ such that $(f_1(\alpha_1) \mid \ldots \mid f_n(\alpha_n)) \neq (0 \mid \ldots \mid 0)$.

The following theorem is direct and hence left for the reader to prove.

**THEOREM 2.1.9**: *Let V = ($V_1$ | ... | $V_n$) be a finite ($n_1$, ..., $n_n$) dimensional super vector space over the field F. For each super vector α = ($\alpha_1$ | ... | $\alpha_n$) in V define*

$$L_\alpha(f) = (L^1_{\alpha_1}(f_1) | \ldots | L^n_{\alpha_n}(f_n))$$
$$= (f_1(\alpha_1) \mid \ldots \mid f_n(\alpha_n)) = f(\alpha)$$
$$f = (f_1 \mid \ldots \mid f_n) \in V^* = (V_1^* | \ldots | V_n^*).$$

*The super mapping α → $L_\alpha$ i.e.*

$$\alpha = (\alpha_1 \mid \ldots \mid \alpha_n) \rightarrow (L^1_{\alpha_1} | \ldots | L^n_{\alpha_n})$$

*is then a super isomorphism of V onto $V^{**}$.*

In view of the above theorem the following two corollaries are direct.

**COROLLARY 2.1.4**: *Let V = ($V_1$ | ... | $V_n$) be a finite ($n_1$, ..., $n_n$) dimensional super vector space over the field F. If L = ($L^1$ | ... |$L^n$) is a super linear functional on the dual super space $V^* = (V_1^* | \ldots | V_n^*)$ of V then there is a unique super vector α = ($\alpha_1$ | ... | $\alpha_n$) in V such that L(f) = f(α); ($L^1(f_1)$| ...| $L^n(f_n)$)= ($f_1(\alpha_1)$| ... | $f_n(\alpha_n)$) for every f = ($f_1$ | ... | $f_n$) in $V^* = (V_1^* | \ldots | V_n^*)$.*

**COROLLARY 2.1.5**: *Let V = ($V_1$ | ... | $V_n$) be a finite ($n_1$, ..., $n_n$) dimensional super vector space over the field F. Each super basis for $V^* = (V_1^* | \ldots | V_n^*)$ is the super dual for some super basis for V.*

**THEOREM 2.1.10:** *If S = ($S_1$ | ... | $S_n$) is a super subset of a finite ($n_1$, ..., $n_n$) dimensional super vector space V = ($V_1$ | ... | $V_n$) then $(S^0)^0 = ((S_1^0)^0 | \ldots | (S_n^0)^0)$ is the super subspace spanned by S = ($S_1$ | ... | $S_n$).*

*Proof:* Let W = ($W_1$| … | $W_n$) be the super subspace spanned by S = ($S_1$ |…| $S_n$). Clearly $W^0 = S^0$ i.e., $(W_1^0|\ldots|W_n^0) = (S_1^0|\ldots|S_n^0)$. Therefore what we have to prove is that

$$W^{00} = W \text{ i.e. } (W_1|\ldots|W_n) = (W_1^{00}|\ldots|W_n^{00}).$$

We have

superdim W + super dim $W^0$ = super dim V.

i.e. (dim $W_1$ + dim $W_1^0, \ldots, \dim W_n + \dim W_n^0$)
= (dim $V_1$, …, dim $V_n$)
= ($n_1$, …, $n_n$).

i.e. ( dim $W_1^0 + \dim W_1^{00}, \ldots, \dim W_n^0 + \dim W_n^{00}$)
= super dim $W^o$ + super dim $W^{oo}$
= super dim $V^*$ = $(\dim V_1^*, \ldots, \dim V_n^*)$
= ($n_1$, …, $n_n$).

Since super dim V = super dim $V^*$, we have super dim W = super dim $W^0$. Since W is a super subspace of $W^{00}$ we see W = $W^{00}$.

Let V = ($V_1$ | … | $V_n$) be a super vector space, a super hyper space in V = ($V_1$ | … | $V_n$) is a maximal proper super subspace of V = ($V_1$ | … | $V_n$).

In view of this we have the following theorem which is left as an exercise for the reader to prove.

**THEOREM 2.1.11**: *If f = ($f_1$ | ... | $f_n$) is a nonzero super linear functional on the super vector space V = ($V_1$ | ... | $V_n$), then the super null space of f = ($f_1$ | ... | $f_n$) is a super hyper space in V. Conversely every super hyper subspace in V is the super null subspace of a non zero super linear functional on V = ($V_1$ | ... | $V_n$).*

The following lemma can easily be proved.

**LEMMA 2.1.1**: *If f = ($f_1$ | ... | $f_n$) and g = ($g_1$ | ... | $g_n$) be linear super functionals on the super vector space then g is a scalar multiple of f if and only if the super null space of g contains the super null space of f that is if and only if f(α) = ($f_1(\alpha_1)$ | ... |*

*$f_n(\alpha_n)) = (0 \mid \ldots \mid 0)$ implies $g(\alpha) = (g_1(\alpha_1) \mid \ldots \mid g_n(\alpha_n)) = (0 \mid \ldots \mid 0)$.*

We have the following interesting theorem on the super null subspaces of super linear functional on V.

**THEOREM 2.1.12**: *Let $g = (g_1 \mid \ldots \mid g_n)$;*

$$f_1 = (f_1^1 \mid \ldots \mid f_n^1), \ldots, f_r = (f_1^r \mid \ldots \mid f_n^r)$$

*be linear super functionals on a super vector space $V = (V_1 \mid \ldots \mid V_n)$ with respective null super spaces $N_1, \ldots, N_r$ respectively. Then $g = (g_1 \mid \ldots \mid g_n)$ is a super linear combination of $f_1, \ldots, f_r$ if and only if N contains the intersection $N_1 \cap \ldots \cap N_r$ i.e., $N_1 = (N^1 \mid \ldots \mid N^n)$ contains $(N_1^1 \cap \ldots \cap N_{r_1}^1 \mid \ldots \mid N_1^n \cap \ldots \cap N_{r_n}^n)$.*

As in case of usual vector spaces in the case of super vector spaces also we have the following :

*Let $V = (V_1 \mid \ldots \mid V_n)$ and $W = (W_1 \mid \ldots \mid W_n)$ be two super vector spaces defined over the field F. Suppose we have a linear transformation $T_s = (T_1 \mid \ldots \mid T_n)$ from V into F. Then $T_s$ induces a linear transformation from $W^*$ into $V^*$ as follows :*

*Suppose $g = (g_1 \mid \ldots \mid g_n)$ is a linear functional on $W = (W_1 \mid \ldots \mid W_n)$ and let $f(\alpha) = g(T(\alpha))$ for each $\alpha = (\alpha_1 \mid \ldots \mid \alpha_n)$ i.e.*

*$(f_1(\alpha_1) \mid \ldots \mid f_n(\alpha_n)) = (g_1(T_1(\alpha_1) \mid \ldots \mid g_n(T_n(\alpha_n))$ ... I*

*for each $\alpha = (\alpha_1 \mid \ldots \mid \alpha_n)$ in $V = (V_1 \mid \ldots \mid V_n)$. Then I defines a function $f = (f_1 \mid \ldots \mid f_n)$ from $V = (V_1 \mid \ldots \mid V_n)$ into $(F \mid \ldots \mid F)$ namely the composition of $T_s$, a super function from V into W with $g = (g_1 \mid \ldots \mid g_n)$ a super function from $W = (W_1 \mid \ldots \mid W_n)$ into $(F \mid \ldots \mid F)$. Since both $T_s$ and g are linear f is also linear i.e. f is a super linear functional on V. Thus $T_s = (T_1 \mid \ldots \mid T_n)$ provides us a rule $T_s^t = (T_1^t \mid \ldots \mid T_n^t)$ which associates with each linear functional $g = (g_1 \mid \ldots \mid g_n)$ on $W = (W_1 \mid \ldots \mid W_n)$ a linear functional $f = T_s^t g$ i.e. $f = (f_1 \mid \ldots \mid f_n) = (T_1^t g_1 \mid \ldots \mid T_n^t g_n)$ on V defined by $f(\alpha) = g(T\alpha)$ i.e. by I; $T_s^t = (T_1^t \mid \ldots \mid T_n^t)$ is actually a linear transformation from $W^* = (W_1^* \mid \ldots \mid W_n^*)$ into*

$V^* = (V_1^* | \ldots | V_n^*)$. *For if $g_1$ and $g_2$ are in* $W^* = (W_1^* | \ldots | W_n^*)$ *and c is a scalar.*

$$
\begin{aligned}
&[T_s^t\ (cg_1 + g_2)]\ (\alpha) \\
&\quad = \left[[T_1^t\ (c_1 g_1^1 + g_2^1)]\ \alpha_1 | \ldots | [T_n^t (c_n g_1^n + g_2^n)]\ \alpha_n\right] \\
&\quad = (cg_1 + g_2)\ (T_s \alpha) \\
&\quad = [(c_1 g_1^1 + g_2^1)\ (T_1 \alpha_1) | \ldots | (c_n g_1^n + g_2^n)\ (T_n \alpha_n)] \\
&\quad = cg_1 (T_s \alpha) + g_2 (T_s \alpha) \\
&\quad = [(c_1 g_1^1 T_1\ \alpha_1 + g_2^1\, T_1 \alpha_1) | \ldots | (c_n g_1^n T_n\, \alpha_n + g_2^n T_n \alpha_n)] \\
&\quad = c(T_s^t g_1)\alpha\ +\ (T_s^t g_2)\alpha \\
&\quad = [(c_1 T_1^t g_1^1)\alpha_1 + (T_1^t g_2^1)\alpha_1 | \ldots | (c_n T_n^t g_1^n)\alpha_n + (T_n^t g_2^n)\alpha_n],
\end{aligned}
$$

*so that* $T_s^t\ (cg_1 + g_2) = cT_s^t\ g_1 + T_s^t g_2$.

This can be summarized by the following theorem.

**THEOREM 2.1.13:** *Let V = ($V_1$ | ... | $V_n$) and W = ($W_1$ | ... | $W_n$) be super vector spaces over the field F. For each linear transformation $T_s$ from V into W there is a unique transformation $T_s^t$ from $W^*$ into $V^*$ such that* $(T_s^t g)\alpha = g(T_s\ \alpha)$ *for every g = ($g_1$ | ... | $g_n$) in* $W^* = (W_1^* | \ldots | W_n^*)$ *and α = ($\alpha_1$ | ... | $\alpha_n$) in V = ($V_1$ | ... | $V_n$).*

*We shall call $T_s^t$ the super transpose of $T_s$. This $T_s^t$ is often called super adjoint of $T_s$ = ($T_1$ | ... | $T_n$).*

The following theorem is more interesting as it is about the associated super matrix.

**THEOREM 2.1.14:** *Let V = ($V_1$ | ... | $V_n$) and W = ($W_1$ | ... | $W_n$) be two finite ($n_1$, ..., $n_n$) dimensional super vector spaces over the field F. Let*

$$B = \{B_1 | \ldots | B_n\} = \{\alpha_1^1 \cdots \alpha_{n_1}^1;\ \alpha_1^2\ \ldots\ \alpha_{n_2}^2; \ldots ; \alpha_1^n \ldots \alpha_{n_n}^n\}$$

*be a super basis for V and* $B^* = \{B_1^* | \ldots | B_n^*\}$ *the dual super basis and let* $B' = \{B_1' | \ldots | B_n'\}$ *be an ordered super basis for W i.e.* $B' = \{\beta_1^1 \ldots \beta_{n_1}^1; \ldots ; \beta_1^n \ldots \beta_{n_n}^n\}$ *be a super ordered super basis for*

*W and let* $B'^{*} = (B_1'^{*} \mid \ldots \mid B_n'^{*})$ *be a dual super basis for B′. Let* $T_s = (T_1 \mid \ldots \mid T_n)$ *be a linear transformation from V into W. Let A be the super matrix of* $T_s$ *which is a super diagonal matrix relative to B and B′; let B be the super diagonal matrix of* $T_s^t$ *relative to* $B'^{*}$ *and* $B^{*}$. *Then* $B_{ij} = A_{ji}$.

*Proof:* Let

$$B = \left\{\alpha_1^1 \ldots \alpha_{n_1}^1; \ldots; \alpha_1^n \ldots \alpha_{n_n}^n\right\}$$

$$B' = \left\{\beta_1^1 \ldots \beta_{m_1}^1; \ldots; \beta_1^n \ldots \beta_{m_n}^n\right\}$$

and

$$B^{*} = \left\{f_1^1 \ldots f_{n_1}^1; \ldots; f_1^n \ldots f_{n_n}^n\right\}.$$

By definition

$$T_s\alpha_j = \sum_{i=1}^{m} A_{ij}\beta_i ;$$

j = 1, 2, …, n.

$$T_s\alpha_j = \left\{\left(\sum_{i_1=1}^{m_1} A_{i_1 j_1}^1 \; \beta_{i_1}^1 \mid \ldots \mid \sum_{i_n=1}^{m_n} A_{i_n j_n}^n \; \beta_{i_n}^n\right)\right\}$$

with $j = (j_1, \ldots, j_n)$; $1 \le j_t \le n_t$ and t = 1, 2, …, n.

$$T_s^t g_j = \left(\sum_{i_1=1}^{n_1} B_{i_1 j_1}^1 f_{i_1}^1 \mid \ldots \mid \sum_{i_n=1}^{n_n} B_{i_n j_n}^n f_{i_n}^n\right)$$

with $j = (j_1, \ldots, j_m)$; $1 \le j_p \le n_p$ and p = 1, 2, …, n.

On the other hand

$$(T_s^t g_j)(\alpha_i) = g_j(T_s\alpha_i) = g_j\left(\sum_{k=1}^{m} A_{ki}\beta_k\right)$$

$$= \left(g_{j_1}^1\left(\sum_{k_1=1}^{m_1} A_{k_1 i_1}^1\right)\beta_{k_1}^1 \mid \ldots \mid g_{j_n}^n\left(\sum_{k_n=1}^{m_n} A_{k_n i_n}^n\right)\beta_{k_n}^n\right)$$

$$= \left(\sum_{k_1=1}^{m_1} A_{k_1 i_1}^1 g_{j_1}^1(\beta_{k_1}^1) \mid \ldots \mid \sum_{k_n=1}^{m_n} A_{k_n i_n}^n g_{j_n}^n(\beta_{k_n}^n)\right)$$

$$= (A_{j_1 i_1}^1 \mid \ldots \mid A_{j_n i_n}^n).$$

For any linear functional $f = (f_1 \mid \ldots \mid f_n)$ on $V = (V_1 \mid \ldots \mid V_n)$.

$$f = (f_1 \mid \ldots \mid f_n) = \left( \sum_{i_1=1}^{m_1} f_1(\alpha^1_{i_1}) f^1_{i_1} \;\Big|\; \ldots \;\Big|\; \sum_{i_n=1}^{m_n} f_n(\alpha^n_{i_n})\, f^n_{i_n} \right).$$

If we apply this formula to the functional $f = T^t_s g_j$ i.e.,

$$(f_1 \mid \ldots \mid f_n) = (T^t_1\, g^1_{j_1} \mid \ldots \mid T^t_n\, g^n_{j_n})$$

and use the fact that

$$\begin{aligned}(T^t_s\; g_j)(\alpha_i) &= (T^t_1\; g^1_{j_1}(\alpha^1_{i_1}) \mid \ldots \mid T^t_n\; g^n_{j_n}\; (\alpha^n_{i_n})) \\ &= (A^1_{j_1 i_1} \mid \ldots \mid A^n_{j_n i_n}),\end{aligned}$$

we have

$$T^t_s g_j = \left( \sum_{i_1=1}^{n_1} A^1_{j_1 i_1}\; f^1_{i_1} \;\Big|\; \ldots \;\Big|\; \sum_{i_n=1}^{n_n} A^n_{j_n i_n}\; f^n_{i_n} \right);$$

from which it immediately follows that $B_{ij} = A_{ji}$ ; by default of notation we have

$$B_{ij} = (B^1_{i_1 j_1} \mid \ldots \mid B^n_{i_n j_n}) = \left(\begin{array}{c|c|c|c} B^1_{i_1 j_1} & 0 & & 0 \\ \hline 0 & B^2_{i_2 j_2} & & 0 \\ \hline & & & \\ \hline 0 & 0 & & B^n_{i_n j_n} \end{array}\right)$$

$$= \left(\begin{array}{c|c|c|c} A^1_{j_1 i_1} & 0 & & 0 \\ \hline 0 & A^2_{j_2 i_2} & & 0 \\ \hline & & & \\ \hline 0 & 0 & & A^n_{j_n i_n} \end{array}\right).$$

We just denote how the transpose of a super diagonal matrix looks like

$$A = \left(\begin{array}{c|c|c|c} A^1_{m_1 \times n_1} & 0 & & 0 \\ \hline 0 & A^2_{m_2 \times n_2} & & 0 \\ \hline & & & \\ \hline 0 & 0 & & A^n_{m_n \times n_n} \end{array}\right).$$

Now $A^t$ the transpose of the super diagonal matrix A is

$$A^t = \left(\begin{array}{c|c|c|c} A^1_{n_1\times m_1} & 0 & & 0 \\ \hline 0 & A^2_{n_2\times m_2} & & 0 \\ \hline & & & \\ \hline 0 & 0 & & A^n_{n_n\times m_n} \end{array}\right)$$

where A is a $(m_1 + m_2 + \ldots + m_n) \times (n_1 + \ldots + n_n)$ super diagonal matrix whereas $A^t$ is a $(n_1 + \ldots + n_n) \times (m_1 + \ldots + m_n)$; super diagonal matrix.

We illustrate by an example.

***Example 2.1.2:*** Let A be a super diagonal matrix, i.e.

$$A = \left(\begin{array}{cccc|ccc|cc|ccccc} 3 & 1 & 0 & 2 & & & & & & & & & & \\ 1 & 0 & 5 & 0 & & 0 & & & 0 & & & 0 & & \\ 0 & 1 & 0 & 1 & & & & & & & & & & \\ \hline & & & & 3 & 4 & 5 & & & & & & & \\ & 0 & & & 1 & 3 & 1 & & 0 & & & 0 & & \\ \hline & & & & & & & 8 & 1 & & & & & \\ & 0 & & & & 0 & & 6 & -1 & & & 0 & & \\ & & & & & & & 2 & 5 & & & & & \\ \hline & & & & & & & & & 1 & 0 & 1 & 2 & 0 \\ & 0 & & & & 0 & & & 0 & 2 & 0 & 2 & 1 & 1 \\ & & & & & & & & & 3 & 5 & 1 & 0 & 0 \\ & & & & & & & & & 4 & 1 & 0 & 3 & 6 \end{array}\right).$$

The transpose of A is again a super diagonal matrix given by

$$A^t = \left(\begin{array}{ccc|cc|ccc|cccc}
3 & 1 & 0 & & & & & & & & & \\
1 & 0 & 1 & & 0 & & 0 & & & 0 & & \\
0 & 5 & 0 & & & & & & & & & \\
2 & 0 & 1 & & & & & & & & & \\
\hline
 & & & 3 & 1 & & & & & & & \\
 & 0 & & 4 & 3 & & 0 & & & 0 & & \\
 & & & 5 & 1 & & & & & & & \\
\hline
 & 0 & & & 0 & 8 & 6 & 2 & & 0 & & \\
 & & & & & 1 & -1 & 5 & & & & \\
\hline
 & & & & & & & & 1 & 2 & 3 & 4 \\
 & & & & & & & & 0 & 0 & 5 & 1 \\
 & 0 & & & 0 & & 0 & & 1 & 2 & 1 & 0 \\
 & & & & & & & & 2 & 1 & 0 & 3 \\
 & & & & & & & & 0 & 1 & 0 & 6
\end{array}\right).$$

Now we proceed onto define the notion of super forms on super inner product spaces.

Let $T_s = (T_1 \mid \ldots \mid T_n)$ be a linear operator on a finite $(n_1, \ldots, n_n)$ dimensional super inner product space $V = (V_1 \mid \ldots \mid V_n)$ the super function $f = (f_1 \mid \ldots \mid f_n)$ is defined on $V \times V = (V_1 \times V_1 \mid \ldots \mid V_n \times V_n)$ by

$$\begin{aligned} f(\alpha, \beta) &= (T_s \alpha \mid \beta) \\ &= ((T_1\alpha_1 \mid \beta_1) \mid \ldots \mid (T_n\alpha_n \mid \beta_n)) \end{aligned}$$

may be regarded as a kind of substitute for $T_s$. Many properties about $T_s$ is equivalent to properties concerning $f = (f_1 \mid \ldots \mid f_n)$. In fact we say $f = (f_1 \mid \ldots \mid f_n)$ determines $T_s = (T_1 \mid \ldots \mid T_n)$. If

$$B = (B_1 \mid \ldots \mid B_n) = \left\{ \alpha_1^1 \ldots \alpha_{n_1}^1 \middle| \ldots \middle| \alpha_1^n \ldots \alpha_{n_n}^n \right\}$$

is an orthonormal super basis for V then the entries of the super diagonal matrix of $T_s$ in B are given by

$$
A_{j_K} = \left(\begin{array}{c|c|c|c}
A^1_{j_1K_1} & 0 & & 0 \\ \hline
0 & A^2_{j_2K_2} & & 0 \\ \hline
 & & & \\ \hline
0 & 0 & & A^n_{j_nK_n}
\end{array}\right)
$$

$$
= \left(\begin{array}{c|c|c|c}
f_1(\alpha^1_{K_1}, \alpha^1_{j_1}) & 0 & & 0 \\ \hline
0 & f_2(\alpha^2_{K_2}, \alpha_{j_2}) & & 0 \\ \hline
 & & & \\ \hline
0 & 0 & & f_n(\alpha^n_{K_n}, \alpha^n_{j_n})
\end{array}\right).
$$

Now we proceed onto define the sesqui linear superform.

**DEFINITION 2.1.5**: *A (sesqui-linear) super form on a real or complex supervector space $V = (V_1 \mid \ldots \mid V_n)$ is a superfunction $f$ on $V \times V = (V_1 \times V_1 \mid \ldots \mid V_n \times V_n)$ with values in the field of scalars such that*

i. $f(c\alpha + \beta, \gamma) = cf(\alpha, \gamma) + f(\beta, \gamma)$
   i.e., $(f_1(c_1\alpha_1 + \beta_1, \gamma_1) \mid \ldots \mid f_n(c_n\alpha_n + \beta_n, \gamma_n))$
   $= (c_1f_1(\alpha_1, \gamma_1) + f_1(\beta_1, \gamma_1) \mid \ldots \mid c_nf_n(\alpha_n, \gamma_n) + f_n(\beta_n, \gamma_n))$

ii. $f(\alpha, c\beta + \gamma) = \bar{c}\, f(\alpha, \beta) + f(\alpha, \gamma)$;

*for all* $\alpha = (\alpha_1 \mid \ldots \mid \alpha_n)$, $\beta = (\beta_1 \mid \ldots \mid \beta_n)$ *and* $\gamma = (\gamma_1 \mid \ldots \mid \gamma_n)$ *in* $V = (V_1 \mid \ldots \mid V_n)$ *and* $c = (c_1 \mid \ldots \mid c_n)$, *with* $c_i \in F$; $1 \leq i \leq n$.

*Thus a sesqui linear superform is a super function $f$ on $V \times V$ such that in $f(\alpha, \beta) = (f_1(\alpha_1, \beta_1) \mid \ldots \mid f_n(\alpha_n, \beta_n))$ is a linear super function of α for fixed β and a conjugate linear super function of $\beta = (\beta_1 \mid \ldots \mid \beta_n)$ for fixed $\alpha = (\alpha_1 \mid \ldots \mid \alpha_n)$. In real case $f(\alpha, \beta)$ is linear as a super function of each argument in other words $f$ is a bilinear superform. In the complex case the sesqui linear super form $f$ is not bilinear unless $f = (0 \mid \ldots \mid 0)$.*

**THEOREM 2.1.15**: *Let $V = (V_1 \mid \ldots \mid V_n)$ be a finite $(n_1, \ldots, n_n)$ dimensional inner product super vector space and $f = (f_1 \mid \ldots \mid f_n)$ a super form on V. Then there is a unique linear operator $T_s$ on V such that $f(\alpha, \beta) = (T_s\alpha \mid \beta)$ for all α, β in V; $\alpha = (\alpha_1 \mid \ldots \mid \alpha_n)$ and $\beta = (\beta_1 \mid \ldots \mid \beta_n)$. $(f_1(\alpha_1, \beta_1) \mid \ldots \mid f_n(\alpha_n, \beta_n)) = ((T_1\alpha_1 \mid \beta_1) \mid \ldots \mid (T_n\alpha_n \mid \beta_n))$; for all α, β in V, the super map $f \to T_s$ (i.e. $f_i \to T_i$ for i = 1, 2, …, n) is super isomorphism of the super space of superforms onto SL (V, V).*

*Proof:* Fix a super vector $\beta = (\beta_1 \mid \ldots \mid \beta_n)$ in V.
Then

$$\alpha \to f(\alpha, \beta) \text{ i.e., } \alpha_1 \to f_1(\alpha_1, \beta_1), \ldots, \alpha_n \to f_n(\alpha_n, \beta_n)$$

is a linear super function on V. By earlier results there is a unique super vector $\beta' = (\beta_1' \mid \ldots \mid \beta_n')$ in $V = (V_1 \mid \ldots \mid V_n)$ such that $f(\alpha, \beta) = (\alpha \mid \beta')$

$$\text{i.e., } (f_1(\alpha_1, \beta_1) \mid \ldots \mid f_n(\alpha_n, \beta_n)) = ((\alpha_1 \mid \beta_1') \mid \ldots \mid (\alpha_n, \beta_n'))$$

for every $\beta = (\beta_1 \mid \ldots \mid \beta_n)$ in V. We define a function $U_s$ from V into V by setting

$$U_s\beta = \beta' \text{ i.e. } (U_1\beta_1 \mid \ldots \mid U_n\beta_n) = (\beta_1' \mid \ldots \mid \beta_n').$$

Then

$f(\alpha, c\beta + \gamma) = (\alpha \mid U(c\beta + \gamma))$

$= (f_1(\alpha_1, c_1\beta_1 + \gamma_1) \mid \ldots \mid f_n(\alpha_n, c_n\beta_n + \gamma_n))$
$= (\alpha_1 \mid U_1(c_1\beta_1 + \gamma_1)) \mid \ldots \mid (\alpha_n \mid U_n(c_n\beta_n + \gamma_n))$
$= \bar{c}\, f(\alpha, \beta) + f(\alpha, \gamma)$
$= (\bar{c}_1 f_1(\alpha_1, \beta_1) \mid \ldots \mid \bar{c}_n f_n(\alpha_n, \beta_n)) + (f_1(\alpha_1 \mid \gamma_1) \mid \ldots \mid f_n(\alpha_n, \gamma_n))$
$= (c_1 f_1(\alpha_1, \beta_1) + f_1(\alpha_1, \gamma_1) \mid \ldots \mid + (c_n f_n(\alpha_n, \beta_n) + f_n(\alpha_n, \gamma_n))$
$= \bar{c}(\alpha \mid U_s\beta) + (\alpha \mid U_s\gamma)$
$= (\bar{c}_1\alpha_1 \mid U_1\beta_1) + (\alpha_1 \mid U_1\gamma_1) \mid \ldots \mid (\bar{c}_n\alpha_n \mid U_n\beta_n) + (\alpha_n \mid U_n\gamma_n))$
$= ((\alpha_1 \mid c_1U_1\beta_1 + U_1\gamma_1) \mid \ldots \mid (\alpha_n \mid c_nU_n\beta_n + U_n\gamma_n)$
$= (\alpha \mid cU\beta + U\gamma)$

for all $\alpha = (\alpha_1 \mid \ldots \mid \alpha_n)$, $\beta = (\beta_1 \mid \ldots \mid \beta_n)$, $\gamma = (\gamma_1 \mid \ldots \mid \gamma_n)$ in V = $(V_1 \mid \ldots \mid V_n)$ and for all scalars $c = (c_1 \mid \ldots \mid c_n)$. Thus $U_s$ is a linear operator on $V = (V_1 \mid \ldots \mid V_n)$ and $T_s = U_s^*$ is an operator such that

$$f(\alpha, \beta) = (T\alpha \mid \beta)$$

i.e.,

$$(f_1(\alpha_1 \mid \beta_1) \mid \ldots \mid f_n(\alpha_n \mid \beta_n)) = ((T_1\alpha_1 \mid \beta_1) \mid \ldots \mid (T_n\alpha_n \mid \beta_n));$$

for all α and β in V.

If we also have

$$f(\alpha, \beta) = (T'_s\alpha \mid \beta)$$

i.e.,

$$(f_1(\alpha_1 \mid \beta_1) \mid \ldots \mid f_n(\alpha_n \mid \beta_n)) = ((T'_1\alpha_1 \mid \beta_1) \mid \ldots \mid (T'_n\alpha_n \mid \beta_n)).$$

Then

$$(T_s\alpha - T'_s\alpha \mid \beta) = (0 \mid \ldots \mid 0)$$

i.e.,

$$((T_1\alpha_1 - T'_1\alpha_1 \mid \beta_1) \mid \ldots \mid (T_n\alpha_n - T'_n\alpha_n \mid \beta_n)) = (0 \mid \ldots \mid 0)$$

$\alpha = (\alpha_1 \mid \ldots \mid \alpha_n)$, $\beta = (\beta_1 \mid \ldots \mid \beta_n)$ so $T_s\alpha = T'_s\alpha$ for all $\alpha = (\alpha_1 \mid \ldots \mid \alpha_n)$. Thus for each superform $f = (f_1 \mid \ldots \mid f_n)$ there is a unique linear operator $T_{sf}$ such that

$$f(\alpha, \beta) = (T_{sf}\, \alpha \mid \beta)$$

i.e.,

$$(f_1(\alpha_1 \mid \beta_1) \mid \ldots \mid f_n(\alpha_n \mid \beta_n)) = ((T_{1f_1}\alpha_1 \mid \beta_1) \mid \ldots \mid (T_{nf_n}\alpha_n \mid \beta_n))$$

for all α, β in V.

If f and g are superforms and c a scalar $f = (f_1 \mid \ldots \mid f_n)$ and $g = (g_1 \mid \ldots \mid g_n)$ then

$$(cf + g)(\alpha, \beta) = (T_{scf+g}\alpha \mid \beta)$$

i.e.,

$$((c_1f_1 + g_1)(\alpha_1, \beta_1) \mid \ldots \mid (c_nf_n + g_n)(\alpha_n, \beta_n)) = ((T_{1c_1f_1+g_1}\alpha_1 \mid \beta_1) \mid \ldots \mid (T_{nc_nf_n+g_n}\alpha_n \mid \beta_n)$$

$= (c_1 f_1 (\alpha_1, \beta_1) + g_1(\alpha_1, \beta_1) \mid \ldots \mid c_n f_n (\alpha_n, \beta_n) + g_n(\alpha_n, \beta_n))$
$= ((c_1 T_{1f_1}\alpha_1 \mid \beta_1) + (T_{1g_1}\alpha_1 \mid \beta_1) \mid \ldots \mid (c_n T_{nf_n}\alpha_n \mid \beta_n) + (T_{ng_n}\alpha_n \mid \beta_n))$
$= ((c_1 T_{1f_1} + T_{1g_1})\, \alpha_1 \mid \beta_1) \mid \ldots \mid (c_n T_{nf_n} + T_{ng_n})\, \alpha_n \mid \beta_n))$
$= ((cT_{sf} + T_{sg})\, \alpha \mid \beta)$

for all $\alpha = (\alpha_1 \mid \ldots \mid \alpha_n)$ and $\beta = (\beta_1 \mid \ldots \mid \beta_n)$ in $V = (V_1 \mid \ldots \mid V_n)$. Therefore

$$T_{scf+g} = cT_{sf} + T_{sg},$$

so $f \to T_{sf}$ i.e.,

$$f = (f_1 \mid \ldots \mid f_n) \to (T_{1f_1} \mid \ldots \mid T_{nf_n})$$

is a linear super map. For each $T_s$ in SL(V,V) the equation

$$f(\alpha, \beta) = (T_s\, \alpha \mid \beta)$$
$$(f_1(\alpha_1, \beta_1) \mid \ldots \mid f_n(\alpha_n, \beta_n)) = ((T_1\alpha_1 \mid \beta_1) \mid \ldots \mid (T_n\alpha_n \mid \beta_n));$$

defines a superform such that $T_{sf} = T_s$ and $T_{sf} = (0 \mid \ldots \mid 0)$ if and only if

$$f = (f_1 \mid \ldots \mid f_n) = (0 \mid \ldots \mid 0)$$
$$(T_1 f_1 \mid \ldots \mid T_n f_n) = (T_1 \mid \ldots \mid T_n).$$

Thus $f \to T_{sf}$ is a super isomorphism.

**COROLLARY 2.1.6:** *The super equation*

$$\begin{aligned}(f \mid g) &= ((f_1 \mid g_1) \mid \ldots \mid (f_n \mid g_n)) \\ &= (T_{sf} \mid T^*_{sg}) = ((T_{1f_1} \mid T^*_{1g_1}) \mid \ldots \mid (T_{nf_n} \mid T^*_{ng_n}));\end{aligned}$$

*defines a super inner product on the super space of forms with the property that*

$$(f \mid g) = ((f_1 \mid g_1) \mid \ldots \mid (f_n \mid g_n))$$
$$= \left( \sum_{j_1 k_1} f_1(\alpha^1_{k_1}, \alpha^1_{j_1})\, \overline{g_1(\alpha^1_{k_1}, \alpha^1_{j_1})} \mid \ldots \mid \sum_{j_n k_n} f_n(\alpha^n_{k_n}, \alpha^n_{j_n})\, \overline{g_n(\alpha^n_{k_n}, \alpha^n_{j_n})} \right)$$

*for every orthonormal superbasis* $\left\{\alpha_1^1 \ldots \alpha_{n_1}^1 \mid \ldots \mid \alpha_1^n \ldots \alpha_{n_n}^n\right\}$ *of* $V = (V_1 \mid \ldots \mid V_n)$.

The proof is direct and is left as an exercise for reader.

**DEFINITION 2.1.6:** *If* $f = (f_1 \mid \ldots \mid f_n)$ *is a super form and* $B = (B_1 \mid \ldots \mid B_n) = (\alpha_1^1 \ldots \alpha_{n_1}^1 \mid \ldots \mid \alpha_1^n \ldots \alpha_{n_n}^n)$ *an ordered super basis of* $V = (V_1 \mid \ldots \mid V_n)$*; the super diagonal matrix with entries*

$$\begin{aligned} A_{jk} &= (A_{j_1k_1}^1 \mid \ldots \mid A_{j_nk_n}^n) \\ &= (f_1(\alpha_{k_1}^1, \alpha_{j_1}^1) \mid \ldots \mid f_n(\alpha_{k_n}^n, \alpha_{j_n}^n)) \\ &= f(\alpha_k, \alpha_j) \end{aligned}$$

*is called the super diagonal matrix of f in the ordered super basis B,*

$$A_{j_k} = \left(\begin{array}{c|c|c|c} A_{j_1k_1}^1 & 0 & & 0 \\ \hline 0 & A_{j_2k_2}^2 & & 0 \\ \hline & & & \\ \hline 0 & 0 & & A_{j_nk_n}^n \end{array}\right).$$

**THEOREM 2.1.16:** *Let* $f = (f_1 \mid \ldots \mid f_n)$ *be a superform on a finite* $(n_1, \ldots, n_n)$ *dimensional complex super inner product space* $V = (V_1 \mid \ldots \mid V_n)$ *in which the super diagonal matrix of f is super upper triangular. We say a super form* $f = (f_1 \mid \ldots \mid f_n)$ *on a real or complex super vector space* $V = (V_1 \mid \ldots \mid V_n)$ *is called super Hermitian if* $f(\alpha,\beta) = \overline{f(\beta,\alpha)}$ *i.e.*

$$(f_1(\alpha_1,\beta_1) \mid \ldots \mid f_n(\alpha_n,\beta_n)) = (\overline{f_1(\beta_1,\alpha_1)} \mid \ldots \mid \overline{f_n(\beta_n,\alpha_n)})$$

*for all* $\alpha = (\alpha_1 \mid \ldots \mid \alpha_n)$ *and* $\beta = (\beta_1 \mid \ldots \mid \beta_n)$ *in* $V = (V_1 \mid \ldots \mid V_n)$.

The following theorem is direct and hence is left for the reader to prove.

**THEOREM 2.1.17:** *Let $V = (V_1 \mid ... \mid V_n)$ be a complex super vector space and $f = (f_1 \mid ... \mid f_n)$ a superform on V such that $f(\alpha, \alpha) = (f_1(\alpha_1, \alpha_1) \mid ... \mid f_n(\alpha_n, \alpha_n))$ is real for every $\alpha = (\alpha_1 \mid ... \mid \alpha_n)$ in V. Then $f = (f_1 \mid ... \mid f_n)$ is a Hermitian superform.*

The following corollary which is a direct consequence of the earlier results is stated without proof.

**COROLLARY 2.1.7:** *Let $T_s = (T_1 \mid ... \mid T_n)$ be a linear operator on a complex finite $(n_1, ..., n_n)$ dimensional super inner product super vector space $V = (V_1 \mid ... \mid V_n)$. Then $T_s$ is super self adjoint if and only if $T_s(\alpha \mid \alpha) = ((T_1 \alpha_1 \mid \alpha_1) \mid ... \mid (T_n \alpha_n \mid \alpha_n))$ is real for every $\alpha = (\alpha_1 \mid ... \mid \alpha_n)$ in $V = (V_1 \mid ... \mid V_n)$.*

However we give sketch of the proof analogous to principal axis theorem for super inner product super spaces.

**THEOREM 2.1.18 (PRINCIPAL AXIS THEOREM):** *For every Hermitian super form $f = (f_1 \mid ... \mid f_n)$ on a finite $(n_1, ..., n_n)$ dimensional super inner product space $V = (V_1 \mid ... \mid V_n)$ there is an orthonormal super basis of V for which $f = (f_1 \mid ... \mid f_n)$ is represented by a super diagonal matrix where each component matrix is also diagonal with real entries.*

*Proof:* Let $T_s = (T_1 \mid \ldots \mid T_n)$ be a linear operator such that $f(\alpha, \beta) = (T_s \alpha \mid \beta)$ for all α and β in V
i.e.

$$(f_1(\alpha_1, \beta_1) \mid \ldots \mid f_n(\alpha_n, \beta_n)) = ((T_1\alpha_1 \mid \beta_1) \mid \ldots \mid (T_n\alpha_n \mid \beta_n)).$$

Then since $f(\alpha,\beta) = \overline{f(\beta,\alpha)}$ i.e.

$$(f_1(\alpha_1, \beta_1) \mid \ldots \mid f_n(\alpha_n, \beta_n)) = (\overline{f_1(\beta_1,\alpha_1)} \mid \ldots \mid \overline{f_n(\beta_n,\alpha_n)})$$

and

$$(T_s\beta \mid \alpha) = (\alpha \mid T_s\beta)$$

i.e.

$$((\overline{T_1\beta_1 \mid \alpha_1}) \mid\mid \ldots \mid (\overline{T_n\beta_n \mid \alpha_n})) = ((\alpha_1 \mid T_1\beta_1) \mid \ldots \mid (\alpha_n \mid T_n\beta_n));$$

it follows that

$$(T_s\alpha \mid \beta) = \overline{f(\beta,\alpha)} = (\alpha \mid T_s\beta); \text{ i.e.,}$$
$$((T_1\alpha_1 \mid \beta_1) \mid \ldots \mid (T_n\alpha_n \mid \beta_n) =$$
$$(\overline{f_1(\beta_1,\alpha_1)} \mid \ldots \mid \overline{f_n(\beta_n,\alpha_n)}) =$$
$$((\alpha_1 \mid T_1\beta_1) \mid \ldots \mid (\alpha_n \mid T_n\beta_n));$$

for all α and β; hence $T_s = T_s^*$

i.e. $(T_1 \mid \ldots \mid T_n) = (T_1^* \mid \ldots \mid T_n^*)$. Thus $T_i = T_i^*$ for every i implies a orthonormal basis for each $V_i$; i = 1, 2, …, n; hence an orthonormal superbasis for $V = (V_1 \mid \ldots \mid V_n)$ which consist of characteristic super vectors for $T_s = (T_1 \mid \ldots \mid T_n)$. Suppose

$$\left\{ \alpha_1^1 \ldots \alpha_{n_1}^1 \mid \ldots \mid \alpha_1^n \ldots \alpha_{n_n}^n \right\}$$

is an orthonormal super basis and that

$$T_s\alpha_j = c_j\alpha_j$$

i.e.

$$(T_1\alpha_{j_1}^1 \mid \ldots \mid T_n\alpha_{j_n}^n) = (c_{j_1}^1\alpha_{j_1}^1 \mid \ldots \mid c_{j_n}^n\alpha_{j_n}^n)$$

for $1 \le j_t \le n_t$; t = 1, 2, …, n.

Then

$$f(\alpha_k, \alpha_j) = (f_1(\alpha_{k_1}^1, \alpha_{j_1}^1) \mid \ldots \mid f_n(\alpha_{k_n}^n, \alpha_{j_n}^n))$$
$$= (T_s\alpha_k \mid \alpha_j) = ((T_1\alpha_{k_1}^1 \mid \alpha_{j_1}^1) \mid \ldots \mid (T_n\alpha_{k_n}^n \mid \alpha_{j_n}^n))$$
$$= \delta_{kj}c_k = (\delta_{k_1j_1}^1 c_{k_1}^1 \mid \ldots \mid \delta_{k_nj_n}^n c_{k_n}^n).$$

Now we proceed onto define the notion of positive superforms.

**DEFINITION 2.1.7:** *A superform $f = (f_1 \mid \ldots \mid f_n)$ on a real or complex super vector space $V = (V_1 \mid \ldots \mid V_n)$ is supernon-negative if it is super Hermitian and $f(\alpha, \alpha) \ge (0 \mid \ldots \mid 0)$ for* *every α in V; i.e.* $(f_1(\alpha_1, \alpha_1) \mid \ldots \mid f_n(\alpha_n, \alpha_n)) \geq (0 \mid \ldots \mid 0)$ *i.e. each* $f_j(\alpha_j, \alpha_j) \geq 0$ *for every* $j = 1, 2, \ldots, n$. *The form is super positive if f is super Hermitian and* $f(\alpha, \alpha) > (0 \mid \ldots \mid 0)$ *i.e.* $(f_1(\alpha_1, \alpha_1) \mid \ldots \mid f_n(\alpha_n, \alpha_n)) > (0 \mid \ldots \mid 0)$ *i.e.* $f_j(\alpha_j, \alpha_j) > 0$ *for every* $j = 1, 2, \ldots, n$.

*The super Hermitian form f is quasi super positive or equivalently quasi super non negative (both mean one and the same) if in* $f(\alpha, \alpha) = (f_1(\alpha_1, \alpha_1) \mid \ldots \mid f_n(\alpha_n, \alpha_n))$ *some* $f_j(\alpha_j, \alpha_j) > 0$ *and some* $f_i(\alpha_i, \alpha_i) \geq 0$; $i \neq j$; $1 \leq i \leq n$.

All properties related with usual non negative and positive Hermitian form can be appropriately extended in case of Hermitian superform.

**THEOREM 2.1.19**: *Let F be the field of real numbers or the field of complex numbers. Let A be a super diagonal matrix of the form*

$$A = \left(\begin{array}{c|c|c|c} A_1 & 0 & & 0 \\ \hline 0 & A_2 & & 0 \\ \hline & & & \\ \hline 0 & 0 & & A_n \end{array}\right)$$

*be a* $(n_1 \times n_1, \ldots, n_n \times n_n)$ *matrix over F. The super function* $g = (g_1 \mid \ldots \mid g_n)$ *defined by* $g(X, Y) = Y^*AX$ *is a positive superform on the super space* $(F^{n_1 \times 1} \mid \ldots \mid F^{n_n \times 1})$ *if and only if there exists an invertible super diagonal matrix*

$$P = \left(\begin{array}{c|c|c|c} P_1 & 0 & & 0 \\ \hline 0 & P_2 & & 0 \\ \hline & & & \\ \hline 0 & 0 & & P_n \end{array}\right).$$

*Each* $P_i$ *is a* $n_i \times n_i$ *matrix* $i = 1, 2, \ldots, n$ *with entries from F such that* $A = P^* P$ *; i.e.,*

$$A = \left(\begin{array}{c|c|c|c} A_1 & 0 & & 0 \\ \hline 0 & A_2 & & 0 \\ \hline & & & \\ \hline 0 & 0 & & A_n \end{array}\right)$$

$$= \left(\begin{array}{c|c|c|c} P_1^* P_1 & 0 & \ldots & 0 \\ \hline 0 & P_2^* P_2 & \ldots & 0 \\ \hline & & \ldots & \\ \hline 0 & 0 & \ldots & P_n^* P_n \end{array}\right).$$

**DEFINITION 2.1.8:** *Let*

$$A = \left(\begin{array}{c|c|c|c} A_1 & 0 & & 0 \\ \hline 0 & A_2 & & 0 \\ \hline & & & \\ \hline 0 & 0 & & A_n \end{array}\right)$$

*be a superdiagonal matrix with each $A_i$ a $n_i \times n_i$ matrix over the field F; i =1, 2, ..., n. The principal super minor of A or super principal minors of A (both mean the same) are scalars*

$$\Delta_k(A) = (\Delta_{k_1}(A_1) \mid \ldots \mid \Delta_{k_n}(A_n))$$

*defined by*

$$\Delta_k(A) = superdet \left\{ \left(\begin{array}{c|c|c} \begin{matrix} A^1_{11} \ldots A^1_{1k_1} \\ \vdots \quad \vdots \\ A^1_{k_1 1} \ldots A^1_{k_1 k_1} \end{matrix} & 0 & 0 \\ \hline 0 & \begin{matrix} A^2_{11} \ldots A^2_{1k_2} \\ \vdots \quad \vdots \\ A^2_{k_2 1} \ldots A^2_{k_2 k_2} \end{matrix} & 0 \\ \hline 0 & 0 & \begin{matrix} A^n_{11} \ldots A^2_{1k_n} \\ \vdots \quad \vdots \\ A^n_{k_n 1} \ldots A^n_{k_n k_n} \end{matrix} \end{array}\right) \right\}$$

$$= \left( \det \begin{pmatrix} A^1_{11} & \cdots & A^1_{1k_1} \\ \vdots & & \vdots \\ A^1_{k_1 1} & \cdots & A^1_{k_1 k_1} \end{pmatrix}, \ \ldots \ , \det \begin{pmatrix} A^n_{11} & \cdots & A^n_{1k_n} \\ \vdots & & \vdots \\ A^n_{k_n 1} & \cdots & A^n_{k_n k_n} \end{pmatrix} \right)$$

*for* $1 \le k_t \le n_t$ *and* $t = 1, 2, \ldots, n$.

Several other interesting properties can also be derived for these superdiagonal matrices.

We give the following interesting theorem and the proof is left for the reader.

**THEOREM 2.1.20:** *Let* $f = (f_1 \mid \ldots \mid f_n)$ *be a superform on a finite* $(n_1, \ldots, n_n)$ *dimensional supervector space* $V = (V_1 \mid \ldots \mid V_n)$ *and let A be a super diagonal matrix of f in an ordered superbasis B =* $(B_1 \mid \ldots \mid B_n)$*. Then f is a positive superform if and only if A = A* and the principal super minor of A are all positive.*
*i.e.*

$$A = \left(\begin{array}{c|c|c|c} A_1 & 0 & & 0 \\ \hline 0 & A_2 & & 0 \\ \hline & & & \\ \hline 0 & 0 & & A_n \end{array}\right)$$

$$= \left(\begin{array}{c|c|c|c} A_1^* & 0 & & 0 \\ \hline 0 & A_2^* & & 0 \\ \hline & & & \\ \hline 0 & 0 & & A_n^* \end{array}\right).$$

***Note**:* The principal minor of $(A_1 \mid \ldots \mid A_n)$ is called as the principal superminors of A or with default of notation the principal minors of $\{A_1, \ldots, A_n\}$ is called the principal super minors of A.

$T_s = (T_1 \mid \ldots \mid T_n)$ a linear operator on a finite $(n_1, \ldots, n_n)$ dimensional super inner product space $V = (V_1 \mid \ldots \mid V_n)$ is said to be super non-negative if $T_s = T_s^*$ i.e.

$$(T_1 \mid \ldots \mid T_n) = (T_1^* \mid \ldots \mid T_n^*)$$

i.e. $T_i = T_i^*$ for i = 1, 2, …, n and

$$(T_s \alpha \mid \alpha) = ((T_1 \alpha_1 \mid \alpha_1) \mid \ldots \mid (T_n \alpha_n \mid \alpha_n)) \geq (0 \mid \ldots \mid 0)$$

for all $\alpha = (\alpha_1 \mid \ldots \mid \alpha_n)$ in V.

A super positive linear operator is one such that $T_s = T_s^*$ and

$$(T \alpha \mid \alpha) = ((T_1 \alpha_1 \mid \alpha_1) \mid \ldots \mid (T_n \alpha_n \mid \alpha_n)) > (0 \mid \ldots \mid 0)$$

for all $\alpha = (\alpha_1 \mid \ldots \mid \alpha_n) \neq (0 \mid \ldots \mid 0)$.

Several properties enjoyed by positive operators and non negative operators will also be enjoyed by the super positive operators and super non negative operators on super vector spaces, with pertinent and appropriate modification. Throughout the related matrix for these super operators $T_s$ will always be a super diagonal matrix A of the form

$$A = \left(\begin{array}{c|c|c|c} A_1 & 0 & & 0 \\ \hline 0 & A_n & & 0 \\ \hline & & & \\ \hline 0 & 0 & & A_n \end{array}\right)$$

where each $A_i$ is a $n_i \times n_i$ square matrix, $1 \leq i \leq n$, $A = A^*$ and the principal minors of each $A_i$ are positive; $1 \leq i \leq n$.

Now we just mention one more property about the super forms.

**THEOREM 2.1.21:** *Let $f = (f_1 \mid \ldots \mid f_n)$ be a super form on a real or complex super vector space $V = (V_1 \mid \ldots \mid V_n)$ and $\left\{\alpha_1^1 \ldots \alpha_{r_1}^1 \mid \ldots \mid \alpha_1^n \ldots \alpha_{r_n}^n\right\}$ a super basis for the finite dimensional super subvector space $W = (W_1 \mid \ldots \mid W_n)$ of $V = (V_1 \mid \ldots \mid V_n)$.*

*Let M be the super square diagonal matrix where each $M_i$ in M; is a $r_i \times r_i$ super matrix with entries $(1 \leq i \leq n)$.* $M^i_{jk} = f_i(\alpha^i_k, \alpha^i_j)$, *i.e.*

$$M = \left(\begin{array}{c|c|c|c} M_1 & 0 & & 0 \\ \hline 0 & M_2 & & 0 \\ \hline & & & \\ \hline 0 & 0 & & M_n \end{array}\right)$$

$$= \left(\begin{array}{c|c|c|c} f^1(\alpha^1_{k_1}, \alpha^1_{j_1}) & 0 & & 0 \\ \hline 0 & f^2(\alpha^2_{k_2}, \alpha^2_{j_2}) & & 0 \\ \hline & & & \\ \hline 0 & 0 & & f^n(\alpha^n_{k_n}, \alpha^n_{j_n}) \end{array}\right)$$

*and* $W' = (W'_1 | \ldots | W'_n)$ *be the set of all super vectors* $\beta = (\beta_1 | \ldots | \beta_n)$ *of V and* $W \cap W' = (W_1 \cap W'_1 | \ldots | W_n \cap W'_n) = (0 | \ldots | 0)$ *if and only if*

$$M = \left(\begin{array}{c|c|c|c} M_1 & 0 & & 0 \\ \hline 0 & M_2 & & 0 \\ \hline & & & \\ \hline 0 & 0 & & M_n \end{array}\right)$$

*is invertible. When this is the case,* $V = W + W'$ *i.e.* $V = (V_1 | \ldots | V_n) = (W_1 + W'_1 | \ldots | W_n + W'_n)$.

The proof can be obtained as a matter of routine.

The projection $E_s = (E_1 | \ldots | E_n)$ constructed in the proof may be characterized as follows.

$$E_s\beta = \alpha;$$
$$(E_1\beta_1 | \ldots | E_n\beta_n) = (\alpha_1 | \ldots | \alpha_n)$$

is in W and β – α belongs $W' = (W_1' \mid \ldots \mid W_n')$. Thus $E_s$ is independent of the super basis of W = ($W_1$ | … | $W_n$) that was used in this construction. Hence we may refer to $E_s$ as the super projection of V on W that is determined by the direct sum decomposition.

$$V = W \oplus W';$$
$$(V_1 \mid \ldots \mid V_n) = (W_1 \oplus W_1' \mid \ldots \mid W_n \oplus W_n').$$

Note that $E_s$ is a super orthogonal projection if and only if W' = $W^{\perp} = (W_1^{\perp} \mid \ldots \mid W_n^{\perp})$. Now we proceed onto develop the analogous of spectral theorem which we call as super spectral theorem.

**THEOREM 2.1.22 (SUPER SPECTRAL THEOREM):** *Let $T_s = (T_1 \mid \ldots \mid T_n)$ be a super normal operator on a finite $(n_1 \mid \ldots \mid n_n)$ dimensional complex super inner product super space $V = (V_1 \mid \ldots \mid V_n)$ or a self-adjoint super operator on a finite super dimensional real inner product super space $V = (V_1 \mid \ldots \mid V_n)$.*

*Let $\{(c_1^1, \ldots c_{k_1}^1 \mid \ldots \mid c_1^n \ldots c_{k_n}^n)\}$ be the distinct characteristic super values of $T_s = (T_1 \mid \ldots \mid T_n)$.*

*Let $W_j = (W_{j_1}^1 \mid \ldots \mid W_{j_n}^n)$ be the characteristic super space associated with $c_{j_t}^t$ of $E_{j_t}^t$, the orthogonal super projection of $V = (V_1 \mid \ldots \mid V_n)$ on $W_j = (W_{j_1}^1 \mid \ldots \mid W_{j_n}^n)$.*

*Then $W_j$ is super orthogonal to $W_i = (W_{i_1}^1 \mid \ldots \mid W_{i_n}^n)$ when $i_t \neq j$, t = 1, 2, …, n; $V = (V_1 \mid \ldots \mid V_n)$ is the super direct sum of $W_1, \ldots, W_k$ and*

$$T_s = (c_1^1 E_1^1 + \ldots + c_{k_1}^1 E_{k_1}^1 \mid \ldots \mid c_1^n E_1^n + \ldots c_{k_n}^n E_{k_n}^n)$$
$$= (T_1 \mid \ldots \mid T_n) \qquad I$$

*This super decomposition I is called the spectral super resolution of $T_s = (T_1 \mid \ldots \mid T_n)$.*

Several interesting results can be derived in this direction.

The following result which is mentioned below would be useful in solving practical problems.

Let $E_s = (E^1 \mid \ldots \mid E^n)$ be a super orthogonal projection where each $E^t = E_1^t \ldots E_{k_t}^t$ ; t = 1, 2, …, n.

If

$$(e_{j_1}^1 \ldots e_{j_n}^n) =$$

$$\left( \prod_{i_1 \neq j_1} \left( \frac{x - c_{i_1}^1}{c_{j_1}^1 - c_{i_1}^1} \right) \middle| \ldots \middle| \prod_{i_n \neq j_n} \left( \frac{x - c_{1_n}^n}{c_{j_n}^n - c_{i_n}^n} \right) \right)$$

then $E_{j_t}^t = e_{j_t}^t (T_t)$ for $1 \leq j_t \leq k_t$ and t = 1, 2, …, n.

$$(E_1^1 \ldots E_{k_1}^1, \ \ldots \ , E_1^n \ldots E_{k_n}^n)$$

are canonically super associated with $T_s$ and

$$I = (I_1 | \ \ldots \ | I_n) = (E_1^1 + \ldots + E_{k_1}^1 | \ldots \ | E_1^n + \ldots + E_{k_n}^n)$$

the family of super projections $\left\{E_1^1 \ldots E_{k_1}^1, \ldots, E_1^n \ldots E_{k_n}^n\right\}$ is called the super resolution of the super identity defined by $T_s$.

Thus we have the following interesting definition about super diagonalizable normal operators.

**DEFINITION 2.1.9**: *Let $T_s = (T_1 \mid \ldots \mid T_n)$ be a super diagonalizable normal operator on a finite $(n_1, \ldots, n_n)$ dimensional inner product super space and*

$$T_s = (T_1 \mid \ldots \mid T_n) = \left( \sum_{j_1=1}^{k_1} c_{j_1}^1 E_{j_1}^1 \middle| \ldots \middle| \sum_{j_n=1}^{k_n} c_{j_n}^n E_{j_n}^n \right)$$

*its super spectral resolution. Suppose $f = (f_1 \mid \ldots \mid f_n)$ is a super function whose super domain includes the super spectrum of $T_s = (T_1 \mid \ldots \mid T_n)$ that has values in the field of scalars F. Then the linear operator $f(T_s) = (f_1(T_1) \mid \ldots \mid f_n(T_n))$ is defined by the equation*

$$f(T_s) = \left( \sum_{j_1=1}^{k_1} f_1(c_{j_1}^1) E_{j_1}^1 \mid \ \ldots \ \mid \sum_{j_n=1}^{k_n} f_n(c_{j_n}^n) E_{j_n}^n \right).$$

Based on this property we have the following interesting theorem.

**THEOREM 2.1.23**: *Let $T_s = (T_1 \mid \ldots \mid T_n)$ be a super diagonalizable normal operator with super spectrum $S = (S_1 \mid \ldots \mid S_n)$ on a finite $(n_1, \ldots, n_n)$ dimensional super inner product super vector space $V = (V_1 \mid \ldots \mid V_n)$. Suppose $f = (f_1 \mid \ldots \mid f_n)$ is a function whose super domain contains S that has super values in the field of scalars. Then $f(T_s) = (f_1(T_1) \mid \ldots \mid f_n(T_n))$ is a super diagonalizable normal operator with super spectrum*

$$f(S_s) = (f_1(S_1) \mid \ldots \mid f_n(S_n)).$$

*If*

$$U_s = (U_1 \mid \ldots \mid U_n)$$

*is a unitary super map of V onto*

$$V' = (V_1' \mid \ldots \mid V_n') \text{ and } T_s' = U_s T_s U_s^{-1}$$
$$= (T_1' \mid \ldots \mid T_n') = (U_1 T_1 U_1^{-1} \mid \ldots \mid U_n T_n U_n^{-1})\,;$$

*then*

$$S = (S_1 \mid \ldots \mid S_n)$$

*is the super spectrum of*

$$T_s' = (T_1' \mid \ldots \mid T_n')$$

*and*

$$f(T') = (f_1(T_1') \mid \ldots \mid f_n(T_n'))$$
$$= (U_1 f_1(T_1)\, U_1^{-1} \mid \ldots \mid U_n f_n(T_n) U_n^{-1}).$$

*Proof:* The normality of $f(T) = (f_1(T_1) \mid \ldots \mid f_n(T_n))$ follows by a simple computation from

$$f(T) = \left( \sum_{j_1=1}^{k_1} f_1(c^1_{j_1}) E^1_{j_1} \;\Big|\; \ldots \;\Big|\; \sum_{j_n=1}^{k_n} f_n(c^n_{j_n}) E^n_{j_n} \right)$$

and the fact that

$$f(T)^* = (f_1(T_1)^* \mid \ldots \mid f_n(T_n)^*)$$
$$= \left( \sum_{j_1} \overline{f_1(c^1_{j_1})}\, E^1_{j_1} \Big| \ldots \Big| \sum_{j_n} \overline{f_n(c^n_{j_n})} E^n_{j_n} \right).$$

Moreover it is clear that for every $\alpha^t = (\alpha^t_1 \mid \ldots \mid \alpha^t_n)$ in $E^t_{j_t}(V_t)$; $t = 1, 2, \ldots, n$;

$$f_t(T_t)\alpha^t = f_t(c^t_{j_t})\, \alpha^t .$$

Thus the superset $f(S) = (f_1(S_1) \mid \ldots \mid f_n(S_n))$ for all $f(c) = (f_1(c_1) \mid \ldots \mid f_n(c_n))$ in $S = (S_1 \mid \ldots \mid S_n)$ is contained in the superspectrum of $f(S) = (f_1(T_1) \mid \ldots \mid f_n(T_n))$. Conversely suppose $\alpha = (\alpha^1 \mid \ldots \mid \alpha^n) \neq (0 \mid \ldots \mid 0)$ and that $f(T)\alpha = b\,\alpha$ i.e.

$$(f_1(T_1)\alpha^1 \mid \ldots \mid f_n(T_n)\alpha^n) = (b_1\alpha_1 \mid \ldots \mid b_n\alpha_n).$$

Then

$$\alpha = \left( \sum_{j_1} E^1_{j_1}\alpha^1 \;\Big|\; \ldots \;\Big|\sum_{j_n} E^n_{j_n}\alpha^n \right)$$

and $\alpha = (f_1(T_1)\alpha^1 \mid \ldots \mid f_n(T_n)\alpha^n)$

$$= \left( \sum_{j_1} f_1(T_1) E^1_{j_1}\alpha^1 \;\Big|\; \ldots \;\Big|\sum_{j_n} f_n(T_n) E^n_{j_n}\alpha^n \right)$$

$$= \left( \sum_{j_1} f_1(c^1_{j_1}) E^1_{j_1}\alpha^1 \Big|\; \ldots \;\Big|\; \sum_{j_n} f_n(c^n_{j_n}) E^n_{j_n}\alpha^n \right)$$

$$= \left( \sum_{j_1} b_1 E^1_{j_1}\alpha^1 \;\Big|\; \ldots \;\Big|\; \sum_{j_n} b_n E^n_{j_n}\alpha^n \right).$$

Hence

$$\left\| \sum_{j} (f(c_{j_1}) - b) E_j \alpha \right\|^2$$

$$= \left( \left\| \sum_{j_1} (f_1(c^1_{j_1}) - b_1) E^1_{j_1} \alpha^1 \right\|^2 \Big| \ldots \Big| \left\| \sum_{j_n} (f_n(c^n_{j_n}) - b_n) E^n_{j_n} \alpha^n \right\|^2 \right)$$

$$= \left( \sum_{j_1} \left| f_1(c^1_{j_1}) - b_1 \right|^2 \left\| E^1_{j_1} \alpha^1 \right\|^2 \Big| \; \ldots \; \Big| \sum_{j_n} \left| f_n(c^n_{j_n}) - b_n \right|^2 \left\| E^n_{j_n} \alpha^n \right\|^2 \right)$$

$$= (0 \mid \ldots \mid 0).$$

Therefore

$$f(c_j) = (f_1(c^1_{j_1}) \big| \ldots \big| f_n(c^n_{j_n})) = (b_1 \mid \ldots \mid b_n)$$

or

$$E_j\alpha = (0 \mid \ldots \mid 0)$$

i.e.,

$$(E^1_{j_1}\alpha^1 \big| \ldots \big| E^n_{j_n}\alpha^n) = (0 \mid \ldots \mid 0).$$

By assumption $\alpha = (\alpha_1 \mid \ldots \mid \alpha_n) \neq (0 \mid \ldots \mid 0)$ so there exists indices $i = (i_1, \ldots, i_n)$ such that $E_i\alpha = (E^1_{i_1}\alpha^1 \big| \ldots \big| E^n_{i_n}\alpha^n) \neq (0 \mid \ldots \mid 0)$. It follows that $f(c_j) = (f_1(c^1_{j_1}) \big| \ldots \big| f_n(c^n_{j_n})) = (b_1 \mid \ldots \mid b_n)$ and hence that $f(S) = (f_1(S_1) \mid \ldots \mid f_n(S_n))$ is the super spectrum of $f(T) = (f_1(T_1) \mid \ldots \mid f_n(T_n))$. In fact that $f(S) = (f_1(S_1) \mid \ldots \mid f_n(S_n)) = \left\{ b^1_1 \ldots b^{r_1}_1, \ldots, b^1_n \ldots b^{r_n}_n \right\}$ where $b^{r_1}_{m_t} \neq b^{r_1}_{n_t}$ when $m_t \neq n_t$ for t = 1, 2, …, n. Let $X_m = (X_{m_1} \big| \ldots \big| X_{m_n})$ indices $i = (i_1, \ldots, i_n)$ such that $1 \leq i_t \leq k_t$; t = 1, 2, …, n and

$$f(c_i) = (f_1(c^1_{i_1}) \big| \ldots \big| f_n(c^n_{i_n})) = (b^1_{m_1} \big| \ldots \big| b^n_{m_n}).$$

Let

$$P_m = \left( \sum_{i_1} E^1_{i_1} \Big| \ldots \Big| \sum_{i_n} E^n_{i_n} \right) = (P^1_{m_1} \big| \ldots \big| P^n_{m_n})$$

the super sum being extended over the indices i = ($i_1$, …, $i_n$) in $X_m = (X_{m_1} | \ldots | X_{m_n})$. Then $P_m = (P^1_{m_1} | \ldots | P^n_{m_n})$ is the super orthogonal projection of V = ($V_1$ | … | $V_n$) on the super subspace of characteristic super vectors belonging to the characteristic super values $b_m = (b^1_{m_1} | \ldots | b^n_{m_n})$ of f(T) = ($f_1(T_1)$ | … | $f_n(T_n)$) and

$$f(T) = \left( \sum_{m_1=1}^{r_1} b^1_{m_1} P^1_{m_1} \; | \ldots | \; \sum_{m_n=1}^{r_n} b^n_{m_n} P^n_{m_n} \right)$$

is the super spectral resolution (or spectral super resolution) of f(T) = ($f_1(T_1)$ | … | $f_n(T_n)$).

Now suppose $U_S$ = ($U_1$ | … | $U_n$) is unitary transformation of V = ($V_1$ | … | $V_n$) onto $V' = (V'_1 | \ldots | V'_n)$ and that

$$T'_s = U_s T_s U_s^{-1}; \; (T'_1 | \ldots | T'_n) = (U_1 T_1 U_1^{-1} | \ldots | U_n T_n U_n^{-1}).$$

Then the equation

$$T_s\alpha = c\alpha; \; (T_1\alpha_1 | \ldots | T_n\alpha_n) = (c_1\alpha_1 | \ldots | c_n\alpha_n)$$

holds good if and only if $T_s' U_s\alpha = cU_s\alpha$ i.e.

$$(T_1' U_1\alpha \,|\, \ldots \,|\, T_n' U_n\alpha_n) = (c_1 U_1\alpha_1 \,|\, \ldots \,|\, c_n U_n\alpha_n).$$

Thus S = ($S_1$ | … | $S_n$) is the super spectrum of $T_s'$ and $U_s$ maps each characteristic super subspace for $T_s$ onto the corresponding super subspace for $T_s'$.
In fact

$$f(T) = \sum_{j=1}^{k} f(c_j) E_j$$

i.e.

$$(f_1(T_1) | \ldots | f_n(T_n)) = \left( \sum_{j_1=1}^{k_1} f_1(c^1_{j_1}) E^1_{j_1} \;|\; \ldots \;|\; \sum_{j_n=1}^{k_n} f_n(c^n_{j_n}) E^n_{j_n} \right)$$

where

$$f_t(T_t) = \sum_{j_t=1}^{k_t} f_t(c^t_{j_t}) E^t_{j_t}$$

for i = 1, 2, …, n.
We see that

$$T_s' = (T_1' \mid \dots \mid T_n')$$
$$= \left( \sum_{j_1} c^1_{j_1} E^1_{j_1} \Big| \dots \Big| \sum_{j_n} c^n_{j_n} E^n_{j_n} \right)$$

$$E'_j = (E'^1_{j_1} \mid \dots \mid E'^n_{j_n}) = \left( U_1 E^1_{j_1} U_1^{-1} \mid \dots \mid U_n E^n_{j_n} U_n^{-1} \right)$$

is the super spectral resolution of $T_s' = (T_1' \mid \dots \mid T_n')$.
Hence

$$f(T') = (f_1(T_1') \mid \dots \mid f_n(T_n'))$$

$$= \left( \sum_{j_1} f_1(c^1_{j_1}) E^1_{j_1} \Big| \dots \Big| \sum_{j_n} f_n(c^n_{j_n}) E^1_{j_n} \right)$$

$$= \left( \sum_{j_1} f_1(c^1_{j_1}) U_1 E^1_{j_1} U_1^{-1} \Big| \dots \Big| \sum_{j_n} f_n(c^n_{j_n}) U_n E^n_{j_n} U_n^{-1} \right)$$

$$= \left( U_1 \sum_{j_1} f_1(c^1_{j_1}) E^1_{j_1} U_1^{-1} \Big| \dots \Big| U_n \sum_{j_n} f_n(c^n_{j_n}) E^n_{j_n} U_n^{-1} \right)$$

$$= \quad U_s \sum_{j} f(c_j) E_j U_s^{-1} = U_s f(T_s) U_s^{-1}.$$

The following corollary is direct and is left as an exercise for the reader to prove.

**COROLLARY 2.1.8:** *With the assumption of the theorem just proved suppose $T_s = (T_1 \mid \dots \mid T_n)$ is represented by the super*

*basis B = (B₁ | … | Bₙ) =* $\left\{\alpha_1^1 \ldots \alpha_{n_1}^1 \mid \ldots \mid \alpha_1^n \ldots \alpha_{n_n}^n\right\}$ *by the superdiagonal matrix*

$$D = \begin{pmatrix} D_1 & 0 & & 0 \\ 0 & D_2 & & 0 \\ & & & \\ 0 & 0 & & D_n \end{pmatrix}$$

*with entries* $(d_1^1 \ldots d_{n_1}^1; \ldots; d_1^n \ldots d_{n_n}^n)$. *Then in the superbasis B, f(T) = (f₁(T₁) | … | fₙ(Tₙ)) is represented by the super diagonal matrix f(D) = (f₁(D₁) | … | fₙ(Dₙ)) with entries* $(f_1(d_1^1) \ldots f_1(d_{n_1}^1); \ldots; f_n(d_1^n) \ldots f_n(d_{n_n}^n))$. *If B' = (B₁' | … | Bₙ') =* $\left\{\beta_1^1 \ldots \beta_{n_1}^1 \mid \ldots \mid \beta_1^n \ldots \beta_{n_n}^n\right\}$ *is another ordered superbasis and P' = (P₁' | … | Pₙ') the super diagonal matrix such that*

$$\beta_{j_t}^t = \sum_{i_t} P_{i_t j_t}^t \alpha_{i_t}^t$$

$$\text{i.e. } (\beta_{j_1}^1 \mid \ldots \mid \beta_{j_n}^n) = \left(\sum_{i_1} P_{i_1 j_1}^1 \alpha_{i_1}^1 \mid \ldots \mid \sum_{i_n} P_{i_n j_n}^n \alpha_{i_n}^n\right)$$

*then*

$$P^{-1} f(D) P = \left(P_1^{-1} f_1(D_1) P_1 \mid \ldots \mid P_n^{-1} f_n(D_n) P_n\right)$$

*is the super diagonal matrix of (f₁(T₁) | … | fₙ(Tₙ)) in the super basis B' = (B₁' | … | Bₙ').*

Thus this enables one to understand that certain super functions of a normal super diagonal matrix. Suppose

$$A = \begin{pmatrix} A_1 & 0 & & 0 \\ 0 & A_2 & & 0 \\ & & & \\ 0 & 0 & & A_n \end{pmatrix}$$

is the normal super diagonal matrix. Then there is an invertible super diagonal matrix

$$P = \left(\begin{array}{c|c|c|c} P_1 & 0 & & 0 \\ \hline 0 & P_2 & & 0 \\ \hline & & & \\ \hline 0 & 0 & & P_n \end{array}\right)$$

in fact superunitary P = ($P_1$ | … | $P_n$) described above as a superdiagonal matrix such that

$$PAP^{-1} = (P_1A_1 P_1^{-1} \mid \ldots \mid P_nA_n P_n^{-1})$$

$$= \left(\begin{array}{c|c|c|c} P_1A_1P_1^{-1} & 0 & & 0 \\ \hline 0 & P_2A_2P_2^{-1} & & 0 \\ \hline & & & \\ \hline 0 & 0 & & P_nA_nP_n^{-1} \end{array}\right)$$

is a super diagonal matrix i.e. each $P_iA_iP_i^{-1}$ is a diagonal matrix say D = ($D_1$ | … | $D_n$) with entries $d_1^1 \ldots d_{n_1}^1, \ldots, d_1^n \ldots d_{n_n}^n$.

Let f = ($f_1$ | … | $f_n$) be a complex valued superfunction which can be applied to $d_1^1 \ldots d_{n_1}^1, \ldots, d_1^n \ldots d_{n_n}^n$ and let f(D) = ($f_1(D_1)$ | … | $f_n(D_n)$) be the superdiagonal matrix with entries

$$f_1(d_1^1) \ldots f_1(d_{n_1}^1), \ldots, f_n(d_1^n) \ldots f_n(d_{n_n}^n).$$

Then

$$P^{-1}f(D)\,P = (P_1^{-1}f_1(D_1)P_1 \mid \ldots \mid P_n^{-1}f_n(D_n)P_n)$$

is independent of D = ($D_1$ | … | $D_n$) and just a super function of A in the following ways.

If

$$Q = \left(\begin{array}{c|c|c|c} Q_1 & 0 & & 0 \\ \hline 0 & Q_2 & & 0 \\ \hline & & & \\ \hline 0 & 0 & & Q_n \end{array}\right)$$

is another super invertible super diagonal matrix such that

$$QAQ^{-1} = \left(\begin{array}{c|c|c|c} Q_1A_1Q_1^{-1} & 0 & & 0 \\ \hline 0 & Q_2A_2Q_2^{-1} & & 0 \\ \hline & & & \\ \hline 0 & 0 & & Q_nA_nQ_n^{-1} \end{array}\right)$$

is a superdiagonal matrix

$$D' = \left(\begin{array}{c|c|c|c} D_1' & 0 & & 0 \\ \hline 0 & D_2' & & 0 \\ \hline & & & \\ \hline 0 & 0 & & D_n' \end{array}\right) = (D_1' \,|\ldots|\, D_n')$$

then f = ($f_1$ | … | $f_n$) may be applied to the super diagonal entries of $D' = P^{-1} f(D) P = Q^{-1} f(D') Q$ under these conditions

$$f(A) = \left(\begin{array}{c|c|c|c} f_1(A_1) & 0 & & 0 \\ \hline 0 & f_2(A_2) & & 0 \\ \hline & & & \\ \hline 0 & 0 & & f_n(A_n) \end{array}\right)$$

is defined as

$$P^{-1}f(D)P = \left(\begin{array}{c|c|c|c} P_1^{-1}f_1(D_1)P_1 & 0 & & 0 \\ \hline 0 & P_2^{-1}f_2(D_2)P_2 & & 0 \\ \hline & & & \\ \hline 0 & 0 & & P_n^{-1}f_n(D_n)P_n \end{array}\right).$$

The superdiagonal matrix

$$f(A) = \left(\begin{array}{c|c|c|c} f_1(A_1) & 0 & & 0 \\ \hline 0 & f_2(A_2) & & 0 \\ \hline & & & \\ \hline 0 & 0 & & f_n(A_n) \end{array}\right)$$

$$= (f_1(A_1) \mid \ldots \mid f_n(A_n))$$

may also be characterized in a different way.

**THEOREM 2.1.24:** *Let*

$$A = \left(\begin{array}{c|c|c|c} A_1 & 0 & & 0 \\ \hline 0 & A_2 & & 0 \\ \hline & & & \\ \hline 0 & 0 & & A_n \end{array}\right)$$

*be a normal superdiagonal matrix and* $\left\{c_1^1 \ldots c_{k_1}^1, \ldots, c_1^n \ldots c_{k_n}^n\right\}$ *be the distinct complex super root of the super*

$$det\ (xI - A) = (\det\ (xI_1 - A_1) \mid \ldots \mid \det\ (xI_n - A_n)) \ .$$

*Let*

$$e_i = (e_{i_1}^1 \mid \ldots \mid e_{i_n}^n) = \left( \prod_{j_1 \neq i_1} \left( \frac{x - c_{j_1}^1}{c_{i_1}^1 - c_{j_1}^1} \right) \Big| \ldots \Big| \prod_{j_n \neq i_n} \left( \frac{x - c_{j_n}^n}{c_{i_n}^n - c_{j_n}^n} \right) \right)$$

*and*

$$E_i = (E_{i_1}^1 \mid \ldots \mid E_{i_n}^n) = e_i\,(A) = (e_{i_1}^1 (A_1) \mid \ldots \mid e_{i_n}^n (A_n))\,;$$

*$1 < i_1 \leq k_t$,*
*then*

$$E_{i_t}^t \, E_{j_t}^t = 0$$

*for t = 1, 2, …, n;* $j_t \neq i_t$; $(E_{i_t}^t)^2 = E_{i_t}^t$, $E_{i_t}^{t*} = E_{i_t}^t$ *and*

$$I = (I_1 \mid \ldots \mid I_n) = (E_1^1 + \ldots + E_{k_1}^1 \mid \ldots \mid E_1^n + \ldots + E_{k_n}^n)\,.$$

*If $f = (f_1 \mid \ldots \mid f_n)$ is a complex valued super function whose super domain includes* $(c_1^1 \ldots c_{k_1}^1 \mid \ldots \mid c_1^n \ldots c_{k_n}^n)$ *then*

$$f(A) = (f_1(A_1) \mid \ldots \mid f_n(A_n)) = f(c_1)E_1 + \ldots + f(c_k)\, E_k$$
$$= ( f_1(c_1^1)E_1^1 + \ldots + f_1(c_{k_1}^1)\, E_{k_1}^1 \mid \ldots \mid f_n(c_1^n)\, E_1^n + \ldots + f_n(c_{k_n}^n)\, E_{k_n}^n ) .$$

*In particular*

$$A = c_1E_1 + \ldots + c_kE_k \text{ i.e., } (A_1 \mid \ldots \mid A_n)$$
$$= (c_1^1E_1^1 + \ldots + c_{k_1}^1 E_{k_1}^1 \mid \ldots \mid c_1^n E_1^n + \ldots + c_{k_n}^n E_{k_n}^n ) .$$

We just recall that an operator $T_s = (T_1 \mid \ldots \mid T_n)$ on an inner product superspace V is super nonnegative if $T_s$ is self adjoint and $(T_s\, \alpha \mid \alpha) \geq (0 \mid \ldots \mid 0)$ i.e.

$$((T_1\alpha_1 \mid \alpha_1) \mid \ldots \mid (T_n\alpha_n \mid \alpha_n) \geq (0 \mid \ldots \mid 0)$$

for every $\alpha = (\alpha_1 \mid \ldots \mid \alpha_n)$ in $V = (V_1 \mid \ldots \mid V_n)$.

We just give a theorem for the reader to prove.

**THEOREM 2.1.25:** *Let $T_s = (T_1 \mid \ldots \mid T_n)$ be a superdiagonalizable normal operator on a finite $(n_1, \ldots, n_n)$ dimensional super inner product super vector space $V = (V_1 \mid \ldots \mid V_n)$. Then $T_s = (T_1 \mid \ldots \mid T_n)$ is self adjoint super non negative or unitary according as each super characteristic value of $T_s$ is real super non negative or of absolute value (1, 1, …, 1).*

*Proof:* Suppose $T_s = (T_1 \mid \ldots \mid T_n)$ has super spectral resolution,

$$T_s = (T_1 \mid \ldots \mid T_n) = (c_1^1E_1^1 + \ldots + c_{k_1}^1 E_{k_1}^1 \mid \ldots \mid c_1^n E_1^n + \ldots + c_{k_n}^n E_{k_n}^n )$$

then

$$T_s^* = (T_1^* \mid \ldots \mid T_n^*)$$
$$= (\bar{c}_1^1E_1^1 + \ldots + \bar{c}_{k_1}^1 E_{k_1}^1 \mid \ldots \mid \bar{c}_1^n E_1^n + \ldots + \bar{c}_{k_n}^n E_{k_n}^n ) .$$

To say $T_s = (T_1 \mid \ldots \mid T_n)$ is super self adjoint is to say $T_s = T_s^*$ or

$$= ((c_1^1 - \overline{c}_1^1)E_1^1 + \ldots + (c_{k_1}^1 - \overline{c}_{k_1}^1)E_{k_1}^1 \big| \ldots \big| (c_1^n - \overline{c}_1^n)E_1^n + \ldots + (c_{k_n}^n - \overline{c}_{k_n}^n)E_{k_n}^n)$$
$= (0 \mid \ldots \mid 0)$.

Using the fact $E_{i_1}^t E_{j_t}^t = 0$, if $i_t \neq j_t$; t = 1, 2, …, n; and the fact that no $E_{j_t}^t$ is a zero operator, we see that $T_s$ is super self adjoint if and if only $c_{j_t}^t = \overline{c}_{j_t}^t$; t = 1, 2, …, n; To distinguish the normal operators which are non negative let us look at

$$(T_s \alpha \mid \alpha) = ((T_1\alpha_1 \mid \alpha_1) \mid \ldots \mid (T_n\alpha_n \mid \alpha_n))$$

$$= \left( \left( \sum_{j_1=1}^{k_1} c_{j_1}^1 E_{j_1}^1 \alpha_1 \middle| \sum_{i_1=1}^{k_1} E_{i_1}^1 \alpha_1 \right) \middle| \ldots \middle| \left( \sum_{j_n=1}^{k_n} c_{j_n}^n E_{j_n}^n \alpha_n \middle| \sum_{j_n=1}^{k_n} E_{j_n}^n \alpha_n \right) \right)$$

$$= \left( \sum_{i_1} \sum_{j_1} c_{j_1}^1 (E_{j_1}^1 \alpha_1 \mid E_{i_1}^1 \alpha_1) \mid \ldots \mid \sum_{i_n} \sum_{j_n} c_{j_n}^n (E_{j_n}^n \alpha_n \mid E_{i_n}^n \alpha_n) \right)$$

$$= \left( \sum_{j_1} c_{j_1}^1 \left\| E_{j_1}^1 \alpha_1 \right\|^2 \mid \ldots \mid \sum_{j_n} c_{j_n}^n \left\| E_{j_n}^n \alpha_n \right\|^2 \right).$$

We have made use of the simple fact

$$(E_{j_t}^t \alpha_t \mid E_{i_t}^t \alpha_t) = 0 \text{ if } i_t \neq j_t;\ 1 \leq i_t, j_t \leq k_t$$

and $t = 1, 2, \ldots, n.$ From this it is clear that the condition

$$(T_s\alpha \mid \alpha) = ((T_1\alpha_1 \mid \alpha_1) \mid \ldots \mid (T_n\alpha_n \mid \alpha_n) \geq (0 \mid \ldots \mid 0)$$

is satisfied if and only if $c_{j_t}^t \geq 0$ for each $j_t$; $1 \leq j_t \leq k_t$ and t = 1, 2, …, n. To distinguish the unitary operators observe that

$$T_s T_s^* = (c_1^1 \overline{c}_1^1 E_1^1 + \ldots + c_{k_1}^1 \overline{c}_{k_1}^1 E_{k_1}^1 \big| \ldots \big| c_1^n \overline{c}_1^n E_1^n + \ldots + c_{k_n}^n \overline{c}_{k_n}^n E_{k_n}^n)$$

$$= (|c_1^1|^2 E_1^1 + \ldots + |c_{k_1}^1|^2 E_{k_1}^1 \mid \ldots \mid |c_1^n|^2 E_1^n + \ldots + |c_{k_n}^n|^2 E_{k_n}^n) .$$

If

$$T_s T_s^* = (I_1 \mid \ldots \mid I_n) = I = (T_1 T_1^* \mid \ldots \mid T_n T_n^*)$$

then

$$(I_1 \mid \ldots \mid I_n) = (|c_1^1|^2 E_1^1 + \ldots + |c_{k_1}^1| E_{k_1}^1 \mid \ldots \mid (|c_1^n|^2 E_1^n + \ldots + |c_{k_n}^n| E_{k_n}^n))$$

and operative with

$$E_{j_t}^t, E_{j_t}^t = |c_{j_t}^t|^2 E_{j_t}^t;$$

$1 \leq j_t \leq k_t$ and $t = 1, 2, \ldots, n$. Since $E_{j_t}^t \neq 0$ we have $|c_{j_t}^t|^2 = 1$ or $|c_{j_t}^t| = 1$. Conversely if $|c_{j_t}^t|^2 = 1$ for each $j_t$ it is clear that

$$T_s T_s^* = (I_1 \mid \ldots \mid I_n) = I = (T_1 T_1^* \mid \ldots \mid T_n T_n^*) .$$

If $T_s = (T_1 \mid \ldots \mid T_n)$ is a general linear operator on the supervector space $V = (V_1 \mid \ldots \mid V_n)$ which has real characteristic super values it does not follow that $T_s$ is super self adjoint. The theorem of course states that if $T_s$ has real characteristic super values and if $T_s$ is super diagonalizable and normal then $T_s$ is super self adjoint. We have yet another interesting theorem.

**THEOREM 2.1.26:** *Let $V = (V_1 \mid \ldots \mid V_n)$ be a finite $(n_1, \ldots, n_n)$ dimensional inner product super space and $T_s$ a super non negative operator on V. Then $T_s = (T_1 \mid \ldots \mid T_n)$ has a unique super non negative square root, that is; there is one and only one non negative super operator $N_s = (N_1 \mid \ldots \mid N_n)$ on V such that $N_s^2 = T_s$ i.e., $(N_1^2 \mid \ldots \mid N_n^2) = (T_1 \mid \ldots \mid T_n)$.*

*Proof:* Let $T_s = (T_1 \mid \ldots \mid T_n)$

$$= (c_1^1 E_1^1 + \ldots + c_{k_1}^1 E_{k_1}^1 \mid \ldots \mid c_1^n E_1^n + \ldots + c_{k_n}^n E_{k_n}^n)$$

be the super spectral resolution of $T_s$. By the earlier results each $c_{j_t}^t \geq 0$; $1 \leq j_t \leq k_t$ and $t = 1, 2, \ldots, n$. If $c^t$ is any non negative real number t = 1, 2, …, n let $\sqrt{c^t}$ denote the non negative square root of c. So if c = ($c_1$, …, $c_n$) then the super square root or square super root of c is equal to, $\sqrt{c} = (\sqrt{c_1}, \ldots, \sqrt{c_n})$. Then according to earlier result $N_s = \sqrt{T_s}$ is a well defined super diagonalizable normal operator on V i.e. $N_s$ = ($N_1$ | … | $N_n$) $\sqrt{T_s} = (\sqrt{T_1} \mid \ldots \mid \sqrt{T_n})$ is a well defined super diagonalizable normal operator on V = ($V_1$ | … | $V_n$). It is super non negative and by an obvious computation $N_s^2 = T_s$ i.e. $(N_1^2 \mid \ldots \mid N_n^2)$ = ($T_1$ | … | $T_n$).

Let $P_s$ = ($P_1$ | … | $P_n$) be a non negative operator V such that $P_s^2 = T_s$ i.e. $(P_1^2 \mid \ldots \mid P_n^2)$ = ($T_1$ | … | $T_n$); we shall prove that $P_s$ = $N_s$. Let

$$P_s = (d_1^1 F_1^1 + \ldots + d_{r_1}^1 F_{r_1}^1 \mid \ldots \mid d_1^n F_1^n + \ldots + d_{r_n}^n F_{r_n}^n)$$

be the super spectral resolution of $P_s$ = ($P_1$ | … | $P_n$). Then $d_{j_t}^t \geq 0$ for 1 ≤ $j_t$ ≤ $k_t$; t = 1, 2, …, n each $j_t$ since $P_s$ is non negative.

From $P_s^2 = T_s$ we have $T_s$ = ($T_1$ | … | $T_n$)

$$= (d_1^{1^2} F_1^1 + \ldots + d_{r_1}^{1^2} F_{r_1}^1 \mid \ldots \mid d_1^{n^2} F_1^n + \ldots + d_{r_n}^{n^2} F_{r_n}^n).$$

Now $(F_1^1 \ldots F_{r_1}^1, \ldots, F_1^n \ldots F_{r_n}^n)$ satisfy the condition

$$(I_1 \mid \ldots \mid I_n) = (F_1^1 + \ldots + F_{r_1}^1 \mid \ldots \mid F_1^n + \ldots + F_{r_n}^n)$$

$F_{i_t}^t F_{j_t}^t = 0$; 1 ≤ t ≤ $r_t$; t = 1, 2, …, n for $i_t \neq j_t$ and no $F_{j_t}^t = 0$. The numbers $d_1^{1^2} \ldots d_{r_1}^{1^2}, \ldots, d_1^{n^2} \ldots d_{r_n}^{n^2}$ are distinct because distinct non negative numbers have distinct squares. By the uniquiness of the super spectral resolution of $T_s$ we must have $r_t = k_t$; t = 1, 2, …, n. $F_{j_t}^t = E_{j_t}^t, \left(d_{j_t}^t\right)^2 = c_{j_t}^t$; t = 1, 2, …, n. Thus $P_s$ = $N_s$.

**THEOREM 2.1.27:** *Let V = ($V_1$ | … | $V_n$) be a finite ($n_1$, …, $n_n$) dimensional super inner product supervector space and let $T_s$ = ($T_1$ | … | $T_n$) be any linear operator on V. Then there exists a*

*unitary operator $U_s = (U_1 \mid \ldots \mid U_n)$ on V and a super non negative operator $N_s = (N_1 \mid \ldots \mid N_n)$ on V such that*

$$T_s = U_sN_s = (T_1 \mid \ldots \mid T_n) = (U_1N_1 \mid \ldots \mid U_nN_n).$$

*The non-negative operators $N_s$ is unique. If $T_s = (T_1 \mid \ldots \mid T_n)$ is invertible, the operator $U_s$ is also unique.*

*Proof:* Suppose we have $T_s = U_sN_s$ where $U_s$ is unitary and $N_s$ is super non negative. Then $T_s^* = (U_sN_s)^* = N_s^* U_s^* = N_sU_s^*$. Thus $T_s^*T_s = N_sU_s^*U_s N_s = N_s^2$. This shows that $N_s$ is uniquely determined as the super non negative square root of $T_s^*T_s$. If $T_s$ is invertible then so is $N_s$ because.

$$(N_s \alpha \mid N_s\alpha) = (N_s^2\alpha|\alpha)$$

i.e. $((N_1\alpha_1 \mid N_1\alpha_1) \mid \ldots \mid (N_n\alpha_n \mid N_n\alpha_n))$

$= \quad (N_1^2\alpha_1 \mid \alpha_1) \mid \ldots \mid (N_n^2 \alpha \mid \alpha_n)$

$= \quad (T_s^*T_s\alpha \mid \alpha) = (T_s^*\alpha \mid T_s\alpha)$

$= \quad ((T_1^*T_1\alpha_1 \mid \alpha_1) \mid \ldots \mid (T_n^*T_n\alpha_n \mid \alpha_n))$

$= \quad ((T_1 \alpha_1|T_1\alpha_1) \mid \ldots \mid (T_n\alpha_n \mid T_n\alpha_n))$ .

In this case we define $U_s = T_sN_s^{-1}$ and prove that $U_s$ is unitary. Now

$$U_s^* = (T_sN_s^{-1})^* = (N_s^{-1})^*T_s^* = (N_s^*)^{-1}T_s^* = N_s^{-1}T_s^*.$$

Thus

$U_sU_s^* = T_sN_s^{-1} N_s^{-1} T_s^*$

$= \quad T_s(N_s^{-1})^2T_s^*$

$= \quad T_s(N_s^2)^{-1} T_s^*$

$= \quad T_s(T_s^*T_s)^{-1} T_*$

$= \quad T_sT_s^{-1}(T_s^*)^{-1} T_s^*$

$= \quad (T_1T_1^{-1}(T_1^*)^{-1} T_1^* \mid \ldots \mid T_nT_n^{-1}(T_n^*)^{-1}T_n^*)$

$= \quad (I_1 \mid \ldots \mid I_n)$,

so $U_s = (U_1 \mid \ldots \mid U_n)$ is unitary.

If $T_s = (T_1 \mid \ldots \mid T_n)$ is not invertible, we shall have to do a bit more work to define $U_s = (U_1 \mid \ldots \mid U_n)$ we first define $U_s$ on the range of $N_s$. Let $\alpha = (\alpha_1 \mid \ldots \mid \alpha_n)$ be a supervector in the superrange of $N_s$ and $\alpha = N\beta$; $(\alpha_1 \mid \ldots \mid \alpha_n) = N_s\beta = (N_1\beta_1 \mid \ldots \mid N_n\beta_n)$.

We define

$$U_s\alpha = T_s\beta \text{ i.e., } (U_1\alpha_1 \mid \ldots \mid U_n\alpha_n) = (T_1\beta_1 \mid \ldots \mid T_n\beta_n),$$

motivated by the fact that we want

$$U_sN_s\beta = T_s\beta.\ (U_1N_1\beta_1 \mid \ldots \mid U_nN_n\beta_n) = (T_1\beta_1 \mid \ldots \mid T_n\beta_n).$$

We must verify that $U_s$ is well defined on the super range of $N_s$; in other words if

$$N_s\beta' = N_s\beta \text{ i.e. } (N_1\beta_1' \mid \ldots \mid N_n\beta_n') = (N_1\beta_1 \mid \ldots \mid N_n\beta_n)$$

then

$$T_s\beta' = T_s\beta;\ (T_1\beta_1' \mid \ldots \mid T_n\beta_n') = (T_1\beta_1 \mid \ldots \mid T_n\beta_n).$$

We have verified above that

$$\|N_s\gamma\|^2 = \left(\|N_1\gamma_1\|^2 \mid \ldots \mid \|N_n\gamma_n\|^2\right)$$

$$= \|T_s\gamma\|^2 = \left(\|T_1\gamma_1\|^2 \mid \ldots \mid \|T_n\gamma_n\|^2\right)$$

for every $\gamma = (\gamma_1 \mid \ldots \mid \gamma_n)$ in V. Thus with $\gamma = \beta - \beta'$ i.e.

$$(\gamma_1 \mid \ldots \mid \gamma_n) = (\beta_1 - \beta_1' \mid \ldots \mid \beta_n - \beta_n'),$$

we see that

$$N_s(\beta - \beta') = (N_1(\beta_1 - \beta_1') \mid \ldots \mid N_n(\beta_n - \beta_n')) = (0 \mid \ldots \mid 0)$$

if and only if

$$T_s(\beta - \beta') = (T_1(\beta_1 - \beta_1') \mid \ldots \mid T_n(\beta_n - \beta_n')) = (0 \mid \ldots \mid 0).$$

So $U_s$ is well defined on the super range of $N_s$ and is clearly linear where defined. Now if $W = (W_1 \mid \ldots \mid W_n)$ is the super range of $N_s$ we are going to define $U_s$ on $W^\perp = (W_1^\perp \mid \ldots \mid W_n^\perp)$.

To do this we need the following observation. Since $T_s$ and $N_s$ have the same super null space their super ranges have the same super dimension. Thus $W^{\perp} = (W_1^{\perp} \mid \ldots \mid W_n^{\perp})$ has the same super dimension as the super orthogonal complement of the super range of $T_s$. Therefore there exists super isomorphism

$$U_s^0 = (U_1^0 \mid \ldots \mid U_n^0) \text{ of } W^{\perp} = (W_1^{\perp} \mid \ldots \mid W_n^{\perp})$$

onto

$$T_s(V)^{\perp} = (T_1(V_1)^{\perp} \mid \ldots \mid T_n(V_n)^{\perp}).$$

Now we have defined $U_s$ on W and we define $U_s W^{\perp}$ to be $U_s^0$.

Let us repeat the definition of $U_s$ since

$$V = W \oplus W^{\perp}$$

i.e.,

$$(V_1 \mid \ldots \mid V_n) = (W_1 \oplus W_1^{\perp} \mid \ldots \mid W_n \oplus W_n^{\perp})$$

each $\alpha = (\alpha_1 \mid \ldots \mid \alpha_n)$ in V is uniquely expressible in the form

$$\alpha = N_s\beta + \gamma \text{ i.e.,}$$

$$\alpha = (\alpha_1 \mid \ldots \mid \alpha_n) = (N_1\beta_1 + \gamma_1 \mid \ldots \mid N_n\beta_n + \gamma_n)$$

where $N_s\beta$ is in the range of $W = (W_1 \mid \ldots \mid W_n)$ of $N_s$ and $\gamma = (\gamma_1 \mid \ldots \mid \gamma_n)$ is in $W^{\perp}$.

We define

$$U_s\alpha = T_s\beta + U_s^0\gamma$$

$$(U_1\alpha_1 \mid \ldots \mid U_n\alpha_n) = (T_1\beta_1 + U_1^0\gamma_1 \mid \ldots \mid T_n\beta_n + U_n^0\gamma_n).$$

This $U_s$ is clearly linear and we have verified it is well defined. Also

$$\begin{aligned}(U_s\alpha \mid U_s\alpha) &= ((U_1\alpha_1 \mid U_1\alpha_1) \mid \ldots \mid (U_n\alpha_n \mid U_n\alpha_n)) \\ &= (T_s\beta + U_s^0\gamma \mid T_s\beta + U_s^0\gamma) \\ &= (T_s\beta \mid T_s\beta) + (U_s^0\gamma \mid U_s^0\gamma) \\ &= ((T_1\beta_1 \mid T_1\beta_1) + (U_1^0\gamma_1 \mid U_1^0\gamma_1) \mid \ldots \mid \\ &\qquad (T_n\beta_n \mid T_n\beta_n) + (U_n^0\gamma_n \mid U_n^0\gamma_n))\end{aligned}$$

$= \quad ((N_1\beta_1 | N_1\beta_1) + (\gamma_1 | \gamma_1) | \ldots | (N_n\beta_n | N_n\beta_n) + (\gamma_n | \gamma_n))$
$= \quad (N_s\beta | N_s\beta) + (\gamma | \gamma) = (\alpha | \alpha)$

and so $U_s$ is unitary. We have $U_sN_s\beta = T_s\beta$ for each $\beta$. Hence the claim.

We call $T_s = U_sN_s$ as in case of usual vector spaces to be the polar super decomposition for $T_s$.

i.e. $T_s = U_s N_s$
i.e. $(T_s | \ldots | T_n) = (U_1N_1 | \ldots | U_n N_n)$.

Now we proceed onto define the notion of super root of the family of operators on an inner product super vector space $V = (V_1 | \ldots | V_n)$.

**DEFINITION 2.1.10:** *Let $F_s$ be a family of operators on an inner product super vector space $V = (V_1 | \ldots | V_n)$. A super function $r = (r_1 | \ldots | r_n)$ on $F_s$ with values in the field F of scalars will be called a super root of $F_s$ if there is a non zero super vector $\alpha = (\alpha_1 | \ldots | \alpha_n)$ in V such that $T_s\alpha = r(T_s)\alpha$ i.e., $(T_1\alpha_1 | \ldots | T_n\alpha_n) = (r_1(T_1)\alpha_1 | \ldots | r_1(T_n)\alpha_n)$ for all $Ts = (T_1 | \ldots | T_n)$ in $F_s$.*

*For any super function $r = (r_1 | \ldots | r_n)$ from $F_s$ to $(F | \ldots | F)$, let $V(r) = (V_1(r_1) | \ldots | V_n(r_n))$ be the set of all $\alpha = (\alpha_1 | \ldots | \alpha_n)$ in V such that $T_s(\alpha) = r(T)\alpha$ for every $T_s$ in $F_s$. Then V(r) is a super subspace of V and $r = (r_1 | \ldots | r_n)$ is a super root of $F_s$ if and only if $V(r) = (V_1(r_1) | \ldots | V_n(r_n)) \neq (\{0\} | \ldots | \{0\})$. Each non zero $\alpha = (\alpha_1 | \ldots | \alpha_n)$ in V(r) is simultaneously a characteristic super vector for every $T_s$ in $F_s$.*

In view of this definition we have the following interesting theorem.

**THEOREM 2.1.28:** *Let $F_s$ be a commuting family of super diagonalizable normal operators on a finite dimensional super inner product space $V = (V_1 | \ldots | V_n)$. Then $F_s$ has only a finite number of super roots. If $r_1^1 \ldots r_{k_1}^1, \ldots, r_1^n \ldots r_{k_n}^n$ are the distinct super roots of $F_s$ then*

i) $V(r_i) = (V_1(r^1_{i_1}) \mid \ldots \mid V_n(r^n_{i_n}))$ is orthogonal to $V(r_j) = (V_1(r^1_{j_1}) \mid \ldots \mid V_n(r^n_{j_n}))$ when i ≠ j i.e., $i_t \neq j_t$; t = 1, 2, …, n and

ii) *V = V(r₁) ⊕ ... ⊕ V(r_k)*
*i.e.,* $V = (V_1(r^1_1) \oplus \ldots \oplus V_1(r^1_{k_1}) \mid \ldots \mid V_n(r^n_1) \oplus \ldots \oplus V_n(r^n_{k_n}))$.

*Proof:* Suppose r = ($r_1$ | … | $r_k$) and s = ($s_1$ | …| $s_k$) distinct super roots of F. Then there is an operator $T_s$ in $F_s$ such that $r(T_s) \neq s(T_s)$; i.e., ($r_1(T_1)$ | … | $r_1(T_n)$) ≠ ($s_1(T_1)$ | … | $s_1(T_n)$) since characteristic super vectors belonging to distinct characteristic super values of $T_s$ are necessarily superorthogonal, it follows that

$$V(r) = (V_1(r^1_1) \mid \ldots \mid V_n(r^n_n))$$

is orthogonal to

$$V(s) = (V_1(s^1_1) \mid \ldots \mid V_n(s^n_n)).$$

Because V is finite ($n_1$, …, $n_n$) dimensional this implies $F_s$ has atmost a finite number of super roots. Let $r^1_1 \ldots r^1_{k_1}, \ldots, r^n_1 \ldots r^n_{k_n}$ be the super roots of F. Suppose $\{T^1_1 \ldots T^1_{m_1} \mid \ldots \mid T^n_1 \ldots T^n_{m_n}\}$ be a maximal linearly independent super subset of $F_s$ and let

$$E^1_{i_1 1}, E^1_{i_1 2} \ldots E^2_{i_2 1}, E^2_{i_2 2}, \ldots, E^n_{i_n 1}, E^n_{i_n 2}, \ldots$$

be the resolution of identity defined by $T^p_{i_p}$; ($1 \leq i_p \leq m_p$); p = 1, 2, …, n; then the super projections $E_{ij} = (E^1_{i_1 j_1} \mid \ldots \mid E^n_{i_n j_n})$ form a commutative super family, for each $E_{ij}$ hence each $E^t_{i_t j_t}$ is a super polynomial in $T^t_{i_t}$ and $T^t_1, \ldots, T^t_{m_t}$ commute with one another. This being true for each p = 1, 2, …, n.
Since

$$I = \left(\left(\sum_{j^1_1} E^1_{i j^1_1}\right)\left(\sum_{j^1_2} E^1_{2 j^1_2}\right)\ldots\left(\sum_{j^1_m} E^1_{m_1 j^1_{m_1}}\right) \Big| \ldots \Big| \left(\sum_{j^n_1} E^n_{i j^n_1}\right)\left(\sum_{j^n_2} E^n_{2 j^n_2}\right)\ldots\right.$$

$$\left(\sum_{j_m^n} E^n_{m_n j^n_{m_n}}\right)\Bigg) = (I_1 \mid \ldots \mid I_n)$$

each super vector $\alpha = (\alpha_1 \mid \ldots \mid \alpha_n)$ in $V = (V_1 \mid \ldots \mid V_n)$ may be written in this form

$$\begin{aligned}\alpha &= (\alpha_1 \mid \ldots \mid \alpha_n)\\ &= \left(\sum_{j_1^1 j_2^1 \ldots j_m^1} E^1_{1 j_1^1} \ldots E^1_{m_1 j^1_{m_1}} \alpha_1 \Big| \ldots \Big| \sum_{j_1^n \ldots j_m^n} E^n_{1 j_1^n} \ldots E^n_{m_n j^n_{m_n}} \alpha_n\right) \qquad \text{(A)}\end{aligned}$$

Suppose $j_1^1 \ldots j^1_{m_1}, \ldots, j_1^n \ldots j^n_{m_n}$ are indices for which

$$\beta = \left(E^1_{1 j_1^1} E^1_{2 j_2^1} \ldots E^1_{m_1 j^1_{m_1}} \alpha_1 \Big| \ldots \Big| E^2_{1 j_1^2} \ldots E^2_{m_2 j^2_{m_2}} \alpha_2 \Big| \ldots \Big| E^n_{1 j_1^n} \ldots E^n_{m_n j^n_{m_n}} \alpha_n\right)$$

$\neq (0 \mid \ldots \mid 0)$,

$$\beta_i = (\beta^1_{i_1} \ldots \beta^1_{1_n}) = \left(\prod_{n_1 \neq i_1} E^1_{n_1 j^1_{n_1}} \alpha_1 \mid \ldots \mid \prod_{n_n \neq i_n} E^n_{n_n j^n_{n_n}} \alpha_n\right).$$

then

$$\beta = \left(E^1_{i_1 j_1} \beta^1_{i_1} \mid \ldots \mid E^n_{i_n j_{in}} \beta^n_{i_n}\right).$$

Hence there is a scalar $c_i$ such that

$$(T_{i_1} \beta^1_{i_1} \mid \ldots \mid T_{i_n} \beta^n_{i_n}) = (c^1_{i_1} \beta^1_{i_1} \mid \ldots \mid c^n_{i_n} \beta^n_{i_n})$$

where $1 \leq i_t \leq m_t$ and $t = 1, 2, \ldots, n$. For each $T_s$ in $F_s$ there exists unique scalars $b^1_{i_1}, \ldots, b^n_{i_n}$ such that

$$T_s = \left(\sum_{i_1=1}^{m_1} b^1_{i_1} T^1_{i_1} \mid \ldots \mid \sum_{i_n=1}^{m_n} b^n_{i_n} T^n_{i_n}\right).$$

Thus

$$T_s\beta = \left(\sum_{i_1} b^1_{i_1} T^1_{i_1} \beta_1 \mid \ldots \mid \sum_{i_n} b^n_{1_n} T^n_{i_n} \beta_n\right)$$

$$= \left(T_1\beta_1 \mid \ldots \mid T_n\beta_n\right) = \left( \sum_{i_1}^{m_1} b^1_{i_1} c^1_{i_1} \beta^1_{i_1} \mid \ldots \mid \sum_{i_n}^{m_n} b^n_{i_n} c^n_{i_n} \beta^n_{i_n} \right).$$

The function

$$T_s = (T_1 \mid \ldots \mid T_n) \rightarrow \left( \sum_{i_1} b^1_{i_1} c^1_{i_1} \mid \ldots \mid \sum_{i_n} b^n_{i_n} c^n_{i_n} \right)$$

is evidently one of the super roots say $r_t = (r^1_{t_t} \mid \ldots \mid r^n_{t_n})$ of $F_s$ and $\beta = (\beta_1 \mid \ldots \mid \beta_n)$ lies in $V(r_t) = (V_1(r^1_{t_1}) \mid \ldots \mid V_n(r^n_{t_n}))$. Therefore each nonzero term in equation (A) belongs to one of the spaces

$$V(r_1) = (V_1(r^1_{1_1}) \mid \ldots \mid V_n(r^n_{1_n})), \ldots, V(r_k) = (V_1(r^1_{k_1}) \mid \ldots \mid V_n(r^n_{k_n})).$$

It follows that $V = (V_1 \mid \ldots \mid V_n)$ is super orthogonal direct sum of $(V(r_1), \ldots, V(r_k))$.

The following corollary is direct and is left as an exercise for the reader to prove.

**COROLLARY 2.1.9:** *Under the assumptions of the theorem, let* $P_j = (P^1_{j_1} \mid \ldots \mid P^n_{j_n})$ *be the super orthogonal projection of* $V = (V_1 \mid \ldots \mid V_n)$ *on* $V(r_j) = (V_1(r^1_{j_1}) \mid \ldots \mid V_n(r^n_{j_n}))$; $1 \le j_t \le k_t$; $t = 1, 2, \ldots, n$. *Then* $P^t_{i_1} P^t_{j_t} = 0$ *when* $i_t \ne j_t$; $t = 1, 2, \ldots, n$.

$$\begin{aligned} I &= (I_1 \mid \ldots \mid I_n) \\ &= (P^1_1 + \ldots + P^1_{k_1} \mid \ldots \mid P^n_1 + \ldots + P^n_{k_n}) \end{aligned}$$

*and every* $T_s$ *in* $F_s$ *may be written in the form*

$$\begin{aligned} T_s &= (T_1 \mid \ldots \mid T_n) \\ &= \left( \sum_{j_1} r^1_{j_1}(T_1) P^1_{j_1} \mid \ldots \mid \sum_{j_n} r^n_{j_n}(T_n) P^n_{j_n} \right). \end{aligned}$$

*The super family of super orthogonal projections $\{P_1^1 \ldots P_{k_1}^1 \mid \ldots \mid P_1^n \ldots P_{k_n}^n\}$ is called the super resolution of the identity determined by $F_s$ and*

$$T_s \;=\; (T_1 \mid \ldots \mid T_n) = \left( \sum_{j_1} r_{j_1}^1 (T_1) P_{j_1}^1 \mid \ldots \mid \sum_{j_n} r_{j_n}^n (T_n) P_{j_n}^n \right)$$

*is the super spectral resolution of $T_s$ in terms of this family of spectral super resolution of $T_s$ (both mean one and the same).*

Although the super projections ($P_1^1 \ldots P_{k_1}^1 \mid \ldots \mid P_1^n \ldots P_{k_n}^n$) in the preceeding corollary are canonically associated with the family $F_s$ they are generally not in $F_s$ nor even linear combinations of operators in $F_s$; however we shall show that they may be obtained by forming certain products of super polynomials in elements of $F_s$.

Thus as in case of usual vector spaces we can say in case of super vector spaces $V = (V_1 \mid \ldots \mid V_n)$ which are inner product super vector spaces the notion of super self adjoint super algebra of operators which is a linear super subalgebra of SL(V, V) which contains the super adjoint of each of its members.

If $F_s$ is the family of linear operators on a finite dimensional inner product super space, the self super adjoint super algebra generated by $F_s$ is the smallest self adjoint super algebra which contains $F_s$.

Now we proceed onto prove an interesting theorem.

***THEOREM 2.1.29**: Let $F_s$ be a commuting family of super diagonalizable normal operators on a finite dimensional inner product super vector space $V = (V_1 \mid \ldots \mid V_n)$ and let $a_s$ be the self adjoint super algebra generated by $F_s$ and the identity operator. Let $\{P_1^1 \ldots P_{k_1}^1 \mid \ldots \mid P_1^n \ldots P_{k_n}^n\}$ be the super resolution of the super identity defined by $F_s$. Then $a_s$ is the set of all operators on $V = (V_1 \mid \ldots \mid V_n)$ of the form*

$$T_s = \left( \sum_{j_1=1}^{k_1} c^1_{j_1} P^1_{j_1} \mid \ldots \mid \sum_{j_n=1}^{k_n} c^n_{j_n} P^n_{j_n} \right) = T = (T_1 \mid \ldots \mid T_n) \qquad I$$

*where ( $(c^1_1\, c^1_2 \ldots c^1_{k_1} \mid \ldots \mid c^n_1\, c^n_2 \ldots c^n_{k_n})$ are arbitrary scalars.*

*Proof:* Let $C_s$ denote the set of all super operators on V of the form given in I of the theorem.
Then $C_s$ contains the super identity operator and the adjoint

$$T^*_s = \left( \sum_{j_1=1}^{k_1} \overline{c}^1_{j_1}\ P^1_{j_1} \mid \ldots \mid \sum_{j_n=1}^{k_n} \overline{c}^n_{j_n}\ P^n_{j_n} \right) = (T^*_1 \mid \ldots \mid T^*_n),$$

of each of its members. If

$$T_s = \left( \sum_{j_1} c^1_{j_1} P^1_{j_1} \mid \ldots \mid \sum_{j_n} c^n_{j_n} P^n_{j_n} \right) = (T_1 \mid \ldots \mid T_n)$$

and

$$U_s = (U_1 \mid \ldots \mid U_n) = \left( \sum_{j_1} d^1_{j_1} P^1_{j_1} \mid \ldots \mid \sum_{j_n} d^n_{j_n} P^n_{j_n} \right)$$

then for every scalar

$$a = (a^1_1 \ldots a^1_n);\ aT_s + U_s = (a^1_1 T_1 + U_1 \mid \ldots \mid a^1_n T_n + U_n)$$

$$= \left( \sum_{j_1} (a^1_1 c^1_1 + d^1_{j_1})\ P^1_{j_1} \mid \ldots \mid \sum_{j_n} (a^n_n c^n_n + d^n_{j_n})\ P^n_{j_n} \right)$$

and

$$T_s U_s = \left( \sum_{i_1 j_1} c^1_{i_1} d^1_{j_1} P^1_{i_1}\ P^1_{j_1} \mid \ldots \mid \sum_{i_n j_n} c^n_{i_n} d^n_{j_n} P^n_{i_n}\ P^n_{j_n} \right)$$

$$= \left( \sum_{j_1} c^1_{j_1} d^1_{j_1} P^1_{j_1} \mid \ldots \mid \sum_{j_n} c^n_{j_n} d^n_{j_n} P^n_{j_n} \right) = U_s T_s .$$

Thus $C_s$ is a self super adjoint commutative super algebra containing $F_s$ and the super identity operator. Thus $C_s$ contains

$a_s$. Now let $r_1^1 \ldots r_{k_1}^1, \ldots, r_1^n \ldots r_{k_n}^n$ be the super roots of $F_s$. Then for each pair of indices ($i_t$, $n_t$), $i_t \neq n_t$, there is an operator $T_{si_tn_t}$ in $F_s$ such that $r_{i_t}^t\ (T_{si_tn_t}) \neq r_{n_t}^t\ (T_{si_tn_t})$. Let

$$a_{i_tn_t}^t = r_{i_t}^t (T_{si_tn_t}) - r_{n_t}^t (T_{si_tn_t})$$

and

$$b_{i_tn_t}^t = r_{n_t}^t\ (T_{si_tn_t})\,.$$

Then the linear operator

$$Q_{si_t} = \prod_{n_t \neq i_t} a_{si_tn_t}^{-1}\ (T_{si_tn_t} - b_{i_tn_t}^t I_t)$$

is an element of the super algebra $a_s$. We will show that $Q_{si_t} = P_{si_t}$ ($1 \leq i_t \leq k_t$). For this suppose $j_t \neq i_t$ and α is an arbitrary super vector in $V(r_j) = (V_1(r_{j_1}^1) \mid \ldots \mid V_n(r_{j_n}^n))$. Then

$$T_{si_tj_t}\alpha = r_j^t (T_{si_tj_t})\alpha = b_{i_tj_t}^t \alpha$$

so that

$$(T_{si_tj_t} - b_{i_tj_t}^t I_t)\ \alpha = (0 \mid \ldots \mid 0).$$

Since the factors in $Q_{si_t}$ all commute it follows that $Q_{si_t}\ \alpha = (0 \mid \ldots \mid 0)$. i.e. $Q_{si_t}\alpha_{i_t} = 0$. Hence $Q_{si_t}$ agrees with $P_{si_t}$ on $V(r_j) = (V_1(r_{j_1}^1) \mid \ldots \mid V_n(r_{j_n}^n))$ whenever $j_t \neq i_t$. Now suppose α is a super vector in $V(r_i)$. Then $T_{ti_tn_t}\alpha_t = r_{i_t}^t(T_{ti_tn_t})\alpha_t$ and $a_{i_tn_t}^{-1}(T_{ti_tn_t} - b_{i_tn_t} I_t)\alpha_t = a_{i_tn_t}^{-1}\ [r_{i_t}(T_{ti_tn_t}) - r_{n_t}(T_{ti_tn_t})]\alpha_t = \alpha_t$. Thus $Q_{ti_t}\alpha_t = \alpha_t$ for t = 1, 2, …, n and $Q_{ti_t}$ agrees with $P_{i_t}$ on $V_t(r_{i_t})$ therefore $Q_{ti_t} = P_{i_t}$ for i = 1, 2, …, $k_t$; t = 1, 2, …, n From which it follows $a_s = c_s$.

The following corollary is left as an exercise for the reader to prove.

**COROLLARY 2.1.10:** Under the assumptions of the above theorem there is an operator $T_s = (T_1 \mid \ldots \mid T_n)$ in $a_s$ such that every member of $a_s$ is a super polynomial in $T_s$.

We now state an interesting theorem on super vector spaces. The proof is left as an exercise for the reader.

**THEOREM 2.1.30:** *Let $T_s = (T_1 \mid \ldots \mid T_n)$ be a normal operator on a finite dimensional super inner product space $V = (V_1 \mid \ldots \mid V_n)$. Let $p = (p_1 \mid \ldots \mid p_n)$ be the minimal polynomial for $T_s$ with $(p_1^1 \ldots p_{k_1}^1, \ldots, p_1^n \ldots p_{k_n}^n)$ its distinct monic prime factors. Then each $p_{j_t}^t$ occurs with multiplicity 1 in the super factorization of p for $1 \leq j_t \leq k_t$ and $t = 1, 2, \ldots, n$ and has super degree 1 or 2. Suppose $W_j = (W_{j_1}^1 \mid \ldots \mid W_{j_n}^n)$ is the null superspace of $p_{j_t}^t(T_t)$, $t = 1, 2, \ldots, n$; $1 \leq j_t \leq k_t$. Then*

i. *$W_j$ is super orthogonal to $W_i$ when $i \neq j$ i.e. $(W_{j_1}^1 \mid \ldots \mid W_{j_n}^n)$ is super orthogonal to $(W_{i_1}^1 \mid \ldots \mid W_{i_n}^n)$; i.e. $W_{j_t}^t$ is orthogonal to $W_{i_t}^t$, $1 \leq i_t, j_t \leq k_t$ and $t = 1, 2, \ldots, n$.*
ii. *$V = (V_1 \mid \ldots \mid V_n) = (W_1^1 \oplus \ldots \oplus W_{k_1}^1 \mid \ldots \mid W_1^n \oplus \ldots \oplus W_{k_n}^n)$*
iii. *$W_j = (W_{j_1}^1 \mid \ldots \mid W_{j_n}^n)$ is super covariant under $T_s$ and $p_j = (p_{j_1}^1 \mid \ldots \mid p_{j_n}^n)$ is the minimal super polynomial for the restriction of $T_s$ to $W_j$.*
iv. *For every $j = (j_1, \ldots, j_n)$ there is a super polynomial $e_j = (e_{j_1}^1, \ldots, e_{j_n}^n)$ with coefficients in the scalar field such that $e_j(T_s) = (e_{j_1}^1(T_1) \mid \ldots \mid e_{j_n}^n(T_n))$ is super orthogonal projection of $V = (V_1 \mid \ldots \mid V_n)$ on $W_j = (W_{j_1}^1 \mid \ldots \mid W_{j_n}^n)$.*

We now prove the following lemma.

**LEMMA 2.1.2:** *Let $N_s = (N_1 \mid \ldots \mid N_n)$ be a normal operator on a super inner product space $W = (W_1 \mid \ldots \mid W_n)$. Then the super null space of $N_s$ is the super orthogonal complement of its super range.*

*Proof:* Suppose

$$(\alpha \mid N_s\beta) = ((\alpha_1 \mid N_1\beta_1) \mid \ldots \mid (\alpha_n \mid N_n\beta_n)) = (0 \mid \ldots \mid 0)$$

for all $\beta = (\beta_1 \mid \ldots \mid \beta_n)$ in W, then

$$(N_s^*\alpha \mid \beta) = ((N_1^*\alpha_1 \mid \beta_1) \mid \ldots \mid (N_n^*\alpha_n \mid \beta_n)) = (0 \mid \ldots \mid 0)$$

for all $\beta$; hence

$$N_s^*\alpha = (N_1^*\alpha_1 \mid \ldots \mid N_n^*\alpha_n) = (0 \mid \ldots \mid 0).$$

By earlier result this implies

$$N_s\alpha = (N_1\alpha_1 \mid \ldots \mid N_n\alpha_n) = (0 \mid \ldots \mid 0).$$

Conversely if

$$N_s\alpha = (N_1\alpha_1 \mid \ldots \mid N_n\alpha_n) = (0 \mid \ldots \mid 0)$$

then

$$N_s^*\alpha = (N_1^*\alpha_1 \mid \ldots \mid N_n^*\alpha_n) = (0 \mid \ldots \mid 0)$$

and

$$\begin{aligned}(N_s^*\alpha \mid \beta) &= (\alpha \mid N_s\beta) \\ &= ((\alpha_1 \mid N_1\beta_1) \mid \ldots \mid (\alpha_n \mid N_n\beta_n)) = (0 \mid \ldots \mid 0)\end{aligned}$$

for all $\beta$ in W. Hence the claim

**LEMMA 2.1.3:** *If $N_s = (N_1 \mid \ldots \mid N_n)$ is a normal operator and $\alpha = (\alpha_1 \mid \ldots \mid \alpha_n)$ is a super vector such that $N_s^2\alpha = (N_1^2\alpha_1 \mid \ldots \mid N_n^2\alpha_n) = (0 \mid \ldots \mid 0)$ then $N_s\alpha = (N_1\alpha_1 \mid \ldots \mid N_n\alpha_n) = (0 \mid \ldots \mid 0)$.*

*Proof*: Suppose $N_s$ is normal and $N_s^2\alpha = (N_1^2\alpha_1 \mid \ldots \mid N_n^2\alpha_n) = (0 \mid \ldots \mid 0)$. Then $N_s\alpha$ lies in the super range of $N_s$ and also lies in the null super space of $N_s$. Just by the above lemma this implies $N_s\alpha = (N_1\alpha_1 \mid \ldots \mid N_n\alpha_n) = (0 \mid \ldots \mid 0)$.

**LEMMA 2.1.4**: *Let $T_s = (T_1 \mid \ldots \mid T_n)$ be a normal operator and $f = (f_1 \mid \ldots \mid f_n)$ be any super polynomial with coefficients in the scalar field. Then $f(T) = (f_1(T_1) \mid \ldots \mid f_n(T_n))$ is also normal.*

*Proof:* Suppose

$$f = (a_0^1 + a_1^1 x + \ldots + a_{n_1}^1 x^{n_1} \mid \ldots \mid a_0^n + a_1^n x + \ldots + a_{n_n}^n x^{n_n})$$
$$= f = (f_1 \mid \ldots \mid f_n);$$

then $f(T_s) = (f_1(T_1) \mid \ldots \mid f_n(T_n))$

$$= (a_0^1 I_1 + a_1^1 T_1 + \ldots + a_{n_1}^1 T_1^{n_1} \mid \ldots \mid a_0^n I_n + a_1^n T_n + \ldots + a_{n_n}^n T_n^{n_n})$$

and

$$f(T_s^*) = (\overline{a}_0^1 I_1 + \overline{a}_1^1 T_1^* + \ldots + \overline{a}_{n_1}^1 (T_1^*)^{n_1} \mid \ldots \mid$$
$$\overline{a}_0^n I_n + \overline{a}_1^n T_n^* + \ldots + \overline{a}_{n_n}^n (T_n^*)^{n_n}).$$

Since $T_s T_s^* = T_s^* T_s$, it follows that $f(T_s)$ commutes with $f(T_s^*)$.

**LEMMA 2.1.5**: *Let $T_s = (T_1 \mid \ldots \mid T_n)$ be a normal operator and $f = (f_1 \mid \ldots \mid f_n)$ and $g = (g_1 \mid \ldots \mid g_n)$, relatively prime super polynomials with coefficients in the scalar field. Suppose $\alpha = (\alpha_1 \mid \ldots \mid \alpha_n)$ and $\beta = (\beta_1 \mid \ldots \mid \beta_n)$ are super vectors such that*

$$f(T_s)\alpha = (f_1(T_1)\alpha_1 \mid \ldots \mid f_n(T_n)\alpha_n) = (0 \mid \ldots \mid 0)$$

*and*

$$g(T_s)\beta = (g_1(T_1)\beta_1 \mid \ldots \mid g_n(T_n)\beta_n) = (0 \mid \ldots \mid 0)$$

*then*

$$(\alpha \mid \beta) = ((\alpha_1|\beta_1) \mid \ldots \mid (\alpha_n|\beta_n)) = (0 \mid \ldots \mid 0).$$

*Proof:* There are super polynomials a and b with coefficients in the scalar field such that $af + bg = (a_1f_1 + b_1g_1 \mid \ldots \mid a_nf_n + b_ng_n) = (1 \mid \ldots \mid 1)$ i.e. for each i, $g_i$ and $f_i$ are relatively prime and we have polynomials $a_i$ and $b_i$ such that $a_if_i + b_ig_i = 1$; $i = 1, 2, \ldots, n$. Thus

$a(T_s)\, f(T_s) + b(T_s)\, g(T_s) = I$ i.e.,
$(a_1(T_1)\, f_1(T_1) + b_1(T_1)\, g_1(T_1) \mid \ldots \mid a_n(T_n)\, f_n(T_n) + b_n(T_n)\, g_n(T_n))$

$= (I_1| \dots | I_n)$

and

$$\begin{aligned} \alpha &= (\alpha_1 | \dots | \alpha_n) \\ &= (g_1(T_1)b_1(T_1)\,\alpha_1 | \dots | g_n(T_n)b_n(T_n)\,\alpha_n) \\ &= g_s(T_s)\, b(T_s)\alpha. \end{aligned}$$

It follows that

$$\begin{aligned} (\alpha\,|\beta) &= ((\alpha_1|\beta_1) | \dots | (\alpha_n | \beta_n) = (g(T_s)\, b(T_s)\, \alpha|\beta) \\ &= ((g_1(T_1)b_1(T_1)\alpha_1|\beta_1 | \dots | g_n(T_n)b_n(T_n)\alpha_n|\beta_n) \\ &= (b_1(T_1)\alpha_1|g_1(T_1)^*\,\beta_1 | \dots | b_n(T_n)\alpha_n | g_n(T_n)^*\,\beta_n)) \\ &= (b(T_s)\,\alpha \mid g\,(T_s)^*\,\beta). \end{aligned}$$

By assumption

$$g(T_s)\beta = ((g_1(T_1)\beta_1 | \dots | g_n(T_n)\beta_n) = (0 | \dots | 0).$$

By earlier lemma

$$g(T) = (g_1(T_1) | \dots | g_n(T_n))$$

is normal. Therefore by earlier result

$$\begin{aligned} g(T)^*\beta &= (g_1(T_1)^*\beta_1 | \dots | g_n(T_n)^*\beta_n) \\ &= (0| \dots |0) \end{aligned}$$

hence

$$\begin{aligned} (\alpha|\,\beta) &= ((\alpha_1 | \beta_1) | \dots | (\alpha_n | \beta_n)) \\ &= (0 | \dots | 0). \end{aligned}$$

We call supersubspaces $W_j = (W^1_{j_1} \mid \dots \mid W^n_{j_n})$ ; $1 \le j_t \le k_t$; t = 1, 2, …, n, the primary super components of V under $T_s$.

**COROLLARY 2.1.11**: Let V = $(T_1 | \dots | T_n)$ be a normal operator on a finite $(n_1 | \dots | n_n)$ dimensional super inner product space V = $(V_1 | \dots | V_n)$ and $W_1, \dots W_k$ where $W_t = (W^1_{t_1} \mid \dots \mid W^n_{t_n})$ ; t = 1, 2, …, n be the primary super components of V under $T_s$; suppose $(W^1 | \dots | W^n)$ is a super subspace of V which is super invariant under $T_s$.

Then $W = \sum_{j} W \cap W_j$

$$= \left( \sum_{j_1} W^1 \cap W^1_{j_1} \mid \ldots \mid \sum_{j_n} \left( W^n \cap W^n_{j_n} \right) \right).$$

The proof is left as an exercise for the reader.

In fact we have to define super unitary transformation analogous to a unitary transformation.

**DEFINITION 2.1.11:** *Let $V = (V_1 \mid \ldots \mid V_n)$ and $V' = (V'_1 \mid \ldots \mid V'_n)$ be super inner product spaces over the same field F. A linear transformation $U_s = (U_1 \mid \ldots \mid U_n)$ from V into V′ is called a super unitary transformation, if it maps V onto V′ and preserves inner products. i.e. $U_i$: $V_i \to V'_i$ and preserves inner products for every i = 1, 2, …, n. If $T_s$ is a linear operator on V and $T'_s$ is a linear operator on V′ then $T_s$ is super unitarily equivalent to $T'_s$ if there exists a super unitary transformation $U_s$ of V onto V′ such that*

$U_s T_s U_s^{-1} = T'_s$ *i.e.* $(U_1 T_1 U_1^{-1} \mid \ldots \mid U_n T_n U_n^{-1}) = (T'_1 \mid \ldots \mid T'_n)$.

**LEMMA 2.1.6**: *Let $V = (V_1 \mid \ldots \mid V_n)$ and $V' = (V'_1 \mid \ldots \mid V'_n)$ be finite $(n_1, \ldots, n_n)$ dimensional super inner product spaces over the same field F. Suppose $T = (T_1 \mid \ldots \mid T_n)$ is a linear operator on $V = (V_1 \mid \ldots \mid V_n)$ and that $T'_s = (T'_1 \mid \ldots \mid T'_n)$ is a linear operator on $V' = (V'_1 \mid \ldots \mid V'_n)$. Then $T_s$ is super unitarily equivalent to $T'_s$ if and only if there is an orthonormal super basis $B = (B_1 \mid \ldots \mid B_n)$ of V and an orthonormal super basis $B' = (B'_1 \mid \ldots \mid B'_n)$ of V′ such that*

$$[T_s]_B = [T'_s]_{B'}$$

$$\text{i.e. } ([T_1]_{B_1} \mid \ldots \mid [T_n]_{B_n} = ([T'_1]_{B'_1} \mid \ldots \mid [T'_n]_{B'_n}).$$

The proof of lemma 2.1.6 and the following theorem are left for the reader.

**THEOREM 2.1.31:** *Let $V = (V_1 \mid \ldots \mid V_n)$ and $V' = (V_1' \mid \ldots \mid V_n')$ be finite $(n_1, \ldots, n_n)$ dimensional super inner product spaces over the same field F. Suppose $T_s$ is a normal operator on V and that $T'_s$ is a normal operator on V′. Then $T_s$ is unitarily equivalent to $T'_s$ if and only if $T_s$ and $T'_s$ have the same characteristic super polynomials.*

## 2.2 Superbilinear Form

Now we proceed onto give a brief description of Bilinear super forms or superbilinear forms before we proceed onto describe the applications of super linear algebra.

**DEFINITION 2.2.1:** *Let $V = (V_1 \mid \ldots \mid V_n)$ be a super vector space over the field F. A bilinear super form on V is a super function $f = (f_1 \mid \ldots \mid f_n)$ which assigns to each ordered pair of super vectors $\alpha = (\alpha_1 \mid \ldots \mid \alpha_n)$ and $\beta = (\beta_1 \mid \ldots \mid \beta_n)$ in V an n-tuple of scalars $f(\alpha, \beta) = (f_1(\alpha_1, \beta_1) \mid \ldots \mid f_n(\alpha_n, \beta_n))$ in F which satisfies:*

i. $f(c\alpha^1 + \alpha^2, \beta) = cf(\alpha^1, \beta) + f(\alpha^2, \beta)$

i.e., $(f_1(c_1\alpha_1^1 + \alpha_1^2, \beta_1) \mid \ldots \mid f_n(c_n\alpha_n^1 + \alpha_n^2, \beta_n))$

$= (c_1 f_1(\alpha_1^1, \beta_1) + f_1(\alpha_1^2, \beta_1) \mid \ldots \mid c_n f_n(\alpha_n^1, \beta_n) + f_n(\alpha_n^2, \beta_n))$

where $\alpha^1 = (\alpha_1^1 \mid \ldots \mid \alpha_n^1)$ and $\alpha^2 = (\alpha_1^2 \mid \ldots \mid \alpha_n^2)$.

ii. $f(\alpha, c\beta^1 + \beta^2) = cf(\alpha, \beta^1) + f(\alpha, \beta^2)$

i.e., $f_1(\alpha_1, c_1\beta_1^1 + \beta_1^2) \mid \ldots \mid f_n(\alpha_n, c_n\beta_n^1 + \beta_n^2))$

$= (c_1 f_1(\alpha_1, \beta_1^1) + f_1(\alpha_1, \beta_1^2) \mid \ldots \mid c_n f_n(\alpha_n, \beta_n^1) + f_n(\alpha_n, \beta_n^2))$

where $\beta^1 = (\beta_1^1 \mid \ldots \mid \beta_n^1)$ and $\beta^2 = (\beta_1^2 \mid \ldots \mid \beta_n^2)$.

*If $V \times V$ denotes the set of all ordered pairs of super vectors in V this definition can be rephrased as follows:*

*A bilinear superform on $V = (V_1 \mid \ldots \mid V_n)$ is a super function $f = (f_1 \mid \ldots \mid f_n)$ from $V \times V = (V_1 \times V_1 \mid \ldots \mid V_n \times V_n)$ into $(F \mid \ldots \mid F)$ which is linear as a superfunction on either of its arguments when the other is fixed. The super zero function (or* *zero super function) from $V \times V$ into (F | … | F) is clearly a bilinear superform. If $f = (f_1 \mid \ldots \mid f_n)$ and $g = (g_1 \mid \ldots \mid g_n)$ then $cf + g$ is also a bilinear superform, for any bilinear superforms f and g where $c = (c_1 \mid \ldots \mid c_n)$ i.e., $cf + g = (c_1f_1 + g_1 \mid \ldots \mid c_nf_n + g_n)$. We shall denote the super space of bilinear superforms on V by SL (V, V, F). SL (V, V, F) = {collection of all bilinear superforms from $V \times V$ into (F | … | F)} = $(L^1(V_1, V_1, F) \mid \ldots \mid L^n(V_n, V_n, F)\}$; where each $L^i(V_i, V_i, F)$ is a bilinear form, from $V_i \times V_i \to F$, i = 1, 2, …, n}.*

**DEFINITION 2.2.2:** *Let $V = (V_1 \mid \ldots \mid V_n)$ be finite dimensional $(n_1, \ldots, n_n)$ super vector space and let $B = (B_1 \mid \ldots \mid B_n) = \{\alpha_1^1 \ldots \alpha_{n_1}^1, \ldots, \alpha_1^n \ldots \alpha_{n_n}^n\}$ be an ordered super basis for V. If $f = (f_1 \mid \ldots \mid f_n)$ is a bilinear superform on V, the super diagonal matrix of f in the ordered super basis B is a $(n_1 \times n_1, \ldots, n_n \times n_n)$ super diagonal matrix*

$$A = \left(\begin{array}{c|c|c|c} A_1 & 0 & & 0 \\ \hline 0 & A_2 & & 0 \\ \hline & & & \\ \hline 0 & 0 & & A_n \end{array}\right)$$

*where each $A_t$ is a $n_t \times n_t$ matrix; t = 1, 2, …, n with entries $A^t_{i_t j_t} = f_t(\alpha^t_{i_t}, \alpha^t_{j_t})$; $1 \le i_t, j_t \le n_t$; t = 1, 2, …, n. At times we shall denote the super diagonal matrix A by $[f]_B = ([f_1]_{B_1} \mid \ldots \mid [f_n]_{B_n})$.*

We now give the interesting theorem on SL(V,V,F).

**THEOREM 2.2.1:** *Let $V = (V_1 \mid \ldots \mid V_n)$ be a finite dimensional super vector space over the field F. For each ordered super basis $B = (B_1 \mid \ldots \mid B_n)$ of V the super function which associates with each bilinear super form on V its super diagonal matrix in the ordered superbasis B is a super isomorphism of the super* *space SL(V, V, F) onto the super space of all ($n_1 \times n_1$, …, $n_n \times n_n$) super diagonal matrix A*

$$= \left(\begin{array}{c|c|c|c} A_1 & 0 & & 0 \\ \hline 0 & A_2 & & 0 \\ \hline & & & \\ \hline 0 & 0 & & A_n \end{array}\right)$$

*where $A_t$' s are $n_t \times n_t$ matrices with entries from F; for t = 1, 2, …, n.*

*Proof:* We observed from above that

$$f = (f_1 \mid \ldots \mid f_n) \to [f]_B = ([f_1]_{B_1} \mid \ldots \mid [f_n]_{B_n})$$

is a one to one correspondence between the set of bilinear superforms on V = ($V_1$ | … | $V_n$) and the set of all ($n_1 \times n_1$, …, $n_n \times n_n$) super diagonal matrices of the forms A with entries over F.

This is a linear transformation for

$$(cf + g)(\alpha_i, \alpha_j) = cf(\alpha_i, \alpha_j) + g(\alpha_i, \alpha_j)$$

$$\text{i.e. } ((c_1f_1 + g_1)(\alpha^1_{i_1}, \alpha^1_{j_1}) \mid \ldots \mid (c_nf_n + g_n)(\alpha^n_{i_n}, \alpha^n_{j_{n'}}))$$

$$= (c_1f_1(\alpha^1_{i_1}, \alpha^1_{j_1}) + g_1(\alpha^1_{i_1}, \alpha^1_{j_1}) \mid \ldots \mid c_nf_n(\alpha^n_{i_n}, \alpha^n_{j_n}) + g_n(\alpha^n_{i_n}, \alpha^n_{j_n}))$$

for each i and j where i = ($i_1$, …, $i_n$) and j = ($j_1$, …, $j_n$).

This simply imply

$$[cf + g]_B = c[f]_B + [g]_B$$

$$= ((c_1f_1 + g_1)_{B_1} \mid \ldots \mid (c_nf_n + g_n)_{B_n})$$

$$= (c_1[f_1]_{B_1} + [g_1]_{B_1} \mid \ldots \mid c_n[f_n]_{B_n} + [g_n]_{B_n}).$$

We now proceed onto give the following interesting corollary.

**COROLLARY 2.2.1**: *If B = (B$_1$ | ... | B$_n$)* $(\alpha_1^1 \ldots \alpha_{n_1}^1 \mid \ldots \mid \alpha_1^n \ldots \alpha_{n_n}^n)$ *is an ordered super basis for V = (V$_1$ | ... | V$_n$) and* $B^* = (B_1^* | \ldots | B_n^*)$ = $(\mathrm{L}_1^1 \ldots \mathrm{L}_{n_1}^1 \mid \ldots \mid \mathrm{L}_1^n \ldots \mathrm{L}_{n_n}^n)$ *is the dual super basis for* $V^* = (V_1^* \mid \ldots \mid V_n^*)$ *then the* $(n_1^2, \ldots, n_n^2)$ *bilinear superforms*

$$f_{ij}(\alpha, \beta) = L_i(\alpha)\, L_j(\beta)$$

i.e.

$$(f_{i_1 j_1}^1(\alpha_{i_1}, \beta_{j_1}) \mid \ldots \mid f_{i_n j_n}^n(\alpha_{i_n}, \beta_{j_n}))$$
$$= ((L_{i_1}^1(\alpha_{i_1}^1)\, L_{j_1}^1(\beta_{j_1}^1) \mid \ldots \mid L_{i_n}^n(\alpha_{i_n}^n)\, L_{j_n}^n(\beta_{j_n}^n));$$

*1 ≤ $i_t$, $j_t$ ≤ $n_t$; t = 1, 2, ..., n; form a super basis for the super space SL(V,V,F). In particular super dimension of SL(V,V,F) is* $(n_1^2, \ldots, n_n^2)$.

*Proof:* The dual super basis $\{\mathrm{L}_1^1 \ldots \mathrm{L}_{n_1}^1 \mid \ldots \mid \mathrm{L}_1^n \ldots \mathrm{L}_{n_n}^n\}$ is essentially defined by the fact that $\mathrm{L}_{i_t}^t(\alpha_{i_t}^t)$ is the $i_t^{th}$ coordinate of α in the ordered super basis B = (B$_1$ | … | B$_n$). Now the superfunction $\mathrm{f}_{ij} = (\mathrm{f}_{i_1 j_1}^1 \mid \ldots \mid \mathrm{f}_{i_n j_n}^n)$ defined by

$$\mathrm{f}_{ij}(\alpha, \beta) = \mathrm{L}_i(\alpha)\, \mathrm{L}_j(\beta)$$

i.e.

$$(\mathrm{f}_{i_1 j_1}^1(\alpha_{i_1}^1, \beta_{j_1}^1) \mid \ldots \mid \mathrm{f}_{i_n j_n}^n(\alpha_{i_n}^n, \beta_{j_n}^n))$$
$$= (\mathrm{L}_{i_1}^1(\alpha_{i_i}^1)\, \mathrm{L}_{j_1}^1(\beta_{j_1}^1) \mid \ldots \mid \mathrm{L}_{i_n}^n(\alpha_{i_n}^n)\, \mathrm{L}_{j_n}^n(\beta_{j_n}^n))$$

are bilinear superforms.
If

$$\alpha = (\mathrm{x}_1^1 \alpha_1^1 + \ldots + \mathrm{x}_{n_1}^1 \alpha_{n_1}^1 \mid \ldots \mid \mathrm{x}_1^n \alpha_1^n + \ldots + \mathrm{x}_{n_n}^n \alpha_{n_n}^n)$$

and

$$\beta = (\mathrm{y}_1^1 \alpha_1^1 + \ldots + \mathrm{y}_{n_1}^1 \alpha_{n_1}^1 \mid \ldots \mid \mathrm{y}_1^n \alpha_1^n + \ldots + \mathrm{y}_{n_n}^n \alpha_{n_n}^n)$$

then

$$\mathrm{f}_{ij}(\alpha, \beta) = \mathrm{x}_i\, \mathrm{y}_j.$$
$$(\mathrm{f}_{i_1 j_1}^1(\alpha_{i_1}^1, \beta_{j_1}^1) \mid \ldots \mid \mathrm{f}_{i_n j_n}^n(\alpha_{i_n}^n, \beta_{j_n}^n))$$
$$= (\mathrm{x}_{i_1}^1\, \mathrm{y}_{j_1}^1 \mid \ldots \mid \mathrm{x}_{i_n}^n\, \mathrm{y}_{j_n}^n).$$

Let f = ($f_1$ | … | $f_n$) be any bilinear superform on V = ($V_1$ | … | $V_n$) and let

$$A = \left(\begin{array}{c|c|c|c} A_1 & 0 & & 0 \\ \hline 0 & A_2 & & 0 \\ \hline & & & \\ \hline 0 & 0 & & A_n \end{array}\right)$$

be the super matrix of f = ($f_1$ | … | $f_n$) in the ordered super basis B = ($B_1$ | … | $B_n$). Then

$$\begin{aligned} f(\alpha, \beta) &= (f_1(\alpha_1, \beta_1) \mid \ldots \mid f_n(\alpha_n, \beta_n)) \\ &= \left( \sum_{i_1 j_1} A^1_{i_1 j_1} x^1_{i_1} y^1_{j_1} \Big| \ldots \Big| \sum_{i_n j_n} A^n_{i_n j_n} x^n_{i_n} y^n_{j_n} \right) \end{aligned}$$

which simply says that

$$\begin{aligned} f &= (f_1 \mid \ldots \mid f_n) \\ &= \left( \sum_{i_1 j_1} A^1_{i_1 j_1} f^1_{i_1 j_1} \Big| \ldots \Big| \sum_{i_n j_n} A^n_{i_n j_n} f^n_{i_n j_n} \right). \end{aligned}$$

It is now clear that the $(n_1^2, \ldots, n_n^2)$ forms $f_{ij} = (f^1_{i_1 j_1} \mid \ldots \mid f^n_{i_n j_n})$ comprise a super basis for SL(V,V,F).

We prove the following theorem.

**THEOREM 2.2.2:** *Let f = ($f_1$ | ... | $f_n$) be a bilinear superform on the finite ($n_1$, ..., $n_n$) dimensional super vector space V = ($V_1$ | ... | $V_n$). Let $L_f = (L^1_{f_1} \mid \ldots \mid L^n_{f_n})$ and $R_f = (R^1_{f_1} \mid \ldots \mid R^n_{f_n})$ be the linear transformation from V into $V^* = (V_1^* \mid \ldots \mid V_n^*)$ defined by*

$$(L_t \alpha)\, \beta = f(\alpha, \beta)$$

*i.e.* $$((L^1_{f_1} \alpha_1) \beta_1 \mid \ldots \mid (L^n_{f_n} \alpha_n) \beta_n)$$

$= (f_1(\alpha_1, \beta_1) \mid \ldots \mid f_n(\alpha_n , \beta_n))$
$= ((R^1_{f_1}\beta_1)\alpha_1 \mid \ldots \mid (R^n_{f_n}\beta_n)\alpha_n)$.

*Then super rank $L_f$ = super rank $R_f$.*

*i.e. super rank ($L_f$)* $= (rank\ L^1_{f_1}, \ldots, rank\ L^n_{f_n})$
$=$ *super rank($R_f$)*
$= (rank\ R^1_{f_1}, \ldots, rank\ R^n_{f_n})$.

The proof is left as an exercise for the reader.

Thus we say if f = ($f_1$ | … | $f_n$) is a bilinear super form on a finite dimensional ($n_1$, …, $n_n$) super vector space V = ($V_1$ | … | $V_n$) the super rank of f = ($f_1$ | … | $f_n$) is the n tuple of integers r = ($r_1$ | … | $r_n$) = super rank of $L_f$ = super rank of $R_f$ i.e. rank of $R^i_{f_i} = L^i_{f_i} = r_i$ for i = 1, 2, …, n..

Based on these results we give the following corollary which is left for the reader to prove.

**COROLLARY 2.2.2:** *The super rank of a bilinear superform is equal to the super rank of the superdiagonal matrix of the super form in the ordered super basis.*

**COROLLARY 2.2.3:** *If f = ($f_1$ | … | $f_n$) is a bilinear super form on the ($n_1$, …, $n_n$) dimensional super vector space V = ($V_1$ | … | $V_n$); the following are equivalent*

(a) *super rank f = (rank $f_1$, …, rank $f_n$) =($n_1$, …, $n_n$).*

(b) *For each nonzero $\alpha = (\alpha_1 \mid \ldots \mid \alpha_n)$ in V there is a $\beta = (\beta_1 \mid \ldots \mid \beta_n)$ in V such that $f(\alpha, \beta) = (f_1(\alpha_1, \beta_1) \mid \ldots \mid f_n(\alpha_n, \beta_n)) \neq (0 \mid \ldots \mid 0)$.*

(c) *For each non zero $\beta = (\beta_1 \mid \ldots \mid \beta_n)$ in V there is an α in V such that $f(\alpha, \beta) = (f_1(\alpha_1, \beta_1) \mid \ldots \mid f_n(\alpha_n , \beta_n)) \neq (0 \mid \ldots \mid 0)$.*

*Proof:* The condition (b) simply says that the super null space of $L_f = (L^1_{f_1} \mid \ldots \mid L^n_{f_n})$ is the zero super subspace. Statement (c) says that super null space of $R_f = (R^1_{f_1} \mid \ldots \mid R^n_{f_n})$ is the super zero subspace. The super linear transformations $L_f$ and $R_f$ have super nullity (0 | … | 0) if and only if they have super rank ($n_1$, …, $n_n$) i.e. if and only if super rank f = ($n_1$, …, $n_n$).

In view of the above conditions we define super non degenerate or non super degenerate or non super singular or super non singular bilinear superform.

**DEFINITION 2.2.3:** *A bilinear superform $f = (f_1 \mid \ldots \mid f_n)$ on a super vector space $V = (V_1 \mid \ldots \mid V_n)$ is called super non degenerate (or super non singular) if it satisfies conditions (b) and (c) of the corollary 2.2.3.*

Now we proceed onto define the notion of symmetric bilinear superforms.

**DEFINITION 2.2.4:** *Let $V = (V_1 \mid \ldots \mid V_n)$ be a super vector space over the field F. A super bilinear form $f = (f_1 \mid \ldots \mid f_n)$ on the super vector space V is super symmetric if $f(\alpha, \beta) = f(\beta, \alpha)$ for all $\alpha = (\alpha_1 \mid \ldots \mid \alpha_n)$ and $\beta = (\beta_1 \mid \ldots \mid \beta_n)$ in V i.e.*
*$f(\alpha, \beta) = (f_1(\alpha_1, \beta_1) \mid \ldots \mid f_n(\alpha_n, \beta_n)) = (f_1(\beta_1, \alpha_1) \mid \ldots \mid f_n(\beta_n, \alpha_n)) = f(\beta, \alpha)$.*

Now interms of the super matrix language we have the following. If $V = (V_1 \mid \ldots \mid V_n)$ be a finite ($n_1$, …, $n_n$) dimensional super vector space over the field F and f is a super symmetric bilinear form if and only if the super diagonal matrix

$$A = \left(\begin{array}{c|c|c|c} A_1 & 0 & & 0 \\ \hline 0 & A_2 & & 0 \\ \hline & & & \\ \hline 0 & 0 & & A_n \end{array}\right)$$

for some super basis B is super symmetric i.e. each $A_i$ is a symmetric matrix of A for i = 1, 2, …, n i.e. $A^t = A$ i.e. f(X, Y) = $X^tAY$ where X and Y super column matrices. This is true if and only if $X^tAY = Y^tAX$ for all supercolumn matrices X and Y, where $X = (X_1 \mid \ldots \mid X_n)^t$ and $Y = (Y_1 \mid \ldots \mid Y_n)^t$ where each $X_i$ and $Y_i$ are row vectors. Now

$$X^tAY = (X_1 \mid \ldots \mid X_n) \times \left(\begin{array}{c|c|c|c} A_1 & 0 & & 0 \\ \hline 0 & A_2 & & 0 \\ \hline & & & \\ \hline 0 & 0 & & A_n \end{array}\right) \times \left(\begin{array}{c} Y_1 \\ \hline \vdots \\ \hline Y_n \end{array}\right)$$

$$= \left(\begin{array}{c|c|c|c} X_1^tA_1Y_1 & 0 & & 0 \\ \hline 0 & X_2^tA_2Y_2 & & 0 \\ \hline & & & \\ \hline 0 & 0 & & X_n^tA_nY_n \end{array}\right)$$

$$= (Y_1 \mid \ldots \mid Y_n) \times \left(\begin{array}{c|c|c|c} A_1 & 0 & & 0 \\ \hline 0 & A_2 & & 0 \\ \hline & & & \\ \hline 0 & 0 & & A_n \end{array}\right) \times \left(\begin{array}{c} X_1 \\ \hline \vdots \\ \hline X_n \end{array}\right)$$

$$= \left(\begin{array}{c|c|c|c} Y_1^tA_1X_1 & & & 0 \\ \hline 0 & Y_2^tA_2X_2 & & 0 \\ \hline & & & \\ \hline 0 & 0 & & Y_n^tA_nX_n \end{array}\right).$$

Since $X^tAY$ is a1 × 1 super matrix we have $X^tAY = Y^tA^tX$. Thus f is super symmetric if and only if $Y^tA^tX = Y^tAX$ for all X, Y. Thus $A = A^t$. If f is a super diagonal, diagonal matrix clearly f is super symmetric as A is also super symmetric.

This paves way for us to define quadratic super form associated with a super symmetric bilinear super form f.

**DEFINITION 2.2.5:** *If $f = (f_1 \mid \ldots \mid f_n)$ is a symmetric bilinear superform the quadratic superform associated with f is the super function $q = (q_1 \mid \ldots \mid q_n)$ from V into $(F \mid \ldots \mid F)$ defined by $q(\alpha) = f(\alpha, \alpha)$ i.e.*

$$q(x) = (q_1(\alpha_1) \mid \ldots \mid q_n(\alpha_n)) = (f_1(\alpha_1, \alpha_1) \mid \ldots \mid f_n(\alpha_n, \alpha_n)).$$

If F is a subfield of the complex number the super symmetric bilinear super form f is completely determined by its associated super quadratic form accordingly the polarization super identity

$$f(\alpha,\beta) = \frac{1}{4} q(\alpha+\beta) - \frac{1}{4} q(\alpha-\beta)$$

i.e.

$$(f_1(\alpha_1, \beta_1) \mid \ldots \mid f_n(\alpha_n, \beta_n)) =$$

$$\left(\left(\frac{1}{4} q_1(\alpha_1 - \beta_1) - \frac{1}{4} q_1(\alpha_1 - \beta_1)\right) \mid \ldots \mid .\right.$$

$$\left.\left(\frac{1}{4} q_n(\alpha_n - \beta_n) - \frac{1}{4} q_n(\alpha_n - \beta_n)\right)\right).$$

If $f = (f_1 \mid \ldots \mid f_n)$ is such that each $f_i$ is the dot product, the associated quadratic superform is given by

$$\begin{aligned} q(x_1, \ldots, x_n) &= (q_1(x_1^1, \ldots, x_{n_1}^1) \mid \ldots \mid q_n(x_1^n, \ldots, x_{n_n}^n)) \\ &= ((x_1^1)^2 + \ldots + (x_{n_1}^1)^2, \ldots, (x_1^n)^2 + \ldots + (x_{n_n}^n)^2) \end{aligned}$$

i.e. $q(\alpha)$ is the super square length of $\alpha$. For the bilinear superform

$$f_A(X,Y) = (f_{A_1}^1(X_1, Y_1) \mid \ldots \mid f_{A_n}^n(X_n, Y_n))$$

$$= X^t A Y = \left(\begin{array}{c|c|c|c} X_1^t A_1 Y_1 & & & 0 \\ \hline 0 & X_2^t A_2 Y_2 & & 0 \\ \hline & & & \\ \hline 0 & 0 & & X_n^t A_n Y_n \end{array}\right)$$

$$= \left( \sum_{i_1, j_1} A^1_{i_1 j_1} x^1_{i_1} y^1_{j_1} \Big| \ldots \Big| \sum_{i_n, j_n} A^n_{i_n j_n} x^n_{i_n} y^n_{j_n} \right).$$

One of the important classes of super symmetric bilinear super forms consists of the super inner products on real vector spaces. If $V = (V_1 \mid \ldots \mid V_n)$ is a real vector super space a super inner product on V is super symmetric bilinear super form f on V which satisfies

$$f(\alpha, \alpha) = (f_1(\alpha_1, \alpha_1) \mid \ldots \mid f_n(\alpha_n, \alpha_n)) > (0 \mid \ldots \mid 0)$$
$$\text{if } \alpha = (\alpha_1 \mid \ldots \mid \alpha_n) \neq (0 \mid \ldots \mid 0). \qquad \text{(I)}$$

A super bilinear superform satisfying I is called super positive definite (or positive super definite). Thus a super inner product on a real super vector space is super positive definite, super symmetric bilinear superform on that space.

We know super inner product is also super non-degenerate i.e. each of its component inner products are non degenerate. Two super vectors $\alpha = (\alpha_1 \mid \ldots \mid \alpha_n)$ and $\beta = (\beta_1 \mid \ldots \mid \beta_n)$ are super orthogonal with respect to the super inner product $f = (f_1 \mid \ldots \mid f_n)$ if

$$f(\alpha, \beta) = (f_1(\alpha_1, \beta_1) \mid \ldots \mid f_n(\alpha_n, \beta_n)) = (0 \mid \ldots \mid 0).$$

The quadratic super forms $q(\alpha) = f(\alpha, \alpha) = (f_1(\alpha_1, \alpha_1) \mid \ldots \mid f_n(\alpha_n, \alpha_n))$ here each $f_i(\alpha_i, \alpha_i)$ takes only non negative values for i = 1, 2, …, n and $q(\alpha) = (q_1(\alpha_1) \mid \ldots \mid q_n(\alpha_n))$ is usually thought of as the super square length of α. i.e. square length of $\alpha_i$ for i = 1, 2, …, n as the orthogonality stems from the dot product.

If $f = (f_1 \mid \ldots \mid f_n)$ is any symmetric bilinear super form on a super vector space $V = (V_1 \mid \ldots \mid V_n)$ it is convenient to apply some terminology of super inner product of f. It is especially convenient to say that α and β are super orthogonal with respect to f if $f(\alpha, \beta) = (f_1(\alpha_1, \beta_1) \mid \ldots \mid f_n(\alpha_n, \beta_n)) = (0 \mid \ldots \mid 0)$. It pertinent to mention here that it is not proper to think of $f(\alpha, \alpha) = (f_1(\alpha_1, \alpha_1) \mid \ldots \mid f_n(\alpha_n, \alpha_n))$ as the super square of the length of α.

We give an interesting theorem for super vector spaces defined over the field of characteristic zero.

**THEOREM 2.2.3:** *Let $V = (V_1 \mid \ldots \mid V_n)$ be a super vector finite dimensional space over the field F of characteristic zero, and let $f = (f_1 \mid \ldots \mid f_n)$ be a super symmetric bilinear super form on V. Then there is an ordered super basis for V in which f is represented by a super diagonal diagonal matrix.*

*Proof*: To find an ordered super basis $B = (B_1 \mid \ldots \mid B_n) = (\alpha_1^1 \ldots \alpha_{n_1}^1 \mid \ldots \mid \alpha_1^n \ldots \alpha_{n_n}^n)$ such that $f(\alpha_i, \alpha_j) = (0 \mid \ldots \mid 0)$ for $i \neq j$ i.e.

$$f(\alpha_i, \alpha_j) = (f_1(\alpha_{i_1}^1, \alpha_{j_1}^1) \mid \ldots \mid f_n(\alpha_{i_n}^n, \alpha_{j_n}^n)) = (0 \mid \ldots \mid 0)$$

for $i_t \neq j_t$; t = 1, 2, …, n.

If $f = (0 \mid \ldots \mid 0)$ or $n = (1, 1, \ldots, 1)$ i.e. each $n_i = 1$ we have nothing to prove as the theorem is true. Thus we suppose a superform $f = (f_1 \mid \ldots \mid f_n) \neq (0 \mid \ldots \mid 0)$ and $n = (n_1, \ldots, n_n) > (1 \mid \ldots \mid 1)$. If $(f_1(\alpha_1, \alpha_1) \mid \ldots \mid f_n(\alpha_n, \alpha_n)) = (0 \mid \ldots \mid 0)$, for every $\alpha = (\alpha_1 \mid \ldots \mid \alpha_n) \in V$, the associated super quadratic form $q = (q_1 \mid \ldots \mid q_n)$ is identically $(0 \mid \ldots \mid 0)$, and the polarization super identity discussed earlier shows that $f = (f_1 \mid \ldots \mid f_n) = (0 \mid \ldots \mid 0)$. Thus there is a super vector $\alpha = (\alpha_1 \mid \ldots \mid \alpha_n)$ in V such that $f(\alpha, \alpha) = q(\alpha)$ i.e.,

$$(f_1 (\alpha_1,\alpha_1) \mid \ldots \mid f_n (\alpha_n,\alpha_n)) = (q_1(\alpha_1) \mid \ldots \mid q_n(\alpha_n))$$
$$= q(\alpha) \neq (0 \mid \ldots \mid 0).$$

Let W be a super dimensional subspace of V spanned by α i.e. $W = (W_1 \mid \ldots \mid W_n)$ is a super subspace of V spanned by $(\alpha_1 \mid \ldots \mid \alpha_n)$; each $W_t$ is spanned by $\alpha_t$, t = 1, 2, …, n. Let $W^{\perp} = (W_1^{\perp} \mid \ldots \mid W_n^{\perp})$ be the set of all super vectors $\beta = (\beta_1 \mid \ldots \mid \beta_n)$ in $V = (V_1 \mid \ldots \mid V_n)$ such that $f(\alpha, \beta) = (f_1(\alpha_1, \beta_1) \mid \ldots \mid f_n(\alpha_n, \beta_n)) = (0 \mid \ldots \mid 0)$.

Now we claim $V = W \oplus W^{\perp}$ i.e.

$$V = (V_1 \mid \ldots \mid V_n) = (W_1 \oplus W_1^{\perp} \mid \ldots \mid W_n \oplus W_n^{\perp}).$$

Certainly the super subspaces W and $W^{\perp}$ are super independent i.e., when we say super independent each $W_t$ and $W_t^{\perp}$ are independent for t = 1, 2, …, n; a typical super vector in $W = (W_1 \mid \ldots \mid W_n)$ is $c\alpha = (c_1\alpha_1 \mid \ldots \mid c_n\alpha_n)$ i.e., each super vector in $W_t$ is only of the form $c_t\alpha_t$; t = 1, 2, …, n where $c = (c_1 \mid \ldots \mid c_n)$ is a scalar n-tuple.

Also each super vector in $V = (V_1 \mid \ldots \mid V_n)$ is the sum of a super vector in W and a super vector in $W^{\perp}$. For let $\gamma = (\gamma_1 \mid \ldots \mid \gamma_n)$ be any super vector in V, and let

$$\beta = \gamma - \frac{f(\gamma,\alpha)}{f(\alpha,\alpha)}\alpha$$

i.e.,

$$\beta = (\beta_1 \mid \ldots \mid \beta_n)$$

$$= \left(\gamma_1 - \frac{f_1(\gamma_1,\alpha_1)}{f_1(\alpha_1,\alpha_1)}\alpha_1 \,\Big|\, \ldots \,\Big|\, \gamma_n - \frac{f_n(\gamma_n,\alpha_n)}{f_n(\alpha_n,\alpha_n)}\alpha_n\right).$$

Then

$$f(\alpha,\beta) = f(\alpha,\gamma) - \frac{f(\gamma,\alpha)f(\alpha,\alpha)}{f(\alpha,\alpha)}$$

i.e.,

$$(f_1(\alpha_1, \beta_1) \mid \ldots \mid f_n(\alpha_n, \beta_n))$$

$$= \left( f_1(\alpha_1,\gamma_1) - \frac{f_1(\gamma_1,\alpha_1)f_1(\alpha_1,\alpha_1)}{f_1(\alpha_1,\alpha_1)} \mid \ldots \mid \right.$$

$$f_n(\alpha_n, \alpha_n) - \frac{f_n(\gamma_n, \alpha_n) f_n(\alpha_n, \alpha_n)}{f_n(\alpha_n, \alpha_n)}\Bigg)$$

and since f is super symmetric f(α, β) = 0. Thus β is in the super subspace $W^{\perp}$. The expression

$$\gamma = \frac{f(\gamma, \alpha)}{f(\alpha, \alpha)} \alpha + \beta$$

i.e.

$$(\gamma_1 \mid \ldots \mid \gamma_n) = \left( \frac{f_1(\gamma_1, \alpha_1)}{f_1(\alpha_1, \alpha_1)} \alpha_1 + \beta_1 \mid \ldots \mid \frac{f_n(\gamma_n, \alpha_n)}{f_n(\alpha_n, \alpha_n)} \alpha_n + \beta_n \right)$$

which shows $V = W + W^{\perp}$
i.e.

$$(V_1 \mid \ldots \mid V_n) = (W_1 + W_1^{\perp} \mid \ldots \mid W_n + W_n^{\perp}).$$

The restriction of f to $W^{\perp}$ i.e. restriction of each $f_i$ to $W_i^{\perp}$ is a symmetric bilinear form, i = 1, 2, …, n; hence f is a symmetric bilinear superform on $W^{\perp}$. Since $W^{\perp}$ is of super dimension ($n_1 - 1, \ldots, n_n - 1$) we may assume by induction $W^{\perp}$ has a super basis $\{\alpha_2^1 \ldots \alpha_{n_1}^1 \mid \ldots \mid \alpha_2^n \ldots \alpha_{n_n}^n\}$ such that $f(\alpha_i, \alpha_j) = 0$; $i \neq j$;
i.e.,

$$f(\alpha_i, \alpha_j) = (f_1(\alpha_{i_t}^1, \alpha_{j_t}^1) \mid \ldots \mid f_n(\alpha_{i_t}^n, \alpha_{j_t}^n)) = (0 \mid \ldots \mid 0);$$

$i_t \neq j_t$; ($i_t \geq 2$; $j_t \geq 2$); $1 \leq i_t, j_t \leq n_t$; for every t = 1, 2, …, n. Putting $\alpha_1 = \alpha = (\alpha_1^1 \mid \ldots \mid \alpha_1^n)$ we obtain a super basis $\{\alpha_1^1 \ldots \alpha_{n_1}^1 \mid \ldots \mid \alpha_1^n \ldots \alpha_{n_n}^n\}$ for $V = (V_1 \mid \ldots \mid V_n)$ such that

$$f(\alpha_i, \alpha_j) = (f_1(\alpha_{i_1}^1, \alpha_{j_1}^1) \mid \ldots \mid f_n(\alpha_{i_n}^n, \alpha_{j_n}^n)) = (0 \mid \ldots \mid 0);$$

for $i \neq j$. i.e. $(i_1, \ldots, i_n) \neq (j_1, \ldots, j_n)$.

**COROLLARY 2.2.4:** *Let F be a field of complex numbers and let A be a super symmetric diagonal matrix over F i.e.*

$$A = \left(\begin{array}{c|c|c|c} A_1 & 0 & & 0 \\ \hline 0 & A_2 & & 0 \\ \hline & & & \\ \hline 0 & 0 & & A_n \end{array}\right)$$

*is super symmetric diagonal matrix in which each $A_i$ is a $n_i \times n_i$ matrix with entries from F, i = 1, 2, …, n. Then there is an invertible super square matrix*

$$P = \left(\begin{array}{c|c|c|c} P_1 & 0 & & 0 \\ \hline 0 & P_2 & & 0 \\ \hline & & & \\ \hline 0 & 0 & & P_n \end{array}\right)$$

*where each $P_i$ is a $n_i \times n_i$ invertible matrix with entries from F such that $P^tAP$ is super diagonal i.e.*

$$\left(\begin{array}{c|c|c|c} P_1^t & 0 & & 0 \\ \hline 0 & P_2^t & & 0 \\ \hline & & & \\ \hline 0 & 0 & & P_n^t \end{array}\right) \times \left(\begin{array}{c|c|c|c} A_1 & 0 & & 0 \\ \hline 0 & A_2 & & 0 \\ \hline & & & \\ \hline 0 & 0 & & A_n \end{array}\right) \times \left(\begin{array}{c|c|c|c} P_1 & 0 & & 0 \\ \hline 0 & P_2 & & 0 \\ \hline & & & \\ \hline 0 & 0 & & P_n \end{array}\right)$$

$$= \left(\begin{array}{c|c|c|c} P_1^tA_1P_1 & 0 & & 0 \\ \hline 0 & P_2^tA_2P_2 & & 0 \\ \hline & & & \\ \hline 0 & 0 & & P_n^tA_nP_n \end{array}\right).$$

*is superdiagonal i.e. each $P_i^tA_iP_i$ is a diagonal matrix, $1 \le i \le n$.*

We give yet another interesting theorem.

**THEOREM 2.2.4:** *Let $V = (V_1 \mid \ldots \mid V_n)$ be a finite $(n_1, \ldots, n_n)$ dimensional super vector space over the field of complex numbers. Let $f = (f_1 \mid \ldots \mid f_n)$ be a symmetric bilinear superform on V which has super rank $r = (r_1, \ldots, r_n)$. Then there is an ordered super basis $B = (B_1, \ldots, B_n) = (\beta_1^1, \ldots, \beta_{n_1}^1; \ldots; \beta_1^n, \ldots, \beta_{n_n}^n)$ for V such that*

*(i) The super diagonal matrix A of f in the basis B is super diagonal, diagonal matrix i.e. if*

$$A = \left(\begin{array}{c|c|c|c} A_1 & 0 & & 0 \\ \hline 0 & A_2 & & 0 \\ \hline & & & \\ \hline 0 & 0 & & A_n \end{array}\right)$$

*each $A_i$ is a diagonal $n_i \times n_i$ matrix, $i = 1, 2, \ldots, n$.*

*ii)* $$f(\beta_j, \beta_j) = \begin{cases} (1 \mid \ldots \mid 1), j = 1, 2, \ldots, r \\ (0 \mid \ldots \mid 0) \; j > r \end{cases}$$

*i.e.* $f(\beta_j, \beta_j) = (f_1(\beta_{j_1}^1, \beta_{j_1}^1) \mid \ldots \mid f_n(\beta_{j_n}^n, \beta_{j_n}^n)) = (1 \mid \ldots \mid 1)$
*if* $j_t = 1, 2, \ldots, r_t;\ 1 \le t \le n$
*and*
$f(\beta_j, \beta_j) = (f_1(\beta_{j_1}^1, \beta_{j_1}^1) \mid \ldots \mid f_n(\beta_{j_n}^n, \beta_{j_n}^n)) = (0 \mid \ldots \mid 0)$
*if* $j_t > r_t$ *for* $t = 1, 2, \ldots, n$.

The proof is left as an exercise for the reader.

**THEOREM 2.2.5:** *Let $V = (V_1 \mid \ldots \mid V_n)$ be a $(n_1, \ldots, n_n)$ dimensional super vector space over the field of real numbers and let $f = (f_1 \mid \ldots \mid f_n)$ be a symmetric bilinear super form on V which has super rank $r = (r_1, \ldots, r_n)$. Then there is an ordered super basis $(\beta_1^1, \ldots, \beta_{n_1}^1, \ldots, \beta_l^n, \ldots, \beta_{n_n}^n)$ for V in which the* *super diagonal matrix of f is a superdiagonal matrix such that the entries are only ±1*

*i.e.* $f(\beta_j, \beta_j) = (f_1(\beta^1_{j_1}, \beta^1_{j_1}) \mid \ldots \mid f_n(\beta^n_{j_n}, \beta^n_{j_n}))$ *= (±1 | … | ±1);*

$j_t$ *= 1, 2, …,* $r_t$*; t = 1, 2, …, n. Further more the number of superbasis vector* $\beta_j = (\beta^1_{j_1}, \ldots, \beta^1_{j_n})$ *for which*

$$f(\beta_j, \beta_j) = (f_1(\beta^1_{j_1}, \beta^1_{j_1}) \mid \ldots \mid f_n(\beta^n_{j_n}, \beta^n_{j_n}))$$
$$= (1 \mid \ldots \mid 1)$$

*is independent of the choice of the superbasis.*

*Proof:* There is a superbasis $\{\alpha^1_1 \ldots \alpha^1_{n_1}; \ldots; \alpha^n_1 \ldots \alpha^n_{n_n}\}$ for V = ($V_1$ | … | $V_n$) i.e. $\{\alpha^t_1 \ldots \alpha^t_{n_t}\}$ is a basis for $V_t$, t = 1, 2, …, n. such that

$$f(\alpha_i, \alpha_j) = (f_1(\alpha^1_{i_1}, \alpha^1_{j_1}) \mid \ldots \mid f_n(\alpha^n_{i_n}, \alpha^n_{j_n})) = (0 \mid \ldots \mid 0)$$

if $i_t \neq j_t$

$$f(\alpha_j, \alpha_j) = (f_1(\alpha^1_{j_1}, \alpha^1_{j_1}) \mid \ldots \mid f_n(\alpha^n_{j_n}, \alpha^n_{j_n})) \neq (0 \mid \ldots \mid 0)$$

for $1 \leq j_t \leq r_t$

$$f(\alpha_j, \alpha_j) = (f_1(\alpha^1_{j_1}, \alpha^1_{j_1}) \mid \ldots \mid f_n(\alpha^n_{j_n}, \alpha^n_{j_n})) = (0 \mid \ldots \mid 0)$$

$j_t > r_t$ for t = 1, 2, …, n.

Let

$$\beta_j = (\beta^1_{j_1} \ldots \beta^n_{j_n}) = |f(\alpha_j, \alpha_j)|^{-1/2} \alpha_j$$
$$= \left(|f_1(\alpha^1_{j_1}, \alpha^1_{j_1})|^{-1/2} \alpha^1_{j_1} \mid \ldots \mid |f_n(\alpha^n_{j_n}, \alpha^n_{j_n})|^{-1/2} \alpha^n_{j_n}\right)$$

$1 \leq j_t \leq r_t$ ; t = 1, 2, …, n.

$$\beta_j = (\beta^1_{j_1} \ldots \beta^n_{j_n}) = \alpha_j = (\alpha^1_{j_1} \ldots \alpha^n_{j_n})$$

$j_t > r_t$ ;   t = 1, 2, …, n;

then $\{\beta^1_1 \ldots \beta^n_{n_1}; \ldots; \beta^n_1 \ldots \beta^n_{n_n}\}$ is a super basis satisfying all the properties.

Let p = ($p_1$ | … | $p_n$) be the number of basis super vectors $\beta_j = (\beta^1_{j_1} \ldots \beta^n_{j_n})$ for which

$$f(\beta_j, \beta_j) = (f_1(\beta^1_{j_1}, \beta^1_{j_1}) \mid \ldots \mid f_n(\beta^n_{j_n}, \beta^n_{j_n})) = (1 \mid \ldots \mid 1) ;$$

we must show the number p is independent of the particular superbasis.

Let $V^+ = (V_1^+ \mid \ldots \mid V_n^+)$ be the super subspace of V = ($V_1$ | … | $V_n$) spanned by the super basis super vectors $\beta_j$ for which f($\beta_j$ $\beta_j$) = (–1 | … | –1). Now p = ($p_1$ | … | $p_n$) = super dim $V^+ = (\dim V_1^+, \ldots, \dim V_n^+)$ so it is the uniqueness of the super dimension of $V^+$ which we must show. It is easy to see that if ($\alpha_1$ | … | $\alpha_n$) is a nonzero super vector in $V^+$ then f(α, α) = $f_1(\alpha_1, \alpha_1)$ | … | $f_n(\alpha_n, \alpha_n)$) > (0 | … | 0) in other words f = ($f_1$, …, $f_n$) is super positive definite i.e. each $f_i$ is positive definite on the subspace $W_i^+$; i = 1, 2, …, n; of $W^+ = (W_1^+ \mid \ldots \mid W_n^+)$; the super subspace of $V^+$. Similarly if α = ($\alpha_1$ | … | $\alpha_n$) is a nonzero super vector in $V^- = (V_1^- \mid \ldots \mid V_n^-)$ then f(α, α) = ($f_1(\alpha_1, \alpha_1)$ | … | $f_n(\alpha_n, \alpha_n)$) < (0 | … | 0) i.e. f is super negative definite on the super subspace $V^-$. Now let $V^\perp = (V_1^\perp \mid \ldots \mid V_n^\perp)$ be super subspace spanned by the super basis of super vectors $\beta_j = (\beta^1_j \mid \ldots \mid \beta^n_{j_n})$ for which

$$f(\beta_j, \beta_j) = (f_1(\beta^1_{j_1}, \beta^1_{j_1}) \mid \ldots \mid f_n(\beta^n_{j_n}, \beta^n_{j_n})) = (0 \mid \ldots \mid 0).$$

If α = ($\alpha_1$ | … | $\alpha_n$) is in $V^\perp$ then f(α,β) = ($f_1$ ($\alpha_1$ , $\beta_1$) | … | $f_n$ ($\alpha_n$, $\beta_n$)) = (0 | 0 | … |0) for all β = ($\beta_1$ | … | $\beta_n$) in V.
Since $(\beta^1_1 \ldots \beta^1_{n_1}; \ldots; \beta^n_1 \ldots \beta^n_{n_n})$ is a super basis for V we have

$$V = V^+ \oplus V^- \oplus V^\perp$$
$$= (V_1^+ \oplus V_1^- \oplus V_1^\perp) \mid \ldots \mid V_n^+ \oplus V_n^- \oplus V_n^\perp).$$

Further if W is any super subspsace of V on which f is super positive definite then the super subspace W, $V^-$ and $V^\perp$ are super independent that is $W_i$, $V_i^-$ and $V_i^\perp$ are independent for i = 1, 2, …, n;

Suppose α is in W, β is in $V^-$ and γ is in $V^\perp$ then

$$\alpha + \beta + \gamma = (\alpha_1 + \beta_1 + \gamma_1 \mid \ldots \mid \alpha_n + \beta_n + \gamma_n) = (0 \mid \ldots \mid 0).$$

Then (0 | … | 0)

$$\begin{aligned}
&= (f_1(\alpha_1, \alpha_1 + \beta_1 + \gamma_1) \mid \ldots \mid f_n(\alpha_n, \alpha_n + \beta_n + \gamma_n)) \\
&= (f_1(\alpha_1, \alpha_1) + f_1(\alpha_1, \beta_1) + f_1(\alpha_1, \gamma_1) \mid \ldots \mid \\
&\quad f_n(\alpha_n, \alpha_n) + f_n(\alpha_n, \beta_n) + f_n(\alpha_n, \gamma_n)) \\
&= f(\alpha, \alpha) + f(\alpha, \beta) + f(\alpha, \gamma) .
\end{aligned}$$

$$\begin{aligned}
(0 \mid \ldots \mid 0) &= f(\beta, \alpha + \beta + \gamma) \\
&= (f_1(\beta_1, \alpha_1 + \beta_1 + \gamma_1) \mid \ldots \mid f_n(\beta_n, \alpha_n + \beta_n + \gamma_n)) \\
&= (f_1(\beta_1, \alpha_1) + f_1(\beta_1, \beta_1) + f_1(\beta_1, \gamma_1) \mid \ldots \mid \\
&\quad f_n(\beta_n, \alpha_n) + f_n(\beta_n, \beta_n) + f_n(\beta_n, \gamma_n)). \\
&= (f(\beta, \alpha) + f(\beta, \beta) + f(\beta, \gamma)).
\end{aligned}$$

Since γ is in $V^\perp = (V_1^\perp \mid \ldots \mid V_n^\perp)$, f(α, γ) = f(β, γ) = (0 | … |0) i.e.

$$\begin{aligned}
(f_1(\alpha_1, \gamma_1) \mid \ldots \mid f_n(\alpha_n, \gamma_n)) &= f_1(\beta_1, \gamma_1) \mid \ldots \mid f_n(\beta_n, \gamma_n)) \\
&= (0 \mid \ldots \mid 0)
\end{aligned}$$

and since f is super symmetric i.e. each $f_i$ is symmetric (i = 1, 2, …, n) we obtain

$$\begin{aligned}
(0 \mid \ldots \mid 0) &= f(\alpha, \alpha) + f(\alpha, \beta) \\
&= (f_1(\alpha_1, \alpha_1) + f_1(\alpha_1, \beta_1) \mid \ldots \mid f_n(\alpha_n, \alpha_n) + f_n(\alpha_n, \beta_n))
\end{aligned}$$

and

$$\begin{aligned}
(0 \mid \ldots \mid 0) &= f(\beta, \beta) + f(\alpha, \beta) \\
&= (f_1(\beta_1, \beta_1) + f_1(\alpha_1, \beta_1) \mid \ldots \mid f_n(\beta_n, \beta_n) + f_n(\alpha_n, \beta_n)).
\end{aligned}$$

Hence

$$f(\alpha, \alpha) = f(\beta, \beta)$$

i.e. $(f_1(\alpha_1, \alpha_1) \mid \ldots \mid f_n(\alpha_n, \alpha_n)) = (f_1(\beta_1, \beta_1) \mid \ldots \mid f_n(\beta_n, \beta_n))$.

Since

$$f(\alpha, \alpha) = (f_1(\alpha_1, \alpha_1) \mid \ldots \mid f_n(\alpha_n, \alpha_n)) \geq (0 \mid \ldots \mid 0)$$

and $f(\beta, \beta) = (f_1(\beta_1, \beta_1) \mid \ldots \mid f_n(\beta_n, \beta_n)) \leq (0 \mid \ldots \mid 0)$

it follows that $f(\alpha, \alpha) = f(\beta, \beta) = (0 \mid \ldots \mid 0)$.
But f is super positive definite on $W = (W_1 \mid \ldots \mid W_n)$ and super negative definite on $V^- = (V_1^- \mid \ldots \mid V_n^-)$. We conclude that $\alpha = (\alpha_1 \mid \ldots \mid \alpha_n) = (\beta_1 \mid \ldots \mid \beta_n) = (0 \mid \ldots \mid 0)$ and hence that $\gamma = (0 \mid \ldots \mid 0)$ as well. Since $V = V^+ \oplus V^- \oplus V^\perp$ and $W, V^-, V^\perp$ are super independent we see that super dim W < super dim $V^+$ i.e. $(\dim W_1, \ldots, \dim W_n) \leq (\dim V_1^+, \ldots, \dim V_n^+)$. That is if $W = (W_1 \mid \ldots \mid W_n)$ is any supersubspace of V on which f is super positive definite, the super dimension of W cannot exceed the superdimension of $V^+$. If $B_1$ is the superbasis given in the theorem, we shall have corresponding supersubspaces $V_I^+, V_I^-$ and $V_I^\perp$ and the argument above shows that superdim $V_I^+ \leq$ superdim $V^+$. Reversing the argument we obtain superdim $V^+ \leq$ super dim $V_I^+$ and subsequently,

$$\text{super dim } V^+ = \text{super dim } V_I^+ .$$

There are several comments we can make about the super basis $\{\beta_1^1 \ldots \beta_{n_1}^1; \ldots; \beta_1^n \ldots \beta_{n_n}^n\}$ and the associated super subspaces $V^+, V^-$ and $V^\perp$. First we have noted that $V^\perp$ is exactly the subsubspace of super vectors which are super orthogonal to all of V. We noted above that $V^\perp$ is contained in the super subspace. But super dim $V^\perp$ = super dim V – (super dim $V^+$ + super dim $V^-$) = super dim V – super rank f; so every super vector α such that

$$f(\alpha, \beta) = (f_1(\alpha_1, \beta_1) \mid \ldots \mid f_n(\alpha_n, \beta_n)) = (0 \mid \ldots \mid 0)$$

for all $\beta = (\beta_1 \mid \ldots \mid \beta_n)$ must be in $V^\perp$. Further the subspace $V^\perp$ is unique. The super dimension of $V^\perp$ is the largest possible super dimension of any subspace on which f is super positive definite. Similarly super dim $V^-$ is the largest super dimension of any supersubspace on which f is super negative definite. Of course super dim $V_I^+$ + super $V^-$ = super rank f.

The super number is the n-tuple superdim $V^+$ - superdim $V^-$ is often called the super signature of f. This is derived because the super dimensions of $V^+$ and $V^-$ are easily determined from the super rank of f and the super signature of f.

This property is worth a good relation of super symmetric bilinear superforms on real vector spaces to super inner products. Suppose V is a finite $(n_1, \ldots, n_n)$ dimensional real super vector space and $W^1$, $W^2$ and $W^3$ are super subspace of V such that

$$V = W^1 \oplus W^2 \oplus W^3$$

i.e.

$$\begin{aligned}(V_1 \mid \ldots \mid V_n) &= (W_1^1 \mid \ldots \mid W_n^1) \oplus (W_1^2 \mid \ldots \mid W_n^2) \oplus (W_1^3 \mid \ldots \mid W_n^3) \\ &= ((W_1^1 \oplus W_1^2 \oplus W_1^3) \mid \ldots \mid (W_n^1 \oplus W_n^2 \oplus W_n^3)).\end{aligned}$$

Suppose that $f^1 = (f_1^1 \mid \ldots \mid f_n^1)$ is an super inner product on $W^1$ and $f^2 = (f_1^2 \mid \ldots \mid f_n^2)$ an super inner product on $W^2$ . We can define a super symmetric bilinear superform $f^1 = (f_1^1 \mid \ldots \mid f_n^1)$ on V as follows. If $\alpha, \beta \in V$ then we write

$$\begin{aligned}\alpha &= (\alpha^1 + \alpha^2 + \alpha^3) \\ &= (\alpha_1^1 \mid \ldots \mid \alpha_n^1) + (\alpha_1^2 \mid \ldots \mid \alpha_n^2) + (\alpha_1^3 \mid \ldots \mid \alpha_n^3) \\ &= (\alpha_1^1 + \alpha_1^2 + \alpha_1^3 \mid \ldots \mid \alpha_n^1 + \alpha_n^2 + \alpha_n^3)\end{aligned}$$

and similarly

$$\beta = (\beta_1^1 + \beta_1^2 + \beta_1^3 \mid \ldots \mid \beta_n^1 + \beta_n^2 + \beta_n^3)$$

with $\alpha_j$ and $\beta_j$ in $V^j$; $1 \leq j \leq 3$. Let $f(\alpha, \beta) = f_1(\alpha^1, \beta^1) - f_2(\alpha^2, \beta^2)$. The super subspace $V^\perp$ for f will be $W^3$, $W^1$ is suitable $V^+$ and $W^2$ is the suitable $V^-$.

Let $V = (V_1 \mid \ldots \mid V_n)$ be a super vector space defined over a subfield F of the field of complex numbers. A bilinear super form $f = (f_1 \mid \ldots \mid f_n)$ on V is called skew super symmetric if $f(\alpha, \beta) = -f(\beta, \alpha)$ for all super vectors $\alpha, \beta$ in V. If V is a finite $(n_1 \mid \ldots \mid n_n)$ dimensional the bilinear super form $f = (f_1 \mid \ldots \mid f_n)$ is skew super symmetric if and only if its super diagonal matrix A in some ordered super basis is skew super symmetric i.e., $A^t = -A$ i.e., if

$$A^t = \left(\begin{array}{c|c|c|c} A_1^t & 0 & & 0 \\ \hline 0 & A_2^t & & 0 \\ \hline & & & \\ \hline 0 & 0 & & A_n^t \end{array}\right)$$

then

$$-A = \left(\begin{array}{c|c|c|c} -A_1 & 0 & & 0 \\ \hline 0 & -A_2 & & 0 \\ \hline & & & \\ \hline o & 0 & & -A_n \end{array}\right)$$

i.e., each $A_i^t = -A_i$ for i = 1, 2, …, n.

Further here f(α, α) = (0 | … | 0) i.e., $(f_1(\alpha_1, \alpha_1) \mid \ldots \mid f_n(\alpha_n, \alpha_n))$ = (0 | … | 0) for every α in V since f(α, α) = – f(α, α). Let us suppose $f = (f_1 \mid \ldots \mid f_n)$ is a non zero super skew symmetric super bilinear form on $V = (V_1 \mid \ldots \mid V_n)$. Since f ≠ (0 | … | 0) there are super vectors α, β in V such that f(α, β) ≠ (0 | … | 0) multiplying α by a suitable scalar we may assume that f(α, β) = (1 | … | 1). Let γ be any super vector in the super subspace spanned by α and β, say

$$\gamma = C\alpha + d\beta$$

i.e.,

$$\begin{aligned} \gamma &= (\gamma_1 \mid \ldots \mid \gamma_n) \\ &= (C_1\alpha_1 + d_1\beta_1 \mid \ldots \mid C_n\alpha_n + d_n\beta_n). \end{aligned}$$

Then

$$\begin{aligned} f(\gamma, \alpha) &= f(C\alpha + d\beta, \alpha) \\ &= df(\beta, \alpha) \\ &= -d \\ &= (-d_1 \mid \ldots \mid -d_n) \end{aligned}$$

and

$$\begin{aligned} f(\gamma, \beta) &= f(C\alpha + d\beta, \beta) \\ &= Cf(\alpha, \beta) \\ &= -C \\ &= (-C_1 \mid \ldots \mid -C_n) \end{aligned}$$

i.e., each $d_i = d_i f_i(\beta_i, \alpha_i)$ and each $C_i = C_i f_i(\alpha_i, \beta_i)$ for i = (1, 2, …, $n_n$). In particular note that α and β are linearly super independent for if $\gamma = (\gamma_1 \mid \ldots \mid \gamma_n) = (0 \mid \ldots \mid 0)$, then

$$f(\gamma, \alpha) = (f_1(\gamma_1, \alpha_1) \mid \ldots \mid f_n(\gamma_n, \alpha_n))$$
$$f(\gamma, \beta) = (f_1(\gamma_1, \beta_1) \mid \ldots \mid f_n(\gamma_n, \beta_n))$$
$$= (0 \mid \ldots \mid 0).$$

Let $W = (W_1 \mid \ldots \mid W_n)$ be a (2, …, 2) dimensional super subspace spanned by α and β i.e., each $W_i$ is spanned by $\alpha_i$ and $\beta_i$ for i = 1, 2, …, $n_n$. Let $W^{\perp} = (W_1^{\perp} \mid \ldots \mid W_n^{\perp})$ be the set of all super vectors $\delta = (\delta_1 \mid \ldots \mid \delta_n)$ such that $f(\delta, \alpha) = f(\delta, \beta) = (0 \mid \ldots \mid 0)$, that is the set of all δ such that $f(\delta, \gamma) = (0 \mid \ldots \mid 0)$ for every γ in the super subspace $W = (W_1 \mid \ldots \mid W_n)$. We claim

$$V = W \oplus W^{\perp} = (W_1 \oplus W_1^{\perp} \mid \ldots \mid W_n \oplus W_n^{\perp}).$$

For let $\in = (\in_1 \mid \ldots \mid \in_n)$ be any super vector in V and $\gamma = f(\in, \beta)\alpha - f(\in, \alpha)\beta$ ; $\delta = \in - \gamma$. Thus γ is in W and δ is in $W^{\perp}$ for

$$\begin{aligned} f(\delta, \alpha) &= f(\in - f(\in, \beta)\alpha + f(\in, \alpha)\beta, \alpha) \\ &= f(\in, \alpha) + f(\in, \alpha)\, f(\beta, \alpha) = (0 \mid \ldots \mid 0) \end{aligned}$$

and similarly $f(\delta, \beta) = (0 \mid \ldots \mid 0)$. Thus every ∈ in V is of the form $\in = \gamma + \delta$, with γ in W and δ in $W^{\perp}$. From earlier results $W \cap W^{\perp} = (0 \mid \ldots \mid 0)$ and so $V = W \oplus W^{\perp}$.

Now restriction of f to $W^{\perp}$ is a skew symmetric bilinear super form on $W^{\perp}$. This restriction may be the zero super form. If it is not, there are super vectors α' and β' in $W^{\perp}$ such that $f(\alpha', \beta') = (1 \mid \ldots \mid 1)$. If we let W' be the two super dimensional i.e., (2, …, 2) dimensional super subspaces spanned by α′ and β′ then we shall have $V = W \oplus W' \oplus W_o$ where $W_o$ is the set of all super vectors δ in $W^{\perp}$ such that $f(\alpha', \delta) = f(\beta', \delta) = (0 \mid \ldots \mid 0)$. If the restriction of f to $W_o$ is not the zero super form, we may select super vectors α'', β'' in $W_o$ such that $f(\alpha'', \beta'') = (1 \mid \ldots \mid 1)$ and continue.

In the finite ($n_1$, …, $n_n$) dimensional case it should be clear that we obtain finite sequence of pairs of super vectors

$\{(\alpha_1^1, \beta_1^1)\ldots(\alpha_{k_1}^1, \beta_{k_1}^1) \mid \ldots \mid (\alpha_1^n, \beta_1^n)\ldots(\alpha_{k_n}^n, \beta_{k_n}^n)\}$ with the following properties

(a) $f(\alpha_j, \beta_j) = (f_1(\alpha_{j_1}^1, \beta_{j_1}^1) \mid \ldots \mid f_n(\alpha_{j_n}^n, \beta_{j_n}^n))$
$= (1 \mid \ldots \mid 1)$ ; $j = 1, 2, \ldots, k$

(b) $f(\alpha_i, \alpha_j) = (f_1(\alpha_{i_1}^1, \alpha_{j_1}^1) \mid \ldots \mid f_n(\alpha_{i_n}^n, \alpha_{j_n}^n))$
$= f(\beta_j, \beta_j)$
$= (f_1(\beta_{i_1}^1, \beta_{j_1}^1) \mid \ldots \mid f_n(\beta_{i_n}^n, \beta_{j_n}^n))$
$= f(\alpha_i, \beta_j)$
$= (f_1(\alpha_{i_1}^1, \beta_{j_1}^1) \mid \ldots \mid f_n(\alpha_{i_n}^n, \beta_{j_n}^n))$
$= (0 \mid \ldots \mid 0)$ ; $i \neq j$
i.e. $i_t \neq j_t$; $1 \leq i_t, j_t \leq k_t$; $t = 1, 2, \ldots, n$.

(c) If $W_j = (W_{j_1}^1 \mid \ldots \mid W_{j_n}^n)$ is the two super dimensional super subspace i.e. super dim $W_j$ is $(2, \ldots, 2)$ and super spanned by $\alpha_j = (\alpha_{j_1}^1 \mid \ldots \mid \alpha_{j_n}^n)$ and $\beta_j = (\beta_{j_1}^1 \mid \ldots \mid \beta_{j_n}^n)$ then
$V = W_1 \oplus \ldots \oplus W_k \oplus W_0$
$= ((W_1^1 \oplus \ldots \oplus W_{k_1}^1 \oplus W_0^1) \mid \ldots \mid (W_1^n \oplus \ldots \oplus W_{k_n}^n \oplus W_0^n))$

where every super vector in $W_0 = (W_0^1 \mid \ldots \mid W_0^n)$ is super orthogonal to all $\alpha_j = (\alpha_{j_1}^1 \mid \ldots \mid \alpha_{j_n}^n)$ and $\beta_j = (\beta_{j_1}^1 \mid \ldots \mid \beta_{j_n}^n)$ and the super restriction of f to $W_0$ is the zero super form.

Next we prove another interesting theorem.

**THEOREM 2.2.6:** *Let $V = (V_1 \mid \ldots \mid V_n)$ be a $(n_1, \ldots, n_n)$ dimensional super vector space over a subfield of the complex numbers and let $f = (f_1 \mid \ldots \mid f_n)$ be a super skew symmetric bilinear superform on V. Then the super rank $r = (r_1, \ldots, r_n)$ of f is even and if $r = (2k_1, \ldots, 2k_n)$ there is an ordered superbasis for V in which the super matrix of f is the super direct sum of the $((n_1 - r_1) \times (n_1 - r_1), \ldots, (n_n - r_n) \times (n_n - r_n))$ zero super diagonal matrix and $(k_1, \ldots, k_n)$ copies of $2 \times 2$ matrix L, where*

$$L = \begin{pmatrix} 0 & 1 \\ -1 & 0 \end{pmatrix}.$$

*Proof:* Let $\alpha_1^1, \beta_1^1 \ldots \alpha_{k_1}^1, \beta_{k_1}^1; \ldots; \alpha_1^n, \beta_1^n, \ldots \alpha_{k_n}^n, \beta_{k_n}^n$ be super vectors satisfying the conditions (a), (b) and (c) mentioned in the page 207. Let $\left\{\gamma_1^1 \ldots \gamma_{s_1}^1, \ldots, \gamma_1^n \ldots \gamma_{s_n}^n\right\}$ be any ordered super basis for the supersubspace $W_0 = (W_0^1 \mid \ldots \mid W_0^n)$.

Then

$$B = \left\{\alpha_1^1, \beta_1^1 \ldots \alpha_{k_1}^1, \beta_{k_1}^1; \gamma_1^1 \ldots \gamma_{s_1}^1, \ldots, \alpha_1^n, \beta_1^n \ldots \alpha_{k_n}^n, \beta_{k_n}^n; \gamma_1^n \ldots \gamma_{s_n}^n\right\}$$

is an ordered super basis for $V = (V_1 \mid \ldots \mid V_n)$. From (a), (b) and (c) it is clear that the super diagonal matrix of $f = (f_1 \mid \ldots \mid f_n)$ in the ordered super basis B is the super direct sum of $((n_1 - 2k_1) \times (n_1 - 2k_1), \ldots, (n_n - 2k_n) \times (n_n - 2k_n))$ zero super matrix and $(k_1 \mid \ldots \mid k_n)$ copies of $2 \times 2$ matrix

$$L = \begin{pmatrix} 0 & 1 \\ -1 & 0 \end{pmatrix}.$$

Further more the super rank of this matrix, hence super rank of f is $(2k_1, \ldots, 2k_n)$.

Several other properties in this direction can be derived, we conclude this section with a brief description of the groups preserving bilinear super forms.

Let $f = (f_1 \mid \ldots \mid f_n)$ be a bilinear super form on a super vector space $V = (V_1 \mid \ldots \mid V_n)$ and $T_s = (T_1 \mid \ldots \mid T_n)$ be a linear operator on V. We say that $T_s$ super preserves f if

$f(T_s\alpha, T_s\beta) = f(\alpha, \beta)$
i.e., $(f_1(T_1\alpha_1, T_1\beta_1) \mid \ldots \mid f_n(T_n\alpha_n, T_n\beta_n))$
$= (f_1(\alpha_1, \beta_1) \mid \ldots \mid f_n(\alpha_n, \beta_n))$
for all $\alpha = (\alpha_1 \mid \ldots \mid \alpha_n)$ and $\beta = (\beta_1 \mid \ldots \mid \beta_n)$ in V.
For any $T_s$ and f the super function $g = (g_1 \mid \ldots \mid g_n)$ defined by

$g(\alpha, \beta) = (g_1(\alpha_1, \beta_1) \mid \ldots \mid g_n(\alpha_n, \beta_n))$
$= f(T_s\alpha, T_s\beta)$
$= (f_1(T_1\alpha_1, T_1\beta_1) \mid \ldots \mid f_n(T_n\alpha_n, T_n\beta_n));$

is easily seen to be a bilinear super form on V. To say that $T_s$ preserves f, is simple (say) g = f. The identity super operator preserves every bilinear super form. If $S_s$ and $T_s$ are linear operators which preserves f the product $S_sT_s$ also preserves f for

$$f(S_sT_s\alpha, S_sT_s\beta) = f(T_s\alpha, T_s\beta) = f(\alpha, \beta)$$
$$\text{i.e., } (f_1(S_1T_1\alpha_1, S_1T_1\beta_1) \mid \ldots \mid f_n(S_nT_n\alpha_n, S_nT_n\beta_n))$$
$$(f_1(T_1\alpha_1, T_1\beta_1) \mid \ldots \mid f_n(T_n\alpha_n, T_n\beta_n))$$
$$= (f_1(\alpha_1, \beta_1) \mid \ldots \mid f_n(\alpha_n, \beta_n))$$

i.e. the collection of all linear operators which super preserve a given bilinear super form is closed under the formation of product.

We have the following interesting theorem.

**THEOREM 2.2.7:** *Let $f = (f_1 \mid \ldots \mid f_n)$ be a super non degenerate bilinear superform of a finite $(n_1, \ldots, n_n)$ dimensional super vector space $V = (V_1 \mid \ldots \mid V_n)$. The set of all super linear operators on V which preserves f is a group called the super group under the operation of composition.*

*Proof:* Let $(G_1 \mid \ldots \mid G_n) = G$ be the super set of all super linear operators preserving the bilinear superform $f = (f_1 \mid \ldots \mid f_n)$ i.e. $G_i$ is a set of linear operators on $V_i$ which preserve $f_i$ ; i = 1, 2, …, n. We see the super identity operator is in $(G_1 \mid \ldots \mid G_n) = G$ and that when ever $S_s$ and $T_s$ are in G the super composition $S_s$ o $T_s = (S_1T_1 \mid \ldots \mid S_nT_n)$ is also in G i.e. each $S_iT_i$ is in $G_i$, i = 1, 2, …, n. Using the fact that f is super non degenerate we shall prove that any super operator $T_s$ in G is super invertible i.e. each component is invertible in every $G_i$, and $T_s^{-1}$ is also in G. Suppose $T_s = (T_1 \mid \ldots \mid T_n)$ preserves $f = (f_1 \mid \ldots \mid f_n)$. Let $\alpha = (\alpha_1 \mid \ldots \mid \alpha_n)$ be a super vector in the super null space of $T_s$. Then for any $\beta = (\beta_1 \mid \ldots \mid \beta_n)$ in V we have

$$f(\alpha, \beta) = f(T_s\alpha, T_s\beta) = f(0, T_s\beta) = (0 \mid \dots \mid 0)$$

i.e.,
$$(f_1(\alpha_1, \beta_1) \mid \dots \mid f_n(\alpha_n, \beta_n))$$
$$= (f_1(T_1\alpha_1, T_1\beta_1) \mid \dots \mid f_n(T_n\alpha_n, T_n\beta_n))$$
$$= (f_1(0, T_1\beta_1) \mid \dots \mid f_n(0, T_n\beta_n)) = (0 \mid \dots \mid 0).$$

Since f is super non degenerate $\alpha = (0 \mid \dots \mid 0)$. Thus $T_s = (T_1 \mid \dots \mid T_n)$ is super invertible i.e. each $T_j$ is invertible; j = 1, 2, …, n Clearly $T_s^{-1} = (T_1^{-1} \mid \dots \mid T_n^{-1})$ also super preserves $f = (f_1 \mid \dots \mid f_n)$ for

$$f(T_s^{-1}\alpha, T_s^{-1}\beta) = f(T_s T_s^{-1}\alpha, T_s T_s^{-1}\beta) = f(\alpha, \beta)$$
$$\text{i.e. } (f_1(T_1^{-1}\alpha_1, T_1^{-1}\beta_1) \mid \dots \mid f_n(T_n^{-1}\alpha_n, T_n^{-1}\beta_n))$$
$$= f_1(T_1T_1^{-1}\alpha_1, T_1T_1^{-1}\beta_1) \mid \dots \mid f_n(T_nT_n^{-1}\alpha_n, T_nT_n^{-1}\beta_n)$$
$$= (f_1(\alpha_1, \beta_1) \mid \dots \mid f_n(\alpha_n, \beta_n))$$

Hence the theorem.

If $f = (f_1 \mid \dots \mid f_n)$ is a super non degenerate bilinear superform on the finite $(n_1, \dots, n_n)$ super space V, then each ordered super basis $B = (B_1 \dots B_n)$ for V determines a super group of super diagonal matrices super preserving f. The set of all super diagonal matrices $[T_s]_B = ([T_1]_{B_1} \mid \dots \mid [T_n]_{B_n})$ where $T_s$ is a linear operator preserving f will be a super group under the super diagonal matrix multiplication. There is another way of description of these super group of matrices.

$$A = \left(\begin{array}{c|c|c|c} A_1 & 0 & & 0 \\ \hline 0 & A_2 & & 0 \\ \hline & & & \\ \hline 0 & 0 & & A_n \end{array}\right) = [f]_B = ([f_1]_{B_1} \mid \dots \mid [f_n]_{B_n})$$

so that if $\alpha$ and $\beta$ are super vectors in V with respective coordinate super matrices X and Y relative to $B = (B_1 \mid \dots \mid B_n)$, we shall have

$$f(\alpha, \beta) = (f_1(\alpha_1, \beta_1) \mid \dots \mid f_n(\alpha_n, \beta_n)) = X^tAY.$$

Suppose $T_s = [T_1 \mid \dots \mid T_n]$ is a linear operator on $V = (V_1 \mid \dots \mid V_n)$ and

$$M = [T_s]_B = ([T_1]_{B_1} \mid \ldots \mid [T_n]_{B_n}).$$

Then

$$f(T_s\alpha, T_s\beta) = (f_1(T_1\alpha_1, T_1\beta_1) \mid \ldots \mid f_n(T_n\alpha_n, T_n\beta_n))$$
$$= (MX)^t A(MY);$$

M, A are super diagonal matrices

$$f(T_s\alpha, T_s\beta) = X^t(M^tAM)Y$$

$$= X^t \left(\begin{array}{c|c|c|c} M_1^tA_1M_1 & 0 & & 0 \\ \hline 0 & M_2^tA_2M_2 & & 0 \\ \hline & & & \\ \hline 0 & 0 & & M_n^tA_nM_n \end{array}\right) Y.$$

Thus $T_s$ preserves f if and only if $M^tAM = A$ i.e. if and only if each $T_i$ preserves $f_i$ i.e. $M_i^t A_i M_i = A_i$ for i = 1, 2, …, n. In the super diagonal matrix language the result can be stated as if A is an invertible super diagonal matrix of $(n_1 \times n_1, \ldots, n_n \times n_n)$ order

$$A = \left(\begin{array}{c|c|c|c} A_1 & 0 & & 0 \\ \hline 0 & A_2 & & 0 \\ \hline & & & \\ \hline 0 & 0 & & A_n \end{array}\right)$$

i.e., each $A_i$ is invertible and $A_i$ is $(n_i \times n_i)$ matrix, i = 1, 2, …, n. $M^tAM = A$ is a super group under super diagonal matrix multiplication. If

$$A = [f]_B = ([f_1]_{B_1} \mid \ldots \mid [f_n]_{B_n});$$

i.e., $$\left(\begin{array}{c|c|c|c} A_1 & 0 & & 0 \\ \hline 0 & A_2 & & 0 \\ \hline & & & \\ \hline 0 & 0 & & A_n \end{array}\right) = \left(\begin{array}{c|c|c|c} [f_1]_{B_1} & 0 & & 0 \\ \hline 0 & [f_2]_{B_2} & & 0 \\ \hline & & & \\ \hline 0 & 0 & & [f_n]_{B_n} \end{array}\right)$$

then M is in this super group of super diagonal matrices if and only if M = $[T_s]_B$

$$\text{i.e.,} \left(\begin{array}{c|c|c|c} M_1 & 0 & & 0 \\ \hline 0 & M_2 & & 0 \\ \hline & & & \\ \hline 0 & 0 & & M_n \end{array}\right) = \left(\begin{array}{c|c|c|c} [T_1]_{B_1} & 0 & & 0 \\ \hline 0 & [T_2]_{B_2} & & 0 \\ \hline & & & \\ \hline 0 & 0 & & [T_n]_{B_n} \end{array}\right)$$

where $T_s$ is a linear operator which preserves f = ($f_1$ | … | $f_n$).

Several properties in this direction can be derived by the reader. We now just prove the following theorem.

**THEOREM 2.2.8:** *Let V = ($V_1$ | … | $V_n$) be a ($n_1$, …, $n_n$) dimensional super vector space over the field of complex numbers and let f = ($f_1$ | … | $f_n$) be a super non-degenerate symmetric bilinear super form on V. Then the super group preserving f is super isomorphic to the complex orthogonal super group O(n, c) = (O($n_1$, c) | … | O($n_n$, c)) where each O($n_i$, c) is a complex orthogonal group preserving $f_i$, i = 1, 2, …, n.*

*Proof:* By super isomorphism between two super groups we mean only isomorphism between the component groups which preserves the group operation.

Let G = ($G_1$ | … | $G_n$) be the super group of linear operators on V = ($V_1$ | … | $V_n$) which preserves the bilinear super form f = ($f_1$ | … | $f_n$). Since f is both super symmetric and super nondegenerate we have an ordered super basis B = ($B_1$ | … | $B_n$). for V in which f is represented by ($n_1 \times n_1$, …, $n_n \times n_n$) super diagonal identity matrix. i.e. each symmetric nondegenerate bilinear form $f_i$ is represented by a $n_i \times n_i$ identity matrix for every i. Therefore the linear operator $T_i$ of $T_s$ preserves f if and only if its matrix in the basis $B_i$ is a complex orthogonal matrix.

Hence $T_i \to [T_i]_{B_i}$ for every i is an isomorphism of $G_i$ onto $O(n_i, c)$; i = 1, 2,…, n. Thus $T_s \to [T_s]_B$ is a super isomorphism of G = ($G_1$ | … | $G_n$) onto O(n, c) = (O($n_1$, c) | … | O($n_n$, c))

We state the following theorem the proof is left as an exercise for the reader.

**THEOREM 2.2.9:** *Let V = ($V_1$ | … | $V_n$) be a ($n_1$, …, $n_n$) dimensional super vector space over the field of real numbers and let f = ($f_1$ | … | $f_n$) be a super non generate bilinear super form on V. Then the super group preserving f is isomorphic to a ($n_1 \times n_1$, …, $n_n \times n_n$) super pseudo orthogonal super group.*

Now we give by an example of a pseudo orthogonal super group.

***Example 2.2.1:*** Let f = ($f_1$ | … | $f_n$) be a symmetric bilinear superform on $(R^{n_1} \,|\, \ldots \,|\, R^{n_n})$ with a quadratic super form q = ($q_1$ | … | $q_n$);

$$q = (x_1, \ldots, x_n) = (q_1(x_1^1 \ldots x_{n_1}^1) \,|\, \ldots \,|\, q_n(x_1^n \ldots x_{n_n}^n))$$

$$= \left( \sum_{j_1=1}^{p_1} (x_{j_1}^1)^2 - \sum_{j_1=p_1+1}^{n_1} (x_{j_1}^1)^2 \,|\, \ldots \,|\, \sum_{j_n=1}^{p_n} (x_{j_n}^n)^2 - \sum_{j_n=p_n+1}^{n_n} (x_{j_n}^n)^2 \right).$$

Then f is a super non degenerate and has super signature

$$2p - n = (2p_1 - n_1 \,|\, \ldots \,|\, 2p_n - n_n).$$

The super group of superdiagonal matrices preserving a super form of this type will be defined as the pseudo-orthogonal super group (or pseudo super orthogonal group or super pseudo orthogonal group) all of them mean the same structure. When each $p_i = n_i$; i = 1, 2, …, n we obtain the super orthogonal group (or orthogonal super group O(n, R) = (O($n_1$, R) | … | O($n_n$, R)) as a particular case of pseudo f orthogonal super group.

## 2.3 Applications

Now we proceed onto give the applications of super matrices, super linear algebras and super vector spaces.

In this section we indicate some of the main applications of super linear algebra / super linear vector spaces / super matrices. For more literature about super matrices please refer [17]. Super linear algebra and super vector spaces have been defined for the first time in this book.

The two main applications we wish to give about these in Markov process and in Leontief economic models.

We first define the new notion of super Markov chain or super Markov process.

A Markov process consists of a set of objects and a set of states such that

i) at any given time each object must be in a state (distinct objects need not be in distinct states).
ii) the probability that an object moves from one state to another (which may be the same as the first state) in one time period depends only on those two states.

If the number of states is finite or countably infinite, the Markov process is a Markov chain. A finite Markov chain is one having a finite number of states we denote the probability of moving from state i to state j in one time period by $p_{ij}$. For an N-state Markov chain where N is a fixed positive integer, the N × N matrix P = $(p_{ij})$ is the stochastic or transition matrix associated with the process.

Denote the $n^{th}$ power of P by $P^n = (p_{ij}^{(n)})$. If P is stochastic then $p_{ij}^{(n)}$ represents the probability that an object moves from state i to state j in n time period it follows that $P^n$ is also a stochastic matrix. Denote the proportion of objects in state i at the end of $n^{th}$ time period by $x_i^{(n)}$ and designate $X^{(n)} \equiv [x_1^{(n)}, \ldots, x_N^{(n)}]$ the distribution super vector for the end of the $n^{th}$ time period.

Accordingly,

$$X^{(0)} = [x_1^{(0)}, \ldots, x_N^{(0)}]$$

represents the proportion of objects in each state at the beginning of the process. $X^{(n)}$ is related to $X^0$ by the equation $X^{(n)} = X^{(0)}P^n$

A stochastic matrix P is erogodic if $\lim_{n \to P^n} p_{ij}^{(n)}$ exists that is if each $p_{ij}^{(n)}$ has a limit as $n \to \infty$. We denote the limit matrix necessarily a stochastic matrix by L. The components of $X^{(\infty)}$ defined by the equation $X^{(\infty)} = X^{(0)}L$ are the limiting state distributions and represent the approximate proportions of objects in various states of a Markov chain after a large number of time periods. Now we define 3 types of Markov chains using 2 types of stochastic or transitive matrix.

Suppose we have some p sets $S_1, \ldots, S_p$ of N objects and a p set of states such that

i) at any given time each set of p objects one object taken from each of the p sets $S_1, \ldots, S_p$ must be in a p-state which denotes at a time, p objects state are considered (or under consideration)

ii) The probability that a set of p objects moves from one to another state in one time period depends only on these two states.

Thus as in case of Markov process these p sets integral numbers of time periods past the moment when the process is started represents the stages of the process, may be finite or infinite. If the number of p set states is finite or countably finite we call that the Markov super row chain i.e. a finite Markov super row chain is one having a finite p set (p-tuple) number of states.

For a N-state Markov super p-row chain we have an associated p-row super N × N square matrix $P = (P_1 \mid \ldots \mid P_p)$ where each $P_t = [p_{ij}^t]$ is the N × N stochastic or transition matrix associated with the process for t = 1, 2, …, p. Thus P =

$(P_1 \mid \ldots \mid P_p) = [[p^1_{ij}] \mid \ldots \mid [p^p_{ij}]]$ is called the stochastic super row square matrix or transition super row square matrix. Necessarily the elements in each row of $P_t$ sum to unity, each $P_t$ is distinct from $P_s$ in its entries if $t \neq s$, $1 \leq t, s \leq p$. Thus we have an N-state p sets of Markov chain defined as super p-row Markov chain or p-row super Markov chain or p-Markov super row chain (all mean one and the same model).

We give an example of a super 5-row Markov chain with two states.

$$P = \left[\begin{array}{cc|cc|cc|cc|cc} 0.19 & 0.81 & 0.31 & 0.69 & 0.09 & 0.91 & 0.18 & 0.82 & 0.73 & 0.27 \\ 0.92 & 0.08 & 0.23 & 0.77 & 0.87 & 0.13 & 0.92 & 0.08 & 0.50 & 0.50 \end{array}\right]$$

$$= (P_1 \mid \ldots \mid P_5) = \left[(p^1_{ij}) \mid (p^2_{ij}) \mid \ldots \mid (p^5_{ij})\right]$$

where the study concerns the economic stability as state 1 of 5 countries and economic depression as state 2 for the same five countries.

Thus this is modeled by the two state super Markov 5-row chain having the super row transition matrix $P = (P_1 \mid \ldots \mid P_p)$. The $n^{th}$ power of a super p-row matrix P is denoted by $P^n = [(p^1_{ij})^n \mid \ldots \mid (p^p_{ij})^n]$.

Denote the proportion of p objects in state i at the end of the $n^{th}$ time period by $x^n_i$ and designate

$$X^{(n)} = [(x^1_1)^{(n)} \ldots (x^1_N)^{(n)} \mid (x^2_1)^{(n)} \ldots (x^2_N)^{(n)} \mid \ldots \mid (x^p_1)^n \ldots (x^p_N)^{(n)}]$$
$$= [X^{(n)}_1 \mid \ldots \mid X^{(n)}_p],$$

the distribution super row vector for the end of the $n^{th}$ time period. Accordingly

$$X^0 = [(x^1_1)^{(0)} \ldots (x^1_N)^{(0)} \mid \ldots \mid [(x^p_1)^{(0)}, \ldots (x^p_N)^{(0)}]$$
$$= [X^{(0)}_1 \mid \ldots \mid X^{(0)}_p]$$
$$\text{i.e., } X^n = X^0 P^n$$
$$\text{i.e. } [X^{(n)}_1 \mid \ldots \mid X^{(n)}_p] = [X^0 P^n_1 \mid \ldots \mid X^0_p P^n_p].$$

A stochastic super row square matrix $P = (P_1 \mid \dots \mid P_n)$ is super ergodic if $\lim_{n\to\infty} P^n \left[\lim_{n\to\infty} P_1^n \mid \dots \mid \lim_{n\to\infty} P_p^n\right]$ exists i.e. if each $(p_{ij}^t)^{(n)}$ has a limit as $n \to \infty$; t = 1, 2, …, p. We denote the limit matrix, necessarily a super row matrix by $L = (L_1 \mid \dots \mid L_p)$. The components of $X^\infty$ defined by the equation

$$X^\infty = X^{(0)}L;\ (X_1^\infty \mid \dots \mid X_p^\infty) = (X_1^{(0)}L_1 \mid \dots \mid X_p^{(0)}L_p)$$

are the limiting super state distribution and represent the approximate proportions.

Thus we see when we have a same set of states to be analyzed regarding p distinct sets of object the Markov super row chain plays a vital role. This method also is helpful in simultaneous comparisons. Likewise when we want to study the outcome of a training program in 5 centres each taking into considerations only 3 states then we can formulate a Markov super row chain with N = 3 and p = 5.

Now when the number of states are the same for all the p sets of objects we can use this Markov super row chain model.

However when we have some p sets of sets of objects and the number of states also vary from time to time among the p sets. Then we have different transition matrix. i.e. if $S_1, \dots, S_p$ are the p sets of objects then each $S_i$ has a $N_i \times N_i$ transition matrix. For the $(N_1, \dots, N_p)$ state Markov chain; if $P_t$ denotes the $[(p_{ij}^{(t)})]$ stochastic matrix then $p_{ij}^n$ represents that, an object moves from state i to state j in $n_t$ time period, this is true for t = 1, 2, …, p.

Thus the matrix which represented the integrated model of the p sets of $S_1, \dots, S_p$ is given by a super diagonal matrix

$$P = \left(\begin{array}{c|c|c|c} P_1 & 0 & & 0 \\ \hline 0 & P_2 & & 0 \\ \hline & & & \\ \hline 0 & 0 & & P_n \end{array}\right)$$

where each $P_t$ is a $N_t \times N_t$ matrix i.e. $P_t = (p_{ij}^t)$.

This is true for every t = 1, 2, …, p. Thus for a $(N_1, \ldots, N_p)$ state Markov chain we have the super diagonal square matrix, or a mixed super diagonal square matrix or a super diagonal square matrix.

Hence

$$P^n = \left(\begin{array}{c|c|c|c} P_1^n & 0 & & 0 \\ \hline 0 & P_2^n & & 0 \\ \hline & & & \\ \hline 0 & 0 & & P_p^n \end{array}\right) = \left(\begin{array}{c|c|c|c} (p_{ij}^1)^{(n)} & 0 & & 0 \\ \hline 0 & (p_{ij}^2)^{(n)} & & 0 \\ \hline & & & \\ \hline 0 & 0 & & (p_{ij}^p)^{(n)} \end{array}\right).$$

Denote the proportion of objects in state i at the end of the $n^{th}$ time period by $(x_i^t)^{(n)}$; t = 1, 2, …, p and designate

$$X^{(n)} = [(x_1^1)^{(n)} \ldots (x_{N_1}^1)^{(n)} \mid \ldots \mid (x_1^p)^{(n)} \ldots (x_{N_p}^p)^{(n)}]$$

$$= [X_1^{(n)} \mid \ldots \mid X_p^{(n)}]$$

here $N_i = N_j$ for $i \neq j$ can also occur.

$$X^0 = [(x_1^1)^{(0)} \ldots (x_{N_1}^1)^{(0)} \mid \ldots \mid (x_1^p)^{(0)} \ldots (x_{N_p}^p)^{(0)}]$$

$$= [X_1^{(0)} \mid \ldots \mid X_p^{(0)}]$$

$$X^n = X^0 P^n$$

i.e. $[X_1^{(n)} \mid \ldots \mid X_p^{(n)}]$

$$= [X_1^{(0)} \mid \ldots \mid X_p^{(0)}] \left(\begin{array}{c|c|c|c} P_1^{(n)} & 0 & & 0 \\ \hline 0 & P_2^{(n)} & & 0 \\ \hline & & & \\ \hline 0 & 0 & & P_p^{(n)} \end{array}\right).$$

i.e. $[X_1^{(n)} \mid \ldots \mid X_p^{(n)}]$

$$= \left(\begin{array}{c|c|c|c} X_1^{(0)}P_1^{(n)} & 0 & & 0 \\ \hline 0 & X_2^{(0)}P_2^{(n)} & & 0 \\ \hline & & & \\ \hline 0 & 0 & & X_p^{(0)}P_p^{(n)} \end{array}\right)$$

i.e. each $X_t^{(n)} = X_t^{(0)} P_t^{(n)}$ true for i = 1, 2, …, n.

This Markov chain model will be know as the super diagonal Markov chain model or equally Markov chain super diagonal model.

Interested reader can apply this model to real world problems and determine the solution. One of the merits of this model is when the expert wishes to study a p-tuple of ($N_1$, …, $N_p$) states $p \geq 2$ this model is handy. Clearly when p = 1 we get the usual Markov chain with $N_1$ state.

## Leontief economic super models

Matrix theory has been very successful in describing the interrelations between prices, outputs and demands in an economic model. Here we just discuss some simple models based on the ideals of the Nobel-laureate Wassily Leontief. Two types of models discussed are the closed or input-output model and the open or production model each of which assumes some economic parameter which describe the inter relations between the industries in the economy under considerations. Using matrix theory we evaluate certain parameters.

The basic equations of the input-output model are the following:

$$\begin{pmatrix} a_{11} & a_{12} & \cdots & a_{1n} \\ a_{21} & a_{22} & \cdots & a_{2n} \\ \vdots & \vdots & & \vdots \\ a_{n1} & a_{n2} & \cdots & a_{nn} \end{pmatrix} \begin{pmatrix} p_1 \\ p_2 \\ \vdots \\ p_n \end{pmatrix} = \begin{pmatrix} p_1 \\ p_2 \\ \vdots \\ p_n \end{pmatrix}$$

each column sum of the coefficient matrix is one

i. $p_i \geq 0$, i = 1, 2, …, n.
ii. $a_{ij} \geq 0$, i , j = 1, 2, …, n.
iii. $a_{ij} + a_{2j} + \ldots + a_{nj} = 1$

for j = 1, 2 , …, n.

$$p = \begin{pmatrix} p_1 \\ p_2 \\ \vdots \\ p_n \end{pmatrix}$$

are the price vector. A = $(a_{ij})$ is called the input-output matrix

$$Ap = p \text{ that is, } (I - A)\, p = 0.$$

Thus A is an exchange matrix, then Ap = p always has a nontrivial solution p whose entries are nonnegative. Let A be an exchange matrix such that for some positive integer m, all of the entries of $A^m$ are positive. Then there is exactly only one linearly independent solution of (I – A) p = 0 and it may be chosen such that all of its entries are positive in Leontief open production model.

In contrast with the closed model in which the outputs of k industries are distributed only among themselves, the open model attempts to satisfy an outside demand for the outputs. Portions of these outputs may still be distributed among the industries themselves to keep them operating, but there is to be some excess some net production with which to satisfy the outside demand. In some closed model, the outputs of the industries were fixed and our objective was to determine the prices for these outputs so that the equilibrium condition that expenditures equal incomes was satisfied.

$x_i$ = monetary value of the total output of the $i^{th}$ industry.

$d_i$ = monetary value of the output of the $i^{th}$ industry needed to satisfy the outside demand.

$\sigma_{ij}$ = monetary value of the output of the $i^{th}$ industry needed by the $j^{th}$ industry to produce one unit of monetary value of its own output.

With these qualities we define the production vector.

$$x = \begin{pmatrix} x_1 \\ x_2 \\ \vdots \\ x_k \end{pmatrix}$$

the demand vector

$$d = \begin{pmatrix} d_1 \\ d_2 \\ \vdots \\ d_k \end{pmatrix}$$

and the consumption matrix,

$$C = \begin{pmatrix} \sigma_{11} & \sigma_{12} & \cdots & \sigma_{1k} \\ \sigma_{21} & \sigma_{22} & \cdots & \sigma_{2k} \\ \vdots & \vdots & & \vdots \\ \sigma_{k1} & \sigma_{k2} & \cdots & \sigma_{kk} \end{pmatrix}.$$

By their nature we have

$$x \geq 0, \; d \geq 0 \text{ and } C \geq 0.$$

From the definition of $\sigma_{ij}$ and $x_j$ it can be seen that the quantity

$$\sigma_{i1} x_1 + \sigma_{i2} x_2 + \ldots + \sigma_{ik} x_k$$

is the value of the output of the $i^{th}$ industry needed by all k industries to produce a total output specified by the production vector x.

Since this quantity is simply the $i^{th}$ entry of the column vector Cx, we can further say that the $i^{th}$ entry of the column vector x – Cx is the value of the excess output of the $i^{th}$ industry available to satisfy the outside demand. The value of the outside demand for the output of the $i^{th}$ industry is the $i^{th}$ entry of the demand vector d; consequently; we are led to the following equation:

$$x - Cx = d \text{ or}$$
$$(I - C)\, x = d$$

for the demand to be exactly met without any surpluses or shortages. Thus, given C and d, our objective is to find a production vector $x \geq 0$ which satisfies the equation $(I - C)x = d$.

A consumption matrix C is said to be productive if $(1 - C)^{-1}$ exists and $(1 - C)^{-1} \geq 0$.

A consumption matrix C is productive if and only if there is some production vector $x \geq 0$ such that $x > Cx$.

A consumption matrix is productive if each of its row sums is less than one. A consumption matrix is productive if each of its column sums is less than one.

Now we will formulate the Smarandache analogue for this, at the outset we will justify why we need an analogue for those two models.

Clearly, in the Leontief closed Input – Output model,
$p_i$ = price charged by the $i^{th}$ industry for its total output in reality need not be always a positive quantity for due to competition to capture the market the price may be fixed at a loss or the demand for that product might have fallen down so badly so that the industry may try to charge very less than its real value just to market it.

Similarly $a_{ij} \geq 0$ may not be always be true. Thus in the Smarandache Leontief closed (Input-Output) model (S-Leontief closed (Input-Output) model) we do not demand $p_i \geq 0$, $p_i$ can be negative; also in the matrix $A = (a_{ij})$,

$$a_{1j} + a_{2j} + \ldots + a_{kj} \neq 1$$

so that we permit $a_{ij}$'s to be both positive and negative, the only adjustment will be we may not have (I – A) p = 0, to have only one linearly independent solution, we may have more than one and we will have to choose only the best solution.
As in this complicated real world problems we may not have in practicality such nice situation. So we work only for the best solution.

Here we introduce a input-output model which has some p number of input-output matrix each of same order say n × n functioning simultaneously. We shall call such models as input – output super row matrix models and describe how it functions. Suppose we have p number of n × n input output matrix given by the super row matrix $A = [A_1 \mid \ldots \mid A_n]$ where each $A_i$ is a n × n input output matrix which are distinct.

$$A = [A_1 \mid \ldots \mid A_n]$$

$$= \left( \begin{pmatrix} a^1_{11} & \ldots & a^1_{1n} \\ a^1_{21} & \ldots & a^1_{2n} \\ \vdots & & \vdots \\ a^1_{n1} & \ldots & a^1_{nn} \end{pmatrix} \Bigg| \ldots \Bigg| \begin{pmatrix} a^p_{11} & \ldots & a^p_{1n} \\ a^p_{21} & \ldots & a^p_{2n} \\ \vdots & & \vdots \\ a^p_{n1} & \ldots & a^1_{nn} \end{pmatrix} \right)$$

where $a^t_{ij} + a^t_{2j} + \ldots + a^t_{nj} = 1$; t = 1, 2, …, p and j = 1, 2, …, n. Suppose

$$P = \left( \begin{array}{c|c|c|c} p^1_1 & p^2_1 & \cdots & p^P_1 \\ \vdots & \vdots & & \vdots \\ p^1_n & p^2_n & \cdots & p^P_n \end{array} \right) = [P_1 \mid \ldots \mid P_p]$$

be the super column price vector then

A * P = P, the (product) * is defined as A * P = P that is

$$[A_1P_1 \mid \ldots \mid A_pP_p] = [P_1 \mid \ldots \mid P_p]$$

$$A * P = P$$

that is

$(I – A) P = (0 | … | 0)$
i.e., $((I – A_1) P_1 | … | (I – A_p) P_p) = (0 | … | 0)$.

Thus A is an super-row square exchange matrix, then AP = P always has a row column vector solution P whose entries are non negative.

Let $A = [A_1 | … | A_n]$ be an exchange super row square matrix such that for some positive integer m all the entries of $A^m$ i.e. entries of each $A_t^m$ are positive for m; m = 1, 2, …, p. Then there is exactly only one linearly independent solution of

$(I – A) P = (0 | … | 0)$
i.e., $((I – A_1) P_1 | … | (I – A_p) P_p) = (0 | … | 0)$

and it may be choosen such that all of its entries are positive in Leontief open production sup model.

Note this super model yields easy comparison as well as this super model can with different set of price super column vectors and exchange super row matrix find the best solution from the p solutions got from the relation

$(I – A) P = (0 | … | 0)$
i.e., $((I – A_1) P_1 | … | (I – A_p) P_p) = (0 | … | 0)$ .

Thus this is also an added advantage of the model. It can study simultaneously p different exchange matrix with p set of price vectors for different industries to study the super interrelations between prices, outputs and demands simultaneously.

Suppose one wants to make a study of interrelation between prices, outputs and demands in an industry with different types of products with different exchange matrix and hence different set of price vectors or of many different industries with same type of products its interrelation between prices, outputs and demands in different locations of the country were the economic status and the education status vary in different locations, how to make a single model to study the situation. In both the cases one can make use of the super input-output model the relation matrix which is a input-output super diagonal mixed square matrix, which will be described presently.

The exchange matrix with p distinct economic models is used to describe the interrelations between prices, outputs and demands. Then the related matrix A will be a super diagonal mixed square matrix

$$A = \left(\begin{array}{c|c|c|c} A_1 & 0 & & 0 \\ \hline 0 & A_2 & & 0 \\ \hline & & & \\ \hline 0 & 0 & & A_p \end{array}\right)$$

$A_1, \ldots, A_p$ are the exchange matrices describing the p-economic models. Now A acts as integrated models in which all the p entities function simultaneously. Now any price vector P will be a super mixed column matrix

$$P = \left(\begin{array}{c} P_1 \\ \hline \vdots \\ \hline P_p \end{array}\right)$$

where each

$$P_t = \left(\begin{array}{c} p_1^t \\ \hline \vdots \\ \hline p_{n_t}^t \end{array}\right);$$

for t = 1, 2, …, p.

Here each $A_t$ is a $n_t \times n_t$ exchange matrix; t = 1, 2, …, p. AP = P is given by

$$A = \left(\begin{array}{c|c|c|c} A_1 & 0 & & 0 \\ \hline 0 & A_2 & & 0 \\ \hline & & & \\ \hline 0 & 0 & & A_p \end{array}\right),$$

$$P = \left(\begin{array}{c} P_1 \\ \hline \vdots \\ \hline P_p \end{array}\right)$$

$$AP = \left(\begin{array}{c|c|c|c} A_1P_1 & 0 & & 0 \\ \hline 0 & A_2P_2 & & 0 \\ \hline & & & \\ \hline 0 & 0 & & A_pP_p \end{array}\right) = \left(\begin{array}{c} P_1 \\ \hline \vdots \\ \hline P_p \end{array}\right)$$

i.e. $A_tP_t = P_t$ for every t = 1, 2, …, p. i.e.

$$\left(\begin{array}{c|c|c|c} (I_1 - A_1)P_1 & 0 & & 0 \\ \hline 0 & (I_2 - A_2)P_2 & & 0 \\ \hline & & & \\ \hline 0 & 0 & & (I_n - A_n)P_n \end{array}\right)$$

$$= \left(\begin{array}{c|c|c|c} 0 & 0 & & 0 \\ \hline 0 & 0 & & 0 \\ \hline & & & \\ \hline 0 & 0 & & 0 \end{array}\right).$$

Thus AP = P has a nontrivial solution

$$P = \left(\begin{array}{c} P_1 \\ \hline \vdots \\ \hline P_p \end{array}\right)$$

whose entries in each $P_t$ are non negative; $1 \leq t \leq p$.
Let A be the super exchange diagonal mixed square matrix such that for some p-tuple of positive integers m = ($m_1$, …, $m_p$), $A_t^{m_t}$

is positive; $1 \leq t \leq p$. Then there is exactly only one linearly independent solution;

$$(I - A)P = \left(\begin{array}{c|c|c|c} 0 & 0 & & 0 \\ \hline 0 & 0 & & 0 \\ \hline & & & \\ \hline 0 & 0 & & 0 \end{array}\right)$$

and it may be choosen such that all of its entries are positive in Leontief open production super model.

Next we proceed on the describe super closed row model (or row closed super model) as the super closed model (or closed super model).

Here we have p sets of K industries which are distributed among themselves i.e. the first set with K industries distributed among themselves, the second set with some K industries distributed among themselves and so on and the p set with some K industries distributed among themselves. It may be that some industries are found in more than one set and some industries in one and only one set and some industries in all the p sets. This open super row model which we choose to call as, when p sets of K industries get distributed among themselves attempts to satisfy an outside demand for outputs. Portions of these outputs may still be distributed among the industries themselves to keep them operating, but there is to be some excess some net production with which they satisfy the outside demand. In some super closed row models the outputs of the industries in those sets which they belong to were fixed and our objective was to determine sets of prices for these outputs so that the equilibrium condition that expenditure equal income was satisfied for each of the p sets individually.

Thus we will have

$x_i^t$ = monetary value of the total output of the $i^{th}$ industry in the $t^{th}$ set $1 \leq i \leq K$ and $1 \leq t \leq p$.

$d_i^t$ = monetary value of the output of the $i^{th}$ industry of the $t^{th}$ set needed to satisfy the outside demand, $1 \le t \le p$, I = 1, 2, …, K.

$\sigma_{ij}^t$ = monetary value of the output of the $i^{th}$ industry needed by the $j^{th}$ industry of the $t^{th}$ set to produce one unit of monetary value of its own output, $1 \le i \le K$; $1 \le t \le p$.

With these qualities we define the production super column vector

$$X = \begin{pmatrix} X_1 \\ \hline \vdots \\ \hline X_t \\ \hline \vdots \\ \hline X_p \end{pmatrix} = \begin{pmatrix} x_1^1 \\ \vdots \\ x_K^1 \\ \hline \vdots \\ \hline x_1^p \\ \vdots \\ x_K^p \end{pmatrix}.$$

The demand column super vector

$$d = \begin{pmatrix} d_1 \\ \hline \vdots \\ \hline d_t \\ \hline \vdots \\ \hline d_p \end{pmatrix} = \begin{pmatrix} d_1^1 \\ \vdots \\ d_K^1 \\ \hline \vdots \\ \hline d_1^P \\ \vdots \\ d_K^P \end{pmatrix}$$

and the consumption super row matrix $C = (C_1 \mid \ldots \mid C_p)$

$$= \left\{ \begin{pmatrix} \sigma^1_{11} & \sigma^1_{12} & \cdots & \sigma^1_{1K} \\ \sigma^1_{21} & \sigma^1_{22} & \cdots & \sigma^1_{2K} \\ \vdots & \vdots & & \vdots \\ \sigma^1_{K1} & \sigma^1_{K2} & \cdots & \sigma^1_{KK} \end{pmatrix} \Big| \ldots \Big| \begin{pmatrix} \sigma^p_{11} & \sigma^p_{12} & \cdots & \sigma^p_{1K} \\ \sigma^p_{21} & \sigma^p_{22} & \cdots & \sigma^p_{2K} \\ \vdots & \vdots & & \vdots \\ \sigma^p_{K1} & \sigma^p_{K2} & \cdots & \sigma^p_{KK} \end{pmatrix} \right\}.$$

By their nature we have

$$X \geq \begin{pmatrix} 0 \\ \hline \vdots \\ \hline 0 \end{pmatrix} ; d > \begin{pmatrix} 0 \\ \hline \vdots \\ \hline 0 \end{pmatrix} \text{ and } C > (0 \mid \ldots \mid 0).$$

For the t$^{th}$ set from the definition of $\sigma^t_{ij}$ and $x^t_j$ it can be seen that the quantity

$$\sigma^t_{i1} x^t_1 + \sigma^t_{i2} x^t_2 + \ldots + \sigma^t_{iK} x^t_K$$

is the value of the i$^{th}$ industry needed by all the K industries (of the set t) to produce a total output specified by the production vector $X_t$. Since this quantity is simply the i$^{th}$ entry of the column vector $C_tX_t$ we can further say that the i$^{th}$ entry of the column vector $X_t – X_tC_t$ is the value of the excess output of the i$^{th}$ industry (from the t$^{th}$ set) available to satisfy the outside demand.

The value of the outside demand for the output of the i$^{th}$ industry (from the t$^{th}$ set) is the i$^{th}$ entry of the demand vector $d_t$; consequently we are lead to the following equation for the t$^{th}$ set $X_t – C_tX_t = d_t$ or $(I – C_t)X_t = d_t$ for the demand to be exactly met without any surpluses or shortages. Thus given $C_t$ and $d_t$ our objective is to find a production vector $X_t \geq 0$ which satisfies the equation

$$(I – C_t)X_t = d_t,$$

so for the all p sets we have the integrated equation to be

$$(I - C)X = d$$
$$\text{i.e., } [(I - C_1)X_1 \mid \dots \mid (I - C_p)X_p]$$
$$= (d_1 \mid \dots \mid d_p) .$$

The consumption super row matrix $C = (C_1 \mid \dots \mid C_p)$ is said to be super productive if

$$(I - C)^{-1} = [(I - C_1)^{-1} \mid \dots \mid (I - C_p)^{-1}]$$

exists and

$$(I - C)^{-1} = [(I - C_1)^{-1} \mid \dots \mid (I - C_p)^{-1}] \geq [0 \mid \dots \mid 0] .$$

A consumption super row matrix is super productive if and only if for some production super vector

$$X = \begin{pmatrix} X_1 \\ \hline \vdots \\ \hline X_n \end{pmatrix} \geq \begin{pmatrix} 0 \\ \hline \vdots \\ \hline 0 \end{pmatrix}$$

such that $X > CX$ i.e. $[X_1 \mid \dots \mid X_p] > [C_1X_1 \mid \dots \mid C_pX_p]$.

A consumption super row matrix is productive if each of its row sums is less than one. A consumption super row matrix is productive if each of its column super sums is less than one. The main advantage of this super model is that one can work with p sets of industries simultaneously provided all the p sets have same number of industries (here K). This super row model will help one to monitor and study the performance of an industry which is present in more than one set and see its functioning in each of the sets. Such a thing may not be possible simultaneously in any other model.

Suppose we have p sets of industries and each set has different number of industries say in the first set output of $K_1$ industries are distributed among themselves. In the second set output of $K_2$ industries are distributed among themselves so on in the $p^{th}$ set output of $K_p$-industries are distributed among themselves the super open model is constructed to satisfy an outside demand for the outputs. Here one industry may find its place in one and only one set or group. Some industries may find its place in several groups. Some industries may find its place in every group. To construct a closed super model to analyze the situation.

Portions of these outputs may still be distributed among the industries themselves to keep them operating, but there is to be some excess some net production with which to satisfy the outside demand.

Let

$X_i^t$ = monetary value of the total output of the $i^{th}$ industry in the $t^{th}$ set (or group).

$d_i^t$ = monetary value of the output of the $i^{th}$ industry of the group t needed to satisfy the outside demand.

$\sigma_{ij}^t$ = monetary value of the output of the $i^{th}$ industry needed by the $j^{th}$ industry to produce one unit monetary value of its own output in the $t^{th}$ set or group, $1 \leq t \leq p$.

With these qualities we define the production super mixed column vector

$$X = \begin{pmatrix} X_1 \\ \hline \vdots \\ \hline X_t \\ \hline \vdots \\ \hline X_p \end{pmatrix} = \begin{pmatrix} x_1^1 \\ \vdots \\ x_{K_1}^1 \\ \hline \vdots \\ \hline x_1^p \\ \vdots \\ x_{K_p}^p \end{pmatrix}$$

and the demand super mixed column vector

$$d = \begin{pmatrix} d_1 \\ \hline \vdots \\ \hline d_p \end{pmatrix} = \begin{pmatrix} d_1^1 \\ \vdots \\ d_{K_1}^1 \\ \hline \vdots \\ \hline d_1^p \\ \vdots \\ d_{K_p}^p \end{pmatrix}$$

and the consumption super diagonal mixed square matrix

$$C = \left(\begin{array}{c|c|c|c} C_1 & 0 & & 0 \\ \hline 0 & C_2 & & 0 \\ \hline & & & \\ \hline 0 & 0 & & C_p \end{array}\right)$$

where

$$C_t = \begin{pmatrix} \sigma_{11}^t & \sigma_{12}^t & \cdots & \sigma_{1K_t}^t \\ \sigma_{21}^t & \sigma_{22}^t & \cdots & \sigma_{2K_t}^t \\ \vdots & \vdots & & \vdots \\ \sigma_{K_t 1}^t & \sigma_{K_t 2}^t & \cdots & \sigma_{K_t K_t}^1 \end{pmatrix};$$

true for t = 1, 2, …, p.

By the nature of the closed model we have

$$X = \begin{pmatrix} X_1 \\ \hline \vdots \\ \hline X_p \end{pmatrix} = \begin{pmatrix} 0 \\ \hline \vdots \\ \hline 0 \end{pmatrix}, \; d = \begin{pmatrix} d_1 \\ \hline \vdots \\ \hline d_p \end{pmatrix} = \begin{pmatrix} 0 \\ \hline \vdots \\ \hline 0 \end{pmatrix}$$

and

$$C = \left(\begin{array}{c|c|c|c} C_1 & 0 & & 0 \\ \hline 0 & C_2 & & 0 \\ \hline & & & \\ \hline 0 & 0 & & C_p \end{array}\right) = \left(\begin{array}{c|c|c|c} 0 & 0 & & 0 \\ \hline 0 & 0 & & 0 \\ \hline & & & \\ \hline 0 & 0 & & 0 \end{array}\right).$$

From the definition of $\sigma_{ij}^{t}$ and $x_j^t$ for every group (set t) it can be seen the quantity $\sigma_{ij}^{t} X_1^t + \ldots \sigma_{i_{K_t}} X_{K_t}^t$ is the value of the output of the i$^{th}$ industry needed by all $K_t$ industries (in the t$^{th}$ group) to produce a total output specified by the production vector $X_t$ ($1 \leq t \leq p$). Since this quantity is simply the i$^{th}$ entry of the super column vector in

$$CX = \left(\begin{array}{c|c|c|c} C_1 & 0 & & 0 \\ \hline 0 & C_2 & & 0 \\ \hline & & & \\ \hline 0 & 0 & & C_p \end{array}\right)_{p\times p} \left(\begin{array}{c} X_1 \\ \hline \vdots \\ \hline X_p \end{array}\right)_{p\times 1}$$

$$= [C_1X_1 \mid \ldots \mid C_pX_p]^t$$

we can further say that the i$^{th}$ entry of the super column vector $X_t - CX_t$ in

$$X - CX = \left(\begin{array}{c} X_1 - \quad C_1X_p \\ \hline \vdots \qquad \vdots \\ \hline X_p - \quad C_pX_p \end{array}\right)$$

is the value of the excess output of the i$^{th}$ industry available to satisfy the output demand.

The value of the outside demand for the output of the i$^{th}$ industry (in t$^{th}$ set / group) is the i$^{th}$ entry of the demand vector $d_t$; consequently we are led to the following equation

$$X_t - C_tX_t = d_t \text{ or } (I_t - C_t)\, X_t = d_t\,, \; (1 \leq t \leq p),$$

for the demand to be exactly met without any surpluses or shortages. Thus given $C_t$ and $d_t$ our objective is to find a production vector $X_t \geq 0$ which satisfy the equation $(I_t - C_t)X_t = d$. The integrated super model for all the p-sets (or groups) is given by X – CX = d i.e.

$$\left(\begin{array}{c} X_1 - C_1X_1 \\ \hline X_2 - C_2X_2 \\ \hline \vdots \\ \hline X_p - C_pX_p \end{array}\right) = \left(\begin{array}{c} d_1 \\ \hline d_2 \\ \hline \vdots \\ \hline d_p \end{array}\right)$$

or

$$\left(\begin{array}{c|c|c|c} (I_1 - C_1) & 0 & & 0 \\ \hline 0 & I_2 - C_2 & & \\ \hline & & & \\ \hline 0 & 0 & & I_p - C_p \end{array}\right) \left(\begin{array}{c} X_1 \\ \hline \vdots \\ \hline X_p \end{array}\right) = \left(\begin{array}{c} d_1 \\ \hline \vdots \\ \hline d_p \end{array}\right)$$

i.e.,

$$\left(\begin{array}{c} (I_1 - C_1)X_1 \\ \hline \vdots \\ \hline (I_p - C_p)X_p \end{array}\right) = \left(\begin{array}{c} d_1 \\ \hline \vdots \\ \hline d_p \end{array}\right)$$

where I is a $K_t \times K_t$ square identity matrix t = 1, 2, …, p.

Thus given C and d our objective is to find a production super column mixed vector

$$X = \left(\begin{array}{c} X_1 \\ \hline \vdots \\ \hline X_p \end{array}\right) \geq \left(\begin{array}{c} 0 \\ \hline \vdots \\ \hline 0 \end{array}\right)$$

which satisfies equation (I – C) X = d

$$\text{i.e.} \quad \left( \begin{array}{c} (I_1 - C_1)X_1 \\ \hline \vdots \\ \hline (I_p - C_p)X_p \end{array} \right) = \left( \begin{array}{c} d_1 \\ \hline \vdots \\ \hline d_p \end{array} \right).$$

A consumption super diagonal matrix C is productive if $(I - C)^{-1}$ exists and i.e.

$$\left( \begin{array}{c|c|c|c} (I_1 - C_1)^{-1} & 0 & & 0 \\ \hline 0 & (I_2 - C_2)^{-1} & & 0 \\ \hline & & & \\ \hline 0 & & & (I_p - C_p)^{-1} \end{array} \right)$$

exists and

$$\left( \begin{array}{c|c|c|c} (I_1 - C_1)^{-1} & 0 & & 0 \\ \hline 0 & (I_2 - C_2)^{-1} & & 0 \\ \hline & & & \\ \hline 0 & & & (I_p - C_p)^{-1} \end{array} \right) \geq$$

$$\left( \begin{array}{c|c|c|c} 0 & 0 & & 0 \\ \hline 0 & 0 & & 0 \\ \hline & & & \\ \hline 0 & & & 0 \end{array} \right).$$

A consumption super diagonal matrix C is super productive if and only if there is some production super vector

$$X = \left( \begin{array}{c} X_1 \\ \hline \vdots \\ \hline X_p \end{array} \right) \geq \left( \begin{array}{c} 0 \\ \hline \vdots \\ \hline 0 \end{array} \right)$$

such that

$$X > CX \text{ i.e. } \begin{pmatrix} X_1 \\ \hline \vdots \\ \hline X_p \end{pmatrix} > \begin{pmatrix} C_1X_1 \\ \hline \vdots \\ \hline C_pX_p \end{pmatrix}.$$

A consumption super diagonal mixed square matrix is productive if each row sum in each of the component matrices is less than one. A consumption super diagonal mixed square matrix is productive if each of its component matrices column sums is less than one.

The main advantage of this system is this model can study different sets of industries with varying strength simultaneously. Further the performance of any industry which is present in one or more group can be studied and also analysed. Such comprehensive and comparative study can be made using these super models.

## Chapter Three

# SUGGESTED PROBLEMS

In this chapter we have given over 160 problems for the reader to understand the subject. Any serious researcher is expected to work out the problems. The complexity of the problems varies.

1. Prove that every m × n simple matrix over the rational Q which is partitioned into a super matrix in the same way is a super vector space over Q.

2. If A = ($A_{ij}$) is the collection of all 4 × 4 matrix with entries from Q all of which are partitioned as

$$\left(\begin{array}{c|ccc} a_{11} & a_{12} & a_{13} & a_{14} \\ a_{21} & a_{22} & a_{23} & a_{24} \\ \hline a_{31} & a_{32} & a_{33} & a_{34} \\ a_{41} & a_{42} & a_{43} & a_{44} \end{array}\right)$$

Prove A is a super vector space over Q.

3. Prove $V = \{(x_1\ x_2 \mid x_3\ x_4\ x_5 \mid x_6\ x_7\ x_8\ x_9) \mid x_i \in R; 1 \le i \le 9\}$ is a super vector space over Q. What is the dimension of V?

4. Let $V = \{(x_1 \mid x_2\ x_3 \mid x_4\ x_5\ x_6) \mid x_i \in R; 1 \le i \le 6\}$ be a super vector space over R. Find dimension of V. Suppose V is a super vector space over Q then what is the dimension of V?

5. Prove

$$V = \left\{ \left( \begin{array}{ccc|cc} a_1 & a_2 & a_3 & a_7 & a_9 \\ a_4 & a_5 & a_6 & a_8 & a_{10} \\ \hline a_{11} & a_{12} & a_{13} & a_{20} & a_{21} \\ a_{14} & a_{15} & a_{16} & a_{22} & a_{23} \\ \hline a_{17} & a_{18} & a_{19} & a_{24} & a_{25} \end{array} \right) \middle| \, a_i \in Q; 1 \le i \le 25 \right\}$$

is a super vector space over Q. Find the dimension of V. Is V a super vector space over R?

6. Let $V = \{(x_1\ x_2 \mid x_3\ x_4 \mid x_5\ x_6\ x_7) \mid x_i \in Q; 1 \le i \le 7\}$ and $W = \{(x_1 \mid x_2\ x_3\ x_4 \mid x_5\ x_6) \mid x_i \in Q, 1 \le i \le 6\}$ be super vector spaces over Q. Define a linear super transformation T from V into W. Find the super null space of T.

7. Let $V = \{(x_1 \mid x_2\ x_3 \mid x_4\ x_5\ x_6 \mid x_7\ x_8) \mid x_i \in Q; 1 \le i \le 8\}$ and $W = \{(x_1\ x_2 \mid x_3 \mid x_4\ x_5 \mid x_6 \mid x_7\ x_8 \mid x_9) \mid x_i \in Q; 1 \le i \le 9\}$ be super vector spaces over Q. Let $T : V \to W$ be defined by $T = (x_1 \mid x_2\ x_3 \mid x_4\ x_5\ x_6 \mid x_7\ x_8) = (x_1 - x_1 \mid x_2 + x_3 \mid x_4 + x_5\ x_6\ x_4 \mid x_7 + x_8 \mid 0\ 0 \mid 0)$. Prove T is a linear super transformation from V into W. Find the super null space of V.

8. Define a different linear transformation $T^1$ from V into W which is different from T defined in problem 7, V and W are taken as given in problem 7. Can a linear super transformation T be defined from V into W so that the super null space of T is just the zero super space?

9. Let V = $\{(x_1\ x_2 \mid x_3\ x_4 \mid x_5\ x_6\ x_7) \mid x_i \in Q;\ 1 \leq i \leq 7\}$ be a super vector space over the field Q. W = $\{(x_1\ x_2\ x_3 \mid x_4\ x_5 \mid x_6\ x_7) \mid x_i \in Q;\ 1 \leq i \leq 7\}$. For any linear super transformation $T_s$ and verify the condition rank $T_s$ + nullity $T_s$ = dim V = 7.

10. Let V = $\{(x_1\ x_2 \mid x_3\ x_4 \mid x_5) \mid x_i \in Q;\ 1 \leq i \leq 5\}$ be a super vector space over Q. W = $\{(x_1\ x_2\ x_3 \mid x_4\ x_5 \mid x_6\ x_7\ x_8 \mid x_9\ x_{10}) \mid x_i \in Q;\ 1 \leq i \leq 10\}$ be a super vector space over Q. Can we have a nontrivial nullity $T_s$; $T_s\ V \to W$ such that rank $T_s$ + nullity $T_s$ = dim V = 5; nullity $T_s \neq 0$.

11. Let V = $\{(x_1\ x_2\ x_3 \mid x_4\ x_5\ x_6\ x_7 \mid x_8\ x_9) \mid x_i \in Q;\ 1 \leq i \leq 9\}$ be a super vector space over Q. W = $\{(x_1\ x_2 \mid x_3\ x_4\ x_5 \mid x_6\ x_7\ x_8) \mid x_i \in Q;\ 1 \leq i \leq 8\}$ a super vector space over Q. Does there exist a linear super transformation $T_s$: V $\to$ W such that nullity $T_s$ = 0? Justify your claim.

12. Let V = $\{(x_1\ x_2 \mid x_3\ x_4\ x_5 \mid x_6\ x_7\ x_8) \mid x_i \in Q;\ 1 \leq i \leq 8\}$ be a super vector space over Q. W = $\{(x_1\ x_2 \mid x_3 \mid x_4\ x_5) \mid x_i \in Q;\ 1 \leq i \leq 5\}$ a super vector space over Q. Does their exist a $T_s$ for which nullity $T_s$ = 0?

13. Let V = $\{(x_1\ x_2 \mid x_3\ x_4\ x_5 \mid x_6\ x_7\ x_8) \mid x_i \in Q;\ 1 \leq i \leq 8\}$ be a super vector space over Q. Find two basis distinct from each other for V which is different from the standard basis.

14. Find a basis for the super vector space

$$V = \left\{ \left( \begin{array}{ccc|cc} a_1 & a_2 & a_3 & a_4 & a_5 \\ \hline a_6 & a_7 & a_8 & a_{15} & a_{16} \\ a_9 & a_{10} & a_{11} & a_{17} & a_{18} \\ a_{12} & a_{13} & a_{14} & a_{19} & a_{20} \end{array} \right) \middle| a_i \in Q;\ 1 \leq i \leq 20 \right\}$$

over Q.

15. Find at least 3 super subspaces of the super vector space

$$V = \left\{ \left( \begin{array}{c|ccc} a_1 & a_7 & a_8 & a_9 \\ a_2 & a_{10} & a_{11} & a_{12} \\ a_3 & a_{13} & a_{14} & a_{15} \\ a_4 & a_{16} & a_{17} & a_{18} \\ \hline a_5 & a_{19} & a_{20} & a_{21} \\ a_6 & a_{22} & a_{23} & a_{24} \end{array} \right) \right.$$

such that $a_i \in Q$; $1 \leq i \leq 24\}$ over Q. Find their dimension show for three other super subspaces $W_1$, $W_2$ and $W_3$ of V we can have $V = W_1 + W_2 + W_3$.

16. Let $V = \{(x_1\ x_2 \mid x_3\ x_4\ x_5 \mid x_6\ x_7\ x_8 \mid x_9\ x_{10}) \mid x_i \in Q; 1 \leq i \leq 10\}$ be a super vector space over the field Q. (1) Find all super subspaces of V. (2) Find two super subspaces $W_1$ and $W_2$ of V such that $W = W_1 \cap W_2$ is not the zero super subspace of V.

17. Let $V = \{(x_1\ x_2\ x_3\ x_4 \mid x_5\ x_6\ x_7\ x_8 \mid x_9\ x_{10}\ x_{11}\ x_{12}) \mid x_i \in Q; 1 \leq i \leq 12\}$ be a super vector space over Q. $W = \{(x_1\ x_2 \mid x_3\ x_4 \mid x_5\ x_6) \mid x_i \in Q; i = 1, 2, \ldots, 6\}$ is a super vector space over Q. Find dimension of SL (V, W).

18. How many super vector subspaces SL (V, W) can be got given V is a super vector space of natural dimension n and W a super vector space of natural dimension m, both defined on the same field F?

19. Given $X = (x_1\ x_2)$ we have only one partition $(x_1 \mid x_2)$. Given $X = (x_1\ x_2\ x_3)$ we have three partitions $(x_1\ x_2 \mid x_3)$, $(x_1 \mid x_2\ x_3)$, $(x_1 \mid x_2 \mid x_3)$. Given $X = (x_1\ x_2\ x_3\ x_4)$ we have $(x_1 \mid x_2 \mid x_3 \mid x_4)$, $(x_1\ x_2 \mid x_3\ x_4)$, $(x_1\ x_2 \mid x_3 \mid x_4)$, $(x_1\ x_2\ x_3 \mid x_4)$ $(x_1 \mid x_2\ x_3\ x_4)$ $(x_1 \mid x_2\ x_3 \mid x_4)$ and $(x_1 \mid x_2 \mid x_3\ x_4)$ seven partitions. Thus given $X = (x_1\ x_2 \ldots x_n)$ how many partitions can we have on X?

20. Let $V = \{(x_1 \mid x_2\ x_3\ x_4\ x_5 \mid x_6\ x_7) \mid x_i \in Q; 1 \leq i \leq 7\}$ be a super vector space over Q. Find SL (V, V). What is the natural dimension of SL (V, V)?

21. Let $V = \{(x_1 \mid x_2\ x_3 \mid x_4\ x_5) \mid x_i \in Q;\ 1 \leq i \leq 5\}$ and $W = \{(x_1\ x_2 \mid x_3 \mid x_4\ x_5) \mid x_i \in Q;\ 1 \leq i \leq 5\}$ be two super vector spaces over Q.

    (a) Find a linear super transformation from V into W which is invertible.

    (b) Is all linear super transformation from V into W in SL (V, W) invertible?

    (c) Suppose SL (W, V) denotes the collection of all linear transformations from W into V. Does their exist any relation between SL (W, V) and SL (V, W)?

    (d) Can we say SL (V, V) and SL (W, W) are identical in this problem?

    (e) Is SL (V, V) any way related with SL (V, W) or SL (W, V)?

    (f) Give a non invertible linear transformations from V into W, W into V, V into V and W into W.

22. Let $V = \{(x_1\ x_2 \mid x_3\ x_4\ x_5 \mid x_6\ x_7\ x_8) \mid x_i \in Q;\ 1 \leq i \leq 8\}$ and $W = \{(x_1\ x_2 \mid x_3\ x_4\ x_5\ x_6 \mid x_7\ x_8 \mid x_9\ x_{10}\ x_{11}) \mid x_i \in Q;\ 1 \leq i \leq 11\}$ be two super vector spaces over the field of rationals. Find SL (V, W). Does SL (V, W) contain a non invertible linear super transformation? Give an example of an invertible super transformation $T_s$: $V \to W$ and verify for $T_s$, rank $T_s$ + nullity $T_s$ = dim V = 8.

23. Let $V = \{(x_1\ x_2 \mid x_3\ x_4 \mid x_5\ x_6\ x_7\ x_8) \mid x_i \in Q;\ 1 \leq i \leq 8\}$ be a super vector space over Q. Will every $T_s$: $V \to V \in$ SL (V, V) satisfy the equality rank $T_s$ + nullity $T_s$ = dim V?

24. Let $V = \{(x_1\ x_2\ x_3 \mid x_4\ x_5 \mid x_6\ x_7) \mid x_i \in Q;\ 1 \leq i \leq 7\}$ be a super vector space over Q. $W = \{(x_1\ x_2\ x_3\ x_4 \mid x_5 \mid x_6\ x_7\ x_8) \mid x_i \in Q;\ 1 \leq i \leq 8\}$ a super vector space over Q. $P = \{(x_1 \mid x_2\ x_3\ x_4 \mid x_5\ x_6) \mid x_i \in Q;\ 1 \leq i \leq 6\}$ be another super vector space over Q. Find SL(V, W), SL (W, P) and SL (V, P). Does then exist any relation between the 3 super spaces SL (V, W), SL (W, P) and SL(V, P)?

25. Let $V = \{(x_1\ x_2\ x_3 \mid x_4\ x_5\ x_6 \mid x_7\ x_8\ x_9 \mid x_{10}\ x_{11}\ x_{12}) \mid x_i \in Q; 1 \leq i \leq 12\}$ be a super vector space of natural dimension 12. Show $12 \times 12$ super diagonal matrix

$$A = \left(\begin{array}{ccc|ccc|ccc|ccc}
3 & 1 & 0 & 0 & 0 & 0 & 0 & 0 & 0 & 0 & 0 & 0 \\
0 & 1 & 1 & 0 & 0 & 0 & 0 & 0 & 0 & 0 & 0 & 0 \\
1 & 1 & 1 & 0 & 0 & 0 & 0 & 0 & 0 & 0 & 0 & 0 \\
\hline
0 & 0 & 0 & 1 & 1 & 1 & 0 & 0 & 0 & 0 & 0 & 0 \\
0 & 0 & 0 & 0 & 1 & 1 & 0 & 0 & 0 & 0 & 0 & 0 \\
0 & 0 & 0 & 0 & 1 & 0 & 0 & 0 & 0 & 0 & 0 & 0 \\
\hline
0 & 0 & 0 & 0 & 0 & 0 & 1 & 0 & 1 & 0 & 0 & 0 \\
0 & 0 & 0 & 0 & 0 & 0 & 0 & 1 & 0 & 0 & 0 & 0 \\
0 & 0 & 0 & 0 & 0 & 0 & 1 & 1 & 0 & 0 & 0 & 0 \\
\hline
0 & 0 & 0 & 0 & 0 & 0 & 0 & 0 & 0 & 1 & 2 & 3 \\
0 & 0 & 0 & 0 & 0 & 0 & 0 & 0 & 0 & 0 & 2 & 1 \\
0 & 0 & 0 & 0 & 0 & 0 & 0 & 0 & 0 & 3 & 0 & 1
\end{array}\right)$$

is associated with a linear operator $T_s$ and find that $T_s$. What is the nullity of $T_s$? Verify rank $T_s$ + nullity $T_s = 12$.

26. Prove any other interesting theorem / results about super vector spaces.

27. Prove all super vector spaces in general are not super linear algebras.

28. Is $W = \{(x_1\ x_2\ x_3 \mid x_4\ x_5 \mid x_6\ x_7\ x_8) \mid x_i \in Q;\ 1 \leq i \leq 8\}$ a super linear algebra over Q. Find a super subspace of W of dimension 6 over Q.

29. Suppose $V = \{(\alpha_1\ \alpha_2\ \alpha_3\ \alpha_4 \mid \alpha_5\ \alpha_6\ \alpha_7\ \alpha_8\ \alpha_9\ \alpha_{10} \mid \alpha_{11}\ \alpha_{12}) \mid x_i \in Q;\ 1 \leq i \leq 12\}$. Prove V is only a super vector space over Q.

Find super subspaces $W_1$ and $W_2$ of V such that $W_1 + W_2 = V$. Can $W_1 \cap W_2 = W$ be a super subspace different from the zero super space?

30. Given

$$V = \left\{ \left( \begin{array}{ccc|cc} \alpha_1 & \alpha_2 & \alpha_3 & & \\ \alpha_4 & \alpha_5 & \alpha_6 & & 0 \\ \alpha_7 & \alpha_8 & \alpha_9 & & \\ \hline & 0 & & \alpha_{10} & \alpha_{11} \\ & & & \alpha_{12} & \alpha_{13} \end{array} \right) \middle| \alpha_i \in Q; 1 \le i \le 13 \right\}.$$

Is V a super linear algebra over Q? Find nontrivial super subspaces of V. Find a nontrivial linear operator $T_s$ on V so that nullity T is not a trivial zero super subspace of V.

31. Show SL (V, W) is a super vector space over F where V and W are super vector spaces of dimension m and n respectively over F. Prove SL (V, W) ≅ {the set of all n × n super diagonal matrices}. Assume $m = m_1 + m_2 + m_3$ and $n = n_1 + n_2 + n_3$ and prove dimension of SL (V, W) is $n_1 \times m_1 + n_2 \times m_2 + n_3 \times m_3$.

32. Let $V = \{(x_1\ x_2\ x_3 \mid x_4\ x_5 \mid x_6\ x_7\ x_8) \mid x_i \in Q;\ 1 \le i \le 8\}$ be a super vector space over Q. Prove SL (V, V) is a super linear algebra of dimension 22. Show

$$SL(V, V) \cong \left\{ \left( \begin{array}{ccc|cc|ccc} \alpha_1 & \alpha_2 & \alpha_3 & 0 & 0 & 0 & 0 & 0 \\ \alpha_4 & \alpha_5 & \alpha_6 & 0 & 0 & 0 & 0 & 0 \\ \alpha_7 & \alpha_8 & \alpha_9 & 0 & 0 & 0 & 0 & 0 \\ \hline 0 & 0 & 0 & \alpha_{10} & \alpha_{11} & 0 & 0 & 0 \\ 0 & 0 & 0 & \alpha_{12} & \alpha_{13} & 0 & 0 & 0 \\ \hline 0 & 0 & 0 & 0 & 0 & \alpha_{14} & \alpha_{15} & \alpha_{16} \\ 0 & 0 & 0 & 0 & 0 & \alpha_{17} & \alpha_{18} & \alpha_{19} \\ 0 & 0 & 0 & 0 & 0 & \alpha_{20} & \alpha_{21} & \alpha_{22} \end{array} \right) \right|$$

$\alpha_i \in Q; 1 \leq i \leq 22\}$.

33. Given $V = \{(x_1\ x_2 \mid x_3 \mid x_4\ x_5 \mid x_6 \mid x_7\ x_8) \mid x_i \in Q; 1 \leq i \leq 8\}$ is a super vector space over Q. $W = \{(x_1 \mid x_2 \mid x_3\ x_4\ x_5 \mid x_6 \mid x_7) \mid x_i \in Q; 1 \leq i \leq 7\}$ is another super vector space over Q. Find SL(V, W) and SL(W, V). Find the dimension of SL (V, W) and SL(W, V). Why does dimension of super vector spaces of linear super transformation decreases in comparison with the vector space of linear transformations?

34. Let $V = \{(x_1 \mid x_2 \mid x_3 \mid x_4 \mid x_5) \mid x_i \in Q; 1 \leq i \leq 5\}$ be a super vector space over Q. Find SL (V, V). Find a basis for V and a basis for SL (V, V). Is SL (V, V) $\cong$ V? Justify your claim.

35. Can we prove if $V = \{(x_1 \mid \ldots \mid x_n) \mid x_i \in Q ; 1 \leq i \leq n\}$ be a super vector space over Q; SL (V, V) the super vector space of super linear operators on V. Is SL (V, V) $\cong$ V?

36. Suppose $V = \{(x_1 \mid \ldots \mid x_n) \mid x_i \in F; 1 \leq i \leq n\}$ be a super vector space over F. Can we prove with the increase in the number of partitions of the row vector $(x_1 \ldots x_n)$, the dimension of SL(V, V) decreases and with the decrease of the number of partition the dimension of SL (V, V) increases?

37. Let $V = \{(x_1\ x_2\ x_3 \mid x_4\ x_5\ x_6) \mid x_i \in Q; 1 \leq i \leq 6\}$ be a super vector space over Q. Prove SL (V, V) is of dimension 18 over Q. If $V = \{(x_1\ x_2 \mid x_3\ x_4 \mid x_5\ x_6) \mid x_i \in Q; 1 \leq i \leq 6\}$ is a super vector space over Q. Prove dimension of SL (V, V) is 12.

If $V = \{(x_1 \mid x_2\ x_3\ x_4\ x_5\ x_6) \mid x_i \in Q; 1 \leq i \leq 6\}$ be a super vector space over Q. Prove dimension of SL (V, V) is 26.

If $V = \{(x_1\ x_2 \mid x_3\ x_4\ x_5\ x_6) \mid x_i \in Q; 1 \leq i \leq 6\}$ be a super vector space over Q. Prove dimension of SL (V, V) is 20.

Prove maximum dimension of same number partition has maximum 26 and minimum is 18.

If $V = \{(x_1\ x_2\ x_3 \mid x_4 \mid x_5\ x_6) \mid x_i \in Q; 1 \leq i \leq 6\}$ is a super vector space over Q. Prove dimension of SL (V, V) is 14.

If $V = \{(x_1\ x_2\ x_3\ x_4 \mid x_5 \mid x_6) \mid x_i \in Q; 1 \leq i \leq 6\}$ is a super vector space over Q. Prove dimension of SL (V, V) is 18.

In this case can we say the minimum of one partition on V is the maximum of 2 partition on V.

38. Let $\{(x_1\ x_2\ x_3 \mid x_4\ x_5\ x_6 \mid x_7) \mid x_i \in Q; 1 \leq i \leq 7\}$ be a super vector space over Q.

    Find a linear super operator $T_s$ on V which is invertible. Give a linear operator $T_s^1$ on V which is non invertible. Obtain the related super matrices of $T_s$ and $T_s^1$.

39. Suppose $V = \{(x_1\ x_2 \mid x_3\ x_4 \mid x_5\ x_6) \mid x_i \in Q; 1 \leq i \leq 6\}$ a super vector space on the field Q. $W = \{(x_1\ x_2 \mid x_3 \mid x_4\ x_5\ x_6) \mid x_i \in Q; 1 \leq i \leq 6\}$ a super vector space over Q of same type as V.

    If $T_s$ is a linear super transformation from V into W and $U_s$ is a linear super transformation from W into V. Is $U_s \circ T_s$ defined? Justify your claim. Can we generalize this result?

40. Let V and W be two super vector spaces of same natural dimension but have the same type of partition. Let U ∈ SL(V, W) such that $U_s$ is an isomorphism. Is $T_s \to U_s T_s U_s^{-1}$ an isomorphism of SL(V, V) onto SL(W, W). Justify your answer.

41. If V and W are super vector spaces over the same field F, when will V and W be isomorphic. Is it enough if natural dimension V = natural dimension W? or it is essential both V and W should have the same dimension and the identical partition?

    Prove or disprove if they have same partition still $V \ncong W$.

42. Let $V = \{(x_1\ x_2\ x_3 \mid x_4\ x_5 \mid x_6\ x_7\ x_8) \mid x_i \in Q; 1 \leq i \leq 8\}$ be a super vector space over Q.

    Let $W = \{(x_1\ x_2 \mid x_3\ x_4 \mid x_5\ x_6\ x_7\ x_8) \mid x_i \in Q; 1 \leq i \leq 8\}$ be a super vector space over Q. Is $V \cong W$? We see V and W are super vector spaces of same dimension and also of same type of partition?

43. Let $V = \{(x_1\ x_2 \mid x_3\ x_4\ x_5 \mid x_6\ x_7\ x_8) \mid x_i \in Q; 1 \leq i \leq 8\}$ and $W = \{(x_1\ x_2\ x_3 \mid x_4\ x_5 \mid x_6\ x_7\ x_8) \mid x_i \in Q; 1 \leq i \leq 8\}$ be super vector

spaces over Q. Find SL(V, W). Find $T_s$ the linear transformation related to the super diagonal matrix.

$$A = \left(\begin{array}{cc|ccc|ccc} 1 & 2 & 0 & 0 & 0 & 0 & 0 & 0 \\ 0 & 1 & 0 & 0 & 0 & 0 & 0 & 0 \\ 1 & 1 & 0 & 0 & 0 & 0 & 0 & 0 \\ \hline 0 & 0 & 1 & 0 & 1 & 0 & 0 & 0 \\ 0 & 0 & 0 & 2 & 1 & 0 & 0 & 0 \\ \hline 0 & 0 & 0 & 0 & 0 & 1 & 2 & 3 \\ 0 & 0 & 0 & 0 & 0 & 0 & -1 & 4 \\ 0 & 0 & 0 & 0 & 0 & 1 & 0 & 2 \end{array}\right).$$

Does A relate to an invertible linear super transformation $T_s$ of V into W. Find nullity of $T_s$. Verify rank $T_s$ + nullity $T_s = 8$.

44. Let $V = \{(x_1\ x_2\ x_3 \mid x_4\ x_5 \mid x_6\ x_7) \mid x_i \in Q; 1 \leq i \leq 7\}$ be a super vector space over Q.
Let

$$A = \left(\begin{array}{ccc|cc|cc} 1 & 0 & 1 & 0 & 0 & 0 & 0 \\ 0 & 1 & 0 & 0 & 0 & 0 & 0 \\ 0 & 0 & 2 & 0 & 0 & 0 & 0 \\ \hline 0 & 0 & 0 & 1 & 0 & 0 & 0 \\ 0 & 0 & 0 & 1 & 2 & 0 & 0 \\ \hline 0 & 0 & 0 & 0 & 0 & 1 & 0 \\ 0 & 0 & 0 & 0 & 0 & 1 & 5 \end{array}\right)$$

be a super diagonal matrix associated with $T_s \in SL(V, V)$. Find the super eigen values of A? Determine the super eigen vectors related with A.

45. Let $V = \{(x_1\ x_2\ x_3 \mid x_4\ x_5\ x_6 \mid x_7\ x_8) \mid x_i \in Q; 1 \leq i \leq 8\}$ be a super vector space over Q. Does their exists a linear operator on V for which all the super eigen values are only imaginary?

46. Find for the above problem a $T_s: V \to V$ such that all the super eigen values are real.

47. Let $V = \{(x_1\ x_2\ x_3 \mid x_4\ x_5\ x_6 \mid x_7\ x_8\ x_9) \mid x_i \in Q; 1 \leq i \leq 9\}$ be a super vector space over Q. Is it ever possible for V to have a linear operator which has all its related eigen super values to be imaginary? Justify your claim.

48. Let $V = \{(x_1\ x_2 \mid x_3\ x_4 \mid x_5\ x_6) \mid x_i \in Q; 1 \leq i \leq 6\}$ be a super vector space over Q. Give a linear super transformation $T_s: V \to V$ which has all its eigen super values to be imaginary. Find $U_s: V \to W$ for which all eigen super values are real?

49. Let $V = \{(x_1\ x_2\ x_3 \mid x_4\ x_5 \mid x_6) \mid x_i \in Q; 1 \leq i \leq 6\}$ be a super vector space over Q. For the super diagonal matrix A associated with a linear operator $T_s$ on V calculate the super characteristic values, characteristic vectors and the characteristic subspace;

$$A = \left(\begin{array}{ccc|cc|c} 0 & 1 & 2 & 0 & 0 & 0 \\ 1 & 0 & 1 & 0 & 0 & 0 \\ 0 & 1 & -2 & 0 & 0 & 0 \\ \hline 0 & 0 & 0 & 1 & -1 & 0 \\ 0 & 0 & 0 & 1 & 2 & 0 \\ \hline 0 & 0 & 0 & 0 & 0 & -1 \end{array}\right).$$

50. Find all invertible linear transformations of V into V where $V = \{(x_1\ x_2 \mid x_3\ x_4 \mid x_5\ x_6 \mid x_7\ x_8) \mid x_i \in Q; 1 \leq i \leq 8\}$ is a super vector space over Q. What is the dimension of SL(V, V)?

51. Let $V = \{(x_1\ x_2\ x_3 \mid x_4\ x_5\ x_6\ x_7 \mid x_8\ x_9) \mid x_i \in Q; 1 \leq i \leq 9\}$ be a super vector space over Q. Is the linear operator $T_s ((x_1\ x_2\ x_3 \mid x_4\ x_5\ x_6\ x_7 \mid x_8\ x_9)) = (x_1 + x_2\ x_2 + x_3\ x_3 - x_1 \mid x_4\ 0\ x_5 + x_7\ x_6 \mid x_8\ x_8 + x_9)$ invertible? Find nullity $T_s$. Find the super diagonal matrix associated with $T_s$. What is the dimension of SL(V, V)?

52. Let $V = \{(x_1\ x_2 \mid x_3 \mid x_4\ x_5\ x_6\ x_7) \mid x_i \in Q;\ 1 \leq i \leq 7\}$ be a super vector space over Q.

If

$$A = \left(\begin{array}{cc|c|cccc} 1 & 2 & 0 & 0 & 0 & 0 & 0 \\ 0 & -3 & 0 & 0 & 0 & 0 & 0 \\ \hline 0 & 0 & 1 & 0 & 0 & 0 & 0 \\ \hline 0 & 0 & 0 & 1 & 2 & 0 & -1 \\ 0 & 0 & 0 & 0 & -1 & 2 & 3 \\ 0 & 0 & 0 & 1 & 0 & 1 & 0 \\ 0 & 0 & 0 & -3 & 1 & 0 & 1 \end{array}\right)$$

find $T_s$ associated with A. Find the characteristic super space associated with $T_s$. Write down the characteristic super polynomial associated with $T_s$.

53. Let $V = \{(x_1\ x_2 \mid x_3 \mid x_4\ x_5\ x_6 \mid x_7\ x_8 \mid x_9) \mid x_i \in Q;\ 1 \leq i \leq 9\}$ be a super vector space over Q. Find a basis for SL(V, V). What is the dimension of SL(V, V)? Find two super subspaces $W_1$ and $W_2$ of V so that $W_1 + W_2 = V$ and $W_1 \cap W_2 = \{0\}$.

54. Let $V = \{(x_1\ x_2\ x_3 \mid x_4\ x_5 \mid x_6\ x_7) \mid x_i \in Q;\ 1 \leq i \leq 7\}$ be a super vector space over Q. $T_s\ (x_1\ x_2\ x_3 \mid x_4\ x_5 \mid x_6\ x_7) = (x_1\ 0\ x_3 \mid 0\ x_5 \mid 0\ x_7)$ be a linear operator on V. Find the associated super diagonal matrix of $T_s$. Is $T_s$ an invertible linear operator? Prove rank $T_s$ + nullity $T_s$ = dim V. Find the associated characteristic super subspace of $T_s$.

55. Define a super hyper space of V, V a super vector space.

56. Give an example of a $10 \times 10$ super square diagonal matrix.

57. Give an example of super diagonal matrix, which is invertible.

58. Give an example of a $17 \times 15$ super diagonal matrix, which is not invertible.

59. Give an example of a 15 × 15 super diagonal matrix whose diagonal matrices are not square matrices.

60. Give an example of a square super diagonal square matrix and find its super determinant.

61. Give an example of a super diagonal matrix which is not a square matrix.

62. Give an example of a square super diagonal matrix whose diagonal entries are not square matrices.

63. Let

$$A = \left(\begin{array}{cc|ccc|cc|ccc} 3 & 1 & & & & & & & & \\ 0 & 0 & & 0 & & & 0 & & 0 & \\ \hline & & 1 & 0 & 2 & & & & & \\ & 0 & 3 & 4 & 5 & & 0 & & 0 & \\ & & 1 & 1 & 1 & & & & & \\ \hline & 0 & & 0 & & 2 & 5 & & 0 & \\ & & & & & -1 & 2 & & & \\ \hline & & & & & & & 1 & 2 & 3 \\ & 0 & & 0 & & & 0 & 4 & 5 & 6 \\ & & & & & & & 7 & 8 & 9 \end{array}\right)$$

be a square super diagonal square matrix. Determine the super determinant of A.

64. Find the characteristic super values associated with the super diagonal matrix A.

$$A = \left(\begin{array}{cccc|cccc|ccc}
3 & 1 & 0 & 1 & & & & & & & \\
0 & 0 & 1 & 2 & & 0 & & & & 0 & \\
-1 & 0 & -1 & 0 & & & & & & & \\
0 & 0 & 0 & 2 & & & & & & & \\
\hline
 & & & & 2 & 3 & 0 & 1 & & & \\
 & 0 & & & 1 & 0 & 1 & 0 & & 0 & \\
 & & & & 0 & 0 & 1 & 0 & & & \\
 & & & & 1 & 0 & 0 & 0 & & & \\
\hline
 & & & & & & & & 0 & 1 & 2 \\
 & 0 & & & & 0 & & & 0 & 1 & 0 \\
 & & & & & & & & 2 & 0 & 1
\end{array}\right)$$

65. Let

$$A = \left(\begin{array}{c|c|c|c}
A_1 & 0 & & 0 \\
\hline
0 & A_2 & & 0 \\
\hline
 & & & \\
\hline
0 & 0 & & A_n
\end{array}\right)$$

be a super diagonal square matrix with characteristic super polynomial $f = (f_1 \mid f_2 \mid \ldots \mid f_n) =$

$$((x - c_1^1)^{d_1^1} \ldots (x - c_{k_1}^1)^{d_1^{k_1}} \mid \ldots \mid (x - c_1^n)^{d_n^1} \ldots (x - c_{k_n}^n)^{d_n^{k_n}}).$$

Show that

$$(c_1^1 d_1^1 + \ldots + c_{k_1}^1 d_1^{k_1} \mid \ldots \mid c_1^n d_1^n + \ldots + c_{k_n}^n d_n^{k_n})$$
$$= (\text{trace } A_1 \mid \ldots \mid \text{trace } A_n).$$

66. Let $V = (V_1 \mid \ldots \mid V_n)$ be a super vector space of $(n_1 \times n_1, \ldots n_n \times n_n)$ super diagonal square matrices over the field F. Let

$$A = \left(\begin{array}{c|c|c|c}
A_1 & 0 & & 0 \\
\hline
0 & A_2 & & 0 \\
\hline
 & & & \\
\hline
0 & 0 & & A_n
\end{array}\right).$$

Let $T_s$ be the linear operator on $V = (V_1 \mid \ldots \mid V_n)$ defined by $T_s(B) = AB$

$$= \left(\begin{array}{c|c|c|c} A_1B_1 & 0 & & 0 \\ \hline 0 & A_2B_2 & & 0 \\ \hline & & & \\ \hline 0 & 0 & & A_nB_n \end{array}\right),$$

show that the minimal super polynomial for $T_s$ is the minimal super polynomial for A.

67. Let $V = (V_1 \mid \ldots \mid V_n)$ be a $(n_1 \mid \ldots \mid n_n)$ dimensional super vector space and $T_s$ be a linear operator on V. Suppose there exists positive integers $(k_1 \mid \ldots \mid k_n)$ so that $T_s^k = (T_1^{k_1} \mid \ldots \mid T_n^{k_n}) = (0 \mid 0 \mid \ldots \mid 0)$. Prove that $T^n = (T_1^{n_1} \mid \ldots \mid T_n^{n_n}) = (0 \mid \ldots \mid 0)$.

68. Let $V = (V_1 \mid \ldots \mid V_n)$ be a $(n_1, \ldots, n_2)$ finite dimensional $(n_1, \ldots, n_n)$ super vector space. What is the minimal super polynomial for the identity operator on V? What is the minimal super polynomial for the zero super operator?

69. Let

$$A = \left(\begin{array}{c|c|c|c} A_1 & 0 & & 0 \\ \hline 0 & A_2 & & 0 \\ \hline & & & \\ \hline 0 & 0 & & A_n \end{array}\right)$$

be a super diagonal square matrix with characteristic super polynomial

$$\left((x - c_1^1)^{d_1^1} \ldots (x - c_{k_1}^1)^{d_{k_1}^1} \mid \ldots \mid (x - c_1^n)^{d_1^n} \ldots (x - c_{k_n}^n)^{d_{k_n}^{n_n}}\right)$$

where $\left((c_1^1 \ldots c_{k_1}^1), \ldots, (c_1^n \ldots c_{k_n}^n)\right)$ are distinct. Let $V = (V_1 \mid \ldots \mid V_n)$ be the super space of $(n_1 \times n_1, \ldots, n_n \times n_n)$ matrices;

$$B = \left(\begin{array}{c|c|c|c} B_1 & 0 & & 0 \\ \hline 0 & B_2 & & 0 \\ \hline & & & \\ \hline 0 & 0 & & B_n \end{array}\right)$$

where $B_i$ is a $n_i \times n_i$ matrix i = 1, 2, …, n such that AB = BA

$$AB = \text{i.e.,} \left(\begin{array}{c|c|c|c} A_1B_1 & 0 & & 0 \\ \hline 0 & A_2B_2 & & 0 \\ \hline & & & \\ \hline 0 & 0 & & A_nB_n \end{array}\right) =$$

$$\left(\begin{array}{c|c|c|c} B_1A_1 & 0 & & 0 \\ \hline 0 & B_2A_2 & & 0 \\ \hline & & & \\ \hline 0 & 0 & & B_nA_n \end{array}\right) = B\,A.$$

Prove that the super dimension of $V = (V_1 \mid \ldots \mid V_n)$ is

$$(d_1^{1^2} + \ldots + d_{k_1}^{1^2} \mid (d_1^2)^2 + \ldots + (d_{k_2}^2)^2 \mid \ldots \mid (d_1^n)^2 + \ldots + (d_{k_n}^n)^2)\,.$$

70. Let $T_s$ be a linear operator on the $(n_1, \ldots, n_n)$ dimensional super vector space $V = (V_1 \mid \ldots \mid V_n)$ and suppose that $T_s$ has a n distinct characteristic super values. Prove that $T_s = (T_1 \mid \ldots \mid T_n)$ is super diagonalizable i.e., each $T_i$ is diagonalizable; i = 1, 2, …, n.

71. Let

$$A = \left(\begin{array}{c|c|c|c} A_1 & 0 & & 0 \\ \hline 0 & A_2 & & 0 \\ \hline & & & \\ \hline 0 & 0 & & A_n \end{array}\right)$$

and

$$B = \left(\begin{array}{c|c|c|c} B_1 & 0 & & 0 \\ \hline 0 & B_2 & & 0 \\ \hline & & & \\ \hline 0 & 0 & & B_n \end{array}\right)$$

be two ($n_1 \times n_1$, …, $n_n \times n_n$) super diagonal square matrices. Prove that if (I – AB) is invertible then (I – BA) is invertible and

$$(I - BA)^{-1} = I + B\,(1 - AB)^{-1}\,A.$$

$$\left(\begin{array}{c|c|c|c} (I - B_1A_1)^{-1} & 0 & & 0 \\ \hline 0 & (I - B_2A_2)^{-1} & & 0 \\ \hline & & & \\ \hline 0 & 0 & & (I - B_nA_n)^{-1} \end{array}\right)$$

$$= I + B\,(I - AB)^{-1}\,A.$$

$$(I_1 \mid \ldots \mid I_n) +$$

$$\left(\begin{array}{c|c|c|c} B_1(I - A_1B_1)^{-1}A_1 & 0 & & 0 \\ \hline 0 & B_2(I - A_2B_2)^{-1}A_2 & & 0 \\ \hline & & & \\ \hline 0 & 0 & & B_n(I - A_nB_n)^{-1}A_n \end{array}\right).$$

$$= \left(\begin{array}{c|c|c|c} I_1 + B_1(I-A_1B_1)^{-1}A_1 & 0 & & 0 \\ \hline 0 & I_2 + B_2(I-A_2B_2)^{-1}A_2 & & 0 \\ \hline & & & \\ \hline 0 & 0 & & I_n + B_n(I-A_nB_n)^{-1}A_n \end{array}\right).$$

Let A and B be two super diagonal square matrices over the field F of same order $(n_1 \times n_1, \ldots, n_n \times n_n)$ where

$$A = \left(\begin{array}{c|c|c|c} A_1 & 0 & & 0 \\ \hline 0 & A_2 & & 0 \\ \hline & & & \\ \hline 0 & 0 & & A_n \end{array}\right)$$

where $A_i$ is a $n_i \times n_i$ matrix and

$$B = \left(\begin{array}{c|c|c|c} B_1 & 0 & & 0 \\ \hline 0 & B_2 & & 0 \\ \hline & & & \\ \hline 0 & 0 & & B_n \end{array}\right)$$

of same order $B_i$ is a $n_i \times n_i$ matrix; $i = 1, 2, \ldots, n$.
The super diagonal square matrices AB and BA have same characteristic super values. Do they have same characteristic super polynomials? Do they have same minimal super polynomial?

73. Let $W = (W_1 \mid \ldots \mid W_n)$ be an invariant super subspace for $T_s = (T_1 \mid \ldots \mid T_n)$ of the super vector space $V = (V_1 \mid \ldots \mid V_n)$. Prove that the minimal super polynomial for the restriction operator $T_w = (T_1/W_1 \mid \ldots \mid T_n/W_n)$ divides the minimal super polynomial for $T_s$, without referring to super diagonal square matrices.

74. Let $T_s = (T_1 \mid \ldots \mid T_n)$ be a diagonalizable super linear operator on the $(n_1, \ldots, n_n)$ dimensional super vector space $V = (V_1 \mid \ldots \mid V_n)$ and let $W = (W_1 \mid \ldots \mid W_n)$ super subspace of V which is super invariant under $T = (T_1 \mid \ldots \mid T_n)$. Prove that the restriction operator $T_W$ is super diagonalizable.

75. Prove that if $T = (T_1 \mid \ldots \mid T_n)$ is a linear super operator on $V = (V_1 \mid \ldots \mid V_n)$, a super vector space. If every super subspace of V is super invariant under $T_s = (T_1 \mid \ldots \mid T_n)$ then $T_s$ is a scalar multiple of the identity operator $I = (I_1 \mid \ldots \mid I_n)$ where each $I_t$ is an identity operator from $V_t$ to itself for t = 1, 2, …, n.

76. Let $V = (V_1 \mid \ldots \mid V_n)$ be a super vector space over the field F. Each $V_t$ is a $n_t \times n_t$ square matrices with entries from F; t = 1, 2, …, n

Let

$$A = \left(\begin{array}{c|c|c|c} A_1 & 0 & & 0 \\ \hline 0 & A_2 & & 0 \\ \hline & & & \\ \hline 0 & 0 & & A_n \end{array}\right)$$

be a super diagonal square matrix where each $A_t$ is of $n_t \times n_t$ order; t = 1, 2, …, n.

Let $T_s$ and $U_s$ be linear super operators on $V = (V_1 \mid \ldots \mid V_n)$ defined by $T_s(B) = AB$

$$U_s(B) = AB - BA.$$

If A is super diagonalizable over F then $T_s$ is diagonalizable; True or false?

If A is super diagonalizable then $U_s$ is also super diagonalizable, prove or disprove.

77. Let $V = (V_1 \mid \ldots \mid V_n)$ be a super vector space over the field F. The super subspace $W = (W_1 \mid \ldots \mid W_n)$ is super invariant under (the family of operators) $\Im_s$; if W is super invariant under each operator in $\Im_s$. Using this prove the following:

Let $\Im_s$ be a commuting family of triangulable linear operators on a super vector space $V = (V_1 \mid \ldots \mid V_n)$. Let $W = (W_1 \mid \ldots \mid W_n)$ be a proper subsuper space of V which is super invariant under $\Im_s$. There exists a super vector $(\alpha_1 \mid \ldots \mid \alpha_n) \in V = (V_1 \mid \ldots \mid V_n)$ such that

(a) $\alpha = (\alpha_1 \mid \ldots \mid \alpha_n)$ is not in $W = (W_1 \mid \ldots \mid W_n)$.

(b) for each $T_s = (T_1 \mid \ldots \mid T_s)$ in $\Im_s$ the super vector $T_s\,\alpha = (T_1\alpha_1 \mid \ldots \mid T_n\alpha_n)$ is the super subspace spanned by α and W.

78. Let V be a finite $(n_1, \ldots, n_n)$ dimensional super vector space over the field F. Let $\Im_s$ be a commuting family of triangulable linear operators on $V = (V_1 \mid \ldots \mid V_n)$. There exists a super basis for V such that every operator in $\Im_s$ is represented by a triangular super diagonal matrix in that super basis.

Hence or other wise prove. If $\Im_s$ is a commuting family of $(n_1 \times n_1, \ldots, n_n \times n_n)$ super diagonal square matrices over an algebraically closed field F. There exists a non singular $(n_1 \times n_1, \ldots, n_n \times n_n)$ super diagonal square matrix P with entries in F such that

$$P^{-1}A\,P =$$

$$\left(\begin{array}{c|c|c|c} P_1^{-1} & 0 & & 0 \\ \hline 0 & P_2^{-1} & & 0 \\ \hline & & & \\ \hline 0 & 0 & & P_n^{-1} \end{array}\right) \left(\begin{array}{c|c|c|c} A_1 & 0 & & 0 \\ \hline 0 & A_2 & & 0 \\ \hline & & & \\ \hline 0 & 0 & & A_n \end{array}\right) \left(\begin{array}{c|c|c|c} P_1 & 0 & & 0 \\ \hline 0 & P_2 & & 0 \\ \hline & & & \\ \hline 0 & 0 & & P_n \end{array}\right)$$

is upper triangular for every super diagonal square matrix A in $\Im_s$.

79. Prove the following theorem. Let $\Im_s$ be a commuting family of super diagonalizable linear operators on a finite $(n_1, \ldots, n_n)$ dimensional super vector space $V = (V_1 \mid \ldots \mid V_n)$. There exists an ordered super basis for V such that every operator in $\Im_s$ is represented in that super basis by a super diagonal matrix.

80. Let F be a field, $(n_1, \ldots, n_n)$ a n tuple of positive integers and let $V = (V_1 \mid \ldots \mid V_n)$ be the super space of $(n_1 \times n_1, \ldots, n_n \times n_n)$ super diagonal square matrices over F; Let $(T_s)_A$ be the linear operator on V defined by

$$(T_s)_A(B) = AB - BA$$

$$\text{i.e., } \left((T_1)_{A_1}(B_1) \mid \ldots \mid (T_n)_{A_n}(B_n)\right)$$

$$= \left(\begin{array}{c|c|c|c} A_1B_1 & 0 & & 0 \\ \hline 0 & A_2B_2 & & 0 \\ \hline & & & \\ \hline 0 & 0 & & A_nB_n \end{array}\right) - \left(\begin{array}{c|c|c|c} B_1A_1 & 0 & & 0 \\ \hline 0 & B_2A_2 & & 0 \\ \hline & & & \\ \hline 0 & 0 & & B_nA_n \end{array}\right)$$

where

$$A = \left(\begin{array}{c|c|c|c} A_1 & 0 & & 0 \\ \hline 0 & A_2 & & 0 \\ \hline & & & \\ \hline 0 & 0 & & A_n \end{array}\right)$$

and

$$B = \left(\begin{array}{c|c|c|c} B_1 & 0 & & 0 \\ \hline 0 & B_2 & & 0 \\ \hline & & & \\ \hline 0 & 0 & & B_n \end{array}\right).$$

Consider the family of linear operators $(T_s)_A$ obtained by letting A vary over all super diagonal square matrices.
Prove that the operators in that family are simultaneously super diagonalizable.

81. Let $E_s = (E_1 \mid \ldots \mid E_n)$ be a super projection on $V = (V_1 \mid \ldots \mid V_n)$ and let $T_s = (T_1 \mid \ldots \mid T_n)$ be a linear operator on V. Prove that super range of $E_s = (E_1 \mid \ldots \mid E_n)$ is super invariant under $T_s$ if and only if $E_s T_s E_s = T_sE_s$ ie $(E_1T_1E_1 \mid \ldots \mid E_nT_nE_n) = (T_1E_1 \mid \ldots \mid T_nE_n)$.

    Prove that both the super range and super null space of E are super invariant under $T_s$ if and only if $E_sT_s = T_sE_s$ i.e., $(E_1T_1 \mid \ldots \mid E_nT_n) = (T_1E_1 \mid \ldots \mid T_nE_n)$.

82. Let $T_s = (T_1 \mid \ldots \mid T_n)$ be a linear operator on a finite $(n_1, \ldots, n_n)$ dimensional super vector space $V = (V_1 \mid \ldots \mid V_n)$.

    Let $R = (R_1 \mid \ldots \mid R_n)$ be the super range of $T_s$ and $N = (N_1 \mid \ldots \mid N_n)$ be the super null space of $T_s$. Prove that R and N are

independent if and only if V = R ⊕ N i.e., $(V_1 \mid \ldots \mid V_n) = (R_1 \oplus N_1 \mid \ldots \mid R_n \oplus N_n)$.

83. Let $T_s = (T_1 \mid \ldots \mid T_n)$ be a linear super operator on $(V = V_1 \mid \ldots \mid V_n)$. Suppose

$$V = W_1 \oplus \ldots \oplus W_k = (W_1^1 \oplus \ldots \oplus W_{n_1}^1 \mid \ldots \mid W_1^n \oplus \ldots \oplus W_{k_n}^n)$$

where each $W_i = (W_{i_1}^1 \mid \ldots \mid W_{i_n}^n)$ is super invariant under $T_s$. Let $T_t^{i_t}$ be the induced restriction operator on $W_{i_t}^t$

Prove:

a. super det T = super det $(T^1)$ … super det $(T^k)$ i.e., $(\det (T_1 \mid \ldots \mid \det T_n) = (\det(T_1^1) \ldots \det(T_1^{k_1}) \mid \ldots \mid \det(T_n^1) \ldots \det(T_n^{k_n}))$.

b. Prove that the characteristic super polynomial for $f = (f_1 \mid \ldots \mid f_n)$ is the product of characteristic super polynomials for $(f_1^1 \ldots f_1^{k_1}), \ldots, (f_n^1, \ldots f_n^{k_n})$.

84. Let $T_s = (T_1 \mid \ldots \mid T_n)$ be a linear operator on $V = (V_1 \mid \ldots \mid V_n)$ which commutes with every projection operator $E_s = (E_1 \mid \ldots \mid E_n)$ i.e., $T_sE_s = E_sT_s$ implies $(T_1E_1 \mid \ldots \mid T_nE_n) = (E_1T_1 \mid \ldots \mid E_nT_n)$. What can be said about $T_s = (T_1 \mid \ldots \mid T_n)$?

85. Let $V = (V_1 \mid \ldots \mid V_n)$ be a super vector space over F, where each $V_i$ is the space of all polynomials of degree less than or equal to $n_i$; i = 1, 2, …, n over F; prove that the differentiation operator $D_s = (D_1 \mid \ldots \mid D_n)$ on V is super nilpotent.
We say $D_s$ is super nilpotent if we can find a n-tuple of positive integers $p = (p_1, \ldots, p_n)$ such that $D_s^p = (D_1^{p_1} \mid \ldots \mid D_n^{p_n}) = (0 \mid \ldots \mid 0)$.

86. Let $T = (T_1 \mid \ldots \mid T_n)$ be a linear super operator on a finite dimensional super vector space $V = (V_1 \mid \ldots \mid V_n)$ with characteristic super polynomial $f = (f_1 \mid \ldots \mid f_n) =$

$$\left((x - c_1^1)^{d_1^1} \ldots (x - c_{k_1}^1)^{d_{k_1}^1} \mid \ldots \mid (x - c_1^n)^{d_1^n} \ldots (x - c_{k_n}^n)^{d_{k_n}^n}\right)$$

and super minimal polynomial $p = (p_1 \mid \ldots \mid p_n)$

$$=\left((x-c_1^1)^{r_1^1}\dots(x-c_{k_1}^1)^{r_{k_1}^1} \mid \dots \mid (x-c_1^n)^{r_1^n}\dots(x-c_{k_n}^n)^{r_{k_n}^n}\right).$$

Let $W_i = (W_{i_1}^1 \mid \dots \mid W_{i_n}^n)$ be the null super subspace of

$$(T - c_i I)^{r_i} = ((T_1 - c_{i_1}^1 I_1)^{r_{i_1}^1} \mid \dots \mid (T_n - c_{i_n}^n I_n)^{r_{i_n}^n}).$$

(a) Prove that $W_i = (W_{i_1}^1 \mid \dots \mid W_{i_n}^n)$ is the set of all super vectors $\alpha = (\alpha_1 \mid \dots \mid \alpha_n)$ in $V = (V_1 \mid \dots \mid V_n)$ such that $(T - c_i I)^m_\alpha$ = $((T_1 - c_{i_1}^1 I_1)^{m_1}_{\alpha_1} \mid \dots \mid (T_n - c_{i_n}^n I_n)^{m_n}_{\alpha_n}) = (0 \mid \dots \mid 0)$ for some n-tuple of positive integers $m = (m_1 \mid \dots \mid m_n)$.

(b) Prove that the super dimension of $W_i = (W_{i_1}^1 \mid \dots \mid W_{i_n}^n)$ is $\left(d_{i_1}^1, \dots, d_{i_n}^n\right)$.

87. Let $V = (V_1 \mid \dots \mid V_n)$ be a finite $(n_1, \dots, n_n)$ dimensional super vector space over the field of complex numbers. Let $T_s = (T_1 \mid \dots \mid T_n)$ be a linear super operator on V and $D_s = (D_1 \mid \dots \mid D_n)$ be the super diagonalizable part of $T_s$. Prove that if $g = (g_1 \mid \dots \mid g_n)$ is any super polynomial with complex coefficients then the diagonalizable part of $g_s(T_s) = (g_1(T_1) \mid \dots \mid g_n(T_n))$ is $g_s(D_s) = (g_1(D_1) \mid \dots \mid g_n(D_n))$.

88. Let $V = (V_1 \mid \dots \mid V_n)$ be a $(n_1, \dots, n_n)$ finite dimensional super vector space over the field F and let $T_s = (T_1 \mid \dots \mid T_n)$ be a linear super operator on V such that rank $(T_s) = (1, 1, \dots, 1)$.
Prove that either $T_s$ is super diagonalizable or $T_s$ is nilpotent, not both.

89. Let $V = (V_1 \mid \dots \mid V_n)$ be a finite $(n_1, \dots, n_n)$ dimensional super vector space over F. $T_s = (T_1 \mid \dots \mid T_n)$ be a linear super operator on V. Suppose that $T_s = (T_1 \mid \dots \mid T_n)$ commutes with every super diagonalizable linear operator on V. Prove that $T_s$ is a scalar multiple of the identity operator.

90. Let $T_s = (T_1 \mid \dots \mid T_n)$ be a linear super operator on $V = (V_1 \mid \dots \mid V_n)$ with minimal super polynomial of the form $p^n = (p_1^{n_1} | \ldots | p_n^{n_n})$ where p is super irreducible over the scalar field. Show that there is a super vector $\alpha = (\alpha_1 | \ldots | \alpha_n)$ in V = $(V_1 | \ldots | V_n)$ such that the super annihilator of $\alpha$ is $p^n = (p_1^{n_1} | \ldots | p_n^{n_n})$. (We say a super polynomial $p = (p_1 | \ldots | p_n)$ is super irreducible if each of the polynomial $p_i$ is irreducible for $i = 1, 2, \ldots, n$).

91. If $N_s = (N_1 | \ldots | N_n)$ is a nilpotent super operator on a $(n_1, \ldots, n_n)$ dimensional vector space $V = (V_1 | \ldots | V_n)$, then the characteristic super polynomial for $N_s = (N_1 | \ldots | N_n)$ is $x^n = (x^{n_1} | \ldots | x^{n_n})$.

92. Let $T_s = (T_1 | \ldots | T_n)$ be a linear super operator on the finite $(n_1, \ldots, n_n)$ dimensional super vector space $V = (V_1 | \ldots | V_n)$ let

$$p = (p^1, \ldots, p^n) = ((p_1^1)^{r_1^1} \ldots (p_{k_1}^1)^{r_{k_1}^1} | \ldots | (p_1^n)^{r_1^n} \ldots (p_{k_n}^n)^{r_{k_n}^n}),$$

be the minimal super polynomial for $T_s = (T_1 | \ldots | T_n)$ and let $V = (V_1 | \ldots | V_n) =$

$$\left(W_1^1 \oplus \ldots \oplus W_{k_1}^1 | \ldots | W_1^n \oplus \ldots \oplus W_{k_n}^n\right)$$

be the primary super decomposition for $T_s$; ie $W_{j_t}^t$ is the null space of $p_{i_t}^t (T_t)^{r_{i_t}^t}$, true for $t = 1, 2, \ldots, n$.

Let $W = (W_1 | \ldots | W_n)$ be any super subspace of V which is super invariant under $T_s$. Prove that $W = (W_1 | \ldots | W_n)$

$$= \left(W_1 \cap W_1^1 \oplus \ldots \oplus W_1 \cap W_{k_1}^1 | \ldots | W_n \cap W_1^n \oplus \ldots \oplus W_n \cap W_{k_n}^n\right).$$

93. Let $V = \{(x_1\ x_2\ x_3 | x_4\ x_5 | x_6\ x_7\ x_8) | x_i \in Q; 1 \le i \le 8\}$ be a super vector space over Q. Find super subspaces $W^1, \ldots, W^5$ in V which are super independent.

94. Find a set of super subspaces $W^1, \ldots, W^k$ of a super vector space $V = (V_1 | \ldots | V_n)$ over the field F which are not super independent.

95. Suppose V = $\{(x_1\ x_2\ x_3\ x_4 \mid x_5\ x_6 \mid x_7 \mid x_8\ x_9\ x_{10}) \mid x_i \in Q; 1 \leq i \leq 10\}$ is a super vector space over Q

(a) Find the maximal number of super subspaces which can be super independent.

(b) Find the minimal number of super subspaces which can be super independent.

(c) Can the collection of all super sub spaces of V be super independent? Justify your claim.

96. Suppose V = $(V_1 \mid \ldots \mid V_n)$ be a super vector space of $(n_1, \ldots, n_n)$ dimension over the field F. Suppose $W^t = (W_1^t \mid \ldots \mid W_n^t)$ be a super subspace of V for t = 1, 2, …, m. Find the number t so that that subset of $\{W^t\}_{t=1}^m$ happens to be super independent super subspaces. If $(m_1^t, \ldots, m_n^t)$ is the dimension of $W^t$ what can be said about $m_i^t$'s?

97. Let V = $\{(x_1\ x_2\ x_3 \mid x_4\ x_5 \mid x_6\ x_7\ x_8\ x_9) \mid x_i \in Q; 1 \leq i \leq 9\}$ be a super vector space over Q. Define $E_s = (E_1 \mid E_2 \mid E_3)$ a projection on V.
If $R_s$ is super range of $E_s$ and $N_s$ the super null space of $E_s$; prove $R_s \oplus N_s = V$ where $R_s = (R_1 \mid R_2 \mid R_3)$ and $N_s = (N_1 \mid N_2 \mid N_3)$.
Show if $T_s = (T_1 \mid T_2 \mid T_3)$ any linear operator on V then

$$T_s^2 = (T_1^2 \mid T_2^2 \mid T_3^2) \neq T_s = (T_1 \mid T_2 \mid T_2).$$

98. Let V = $(V_1 \mid \ldots \mid V_n)$ be a super vector space of finite $(n_1, \ldots, n_n)$ dimension over a field F. Suppose $E_s$ is any projection on V, prove $E_s = (E_1 \mid \ldots \mid E_n)$ is super diagonalizable.

99. Let V = $(V_1 \mid \ldots \mid V_n)$ be a super vector space over a field F. Let $T_s = (T_1 \mid \ldots \mid T_n)$ a linear operator V. Let $E_s = (E_1 \mid \ldots \mid E_n)$ be any projection on V. Is $T_sE_s = E_sT_s$? Will $(T_1 E_1 \mid \dots \mid T_n E_n) = (E_1 T_1 \mid \dots \mid E_nT_n)$? Justify your claim.

100. Derive primary decomposition theorem for super vector space $V = (V_1 \mid \dots \mid V_n)$ over F of finite $(n_1, \dots, n_n)$ dimension.

101. Define super diagonalizable part of a linear super operator $T_s$ on V ($T_s = (T_1 \mid \dots \mid T_n)$ and $V = (V_1 \mid \dots \mid V_n)$).

102. Define the notion of super nil potent linear super operator on a super vector space $V = (V_1 \mid \dots \mid V_n)$ over a field F.

103. Let $T_s = (T_1 \mid \dots \mid T_n)$ be a linear operator on $V = (V_1 \mid \dots \mid V_n)$ over the field F. Suppose that the minimal super polynomial for $T_s = (T_1 \mid \dots \mid T_n)$ decomposes over F into product of linear super polynomial, then prove there is a super diagonalizable super operator $D_s = (D_1 \mid \dots \mid D_n)$ on V and nilpotent super operator $N_s = (N_1 \mid \dots \mid N_n)$ on V such that (i) $T_s = D_s + N_s$ i.e., $T_s = (T_1 \mid \dots \mid T_n) = D_s + N_s = (D_1 + N_1 \mid \dots \mid D_n + N_n)$ (ii) $D_sN_s = N_sD_s$ ie $(D_1N_1 \mid \dots \mid D_nN_n) = (N_1D_1 \mid \dots \mid N_nD_n)$.

104. Does their exists a linear operator $T_s = (T_1 \mid \dots \mid T_n)$ on the super vector space $V = (V_1 \mid \dots \mid V_n)$ such that $T_s \neq D_s + N_s$?

105. Let $T_s = (T_1 \mid \dots \mid T_s)$ be a linear operator on a finite $(n_1, \dots, n_n)$ dimensional super vector space $V = (V_1 \mid \dots \mid V_n)$.
If $T_s = (T_1 \mid \dots \mid T_s)$ is super diagonalizable and if

$$c = (c_1 \dots c_k) = \left\{(c_1^1 \dots c_{k_1}^1), \dots, (c_1^n, \dots c_{k_n}^n)\right\}$$

are distinct characteristic super values of $T_s$ then there exists linear operators $E_s^1, \dots E_s^k$ on V. Prove that

a. $T_s = c_1E^1 + \dots + c_kE^k$ i.e., $(T_1 \mid \dots \mid T_n) = (c_1^1E_1^1 + \dots + c_{k_1}^1 E_{k_1}^1 \mid \dots \mid c_1^n E_1^n + \dots + c_{k_n}^n E_{k_n}^n)$
i.e., each $T_p = c_1^p E_1^p + \dots + c_{k_p}^p E_{k_p}^p$.

b. $I = E_s^1 + \dots + E_s^k = (I_1 \mid \dots \mid I_n)$ i.e., $I_t = E_1^t + \dots + E_{k_t}^t$; t = 1, 2, …, n

c. $E_s^i E_s^t = (0 \mid \ldots \mid 0)$ if $i \neq j$

i.e., $E_s^i\, E_s^j = (E_1^i\, E_1^j \mid \ldots \mid E_n^i\, E_n^j) = (0 \mid \ldots \mid 0)$ if $i \neq j$

d. $(E_s^i)^2 = E_s^i$ i.e., $(E_1^i \mid \ldots \mid E_n^i)^2 = (E_1^i \mid \ldots \mid E_n^i)$; i = 1, 2, …, n.

e. The super range of $E_s^i$ is the characteristic super space for $T_s$ associated with $c_i = (c_{i_1}^1, \ldots, c_{i_n}^n)$ where $E_s^i = (E_1^i \mid \ldots \mid E_n^i)$.

106. Define for a linear transformation $T_s$: V → W; V and W super inner product spaces a super isomorphism $T_s$ of V on to W.

107. Give an example of a complex inner product super vector space of (3, 9, 6) dimension.

108. Let $T_s$ = V → V be a linear super operator of a super complex inner product space. When will $T_s$ be a super self adjoint on V.

109. Can the notion of "super normal" be defined for any super matrix A? Justify your answer.

110. Let

$$A = \left(\begin{array}{c|c|c|c} A_1 & 0 & & 0 \\ \hline 0 & A_2 & & 0 \\ \hline & & & \\ \hline 0 & 0 & & A_n \end{array}\right)$$

be a super diagonal square complex matrix. Can A be defined to be super normal if $A_i A_i^* = A_i^* A_i$ for i = 1, 2, …, n.

111. Give an example of a super normal super diagonal square matrix.

112. Let $T_s = (T_1 \mid \ldots \mid T_n)$ be a linear super operator on a super vector space $V = (V_1 \mid \ldots \mid V_n)$ over the field F.

Define the super normal linear super operator $T_s$ on V and illustrate it by an example.

113. Prove only super diagonal square matrices can be super invertible matrices.

114. Is

$$A = \left(\begin{array}{ccc|ccc} 3 & 4 & 5 & 0 & 0 & 1 \\ \hline 0 & 1 & 1 & 3 & 2 & 1 \\ 9 & 2 & 0 & 1 & 1 & 1 \\ 1 & 1 & 1 & 1 & 0 & 1 \\ \hline 3 & 7 & 1 & 8 & 0 & 5 \\ 5 & 0 & 1 & 9 & 9 & 2 \end{array}\right),$$

an invertible matrix?

Justify your answer.

115. Let

$$A = \left(\begin{array}{ccc|cc|ccc} 3 & 1 & 2 & & & & & \\ 5 & 0 & 1 & 0 & & 0 & & \\ 1 & 2 & 3 & & & & & \\ \hline & 0 & & 3 & 4 & & 0 & \\ & & & 7 & 2 & & & \\ \hline & & & & & 5 & 1 & 0 \\ & 0 & & 0 & & 2 & 3 & 1 \\ & & & & & 0 & 1 & 5 \end{array}\right)$$

be a super diagonal square matrix. Is A a super invertible matrix?

116. Can every super diagonal matrix be an invertible matrix?

117. Let

$$A = \left(\begin{array}{cccc|cccc|ccc}
3 & 6 & 7 & 2 & & & & & & & \\
0 & 1 & 1 & 1 & & 0 & & & & 0 & \\
5 & 0 & 2 & 1 & & & & & & & \\
\hline
& & & & 3 & 2 & 1 & 0 & & & \\
& & & & 1 & 1 & 1 & 0 & & & \\
& 0 & & & 0 & 1 & 1 & 1 & & 0 & \\
& & & & 5 & 7 & 2 & 1 & & & \\
\hline
& 0 & & & & 0 & & & 3 & 7 & 5 \\
& & & & & & & & 4 & 2 & 1
\end{array}\right)$$

be super diagonal matrix. If A invertible?

118. Give an example of a super symmetric matrix.

119. Will the partition of a symmetric matrix always be a super symmetric matrix?

120. Let $T_s$ be a linear super operator on a super inner product space $V = (V_1 \mid \ldots \mid V_n)$ on a field F. When will $T_s = T_s^*$ ?

121. Suppose A and B are super square matrices of same natural order can we ever make A unitarily super equivalent to B. Justify your claim.

122. Suppose

$$A = \left(\begin{array}{c|c|c|c}
A_1 & 0 & & 0 \\
\hline
0 & A_2 & & 0 \\
\hline
 & & & \\
\hline
0 & 0 & & A_n
\end{array}\right)$$

is a super diagonal square matrix. Can we define for any

$$B = \left(\begin{array}{c|c|c|c} B_1 & 0 & & 0 \\ \hline 0 & B_2 & & 0 \\ \hline & & & \\ \hline 0 & 0 & & B_n \end{array}\right)$$

a super diagonal square matrix of same order. When can we say B is unitarily super equivalent to A.

123. Let

$$A = \left(\begin{array}{ccc|cccc|cc} 3 & 5 & 1 & & & & & & \\ 0 & 1 & 2 & & 0 & & & 0 & \\ 0 & 0 & 1 & & & & & & \\ \hline & & & 2 & 1 & 1 & 1 & & \\ & 0 & & 0 & 0 & 1 & 2 & 0 & \\ & & & 0 & 0 & 1 & 1 & & \\ & & & 0 & 0 & 0 & -1 & & \\ \hline & 0 & & & 0 & & & 9 & 2 \\ & & & & & & & 0 & 1 \end{array}\right)$$

and

$$B = \left(\begin{array}{ccc|cccc|cc} 5 & 0 & 1 & & & & & & \\ 0 & 1 & 2 & & 0 & & & 0 & \\ 0 & 0 & 1 & & & & & & \\ \hline & & & 1 & 0 & 1 & 1 & & \\ & 0 & & 0 & 1 & 2 & 1 & 0 & \\ & & & 0 & 0 & 1 & 2 & & \\ & & & 0 & 0 & 0 & 1 & & \\ \hline & 0 & & & 0 & & & 2 & 1 \\ & & & & & & & 0 & 1 \end{array}\right)$$

be two super diagonal square matrix. Is A super unitarily equivalent to B?

124. Define super ring of polynomials over Q. Is it a super vector space over Q?

125. Let V = [Q[x] | Q[x] | Q[x]] be a super vector space over Q. Find a super ideal of V which is a super minimal ideal of V.

126. Let $V = (V_1 \mid \ldots \mid V_n)$ be a super vector space over a field F. $T_s$ a linear super transformation from V into V.
Prove that the following two statements about $T_s = (T_1 \mid \ldots \mid T_n)$ are equivalent.

a. The intersection of the super range of $T_s = (T_1 \mid \ldots \mid T_n)$ and super null space of $T_s = (T_1 \mid \ldots \mid T_s)$ is a zero super subspace of V.

b. If $T_s (T_s(\alpha)) = (T_1(T_1(\alpha_1)) \mid \ldots \mid T_n (T_n(\alpha_n))] = (0 \mid \ldots \mid 0)$ then $T_s\alpha = (T_1\alpha_1 \mid \ldots \mid T_n\alpha_n) = (0 \mid \ldots \mid 0)$.

127. Define super linear functional? Give an example.

128. Define the concept of dual super space of a super space $V = (V_1 \mid \ldots \mid V_n)$ over the field F.

129. Can polarization identities be derived for super norms defined over super vector spaces?

130. Define the super matrix of the super inner product for a given super basis for a super vector space $V = (V_1 \mid \ldots \mid V_n)$ over a field F.

131. Verify the super standard inner product on $V = (F^{n_1} \mid \ldots \mid F^{n_n})$ over the field F is an super inner product on V.

132. Can Cauchy Schwarz super inequality for super vector spaces be oblained?

133. Can Bessels inequality of super vector spaces with super inner product be derived?

134. Can we have a polar decomposition in case of linear operators $T_s$. $U_s$ and $W_s$ on a super vector space V such that $T_s = U_sN_s$?

135. Give a proper definition of a non-negative super diagonal square matrix.

$$A = \left(\begin{array}{c|c|c|c} A_1 & 0 & & 0 \\ \hline 0 & A_2 & & 0 \\ \hline & & & \\ \hline 0 & 0 & & A_n \end{array}\right)$$

where each $A_i$ is a $n_i \times n_i$ matrix i = 1, 2, …, n.
Then prove that such a super diagonal square matrix has a unique non negative super square root. Illustrate this by an example.

136. If $U_s$ and $T_s$ are normal operators in SL (V, V) which commute prove $T_s + U_s$ and $U_s T_s$ are also normal.

137. Let SL (V, V) be the set of operators on a super vector space V $= (V_1 \mid \ldots \mid V_n)$ over a complex field i.e., V itself is finite $(n_1, \ldots, n_n)$ dimensional complex super inner product space. Prove that the following statements about $T_s$ are equivalent.

a. $T_s = (T_1 \mid \ldots \mid T_n)$ is (super) normal

b. $\| T_s\alpha \| = (\| T_1\alpha_1 \| \mid \ldots \mid \| T_n\alpha_n \|)$
$= (\| T_1^*\alpha_1 \| \mid \ldots \mid \| T_n^*\alpha_n \|)$
$= \| T_s^*\alpha \|$.
for every $\alpha = (\alpha_1 \mid \ldots \mid \alpha_n) \in V = (V_1 \mid \ldots \mid V_n)$.

c. $T_s = T_s^1 + i\, T_s^2$ where $T_s^1$ and $T_s^2$ are super self adjoint and $T_s^1\, T_s^2 = T_s^2\, T_s^1$ where $T_s = (T_1 \mid \ldots \mid T_n)$
$= (T_1^1 + i\, T_1^2 \mid \ldots \mid T_n^1 + i\, T_n^2)$
and

$$T_s^1\, T_s^2 = (T_1^1\, T_1^2 \mid \ldots \mid T_n^1\, T_n^2)$$
$$= (T_1^2\, T_1^1 \mid \ldots \mid T_n^2\, T_n^1)$$
$$= T_s^2\, T_s^1.$$

d. If $\alpha = (\alpha_1 \mid \ldots \mid \alpha_n)$ is a super vector and $c = (c_i, \ldots, c_n)$ any scalar n–tuple then $T_s\, \alpha = c\alpha$ i.e., $(T_1\alpha_1 \mid \ldots \mid T_n\alpha_n) = (c_1\alpha_1 \mid \ldots \mid c_n\alpha_n)$ then $T_s^* \alpha = \bar{c}\, \alpha$ i.e., $(T_s^*\, \alpha_1 \mid \ldots \mid T_n^* \alpha_n) = (\bar{c}\, \alpha_1 \mid \ldots \mid \bar{c}_n\, \alpha_n)$,

e. There is an orthonormal super basis for $V = (V_1 \mid \ldots \mid V_n)$ consisting of characteristic super vectors for $T_s = (T_1 \mid \ldots \mid T_n)$.

f. There is an orthonormal super basis $B = (B_1 \mid \ldots \mid B_n)$. $B_i$ a basis for $V_i$; $i = 1, 2, \ldots, n$ such that $[T_s]_B$ is a super diagonal matrix A.
i.e., $\left[\left([T_1]_{B_1} \mid \ldots \mid [T_n]_{B_n}\right) = \left(A_1 \mid \ldots \mid A_n\right)\right]$ where each $A_i$ is a diagonal matrix, $i = 1, 2, \ldots, n$.
i.e.,

$$A = \left(\begin{array}{c|c|c|c} A_1 & 0 & & 0 \\ \hline 0 & A_2 & & 0 \\ \hline & & & \\ \hline 0 & 0 & & A_n \end{array}\right).$$

g. There is a super polynomial $g = (g_1 \mid \ldots \mid g_n)$ with complex coefficients such that $T_s^* = g(T_s)$ i.e., $(T_1^* \mid \ldots \mid T_n^*) .= (g_1(T_1) \mid \ldots \mid g_n(T_n))$.

h. Every super subspace which is super invariant under $T_s$ is also super invariant under $T_1^*$.

i. $T_s = N_s\, U_s$ where $N_s$ is super non negative, $U_s$ is super unitary and $N_s$ super commutes with $U_s$ ie $(N_1\, U_1 \mid \ldots \mid N_n\, U_n) = N_s\, U_s = (U_1\, N_1 \mid \ldots \mid U_n\, N_n) = U_s\, N_s$.

j. $T_s = (C_1^1 E_1^1 + \ldots + C_{k_1}^1 E_{k_1}^1 | \ldots | C_1^n E_1^n + \ldots + C_{k_n}^n E_{k_n}^n)$ where I $= (I_1 | \ldots | I_n)$
$= (E_1^1 + \ldots + E_{k_1}^1 | \ldots | E_1^n + \ldots + E_{k_n}^n)$
with
$E_{i_t}^t E_{j_t}^t = 0$ if $i_t \neq j_t$; t = 1, 2, …, n and $(E_{i_t}^t)^2 = E_{i_t}^t = E_{i_t}^{t*}$ for $1 \leq i_t \leq k_t$ and t = 1, 2, …, n.

138. Let $V = (V_1 | \ldots | V_n)$ be a super complex $(n_1 \times n_1, \ldots, n_n \times n_n)$ super diagonal matrices equipped with a super inner product $(A | B) = \text{trace}(AB^*)$
i.e., $((A_1 | B_1) | \ldots | (A_n | B_n))$
$= (\text{tr}(A_1 B_1^*) | \ldots | \text{tr}(A_n, B_n^*))$
where

$$A = \left(\begin{array}{c|c|c|c} A_1 & 0 & & 0 \\ \hline 0 & A_2 & & 0 \\ \hline & & & \\ \hline 0 & 0 & & A_n \end{array}\right)$$

and

$$B = \left(\begin{array}{c|c|c|c} B_1 & 0 & & 0 \\ \hline 0 & B_2 & & 0 \\ \hline & & & \\ \hline 0 & 0 & & B_n \end{array}\right).$$

If B is a super diagonal $(n_1 \times n_1, \ldots n_n \times n_n)$ matrix of V, let

$$L_B = (L_{B_1}^1 | \ldots | L_{B_n}^n),$$
$$R_B = (R_{B_1}^1 | \ldots | R_{B_n}^n)$$

and

$$T_B = (T_{B_1}^1 | \ldots | T_{B_n}^n),$$

denote the linear super operators on $V = (V_1 | \ldots | V_n)$ defined by

(a) $L_B(A) = BA$
i.e., $(L^1_{B_1}(A_1) | \ldots | L^n_{B_n}(A_n)) = (B_1 A_1 | \ldots | B_n A_n)$.

(b) $R_B(A) = AB$ i.e.,
$(R^1_{B_1}(A_1) | \ldots | R^n_{B_n}(A_n)) = (A_1B_1 | \ldots | A_nB_n)$.

(c) $T_B(A) = (T_{B_1}(A_1) | \ldots | T_{B_n}(A_n))$
$= ((B_1A_1 - A_1B_1) | \ldots | (B_nA_n - A_nB_n)) = BA - AB$.

139. Let $\mathfrak{I}_s$ be a commuting family of super diagonalizable normal operators on a finite $(n_1, \ldots, n_n)$ dimensional super inner product space $V = (V_1 | \ldots | V_n)$ and $A_0$ the self adjoint super algebra generated by $\mathfrak{I}_s$. Let $a_s$ be the self adjoint super algebra generated by $\mathfrak{I}_s$ and the super identity operator $I = (I_1 | \ldots | I_n)$. Show that

a. $a_s$ is the set of all operators on V of the form $cI + T_s$ i.e., $(c_1 I_1 + T_1 | \ldots | c_n I_n + T_n)$ where $c = (c_1, \ldots, c_n)$ is a scalar n tuple and $T_s = (T_1 | \ldots | T_n)$ is a super operator in $a_s$ and $T_s$ an operator in $a_{s_0}$

b. $a_s = a_{s_0}$ if and only if for each super root $r = (r_1, \ldots, r_n)$ of $a_s$ there exists an operator $T_s$ in $a_{s_0}$ such that $r(T_s) = (r_1(T_1) | \ldots | r_n(T_n)) \neq (0 | \ldots | 0)$.

140. Find all linear super forms on the super space of column super vectors $V = (n_1 \times 1 | \ldots | n_n \times 1)$, super diagonal matrices over C which are super invariant under $o(n, c) = (o(n_1, c) | \ldots | o(n_n, c))$

141. Find all bilinear super forms on the super space of column super vector $V = (n_1 \times 1 | \ldots | n_n \times 1)$, super diagonal matrices over R which are super invariant under o(n, R).

142. Does their exists any relation between the problems 140 and 141.

143. Let m = $(m_1 \mid \ldots \mid m_n)$ be a member of the complex orthogonal super group $(o(n_1, c) \mid \ldots \mid o(n_n, c)$ Show that

$$m^t = (m_1^t \mid \ldots \mid m_n^t) = \bar{m} = (\bar{m}_1 \mid \ldots \mid \bar{m}_n)$$

and

$$m^* = (m_1^* \mid \ldots \mid m_n^*) = \bar{m}^t = (\bar{m}_1^t \mid \ldots \mid \bar{m}_n^t)$$

also belong to $o(n, c) = (o(n_1, c) \mid \ldots \mid o(n_n, c))$.

144. Suppose $m = (m_1 \mid \ldots \mid m_n)$ belongs to $o(n, c) = (o(n_1, c) \mid \ldots \mid o(n_n, c))$ and that $m' = (m_1' \mid \ldots \mid m_n')$ similar to m. Does m' also belong to o (n, c).

145. Let

$$y_1 = (y_{i_1} \mid \ldots \mid y_{i_n}) = \left( \sum_{k_1=1}^{n_1} m^1_{j_1 k_1} x^1_{k_1} \mid \ldots \mid \sum_{k_n=1}^{n_n} m^n_{j_n k_n} x^n_{k_n} \right)$$

where $m = (m_1 \mid \ldots \mid m_n)$ is a member of $o(n, c) = (o(n_1, c) \mid \ldots \mid o(n_n, c))$.

Show that

$$\sum_j y_c^2 = \left( \sum_{j_1} (y^1_{i_1})^2 \mid \ldots \mid \sum_{j_n} (y^n_{i_n})^2 \right)$$

$$= \left( \sum_{j_1} (x^1_{i_1})^2 \mid \ldots \mid \sum_{j_n} (x^n_{i_n})^2 \right)$$

$$= \sum_j x_i^2 \ .$$

146. Let $m = (m_1 \mid \ldots \mid m_n)$ be an $(n_1 \times n_1, \ldots, n_n \times n_n)$ super diagonal matrix over C with columns

$$m_1^1 \ldots m_{n_1}^1, m_1^2 \ldots m_{n_2}^2, \ldots, m_1^n, \ldots, m_{n_n}^n \ ,$$

show that m belongs to $o(n, c) = (o(n_1, c) \mid \ldots \mid o(n_n, c))$ if and only if

$$m_j^t m_k = \delta_{jk} \text{ i.e., } ((m^1_{j_1})^t \ m^1_{k_1} \mid \ldots \mid (m^n_{j_n})^t \ m^n_{k_n}) = (\delta_{j_1 k_1} \mid \ldots \mid \delta_{j_n k_n}) \ .$$

147. Let x = ($x_1$ | … | $x_n$) be an ($n_1 \times 1$, …, $n_n \times 1$) super diagonal matrix over C. Under what condition o (n, c) = (o ($n_1$, $c_1$) | … | o($n_n$, c)) contain a super diagonal matrix m whose first super column is X i.e., if

$$m = \left(\begin{array}{c|c|c|c} m_1 & 0 & & 0 \\ \hline 0 & m_2 & & 0 \\ \hline & & & \\ \hline 0 & 0 & & m_n \end{array}\right)$$

i.e., o(n, c) has a super diagonal matrix m such that the matrix $m_i$ whose first column is $x_i$; i = 1, 2, …, n.

148. Let V = ($V_1$ | … | $V_n$) be the space of all n × 1 = ($n_1 \times 1$ | … | $n_n \times 1$) matrices over C and f = ($f_1$ | … | $f_n$) the bilinear super form on V given by

$$\begin{aligned} f(x,y) &= (f_1(x_1, y_1) \mid \ldots \mid f_n(x_n, y_n)) \\ &= (x_1^t y_1 \mid \ldots \mid x_n^t y_n). \end{aligned}$$

Let m belong to o(n c) = (o($n_1$, c) | … | o($n_n$, c)). What is the super diagonal matrix of f in the super basis of V containing super columns $m_1^1 \ldots m_{n_1}^1, \ldots, m_1^n, \ldots, m_{n_n}^n$ of m?

149. Let x = ($x_1$ | … | $x_n$) be a ($n_1 \times 1$ | … | $n_n \times 1$) super matrix over C such that $x^t x = (x_1^t x_1 \mid \ldots \mid x_n^t x_n) = (1 \mid \ldots \mid 1)$ and $I_j = (I_{j_1} \mid \ldots \mid I_{j_n})$ be the $j^{th}$ super column of the identity super diagonal matrix. Show there is a super diagonal matrix

$$m = \left(\begin{array}{c|c|c|c} m_1 & 0 & & 0 \\ \hline 0 & m_2 & & 0 \\ \hline & & & \\ \hline 0 & 0 & & m_n \end{array}\right)$$

in o(n, c) = (o($n_1$, c) | … | o($n_n$, c)) such that m x = $I_j$; i.e.,

$$\begin{pmatrix} m_1x_1 & 0 & & 0 \\ 0 & m_2x_2 & & 0 \\ & & & \\ 0 & 0 & & m_nx_n \end{pmatrix}$$

$= [I_{j_n} | \ldots | I_{j_n}]$. If $x = (x_1 | \ldots | x_n)$ has real entries show there is a m in o (n, R) with the property that $mx = I_j$.

150. Let $V = (V_1 | \ldots | V_n)$ be a super space of all $(n_1 \times 1 | \ldots | n_n \times 1)$ super diagonal matrices over C.

$$A = \begin{pmatrix} A_1 & 0 & & 0 \\ 0 & A_2 & & 0 \\ & & & \\ 0 & 0 & & A_n \end{pmatrix}$$

an $(n_1 \times n_1, \ldots, n_n \times n_n)$ super diagonal matrix over C, here each $A_i$ is a $n_i \times n_i$ matrix; $i = 1, 2, \ldots, n$ and $f = (f_1 | \ldots | f_n)$ the bilinear super form on V given by

$$f(x, y) = (f_1 (x_1, y_1) | \ldots | f_n(x_n, y_n))$$
$$= x^t A\, y = (x_1^t A_1 y_1 | \ldots | x_n^t A_n y_n)$$

$$= \begin{pmatrix} x_1^t A_1 Y_1 & 0 & & 0 \\ 0 & x_2^t A_2 Y_2 & & 0 \\ & & & \\ 0 & 0 & & x_n^t A_n Y_n \end{pmatrix}.$$

Show that f is super invariant under o (n c) = (o $(n_1, c)$ | … | o$(n_n, c)$) i.e., f(mx; my) = f(x, y) i.e., $(f_1(m_1x_1, m_1y_1) | \ldots | f_n(m_nx_n, m_ny_n)) = (f_1(x_1, y_1) | \ldots | f_n(x_n, y_n))$ for all $x = (x_1 | \ldots | x_n)$ and $y = (y_1 | \ldots | y_n)$ in V and

$$m = \left(\begin{array}{c|c|c|c} m_1 & 0 & & 0 \\ \hline 0 & m_2 & & 0 \\ \hline & & & \\ \hline 0 & 0 & & m_n \end{array}\right)$$

in (o ($n_1$, c) | … | o ($n_n$, c)) = o (n, c) if and only if

$$A = \left(\begin{array}{c|c|c|c} A_1 & 0 & & 0 \\ \hline 0 & A_2 & & 0 \\ \hline & & & \\ \hline 0 & 0 & & A_n \end{array}\right)$$

commutes with each member of o (n, c).

151. Let F be a subfield of C, V be the super space of ($n_1 \times 1$ | … | $n_n \times 1$) matrices over F i.e.,

$$V = \left\{\left(x_1^1 \ldots x_{n_1}^1 \,|\ldots|\, x_1^n \ldots x_{n_n}^n\right)\right\}^t$$

is the collection of all super column vectors.

$$A = \left(\begin{array}{c|c|c|c} A_1 & 0 & & 0 \\ \hline 0 & A_2 & & 0 \\ \hline & & & \\ \hline 0 & 0 & & A_n \end{array}\right)$$

is a super diagonal matrix where each $A_i$ is a $n_i \times n_i$ matrix over F, and f = ($f_1$ | … | $f_n$) the bilinear super form on V given by f(x, y) = $x^t Ay$ i.e., ($f_1$ ($x_1$, $y_1$) | … | $f_n$ ($x_n$, $y_n$))

$$= \left(\begin{array}{c|c|c|c} x_1^t A_1 y_1 & 0 & & 0 \\ \hline 0 & x_2^t A_2 y_2 & & 0 \\ \hline & & & \\ \hline 0 & 0 & & x_n^t A_n y_n \end{array}\right).$$

If m is a super diagonal matrix

$$\left(\begin{array}{c|c|c|c} m_1 & 0 & & 0 \\ \hline 0 & m_2 & & 0 \\ \hline & & & \\ \hline 0 & 0 & & m_n \end{array}\right)$$

where each $m_i$ is a $n_i \times n_i$ matrix over F; show that m super preserves f if and only if $A^{-1}m^tA =$

$$\left(\begin{array}{c|c|c|c} A_1^{-1} & 0 & & 0 \\ \hline 0 & A_2^{-1} & & 0 \\ \hline & & & \\ \hline 0 & 0 & & A_n^{-1} \end{array}\right) \times \left(\begin{array}{c|c|c|c} m_1^t & 0 & & 0 \\ \hline 0 & m_2^t & & 0 \\ \hline & & & \\ \hline 0 & 0 & & m_n^t \end{array}\right)$$

$$\times \left(\begin{array}{c|c|c|c} A_1 & 0 & & 0 \\ \hline 0 & A_2 & & 0 \\ \hline & & & \\ \hline 0 & 0 & & A_n \end{array}\right)$$

$$= \left(\begin{array}{c|c|c|c} A_1^{-1}m_1^tA_1 & 0 & & 0 \\ \hline 0 & A_2^{-1}m_2^tA_2 & & 0 \\ \hline & & & \\ \hline 0 & 0 & & A_n^{-1}m_n^tA_n \end{array}\right)$$

$$= \left(\begin{array}{c|c|c|c} m_1^{-1} & 0 & & 0 \\ \hline 0 & m_2^{-1} & & 0 \\ \hline & & & \\ \hline 0 & 0 & & m_n^{-1} \end{array}\right)$$

$= m^{-1}$.

152. Let g = ($g_1$ | … | $g_n$) be a non singular bilinear super form on a finite ($n_1$, …, $n_n$) dimensional super vector space V = ($V_1$ | … | $V_n$). Suppose $T_s$ = ($T_1$ | … | $T_n$) is a linear operator on V and that f = ($f_1$ | … | $f_n$) be a bilinear super form on V given by f (α, β ) = g (α, $T_s$β). i.e.,

$$(f_1 (\alpha_1, \beta_1) | \ldots | f_n(\alpha_n, \beta_n)) = (g_1(\alpha_1, T_1\beta_1) | \ldots | g_n(\alpha_n, T_n\beta_n)).$$

If $U_s$ = ($U_1$ | … | $U_n$) is a linear operator on V find necessary and sufficient condition for $U_s$ to preserve f.

153. Let q = ($q_1$ | $q_2$) be the quadratic super form on ($R^2$ | $R^3$) given by

$$\begin{aligned} q (x, y) &= (q_1(x_1^1, x_2^1) | q_2(x_1^2\ x_2^2\ x_3^2)) \\ &= (2bx_1^1\ x_2^1 \mid x_1^2\ x_2^2 + 2x_1^2\ x_3^2 + (x_3^2)^2) \end{aligned}$$

Find a super invertible linear operator $U_s$ = ($U_1$ | $U_2$) on ($R^2$ | $R^3$) such that

$$\begin{aligned} &\left((U_1^t q_1)(x_1^1, x_2^1) \mid (U_2^t\ q_2)(x_1^2\ x_2^2\ x_3^2)\right) \\ &\quad = (2b(x_1^1)^2 - 2b(x_2^1)^2 \mid (x_1^2)^2 - (x_2^2)^2 + (x_3^2)^2). \end{aligned}$$

154. Let V = ($V_1$ | … | $V_n$) be a finite ($n_1$, …, $n_n$) dimensional super vector space and f = ($f_1$ | … | $f_n$) a super non degenerate symmetric bilinear super form on V associated with f is a natural super homomorphism of V into the dual super space $V^* = (V_1^* | \ldots | V_n^*)$, this super isomorphism being the transformation $L_f = (L_{f_1}^1 | \ldots | L_{f_n}^n)$. Using $L_f$ show that for each super basis B = [$B_1$ | … | $B_n$] = $(\alpha_1^1 \ldots \alpha_{n_1}^1 | \ldots | \alpha_1^n \ldots \alpha_{n_n}^n)$ on V there exists a unique super basis $B^1 = (\beta_1^1 \ldots \beta_{n_1}^1 | \ldots | \beta_1^n \ldots \beta_{n_n}^n)$ = $(B_1^1 | \ldots | B_n^1)$ of V such that

$$f(\alpha_i, \beta_i) = (f_1(\alpha_{i_1}, \beta_{j_1}^1) | \ldots | f_n(\alpha_{i_n}^n, \beta_{j_n}^n)) = (\delta_{i_1 j_1} | \ldots | \delta_{i_n j_n}).$$

Then show that for every super vector α = ($\alpha_1$ | … | $\alpha_n$) in V we have

$$\alpha = \left( \sum_{i_1} f_1(\alpha_1, \beta_{i_1}^1)\, \alpha_{i_1}^1 \,|\ldots|\ldots|\, \sum_{i_n} f_n(\alpha_n, \beta_{i_n}^n)\, \alpha_{i_n}^n \right)$$

$$= \left( \sum_{i_1} f_1(\alpha^1_{i_1}, \alpha_1)\, \beta^1_{i_1} \mid \ldots \mid \sum_{i_n} f_n(\alpha^n_{i_n}, \alpha_n)\, \beta^n_{i_n} \right).$$

154. Let V, f, B and $B_1$ be as in problem (153); suppose $T_s = (T_1 \mid \ldots \mid T_n)$ is a linear super operator on V and $T'_s$ is the linear operator which f associates with $T_s$ given by $f(T_s\alpha, \beta) = f(\alpha, T'_s\beta)$ i.e., $(f_1(T_1\alpha_1, \beta_1) \mid \ldots \mid f_n(T_n\alpha_n, \beta_n)) = (f_1(\alpha_1, T'_1\beta_1) \mid \ldots \mid f_n(\alpha_n, T'_n\beta_n))$;

(a) Show that $[T'_s]_{B'} = [T]^t_B$
i.e., $[[T'_1]_{B_1} \mid \ldots \mid [T'_n]_{B_n}] = [[T_1]^t_{B_1} \mid \ldots \mid [T_n]^t_{B_n}]$.

(b) super tr $(T_s)$

$$= \text{super trace } (T^t) = \sum_i f\,(T_s\alpha_i, \beta_i) \text{ i.e.,}$$

$$(\text{tr}\,(T_1) \mid \ldots \mid \text{tr}\,(T_n)) = (\text{tr}(T'_1) \mid \ldots \mid \text{tr}\,(T'_n))$$

$$= \left( \sum_{i_1} f_1(T_1\alpha^1_{i_1}, \beta^1_{i_1}) \mid \ldots \mid \sum_{i_n} f_n(T_n\alpha^n_{i_n}, \beta^n_{i_n}) \right).$$

155. Let V, f, B and $B^1$ be as in problem (153) suppose $[f]_B = A$ i.e., $((f_1)_{B_1} \mid \ldots \mid (f_n)_{B_n}) =$

$$\left( \begin{array}{c|c|c|c} A_1 & 0 & & 0 \\ \hline 0 & A_2 & & 0 \\ \hline & & & \\ \hline 0 & 0 & & A_n \end{array} \right).$$

Show that $\beta_i = (\beta^1_{i_1} \mid \ldots \mid \beta^n_{i_n}) = \sum_j (A^{-1})_{ij}\alpha_j$

$$= \left( \sum_{j_1} (A_1^{-1})_{i_1 j_1}\, \alpha^1_{j_1} \mid \ldots \mid \sum_{j_n} (A_n^{-1})_{i_n j_n}\, \alpha^n_{j_n} \right)$$

$$= \left( \sum_{j_1} (A_1^{-1})_{j_1 i_1} \alpha_{j_1}^1 | \ldots | \sum_{j_n} (A_n^{-1})_{j_n i_n} \alpha_{j_n}^n \right)$$

$$= \sum_j (A^{-1})_{ji} \alpha_j .$$

156. Let $V = (V_1 | \ldots | V_n)$ be a finite $(n_1, \ldots, n_n)$ dimensional super vector space over the field F and $f = (f_1 | \ldots | f_n)$ be a symmetric bilinear super form on V. For each super subspace $W = (W_1 | \ldots | W_n)$ of V let $W^\perp = (W_1^\perp | \ldots | W_n^\perp)$ be the set of all super vector $\alpha = (\alpha_1 | \ldots | \alpha_n)$ in V such that $f(\alpha, \beta) = (f_1(\alpha_1, \beta_1) | \ldots | f_n(\alpha_n, \beta_n)$ in W show that

a. $W^\perp$ is a super subspace

b. $V = \{0\}^\perp$ i.e., $(V_1 | \ldots | V_n) = \{\{0\}^\perp | \ldots |\{0\}^\perp\}$.

c. $V^\perp = (V_1^\perp | \ldots | V_n^\perp) = \{0 | \ldots | 0\}$ if and only if $f = (f_1 | \ldots | f_n)$ is super non degenerate i.e., if and only if each $f_i$ is non degenerate for i = 1, 2, …, n.

d. super rank f = (rank $f_1$, …, rank $f_n$) = super dim V - super dim $V^\perp$ i.e., (dim $V_1$ – dim $V_1^\perp$, …, dim $V_n$ – dim $V_n^\perp$).

e. If super dim V = (dim$V_1$, …, dim $V_n$) and super dim W = (dim $W_1$, …, dim $W_n$) = $(m_1, \ldots, m_n)$ ($m_i < n_i$ for i = 1, 2, …, n) then super dim $W^\perp$ = (dim $W_1^\perp$, …, dim $W_n^\perp$) $\geq (n_1 - m_1, \ldots, n_n - m_n)$.
(Hint: If $(\beta_1^1 \ldots \beta_{m_1}^1; \ldots; \beta_1^n \ldots \beta_{m_n}^n)$ is a super basis of $W = (W_1 | \ldots | W_n)$, consider the super map;

$$(\alpha_1 | \ldots | \alpha_n) \rightarrow$$
$$(f_1(\alpha^1, \beta_1^1), \ldots, f_1(\alpha^1, \beta_{m_1}^1) | \ldots | (f_n(\alpha^n, \beta_1^n), \ldots, f_n(\alpha^n, \beta_{m_n}^n))$$

of V into $(F^{m_1} | \ldots | F^{m_n})$.

f. The super restriction of f to W is super non degenerate if and only if

$W \cap W^{\perp} = (W_1 \cap W_1^{\perp} | \ldots | W_n \cap W_n^{\perp})$

$= (0 | \ldots | 0)$.

g. $V = W \oplus W^{\perp} = (V_1 | \ldots | V_n) = (W_1 \oplus W_1^{\perp} | \ldots | W_n \oplus W_n^{\perp})$

if and only if the super restriction of $f = (f_1 | \ldots | f_n)$ to $W = (W_1 | \ldots | W_n)$ is super non generate ie each $f_i$ to $W_i$ is non generate for $i = 1, 2, \ldots, n$.

157. Let $S_s$ and $T_s$ be super positive operators. Prove that every characteristic super value of $S_s T_s$ is super positive.

158. Prove that the product of two super positive linear operators $T_sU_s = (T_1U_1 | \ldots | T_nU_n)$ is positive if and only if they super commute i.e., if and only if $T_iU_i = U_iT_i$ for every $i = 1, 2, \ldots, n$.

159. If

$$A = \left(\begin{array}{c|c|c|c} A_1 & 0 & & 0 \\ \hline 0 & A_2 & & 0 \\ \hline & & & \\ \hline 0 & 0 & & A_n \end{array}\right)$$

is a super self adjoint $(n_1 \times n_1, \ldots, n_n \times n_n)$ super diagonal matrix i.e., each $A_i$ is a $n_i \times n_i$ matrix; $i = 1, 2, \ldots, m$.

Prove that there is a real n-tuple of numbers $c = (c_1, \ldots, c_n)$ such that the super diagonal matrix $cI + A$

$$= \left(\begin{array}{c|c|c|c} c_1I_1 + A_1 & 0 & & 0 \\ \hline 0 & c_2I_2 + A_2 & & 0 \\ \hline & & & \\ \hline 0 & 0 & & c_nI_n + A_n \end{array}\right)$$

is super positive.

160. Obtain some interesting results on super linear algebra $A = (A_1 \mid \ldots \mid A_n)$ over the field of reals.

161. Let $V = (V_1 \mid \ldots \mid V_n)$ be a finite $(n_1, \ldots, n_n)$ dimensional super inner product space. If $T_s = (T_1 \mid \ldots \mid T_n)$ and $U_s = (U_1 \mid \ldots \mid U_n)$ are linear operators on V we write $T_s < U_s$ if $U - T = (U_1 - T_1 \mid \ldots \mid U_n - T_n)$ is a super positive operator ie each $U_i - T_i$ is a positive operator on $V_i$; $i = 1, 2, \ldots, n$.
Prove the following

a. $T_s < U_s$ then $U_s < T_s$ is impossible.
b. If $T_s < U_s$ and $U_s < P_s$ then $T_s < P_s$.
c. If $T_s < U_s$ and $0 < P_s$ ; it need not imply that $P_s T_s < P_s U_s$.

i.e., each $P_i T_i < P_i U_i$ may not hold good for each i even if $T_i < U_i$ and $0 < P_i$ for $i = 1, 2, \ldots, n$.

# FURTHER READING

# INDEX

**B**

**C**

**D**

**N**

**O**

**P**

**Q**

**T**

# ABOUT THE AUTHORS

**Dr.W.B.Vasantha Kandasamy** is an Associate Professor in the Department of Mathematics, Indian Institute of Technology Madras, Chennai. In the past decade she has guided 12 Ph.D. scholars in the different fields of non-associative algebras, algebraic coding theory, transportation theory, fuzzy groups, and applications of fuzzy theory of the problems faced in chemical industries and cement industries.

She has to her credit 640 research papers. She has guided over 64 M.Sc. and M.Tech. projects. She has worked in collaboration projects with the Indian Space Research Organization and with the Tamil Nadu State AIDS Control Society. This is her $35^{th}$ book.

On India's 60th Independence Day, Dr.Vasantha was conferred the Kalpana Chawla Award for Courage and Daring Enterprise by the State Government of Tamil Nadu in recognition of her sustained fight for social justice in the Indian Institute of Technology (IIT) Madras and for her contribution to mathematics. (The award, instituted in the memory of Indian-American astronaut Kalpana Chawla who died aboard Space Shuttle Columbia). The award carried a cash prize of five lakh rupees (the highest prize-money for any Indian award) and a gold medal.

She can be contacted at vasanthakandasamy@gmail.com
You can visit her on the web at: http://mat.iitm.ac.in/~wbv

---

**Dr. Florentin Smarandache** is a Professor of Mathematics and Chair of Math & Sciences Department at the University of New Mexico in USA. He published over 75 books and 150 articles and notes in mathematics, physics, philosophy, psychology, rebus, literature.

In mathematics his research is in number theory, non-Euclidean geometry, synthetic geometry, algebraic structures, statistics, neutrosophic logic and set (generalizations of fuzzy logic and set respectively), neutrosophic probability (generalization of classical and imprecise probability). Also, small contributions to nuclear and particle physics, information fusion, neutrosophy (a generalization of dialectics), law of sensations and stimuli, etc. He can be contacted at smarand@unm.edu